\documentclass[singlespacing,manychapters,school=upenn,oneside,allcapstitle]{puclassicthesis}

\usepackage{thesis_preamble}

\begin{document}

\makefrontmatter

\tocundopartskip

\ctparttext{This part serves multiple groups of people and can be used in different ways:
\begin{itemize}
\item For those who are category theory neophytes, a reading of Chapter~\ref{sec:categories} is advised, after which they should move on to Part~\ref{part:easy_sheaves}, with particular emphasis on the beginning of Chapter~\ref{sec:cell_sheaves} and Chapter~\ref{sec:homology}.

\item Chapter~\ref{sec:abstract_sheaves} is designed for those who want a general definition of sheaves and cosheaves on a topological space. After looking at the definition, one should proceed as quickly as possible to Chapter~\ref{sec:examples} to get some simple examples. 

\item Section~\ref{subsec:cover_structure} is meant for people who have always found the expression of the sheaf axiom as an exact sequence a little opaque. Such people are usually frustrated by the notation used in \v{C}ech homology, which is the subject of Section~\ref{subsec:cech}. 

\item Sections~\ref{subsec:generalities} and~\ref{subsec:cosheafification} are for those who think of cosheaves simply as sheaves valued in the opposite category.
\end{itemize}
}
\part{A Mathematical Introduction}
\label{part:abstract_intro}

%
%

\chapter{A Primer on Category Theory}
\label{sec:categories}

\begin{quote}
{\em ``A healthy new seed was planted some twenty odd years ago in the well fertilized soil of the mathematical periodical literature --- the notion of a category. It sprouted, took root, flowered, attracted bees, and by now the landscape is dotted with its progeny. It is a beneficent plant: mathematical gardeners have come to appreciate its usefulness in holding down the topsoil and preventing dust storms; indeed, some half dozen books have appeared within the past dozen years putting it to this use. It is a beautiful plant too, whose rapid proliferation has produced many unique and exotic variants; but, perhaps because of its increasingly multiform variety, the book extolling all its loveliness has not yet been written.''}
\begin{flushright} --- F.E.J. Linton~\cite{linton1965} \end{flushright}
\end{quote}

Categories emerged out of the study of functors, which were originally conceived as a principled way of assigning algebraic invariants to topological spaces. Thus, category theory is part and parcel of the study of algebraic topology. However, from its conception in Samuel Eilenberg and Saunders Mac Lane's 1945 paper on a ``General Theory of Natural Equivalence''~\cite{em-cat45}, it was realized that the language of categories provides a way of identifying formal similarities throughout mathematics. The success of this perspective is largely due to the fact that category theory --- as opposed to set theory --- emphasizes understanding the relationships between objects rather than the objects themselves.

In this section, we provide a brief review of the parts of category theory needed to understand the abstract definitions of a sheaf and cosheaf in Chapter~\ref{sec:abstract_sheaves}. Most importantly, the reader should be able to do the following before moving onto that section:
\begin{itemize}
	\item Think of the set of open sets of a topological space $X$ as a category.
	\item Understand how to summarize the behavior of various functors via limits and colimits.
\end{itemize}

We have tried to provide a self-contained introduction to category theory, but the reader is urged to consult Mac Lane's ``Categories for the Working Mathematician''~\cite{catwork} for a book that very well may be the book anticipated by the quote above.

\section{Categories}

One should visualize categories as graphs with objects corresponding to vertices and maps as edges between vertices, subject to relations that specify when following one sequence of edges is equivalent to another sequence. One can think of some of the axioms of a category as gluing in triangles and tetrahedra to witness these relations.
\[
	\xymatrix{ & & \\ \bullet \ar[rr] & & \bullet }  \qquad
	\xymatrix{ & \bullet \ar[dr] & \\ \bullet \ar[ru] \ar[rr] & & \bullet } \qquad
	\xymatrix{ & \bullet \ar[dr] \ar[r] & \bullet \\ \bullet \ar[ru] \ar@{.>}[rru] \ar[rr] & & \bullet \ar[u] }
\]

\begin{defn}[Category]\index{category}
A \textbf{category} $\cat$ consists of a class of objects denoted $\obj(\cat)$ and a set of morphisms $\Hom_{\cat}(a,b)$ between any two objects $a,b\in\obj(\cat)$. An individual morphism $f:a\to b$ is also called an arrow since it points (maps) from $a$ to $b$. We require that the following axioms hold:
\begin{itemize}
	\item Two morphisms $f\in \Hom_{\cat}(a,b)$ and $g\in\Hom_{\cat}(b,c)$ can be composed to get another morphism $g\circ f\in \Hom_{\cat}(a,c)$.
	\item Composition is associative, i.e. if $h\in\Hom(c,d)$, then $(h\circ g)\circ f=h \circ (g\circ f)$.
	\item For each object $x$ there is an identity morphism $\id_x\in\Hom_{\cat}(x,x)$ that satisfies $f\circ \id_a=f$ and $\id_b\circ f=f$.
\end{itemize}
When the category $\cat$ is understood, we will sometimes write $\Hom(a,b)$ to mean $\Hom_{\cat}(a,b)$.
\end{defn}

One can usually ignore the technicality that the collection of objects forms a class rather than a set. A class is a collection of sets that one can refuse to quantify over in a logical sense. This prohibits Russell-type paradoxes gotten by considering the category of all categories that do not contain themselves. Colloquially, one says a proper class is ``bigger'' than a set. In order to avoid certain machinery that accompanies the use of classes, we will often consider categories that are ``small'' in a precise sense.\footnote{The machinery we are referring to is that of Grothendieck universes.}

\begin{defn}[Small and Finite Categories]\index{category!small}\index{category!finite}
	A category is \textbf{small} if its class of objects is actually a set. A category is \textbf{finite} if its set of objects has finite cardinality.
\end{defn}

\begin{ex}[Discrete Category]\index{category!discrete}
	Any set $X$ can be regarded as a \textbf{discrete category} $\bar{X}$ with only the identity morphism $\id_x$ sitting over each object. There are no non-identity morphisms.
\end{ex}

Recall that a \textbf{relation}\index{relation} $R$ on a set $X$ is a subset of the product set $X\times X$. If two elements are related by $R$, one writes $x R y$ to mean that $(x,y)\in R$. We now give an example of some relations on a set that endow that set with the structure of a category.

\begin{ex}[Posets and Preorders]
	A \textbf{preordered set}\index{preordered set} is a set $X$ along with a relation $\leq$ that satisfies the following two axioms:
\begin{description}
 \item[Reflexivity] --- $x\leq x$ for all $x\in X$
 \item[Transitivity] --- $x\leq y$ and $y\leq z$ implies $x\leq z$
\end{description}
	A \textbf{partially ordered set}, or \textbf{poset}\index{poset} for short, is a preordered set that additionally satisfies the following third axiom:
	\begin{description}
		\item[Anti-Symmetry] --- $x\leq y$ and $y\leq x$ implies $x=y$
	\end{description}
Any preordered set $(X,\leq)$ defines a category by letting the objects be the elements of $X$ and by declaring each Hom set $\Hom(x,y)$ to either have a unique morphism if $x\leq y$ or to be empty if $x\nleq y$. 
\end{ex}

We now reach our example of fundamental importance.

\begin{ex}[Open Set Category]\index{category!of open sets}
 The \textbf{open set category} associated to a topological space $X$, denoted $\Open(X)$, has as objects the open sets of $X$ and a unique morphism $U\to V$ for each pair related by inclusion $U\subseteq V$.
\end{ex}

 There is an example very closely aligned with the category of open sets that is allegedly due to Raoul Bott, who gave it as an example of a topological category~\cite{bott-mexico,bott-legacy}.

\begin{ex}[Pointed Open Set Category]\index{category!of points and open sets}
The \textbf{pointed open set category} $\Open^*(X)$ associated to a topological space $X$ has pairs $(U,x)$, where $U$ is an open set and $x$ is a point in $U$, for objects and a unique morphism $(U,x)\to (V,y)$ if $U\subset V$ and $x=y$.
\end{ex}

This pointed open set category takes us nicely over to a category whose objects are points of a topological space. First, we introduce some terminology.

\begin{defn}[Groupoid]\label{defn:groupoid}\index{category!groupoid}\index{groupoid}
A \textbf{groupoid} is a category where every morphisms is invertible. In other words if $\bsG$ is a groupoid, then for every pair of objects $x,y\in\obj(\bsG)$ and every morphism $\alpha\in\Hom_{\bsG}(x,y)$ there exists a morphism $\beta\in \Hom_{\bsG}(y,x)$ such that $\alpha\circ\beta=\id_y$ and $\beta\circ\alpha=\id_x$.
\end{defn}

\begin{exr}
Let $\bsG$ be a groupoid with only one object. Show that the structure axioms of a category along with the property of being a groupoid guarantees that $\bG$ is a group. Observe that this gives us a way of treating every group as a category, where multiplication in the group corresponds to composition of morphisms.
\end{exr}

\begin{defn}[The Fundamental Groupoid]\label{defn:fundamental_groupoid}\index{category!fundamental groupoid $\pijuan(X)$}\index{fundamental groupoid $\pijuan(X)$}
Let $X$ be a topological space. The \textbf{fundamental groupoid} $\pijuan(X)$ has points $x\in X$ for objects and homotopy classes of paths relative endpoints for morphisms. Specifically,
\[
	\Hom_{\pijuan(X)}(x,y):=\{\gamma:[0,1]\to X\,|\,\gamma(0)=x, \gamma(1)=y\}/\sim
\]
where $\gamma\sim\gamma'$ if there exists a third continuous map $h:[0,1]^2\to X$ such that $h(0,t)=\gamma(t)$, $h(1,t)=\gamma'(t)$, $h(s,0)=x$ and $h(s,1)=y$.
\end{defn} 

\begin{rmk}[Poincar\'e $\infty$-Groupoid]\index{fundamental groupoid $\pijuan(X)$!Poincar\'e $\infty$-groupoid}
	To a topological space $X$, one can consider a generalization of the fundamental groupoid, called the \textbf{Poincar\'e $\infty$-groupoid} $\piinfty(X)$, which has an object for each point of $X$, a morphism for every path $\gamma:[0,1]\to X$ , a ``2-morphism'' for every continuous map $\sigma:\Delta^2\to X$, and so on for higher $\Delta^n$. The 2-morphisms should be regarded as providing a homotopy between $\sigma|_{0,2}$ and $\sigma|_{1,2}\circ\sigma|_{0,1}$, i.e. a morphism between morphisms. Here $\sigma|_{i,j}$ is the restriction of the map $\sigma$ to the edge going from vertex $i$ to $j$. As stated, this is an example of an $\infty$-category, which is currently vying to replace ordinary category theory as the foundation for mathematics~\cite{lurie-htt}.  
\end{rmk}

The above examples of categories are quite small when compared to the categories that Eilenberg and Mac Lane first introduced. The categories considered there correspond to \textbf{data types}\index{category!data}\index{data types} and we will usually refer to them with the letter $\dat$. For this paper $\dat$ will usually mean one of the following: 
\begin{description}
	\item[$\Set$] --- the category whose objects are sets and whose morphisms are all set maps (multi-valued maps are prohibited as are partially defined maps)
	\item[$\Ab$] --- the category whose objects are abelian groups and whose morphisms are group homomorphisms
	\item[$\Vect$] --- the category whose objects are vector spaces and whose morphisms are linear transformations
	\item[$\vect$] --- the category whose objects are \emph{finite-dimensional} vector spaces and linear transformations
	\item[$\Top$] --- the category whose objects are topological spaces and whose morphisms are continuous maps
\end{description}
The category $\vect$ is an example of a subcategory, which we now define.
\begin{defn}[Subcategories]\index{category!subcategory}
	Let $\cat$ be a category. A \textbf{subcategory} $\bat$ of $\cat$ consists of a subcollection of objects from $\cat$ and a choice of subset of the morphism set $\Hom_{\cat}(x,y)$ for each pair $x,y\in\obj(\bat)$. We require that these morphism sets have the identity and be closed under composition so as to guarantee that $\bat$ is a category. We say that a subcategory is \textbf{full} if $\Hom_{\bat}(x,y)=\Hom_{\cat}(x,y)$. 
\end{defn}

Categories have a built-in notion of directionality. For example, in $\Set$ every object $X$ has a unique map from the empty set $\emptyset$, but there are no maps to the empty set. We can abstract out this property, so as to make it apply in other situations.

\begin{defn}[Initial and Terminal Objects]\index{initial object}\index{terminal object}
	An object $x\in\obj(\cat)$ is said to be \textbf{initial} if for any other object $y\in\obj(\cat)$ there is a unique morphism from $x$ to $y$. Dually, an object $y$ is said to be \textbf{terminal} if for any object $x$ there is a unique morphism from $x$ to $y$. 
\end{defn}

As already mentioned, in $\Set$ the empty set is initial, but it is not terminal. On the contrary, the terminal object is the one point set $\{\star\}$ since there is only one constant map. Similarly, for $\Open(X)$ the empty set is initial, but the whole space $X$ is terminal. In $\Vect$ the initial and terminal objects coincide with the zero vector space. In some sense, the difference between the initial and terminal objects in a category measure how different it is from its reflection. We now say what we mean by a category's reflection. 

\begin{ex}[Opposite Category]\index{category!opposite}
 For any category $\cat$ there is an \textbf{opposite category} $\cat^{op}$ where all the arrows have been turned around, i.e. $\Hom_{\cat^{op}}(x,y)=\Hom_{\cat}(y,x)$.
\end{ex}
\begin{rmk}[Duality and Terminology]
Because one can always perform a general categorical construction in $\cat$ or $\cat^{op}$ every concept is really two concepts. As we shall see, this causes a proliferation of ideas and is sometimes referred to as the \textbf{mirror principle}. The way this affects terminology is that a construction that is dualized is named by placing a ``co'' in front of the name of the un-dualized construction. Thus, as we will see shortly, there are limits and colimits, products and coproducts, equalizers and coequalizers, among other things.
\end{rmk}

Now we introduce the fundamental device that assigns objects and morphisms in one category to objects and morphisms in another category. Historically, this device was introduced first and categories were summoned into existence to provide a domain and range for this assignment. 

\begin{defn}[Functor]\index{functor}
	A \textbf{functor} $F:\cat\to\dat$ consists of the following data: To each object $a\in\cat$ an object $F(a)\in\dat$ is associated, i.e. $a\rightsquigarrow F(a)$. To each morphism $f:a\to b$ a morphism $F(f):F(a)\to F(b)$ is likewise associated. We require that the functor respect composition and preserve identity morphisms, i.e. $F(f\circ g)=F(f)\circ F(g)$ and $F(\id_a)=\id_{F(a)}$. For such a functor $F$, we say $\cat$ is the \textbf{domain} and $\dat$ is the \textbf{codomain} of $F$.
\end{defn}
\begin{rmk}
We can phrase the definition of a functor differently by saying that we have a \emph{function} $F:\obj(\cat)\to\obj(\dat)$ and \emph{functions} $F(a,b):\Hom_{\cat}(a,b)\to\Hom_{\dat}(F(a),F(b))$ for every pair of objects $a,b\in\obj(\cat)$. We require that these functions preserve identities and composition. When $F(a,b):\Hom_{\cat}(a,b)\to\Hom_{\dat}(F(a),F(b))$ is injective for every pair of objects we say $F$ is \textbf{faithful}. When $F(a,b)$ is surjective for every pair of objects we say $F$ is \textbf{full}. When a functor is both full and faithful, we say it is \textbf{fully faithful}.
\end{rmk}

\begin{exr}
Check that the definition of a subcategory guarantees that the inclusion $\bat\hookrightarrow \cat$ is a functor.
\end{exr}

An example familiar to every topologist is that of homology and cohomology with field coefficients. In every non-negative degree $i$, these invariants define functors
\[
	H_i(-;k):\Top\to\Vect \qquad \mathrm{and} \qquad H^i(-;k):\Top^{op}\to\Vect
\]
respectively. Here we have used the opposite category as an alternative way of saying cohomology is \textbf{contravariant}.
     
Historically, there was a plethora of different homology theories --- simplicial, singular, \v{C}ech, Vietoris, Alexander, et al --- and every time one was introduced a long repetition of the basic properties of that homology theory ensued. Understanding the precise relationships between these motivated the notion of a map between functors, which led in turn to the Eilenberg-Steenrod axioms~\cite[p.335]{maclane-concepts}.

\begin{defn}[Natural Transformation]\index{natural transformation}
	Given two functors $F,G:\cat\to\dat$ a \textbf{natural transformation}, sometimes written $\eta:F \Rightarrow G$, consists of the following information: to each object $a\in\cat$, a morphism $\eta(a):F(a)\to G(a)$ is assigned such that for every morphism $f:a\to b$ in $\cat$ the following diagram \textbf{commutes}:
	\begin{displaymath}
		\xymatrix{F(a) \ar[r]^{\eta(a)} \ar[d]_{F(f)} & G(a)\ar[d]^{G(f)}\\F(b) \ar[r]^{\eta(b)} & G(b)}
	\end{displaymath}
	By commutes, we mean $G(f)\circ \eta(a)=\eta(b)\circ F(f)$.
\end{defn}

\begin{defn}\index{naturally isomorphic}
	Two functors $F,G:\cat\to\dat$ are said to be \textbf{naturally isomorphic} if there is a natural transformation $\eta:F\Rightarrow G$ such that for every object $a\in \cat$ the morphism $\eta(a)$ is an isomorphism, i.e. it is invertible. These inverse maps $\eta(a)^{-1}$ define an inverse natural transformation $\eta^{-1}:G\Rightarrow F$.
\end{defn}

Functors and natural transformations assemble themselves into a category in their own right. Since an arrow is an arrow by any other symbol, we will sometimes use the notation $F\to G$ to denote a natural transformation, instead of $F\Rightarrow G$. In the functor category, we will see that naturally isomorphic functors are isomorphic objects. This demonstrates again the linguistic efficiency of category theory.

\begin{ex}[Functor Category]\index{category!of functors}
 $\Fun(\cat,\dat)$ denotes the category whose objects are functors from $C$ to $D$ and whose morphisms are natural transformations.
\end{ex}

Certain functors deserve special attention. These are the ones that allow us to identify two different categories. One approach to identifying categories is to say that two categories $\cat$ and $\dat$ are \textbf{isomorphic} if there are functors $F:\cat\to\dat$ and $G:\dat\to\cat$ such that $G\circ F=\id_{\cat}$ and $F\circ G=\id_{\dat}$. This definition is so restrictive that it rarely occurs. Thus, we have a looser notion that includes isomorphism as a special case. Instead of asking that $F\circ G$ be equal to $\id_{\dat}$, we only require that they be isomorphic as objects in $\Fun(\dat,\dat)$ and similarly for $G\circ F$ and $\id_{\cat}$ in $\Fun(\cat,\cat)$. The reader should compare this with the notion of homotopy equivalence. 

\begin{defn}\index{equivalence!adjoint}
A pair of functors $F:\cat\to\dat$ and $G:\dat\to\cat$ together define an \textbf{adjoint equivalence} of categories if there are two natural isomorphisms of functors $\epsilon:F\circ G\to \id_{\dat}$ and $\eta:\id_{\cat}\to G\circ F$.
\end{defn}
We will see that this notion of a equivalence is a special instance of an adjunction, which is taken up in Section~\ref{sec:adjunctions}

Equivalence can also be phrased in a way that doesn't require us to construct $G$ as a ``weak inverse'' of $F$.
\begin{defn}[Fully Faithful and Essentially Surjective]\index{equivalence!fully faithful and essentially surjective}
	A functor $F:\cat\to\dat$ induces an \textbf{equivalence} of categories if it is bijective on $\Hom$ sets (fully faithful) and is \textbf{essentially surjective}. This last property means that for every object $d\in\dat$ there is an object $c\in\cat$ such that $F(c)$ is isomorphic to $d$, i.e. $F$ is bijective on isomorphism classes of $\cat$ and $\dat$.
\end{defn}

The notion of equivalence allows us to find compressed presentations of a category.

\begin{defn}[Skeletal Subcategory]\index{category!skeletal subcategory}\index{skeletal subcategory}
	Suppose $\cat$ is a category, then a subcategory $\covS$ is \textbf{skeletal} if the inclusion functor is an equivalence, and no two objects of $\covS$ are isomorphic.
	
	If $\cat$ is small, then we can describe explicitly how to construct a skeletal subcategory $\covS$. On the objects of $\cat$ we define an equivalence relation that says $x\sim x'$ if and only if $x$ and $x'$ are isomorphic. To define a skeletal subcategory we pick one object $x\in\bar{x}$ from each equivalence class and define the morphisms to be $\Hom_{\covS}(\bar{x},\bar{y}):=\Hom_{\cat}(x,y)$.
\end{defn}

\begin{exr}[Fundamental groupoid]
Suppose $X$ is a path connected space. Show that for any point $x_0\in X$, the fundamental group $\pijuan(X,x_0)$ is a skeletal subcategory of $\pijuan(X)$.
\end{exr}

Finally, let's analyze how working in the opposite category impacts functors and natural transformations. Observe, first and foremost, that formality allows us to take a functor $F:\cat\to\dat$ and define a functor $F^{op}:\cat^{op}\to\dat^{op}$. Moreover, a natural transformation $\eta:F\Rightarrow G$ translates to a natural transformation $\eta^{op}:G^{op}\Rightarrow F^{op}$. This observation allows us to state the equalities
\[
	\Fun(\cat^{op},\dat^{op})=\Fun(\cat,\dat)^{op} \qquad \mathrm{or} \qquad \Fun(\cat^{op},\dat^{op})^{op}=\Fun(\cat,\dat)
\]
since $(\cat^{op})^{op}$ is isomorphic to $\cat$ (not just equivalent). See the wonderful work ``Abstract and Concrete Categories: The Joy of Cats''~\cite{joyofcats} for more on duality and category theory more generally.

\section{Diagrams and Representations}
\label{subsec:reps}
Categories and functors allow us to develop an \emph{algebra of shape}, the shapes being modeled on the domain category of a functor. For example, we will be interested in studying data arranged in the following forms:
\[
	\xymatrix{\bullet \ar[r] \ar[d] & \bullet \\ \bullet &} \qquad \qquad
	\xymatrix{ \bullet \ar@<-.5ex>[rd] \ar@<.5ex>[rd] & \bullet \\ \bullet & \bullet} \qquad \qquad
	\xymatrix{& \bullet\ar[d] \\ \bullet \ar[r] &\bullet}
\]
If we imagine the identity arrows in a category as being the vertices themselves, and thus not drawn independently of the objects, each of these shapes gives an example of a finite category.

\begin{defn}[Diagram]\index{diagram}
	Suppose $\Iat$ is a small category and $\cat$ is an arbitrary category. A \textbf{diagram} is simply a functor $F:\Iat\to \cat$.
\end{defn}

\begin{ex}[Constant Diagram]
	For any category $\Iat$ there is always a diagram for each object $O\in \cat$, called the \textbf{constant diagram}, $\mathrm{const}_O:\Iat\to\cat$ where $\mathrm{const}_O(x)=\mathrm{const}_O(y)=O$ for all objects $x,y\in \Iat$. Every morphism in $\Iat$ goes to the identity morphism.
\end{ex}

\begin{defn}[Representation]\index{representation}
	A \textbf{representation} of a category $\cat$ is a functor $F:\cat\to\Vect$. 
\end{defn}

One should note that this definition generalizes the notion of a representation of a group. Every group, say $\ZZ$ for example, can be considered as a small category with a single object $\star$ and $\Hom(\star,\star)=\ZZ$. A representation of $\ZZ$ then corresponds to picking a vector space $V$ and assigning an endomorphism of $V$ for each element of $\ZZ$, i.e. it is a functor.
\[
	\xymatrix{\star \ar@{~>}[r] \ar[d]_g & V \ar[d]^{\rho(g)} \\
	\star \ar@{~>}[r] & V}
\]
Maps of representations correspond precisely with natural transformations of such functors. Isomorphic representations are naturally isomorphic functors.\footnote{Confusingly, the term ``equivalent representations'' is often used.} These basic notions carry over to the representation theory of arbitrary categories, which allows us to compare different situations in one language.

\section{Cones and Limits}
\label{subsec:limits}
The next two sections are devoted to studying one way (and a dual way) of summarizing a functor's behavior. This gives a way of compressing the data of a functor into a single object. These concepts are fundamental to the study of sheaves and cosheaves.

\begin{defn}[Cone]\index{cones}
	Suppose $F:\Iat\to\cat$ is a diagram. A \textbf{cone} on $F$ is a natural transformation from a constant diagram to $F$. Specifically, it is a choice of object $L\in\cat$ and a collection of morphisms $\psi_{x}:L\to F(x)$, one for each $x$, such that if $g:x\to y$ is a morphism in $\Iat$, then $F(g)\circ 	\psi_x=\psi_y$, i.e. the following diagram commutes:
\[
\xymatrix{F(x) \ar[rr]^{F(g)} & & F(y)\\ & L \ar[ul]^{\psi_x} \ar[ru]_{\psi_y} &}
\]
In other words, $\psi_y=F(g)\circ \psi_x$.
\end{defn}

\begin{defn}\index{cones!category of}\index{category!of cones}
	The collection of cones on a diagram $F$ form a category, which we will call $\Cone(F)$. The objects are cones $(L,\psi_x)$ and a morphism between two cones $(L',\psi_x')$ and $(L,\psi_x)$ consists of a map $u:L'\to L$ such that $\psi_x'=\psi_x\circ u$ for all $x$
\end{defn}

A limit is simply a distinguished or universal object in the category of cones on $F$.

\begin{defn}[Limit]\index{limit}
	The \textbf{limit} of a diagram $F:\Iat\to\cat$, denoted $\varprojlim F$ is the terminal object in $\Cone(F)$. This means that a limit is an object $\varprojlim F\in\cat$ along with a collection of morphisms $\psi_x:L\to F(x)$ that commute with arrows in the diagram such that whenever there is another object $L'$ and morphisms $\psi'_x$ that also commute there then exists is a unique morphism $u:L'\to \varprojlim F$ that additionally commutes with everything in sight, i.e. $\psi_x'=\psi_x\circ u$ for all $x$.
	\[
	\xymatrix{F(x) \ar[rr]^{F(g)} & & F(y)\\ & \varprojlim F \ar[ul]_{\psi_x} \ar[ru]^{\psi_y} & \\ & L' \ar[uul]^{\psi'_x} \ar@{.>}[u]^{\exists!}_u  \ar[uur]_{\psi'_y}  &}
	\]
\end{defn}

\begin{rmk}[Glossary]\index{limit!synonyms for}
	Quite confusingly, the following terms are synonyms for limits: inverse limits, projective limits, left roots, $\lim$ and $\varprojlim$ are all common.
\end{rmk}

We now consider some examples of limits over discrete categories.

\begin{ex}[Products]\label{ex:products}\index{product}
	Consider the following index category and diagram:
	\[
		\xymatrix{\bullet & \bullet} \qquad \qquad \qquad \xymatrix{F(i) & F(j)}
	\]
	The limit of this diagram is called the \textbf{product} and is usually written
	\[
		F(i) \prod F(j).
	\]
	More generally, we define the product to be the limit of any diagram $F:\Iat\to\cat$ indexed by a discrete category and write $\prod_i F(i)$. Sometimes one writes $\times_i F(i)$ for the product.
\end{ex}

We give an unusual example of a product that will prepare the reader for thinking about the category of open sets.

\begin{ex}[Open Sets: Limits are Intersections]\index{limit!is an intersection}
	Suppose $\Lambdat=\{1,\ldots, n\}$ is a finite discrete category, i.e. it has $n$ objects and the only morphisms are the identity morphisms. Now let $X$ be a topological space and let $\cat=\Open(X)$ be the category of open sets in $X$. This is a category that has an object for each open set and a single morphism $U\to V$ if $U\subset V$. A functor $F:\Lambdat\to\Open(X)$ is nothing more than a choice of $n$ not necessarily distinct open sets. A cone to $F$ is an open set that includes into all the open sets picked out by $F$. The limit of $F$ is the largest possible open set that includes into all the open sets picked out by $F$, i.e. $$\varprojlim F=\cap_{i=1}^n F(i).$$
\end{ex}

\begin{ex}
	Consider the following small category $\Iat$ along with some representation $F:\Iat\to\Vect$.
	\[
	\xymatrix{\bullet \ar[r] \ar[d] & \bullet \\ \bullet &} \qquad \qquad \qquad
	\xymatrix{U \ar[r]^A \ar[d]_B & V \\ W &}
	\]
	By thinking about the definition, one can see that
	\[
		\varprojlim F\cong U.
	\]
\end{ex}

\begin{ex}[Pullbacks]\label{ex:pullback}\index{pullback}	
	Consider the category $\Jat=\Iat^{op}$ and a representation $F:\Jat\to\Vect$.
	\[
	\xymatrix{& \bullet\ar[d] \\ \bullet \ar[r] &\bullet}
	\qquad \qquad \qquad
	\xymatrix{ & V\ar[d]^A \\ W \ar[r]_B & U}
	\]
	With some thought one can describe the limit set-theoretically as
	\[
		\varprojlim F\cong \{(v,w)\in V\times W | Av=Bw\},
	\]
	which is called the \textbf{pullback}. If $U=0$, then we re-obtain the product of $V$ and $W$ and one usually writes $V\times W$. 
\end{ex}

\begin{ex}[Equalizers and Kernels]\label{ex:equalizer}\index{equalizer}\index{kernel}
	Consider the following category $\Kat$ and an arbitrary functor $F:\Kat\to\dat$.
	\[
	 \xymatrix{ \bullet \ar@<-.5ex>[r] \ar@<.5ex>[r] & \bullet}
	\qquad \qquad \qquad
	 \xymatrix{ X \ar@<-.5ex>[r]_{g} \ar@<.5ex>[r]^{f} & Y}
	\]
	The limit of this diagram, which is also called the \textbf{equalizer}, is an object $E$ along with a map $h$ that satisfies $f\circ h=g\circ h$.
	\[
	 \xymatrix{E \ar[r]^h & X \ar@<-.5ex>[r]_g \ar@<.5ex>[r]^f & Y}
	\]
	If $\dat=\Vect$ and one sets $g=0$, then the equalizer is the \textbf{kernel}. Thus, if one wants to mimic kernels in data types lacking of zero maps and objects, equalizers can be substituted.
\end{ex}

Finally, we finish with an example from representation theory.

\begin{ex}[Invariants]\label{ex:invariants}\index{invariants}
	Suppose that $V$ is a vector space with an endomorphism $T:V\to V$, i.e. a $k[x]$-module. Just as a group can be viewed as a category with one object, a ring can be viewed as a category with multiplication corresponding to composition of morphisms and addition corresponding to addition of morphisms, thus such a category has extra structure. Thus the $k[x]$-module determined by $V$ and $T$ is equivalent to a functor $k[x]\to\Vect$ that sends the unique object $\star$ to $V$ and sends $x$ to $T$. The limit of such a functor is called the \textbf{invariants} of the action, i.e.
	\[
		I=\{v\in V \, | \, T(v)=v\}.
	\]
\end{ex}

\section{Co-Cones and Colimits}
\label{subsec:colimits}
Here we invoke the mirror principle to dualize the theory of cones and limits. In accordance with usual terminology, we refer to these as \emph{co}cones and \emph{co}limits.

\begin{defn}[Co-Cone]\index{cocones}
	Given a diagram $F:\Iat\to \cat$, a \textbf{cocone} is a natural transformation from $F$ to a constant diagram. In other words, it consists of an object $C\in\cat$ along with a collection of maps $\phi_x:F(x)\to C$ such that these maps commute with the ones internal to the diagram.
	\[
		\xymatrix{ & C & \\ F(x) \ar[rr]_{F(g)} \ar[ru]^{\phi_x} & & F(y) \ar[ul]_{\phi_y} }
	\]
\end{defn}

Similarly, there is a category of cocones\index{cocones!category of}\index{category!of cocones} to a diagram $F$, denoted $\CoCone(F)$. A colimit is a distinguished object in this category.

\begin{defn}[Colimit]\index{colimit}
	The \textbf{colimit} of a diagram $F$ is the initial object in the category $\CoCone(F)$. One should practice dualizing the explicit description of the limit in order to understand the following diagram: 
	\[
		\xymatrix{
		& C' & \\
		& \varinjlim F \ar@{.>}[u]^{\exists!}_u & \\
		F(x) \ar[ur]_{\phi_x} \ar[uur]^{\phi'_x} \ar[rr]_{F(g)} & & F(y) \ar[ul]^{\phi_y} \ar[uul]_{\phi'_y}
		}
	\]
\end{defn}

\begin{rmk}[Glossary]\index{colimit!synonyms for}
	The following terms are synonyms for colimits: direct limits, inductive/injective limits, right roots, $\colim$ and $\varinjlim$ are all used.
\end{rmk}

To better understand the similarities and differences between limits and colimits, let us re-examine the same examples in the previous section.

\begin{ex}[Coproducts]\label{ex:coproducts}\index{coproduct}
	Consider the following index category and diagram:
	\[
		\xymatrix{\bullet & \bullet} \qquad \qquad \qquad \xymatrix{F(i) & F(j)}
	\]
	The colimit of this diagram is called the \textbf{coproduct} and is usually written
	\[
		F(i) \coprod F(j).
	\]
	More generally, we define the product to be the limit of any diagram $F:\Iat\to\cat$ indexed by a discrete category and write $\coprod_i F(i)$. Alternative notations for the coproduct, depending usually on whether the target category is $\Set,\Vect$, $\Ab$ or $\Top$ include
	\[
		\bigoplus_i F(i) \qquad \mathrm{and} \qquad \sum_i F(i) \qquad \mathrm{and} \qquad \bigsqcup_i F(i).
	\]
\end{ex}

\begin{ex}[Open Sets: Colimits are Unions]\index{colimit!is a union}
	Suppose $\Lambdat=\{1,\ldots, n\}$ is a finite discrete category. Let $\cat=\Open(X)$ be the category of open sets in $X$. A functor $F:\Lambdat\to\Open(X)$ is a choice of $n$ not necessarily distinct open sets. A cocone to $F$ is an open set that contains all the open sets picked out by $F$. The colimit of $F$ is the smallest possible open set containing all the open sets picked out by $F$, i.e. the union: $$\varinjlim F=\cup_{i=1}^n F(i)$$
	One should note that since the arbitrary union of open sets is still open one could have worked over a larger indexing category $\Lambdat$.
\end{ex}

\begin{ex}[Pushouts]\label{ex:pushout}\index{pushout}
	Consider the following small category $\Iat$ and a representation $F:\Iat\to\Vect$.
	\[
	\xymatrix{\bullet \ar[r] \ar[d] & \bullet \\ \bullet &}
	\qquad \qquad \qquad
	\xymatrix{U \ar[r]^A \ar[d]_B & V \\ W &}
	\]
	Contrary to the case of the limit, this one requires a bit more thought. Let's start with something that is \emph{not} a cocone, but is nevertheless naturally built out of pieces of the diagram.
	\[
	\xymatrix{U \ar[r]^A \ar[d]_B \ar[rd]^{B\oplus A} & V \ar[d]^{\iota_V} \\ W \ar[r]_{\iota_W} & W\oplus V}
	\]
	This is not a cocone because the diagram does not commute since $(Bu,0)\neq (Bu,Au)\neq (0,Au)$. We can force commutativity by forcing the equivalence relation $[(Bu,0)]\sim [(0,Au)]$ or equivalently $[(Bu,-Au)]\sim [(0,0)]$. We thus conclude that
	\[
	\varinjlim F=W\oplus V/\im(B\oplus -A) \quad \phi_U=q\circ \iota_WB=q\circ\iota_VA \quad \phi_W=q\circ \iota_W \quad \phi_V=q\circ \iota_W
	\]
	where $q$ is the quotient map. One should note that this is clearly dual to the limit computation in \ref{ex:pullback} with the added complication that whereas the limit is a sub-object, the colimit is a quotient object.
	
	Like before, if $U=0$ then the pushout reduces to the coproduct of $V$ and $W$ and one writes it as $V\oplus W$.
\end{ex}

\begin{ex}
	Consider the example $\Jat=\Iat^{op}$ and corresponding representation $F:\Jat\to\Vect$.
	\[
	\xymatrix{& \bullet\ar[d] \\ \bullet \ar[r] &\bullet}
	\qquad \qquad \qquad
	\xymatrix{ & V\ar[d]^A \\ W \ar[r]_B & U}
	\]
	One can see that
	\[
		\varinjlim F\cong U.
	\]
\end{ex}

\begin{ex}[Coequalizers and Cokernels]\label{ex:coequalizer}\index{cokernel}\index{coequalizer}
	Consider the same category $\Kat$ as before and a functor $F:\Kat\to\dat$.
	\[
	 \xymatrix{ \bullet \ar@<-.5ex>[r] \ar@<.5ex>[r] & \bullet}
	\qquad \qquad \qquad
	 \xymatrix{ X \ar@<-.5ex>[r]_g \ar@<.5ex>[r]^f & Y}
	\]
	The colimit, which is called the \textbf{coequalizer}, is an object $E$ and map $h$ such that $h \circ f=h\circ g$.
	\[
	 \xymatrix{ X \ar@<-.5ex>[r]_g \ar@<.5ex>[r]^f & Y \ar[r]^h & E}
	\]
	If $\dat=\Vect$ and one sets $g=0$, then the coequalizer is the \textbf{cokernel}. Thus if one wants to mimic cokernels in data types lacking of zero maps and objects, coequalizers can be substituted.
\end{ex}

\begin{ex}[Co-invariants]\index{coinvariants}
	As described in Example~\ref{ex:invariants}, a vector space $V$ with an endomorphism $T$ is equivalent to a functor $k[x]\to\Vect$. The colimit of this functor is called the \textbf{coinvariants} of $T$, i.e.
	\[
		C=V/<Tv-v>.
	\]
\end{ex}

\section{Adjunctions}
\label{sec:adjunctions}

Adjunctions allow us to derive interesting relationships with almost no effort; they are in essence dualities. For the individual interested in using category theory to model the world, facile manipulations of adjunctions is essential. One often can transform a complicated problem into a simpler one via an adjunction, thereby gaining a computational payoff at the cost of abstraction. This is why using adjunctions between the functors defined in Section~\ref{subsec:map_functors} is one of the key technical skills every sheaf theorist must master. Adjunctions also have played an essential role in the development of sheaf theory. Finding an adjoint to the functor $f_!$ was one of the primary reasons that the notion of a derived category was invented. Only by enlarging the domain could a new, adjoint functor $f^!$ be defined. Here we introduce the general theory.

\begin{defn}\index{adjoint}
 Suppose $F:\cat\to\dat$ and $G:\dat\to\cat$ are functors. We say that $(F,G)$ is an \textbf{adjoint pair} or that \textbf{$F$ is left adjoint to $G$} (or equivalently $G$ is right adjoint to $F$) if we have a natural transformation $\eta:\id_{\cat}\to G\circ F$ and a natural transformation to $\epsilon:F\circ G\to \id_{\dat}$ such that
\[
	\xymatrix{G\ar[r]^{\eta G} & GFG \ar[r]^{G\epsilon} & G},\qquad \xymatrix{F\ar[r]^{F\eta} & FGF \ar[r]^{\epsilon F}& F}
\]
We call $\eta$ the \textbf{unit} of the adjunction and $\epsilon$ the \textbf{counit} of the adjunction.
\end{defn}

There are about a half-dozen different, but equivalent, ways of defining an adjunction; see~\cite[p.81]{catwork} for a list. One can just specify $\eta$ and ask that it is universal,\footnote{In other words, initial in a particular comma category; see~\cite[p.56]{catwork}} i.e. for each $x\in\cat$ and for every $y\in\dat$ there is a map $\eta_x:x\to GF(x)$ such that if we have $f:x\to G(y)$, then there exists a unique map $f':F(x)\to y$ with $G(f')\circ \eta_x=f$. 
\[
\xymatrix{x \ar[r]^{\eta_x} \ar[dr]_{f} & GF(x) \ar@{.>}[d] \\ & G(y)}
\]
Of course we could have just defined $\epsilon$ and asked that it is universal in a dual sense.\footnote{It is final in a different comma category.} The point is this: an adjunction is equivalent to specifying for every $x\in\cat$ and $y\in\dat$ a natural bijection $\varphi_{x,y}$
\[
	\Hom_{\dat}(F(x),y)\cong\Hom_{\cat}(x,G(y)).
\]

The following theorem gives us an abstract criterion for determining when a functor has an adjoint.

\begin{thm}[Freyd's Adjoint Functor Theorem]\label{thm:freyd}\index{adjoint!functor theorem, Freyd's}
	Let $\dat$ be a complete category and $G:\dat\to\cat$ a functor, then $G$ has a left adjoint $F$ if and only if $G$ preserves all limits and satisfies the \textbf{solution set condition}. This condition states that for each object $x\in\cat$ there is a set $\Iat$ and an $\Iat$-indexed family of arrows $f_i:c\to G(a_i)$ such that every arrow $f:x\to G(a)$ can be factored as $x\to G(a_i)\to G(a)$, where the first map is $f_i:x\to G(a_i)$ and the second is $G$ applied to to some $t:a_i\to a$.
\end{thm}

The solution set condition holds nearly all the time, so in practice one only needs to check that $G$ preserves limits, in which case $G$ is a right adjoint (has a left adjoint). Dually, for a functor to be a left adjoint it needs to preserve colimits.

%
%

\chapter{The Theory of Sheaves and Cosheaves}
\label{sec:abstract_sheaves}

\begin{quote}
{\em ``Nous nous proposons d'indiquer sommairement comment les m\'ethodes par lesquelles nous avons etudi\'e la topologie d'un espace peuvent \^{e}tre adapt\'ees \`a l'\'etude de la topologie d'une repr\'esentation.''}\footnote{``We propose to state briefly how the methods by which we have studied the topology of a space can be adapted to the study of the topology of maps.''}
\begin{flushright} --- Jean Leray~\cite{leray46representation} \end{flushright}
\end{quote}

In its most general form, the subject of this thesis involves the assignment of data to subsets of a space $X$. This should sound like a very useful thing to do. After all, we have in both pure and applied mathematics many an occasion to record data or solutions in a local, spatially distributed way. Immediate questions arise: To which subsets should we assign data? What should these assignments be used for? What are they to be called?

The author believes such assignments are to be called sheaves or cosheaves depending on whether it is natural to restrict data from larger spaces to smaller spaces or by extending data from smaller spaces to larger ones. The evolution of these ideas deserves some discussion and the eager historian should consult John Gray's ``Fragments of the History of Sheaf Theory,''~\cite{gray} for a more thorough account. However, we outline three basic opinions on what a sheaf (or cosheaf) is really:
\begin{itemize}\index{sheaf!three perspectives on}\index{cosheaf!three perspectives on}
	\item A sheaf is a \textbf{system of coefficients} for computing cohomology that weighs and measures parts of the space differently. A cosheaf, in like manner, is a system of coefficients for homology that varies throughout the space.
	\item A sheaf is an \textbf{\'etal\'e space} $E$ along with a local homeomorphism $\pi:E\to X$. Analogously, a cosheaf is a locally-connected space $D$, called the \textbf{display locale}, that maps to $X$~\cite{funk-display}.
	\item A sheaf (or a cosheaf) is an \textbf{abstract assignment of data} --- a functor --- that further satisfies a gluing axiom expressed by limits (or colimits).
\end{itemize}

Historically, the system of coefficients perspective came first. In a 1943 paper Norman Steenrod defined a new homology theory determined by assigning abelian groups directly to \emph{points} of a space $X$ and group isomorphisms to (homotopy classes of) paths between points~\cite{steenrod}. This theory was vastly generalized in 1946 by Jean Leray where a \textbf{faisceau} (or sheaf) was defined to be a way of assigning modules to \emph{closed sets} in an inclusion-reversing way.\index{sheaf!on closed sets}
\[
	\xymatrix{V \ar[rr] \ar[rd] & & X \\ & W \ar[ru] &} \qquad \qquad \qquad \xymatrix{F(V) & & \ar[ll] \ar[ld] F(X) \\ & F(W) \ar[lu] &}
\]
Although this strengthened the abstract assignment perspective, Leray was still concerned with the cohomological ideas developed by Georges de Rham, Kurt Reidemeister and Hassler Whitney. 

By the early 1950s, Henri Cartan and his seminar revised Leray's definition of a sheaf to consist of a local homeomorphism $\pi:E\to X$. One could re-obtain the assignment perspective by attaching to each \emph{open set} $U$ the set of sections of this map over $U$: 
\[
	U \rightsquigarrow \{s:U\to E\, |\, \pi\circ s(x)=x\}
\]
One plausible explanation for using open sets is provided by the \textbf{open pasting lemma}, which states\footnote{Munkres calls this the ``local formulation of continuity'' in theorem 18.2(f)~\cite{munkres}. Munkres reserves the term ``pasting lemma'' for the closed set version, which is stated directly afterwards as theorem 18.3.} that if $X=\cup U_i$ is a (potentially infinite) union of open sets equipped with continuous sections $s_i:U_i\to E$ that agree on overlaps, then the set-theoretically defined section $s:X\to E$ will also be continuous. If closed sets are used, then this gluing argument only works for covers consisting of finitely many closed sets.

Finally, the Weil conjectures in algebraic geometry motivated the introduction of a more general notion of a topology and cohomology. Following suggestions of Jean-Pierre Serre, the domain of a sheaf was abstracted by Alexander Grothendieck from subsets $U\subseteq X$ to collections of mappings $U\to X$ that satisfy certain conditions reminiscent of an open cover~\cite{mm-topos}. Defining a sheaf on a Grothendieck topology ushered in the abstract formulation of sheaves using categories, functors and equalizers (limits) found in Michael Artin's 1962 Harvard notes on the subject~\cite{artin-gt}.

All three of these models are useful for thinking about sheaves and cosheaves, but the abstract assignment model is powerful and elegant enough to capture the other two. Moreover, whereas the \'etal\'e space perspective can be adapted from sheaves of sets to sheaves of more general data types, the display space perspective on cosheaves appears to only be valid for set-valued cosheaves and cannot be adapted more generally. In particular, since homology requires working with abelian groups or vector spaces, the display space model and the homology perspective describe different types of cosheaves. Thus, the only vantage point capable of reasoning about cosheaves in a unified way is the functorial perspective, where the dualities of category theory can be employed.

In this section, we provide the general definition of sheaves and cosheaves, but restrict ourselves to considering open sets and covers in a topological space. We phrase things using limits and colimits that take the shape of a simplicial complex: the nerve of a cover. The sheaf or cosheaf condition says that the value of this limit or colimit is independent of the cover chosen. To make the limits and colimits over covers more computable, we reduce to equalizers and coequalizers. We then specialize to the data type of vector spaces, where \v{C}ech homology for a cover is introduced. This evolves into a discussion of why singular zeroth homology defines a cosheaf. As set up for the discussion on general differences between sheaves and cosheaves, we consider how refinement of covers plays with the sheaf and cosheaf property.

\section{The General Definition}
\label{subsec:abstract_defn}
In elementary mathematics one learns that functions are devices for assigning points in one set to points in another. Motivated by differential calculus, one learns properties of functions on metric and topological spaces such as continuity. In its simplest form, continuity of a function states that if $f:X\to Y$ is a function and $\{x_n\}_{n=1}^{\infty}$ is a sequence of points in $X$ converging to some point $x$, then
\[
	\lim_{n\to\infty}f(x_n)=f(\lim_{n\to\infty} x_n)=f(x),
\] 
i.e. $f$ commutes with the limits one learns in analysis. Moreover, there is an independence result: The value $f(x)$ is independent of which sequence one used to approximate the point $x$.

The exact analogous situation occurs in category theory. A functor assigns objects and morphisms of one category to objects and morphisms in another. If a functor commutes with the categorical notion of a limit, then we also say that the functor is continuous. However, since there are so many different shapes of limits in arbitrary categories, this notion is too restrictive. A sheaf is a functor that commutes with limits coming from open covers. Applying the duality principle in category theory, a cosheaf is a functor that preserves colimits coming from open covers. 

\begin{defn}\index{cover}
Let $X$ be a topological space and $U$ an open set in $X$. An \textbf{open cover} of $U$ is a collection of open sets $\cU:=\{U_i\}_{i\in\Lambda}$ whose union is $U$.
\end{defn}

Pavel Alexandrov introduced in 1928 a method\footnote{In Definition~\ref{defn:MV_blowup} we consider the ``correct'' generalization of the nerve.} for associating to every open cover an abstract simplicial complex~\cite{nerve}. We will use these shapes to model our limits and colimits of interest. 
\begin{defn}\label{defn:nerve}\index{nerve}
	Suppose $\covU:=\{U_i\}_{i\in \Lambda}$ is an open cover of $U$. We can take the \textbf{nerve} of the cover to get an abstract simplicial complex $N(\covU)$, whose elements are subsets $I=\{i_0,\ldots,i_n\}$ for which $U_I:=U_{i_0}\cap\cdots\cap U_{i_n}\neq\emptyset$. We can regard $N(\covU)$ as a category whose objects are the finite subsets $I$ such that $U_I\neq\emptyset$ with a unique arrow from $I\to J$ if $J\subseteq I$. Since our intersections are only finite, and the finite intersection of open sets is open, we get natural functors
	\[
	\iota_{\covU}:N(\covU)\to\Open(X) \qquad \mathrm{or} \qquad \iota^{op}_{\covU}:N(\covU)^{op}\to \Open(X)^{op}.
	\]
\end{defn}

\begin{rmk}
Sometimes we will use the notation $N(\cU)$, $N_{\cU}$ and $N$ interchangeably, depending on the context.
\end{rmk}

\begin{figure}[ht]
\centering
\includegraphics[width=\textwidth]{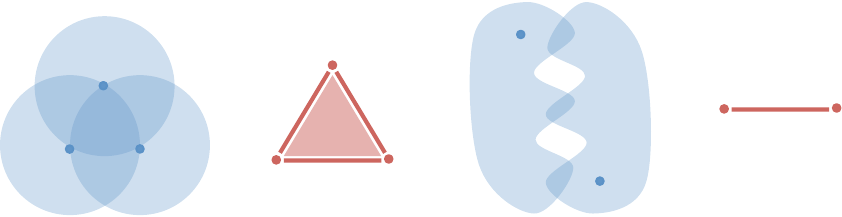}
\caption{Covers and Their Nerves}
\label{fig:nerves}
\end{figure}

In Figure~\ref{fig:nerves} we have drawn two different arrangements of open sets and their corresponding nerves, which we have represented graphically to the right. We have added points to each open set to make it clear how many open sets are in the cover. Note that in general, there is nothing to prevent a disconnected open set from being marked by a single label. 

The nerve is purely an algebraic and combinatorial model for the cover --- it need not respect the topology of the union. However, the \textbf{nerve theorem}\index{nerve theorem} of Leray and Borsuk~\cite{leray45,borsuk-nerve} states that if the intersections are contractible then the nerve and the union have the same homotopy type. The example on the left in Figure~\ref{fig:nerves} gives a positive example of the nerve lemma, whereas the example on the right gives a negative one.

The definition of a sheaf or cosheaf requires the synthesis of covers and data. We now introduce the functor that assigns data to open sets.
\begin{defn}[Pre-Sheaf and Pre-Cosheaf]\index{presheaf}\index{pre-cosheaf}
	A \textbf{pre-sheaf} is a functor $F:\Open(X)^{op}\to\dat$ and a \textbf{pre-cosheaf} is a functor $\hF:\Open(X)\to\dat$. If $V\subset U$, then we usually write the \textbf{restriction map} as $\rho_{V,U}^F:F(U)\to F(V)$ and the \textbf{extension map} as $r^{\hF}_{U,V}:\hF(V)\to\hF(U)$. Often we omit the superscript $F$ or $\hF$.
\end{defn}

If one imagines the pre-cosheaf that associates a copy of the field $k$ to every connected component of an open set, then the following diagrams of vector spaces emerge from Figure~\ref{fig:nerves}:
\[
\xymatrix@R=1em{ & k & \\ k \ar[ru] \ar[dd] & & k \ar[lu] \ar[dd]\\ & \ar[lu] \ar[uu] \ar[ru] k \ar[dl] \ar[dd] \ar[dr] & \\ k & & k \\ & \ar[lu] \ar[ru] k &}
\qquad \qquad \qquad
\xymatrix@R=1em{ & & \\ & & \\ k & k^3 \ar[l] \ar[r] & k \\ & & \\ & &}
\]
We will examine various ways for computing the colimits of these diagrams explicitly. Since the colimits occur over simplicial complexes, we introduce a structure theorem that allows us to use coequalizers. In the vector space case, this reduces to linear algebra --- the colimit will be $H_0$ of a suitable chain complex.

We want to express the fact that since the colimit of a cover $N(\covU)\to \Open(X)$ is just the union $U=\cup U_i$, the data associated to $U$ should be expressible as the colimit of data assigned to the nerve. Moreover, this should be independent of which cover we take. Examples where this does not occur are given in Example~\ref{ex:nonlocalpresheaf} and Example~\ref{ex:inconsistentpresheaf}.
\begin{defn}[Sheaves and Cosheaves]\label{defn:sheaf}\index{sheaf!definition of}\index{cosheaf!definition of}
	Suppose $F$ is a pre-sheaf and $\hF$ is a pre-cosheaf, both of which are valued in $\dat$. Suppose $\covU=\{U_i\}$ is an open cover of $U$. We say that $F$ is a \textbf{sheaf on $\covU$} if the unique map from $F(U)$ to the limit of $F\circ \iota^{op}_{\covU}$, written 
	\[
	F(U)\to \varprojlim_{I\in N(\covU)} F(U_I)=:F[\covU],
	\]
	is an isomorphism. Similarly, we say $\hF$ is a \textbf{cosheaf on $\covU$} if the unique map from the colimit of $\hF\circ\iota_{\covU}$ to $\hF(U)$, written
	\[
	\hF[\covU]:=\varinjlim_{I\in N(\covU)} \hF(U_I) \to \hF(U),
	\]
	is an isomorphism. We say that $F$ is a \textbf{sheaf} or $\hF$ is a \textbf{cosheaf} if for every open set $U$ and every open cover $\covU$ of $U$, $F(U)\to F[\covU]$ or $\hF[\covU]\to \hF(U)$ is an isomorphism. For a catchy slogan, we say 
\begin{quote}
	\begin{center}
	\emph{On an open set (co)sheaves turn different covers into isomorphic (co)limits}.
\end{center}
\end{quote}
\end{defn}
\begin{rmk}[Stable Under Finite Intersection]
	Most authors do not introduce the nerve as any part of the definition of a sheaf or cosheaf. Instead, some will require that the cover $\covU$ is ``stable under finite intersection,'' i.e. if $U_i,U_j\in\covU$, then $U_i\cap U_j\in\covU$. This allows those authors to just consider the limit or colimit over the cover and not over some auxiliary construction, like we have done. This works because one can take any cover and then add the intersections after the fact, but this tends to be done unconsciously and without any warning to the reader. Our approach is equivalent to that approach, but we believe it has some added benefits.
\end{rmk}

We have not stated any requirements on the data category $\dat$, but in order to even parse the statement of the (co)sheaf axiom we require that the (co)limits coming from such covers exist. For the most part, we will work in categories where all limits and colimits exist. In analogy with analysis, a category where the limit of any diagram $F:I\to\dat$ exists is called \textbf{complete}. Similarly, if the colimit of an arbitrary diagram exists, we say $\dat$ is \textbf{co-complete}. The category $\Vect$ is both complete and co-complete.
	
A particular consequence of the axiom is that for a sheaf, $F(\emptyset)$ must be the limit over covers of the empty set, but since there are no such covers\footnote{Alternatively one argues that the empty set covers itself and hence the value there is chosen to be the initial/terminal object of the category $\dat$.} this is the limit over the empty diagram, i.e. $\Cone(\emptyset)=\dat$, whose terminal object is the terminal object of $\dat$. Similarly, for a cosheaf $\hF(\emptyset)$ must be the initial object in $\dat$. For $\dat=\Vect$ the initial and terminal objects coincide with the zero vector space.
	
It is true that if $\dat$ has pullbacks (see Example~\ref{ex:pullback} in Section~\ref{sec:categories}) and a terminal object then it has all finite limits. The dual statement that having an initial object and pushouts (see Example~\ref{ex:pushout}) implies finitely co-complete is also true. Thus, if one focuses on sheaves and cosheaves valued in $\vect$ (the category of finite dimensional vector spaces and linear maps), then the collection of covers of $U$ one can consider must be restricted. In particular, if the sheaf or cosheaf axiom holds for open covers with two sets, then we can only guarantee that it holds for covers with finitely many open sets. As a purely philosophical point, one wonders whether working with the cover of the complement of the Cantor set given by
\[
	\covU=\{(\frac{3k+1}{3^n},\frac{3k+2}{3^n})\subset[0,1] | 0\leq k\leq 3^{n-1}-1, \, 0\leq n<\infty\}
\]
would ever be computationally tractable. One might wish to systematically revise the notion of a ``cover,'' and this would lead to the notion of a \textbf{Grothendieck site}, which we do not address here.

We now examine the axioms just for covers with only two open sets.

\begin{figure}[ht]
\centering
\includegraphics[width=\textwidth]{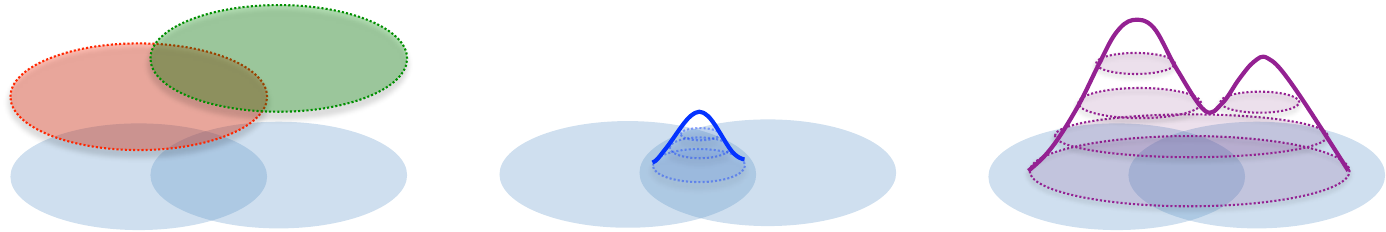}
\caption{Sheaves and Cosheaves of Functions}
\label{fig:fun_sheaves}
\end{figure}

\begin{ex}[Cover by Two Sets]\index{sheaf!axiom for two opens}\index{cosheaf!axiom for two opens}\index{cosheaf!of compactly supported functions}
Suppose $\dat=\Set$, and suppose $\covU=\{U_1,U_2\}$ is a cover of $U$. The sheaf condition says that
\[
	F(U)\cong\{(s_1,s_2)\in F(U_1)\prod F(U_2) | \rho_{U_{12},U_1}(s_1)=\rho_{U_{12},U_2}(s_2)\}=:F[\covU],
\]
i.e. $F(U)$ lists the set of consistent choices of elements from $F(U_1)$ and $F(U_2)$. In particular, $F[\covU]$ is a sub-object of the product of $F(U_1)$ and $F(U_2)$. For an example, one can let $F$ be the assignment
\[
	U \rightsquigarrow \{f:U\to\RR \,|\, \mathrm{continuous} \}.
\]
The sheaf axiom then says in order for two functions (or sections) $s_1=f_1:U_1\to\RR$ and $s_2=f_2:U_2\to\RR$ to determine an element in $U=U_1\cup U_2$ it is necessary and sufficient that the functions $f_1(x)$ and $f_2(x)$ agree on the overlap $U_{12}=U_1\cap U_2$.

The cosheaf condition for $\dat=\Set$ is slightly strange. It says that
\[
	\hF(U)\cong (\hF(U_1)\coprod \hF(U_2))/\sim
\]
\[
	s_1\sim s_2 \, \Leftrightarrow \, \exists s_{12} \quad s_1=r_{U_1,U_{12}}(s_{12}) \quad  s_2=r_{U_2,U_{12}}(s_{12}).
\]
	In contrast to the sheaf case, the notion of consistent choices no longer applies for cosheaves, because it requires thinking in terms of quotient objects --- something human beings are not accustomed to. However, a useful analogy is that one must subtract out or identify those elements that might be counted twice because they come from the intersection. For an example similar in spirit to the sheaf of real-valued functions, we begin by considering the pre-cosheaf of compactly supported functions gotten by the assignment
	\[
		U \rightsquigarrow \{f:U\to\RR \,|\, \mathrm{continuous}\,\,\mathrm{and}\,\,\mathrm{compactly}\,\,\mathrm{supported} \}.
	\]
	Extending by zero provides the extension map and identifying the two copies of a function whose support is contained in $U_{12}=U_1\cap U_2$ prevents double counting on $U$. However, this is not all that the cosheaf axiom requires. Any compactly supported function should appear as one supported in $U_1$ or $U_2$, but this is not always true. Some compactly supported functions are not compact when restricted to any particular open set in a cover. Thus, this pre-cosheaf is \emph{not} a cosheaf.
	
	The reader familiar with partitions of unity will realize that if $X$ a paracompact Hausdorff space then we can express any compactly supported function $f(x)$ defined on all of $U$ as a \emph{sum} of compactly supported functions on $U_1$ and $U_2$. By taking a partition of unity subordinate to the cover $\covU$ we get two functions $\lambda_1(x)$ and $\lambda_2(x)$ such that
	\[
		f(x)=f_1(x)+f_2(x) \quad \mathrm{where} \quad f_1(x):=\lambda_1(x)f(x) \quad \mathrm{and} \quad f_2(x):=\lambda_2(x)f(x).
	\]
	By carrying out the colimit in a data category equipped with sums, such as $\dat=\Vect$ of $\Ab$, then compactly supported functions do define a cosheaf valued there.

	More generally, if $\dat=\Vect$, then the cosheaf axiom for the cover says the sequence
	\[
	\hF(U_{12}) \to \hF(U_1)\oplus\hF(U_2) \to \hF(U) \to 0
	\]
	is exact, where the maps are $(-r_{U_1,U_{12}},r_{U_2,U_{12}})$ and $r_{U,U_2}+r_{U,U_1}$. Dually, the sheaf axiom says the dual sequence
	\[
	0 \to F(U)\to F(U_1)\times F(U_2) \to F(U_{12})
	\]
	is exact, where the second map is $\rho_{U_{12},U_2}+\rho_{U_{12},U_1}$ and the first map is $(-\rho_{U_1,U},\rho_{U_2,U})$.
\end{ex}

\section{Limits and Colimits over Covers: a Structure Theorem}
\label{subsec:cover_structure}\index{limit!over a cover}\index{colimit!over a cover}

The sheaf and cosheaf axioms as stated are meant to emphasize that if one is comfortable with the operations of limits and colimits, then one is already comfortable with sheaves and cosheaves. However, the limits and colimits considered in Definition~\ref{defn:sheaf} have a special structure. This structure comes from the fact that the indexing category --- the nerve --- is a simplicial complex.

The first observation one can make is that for any functor $F:N(\covU)^{op}\to\dat$ the limit can be thought of as ``sitting inside'' the product over the vertices --- the vertices corresponding to the elements of the cover through the nerve construction. Dually, the colimit of a functor $\hF:N(\covU)\to\dat$ can be thought of as a quotient of the coproduct of the functor over the vertices. Said using formulas, this is
\begin{equation*}
	\varprojlim F \hookrightarrow \prod F(i) \qquad \coprod \hF(i) \twoheadrightarrow \varinjlim \hF.
\end{equation*}
The way to see this is to note that any cone or cocone's morphism must factor through a vertex. However, the difference between the limit or colimit from the functor's aggregate value on vertices is measured by edges in the nerve. This is a reflection of a more general theorem, which we now state.

\begin{thm}\label{thm:limits_equal}\index{limit!as an equalizer}\index{colimit!as a coequalizer}
	A category $\dat$ has all (co)limits of an appropriate size if it has all (co)products and (co)equalizers of same such size. Here ``size'' corresponds to the cardinality of the indexing category of the (co)limit in question.
\end{thm}
\begin{proof}[Proof (Sketch)]
One should consult~\cite[Prop. 5.22-3]{awodey} for a complete proof. To give the reader the idea, one can compute the limit of $F:I\to\dat$ by taking the product over all the objects $x\in I$ and separately the product over all morphisms in the indexing category $I$. The limit is isomorphic to the equalizer going from the first product to the latter, i.e.
\begin{equation*}
 \xymatrix{ \varprojlim F \ar[r] & \prod_{x\in I} F(x) \ar@<-.5ex>[r] \ar@<.5ex>[r] & \prod_{x\to x'} F(x')}.
\end{equation*}
By dualizing, one can prove the analogous result for colimits. 
\end{proof}

This theorem gives us effective means for computing limits and colimits for general data types. We now specialize this result to the limits and colimits pertinent to sheaves and cosheaves.

\subsection{Rephrased as Equalizers or Co-equalizers}
\label{subsubsec:equalizers}
\index{sheaf!axiom as an equalizer}\index{cosheaf!axiom as a coequalizer}

The method outlined in Theorem~\ref{thm:limits_equal} for computing limits and colimits contains too much redundant information for the case $I=N(\covU)^{op}$. As such, we state the precise, simplified formulation here. The sheaf and cosheaf axioms can be rephrased as saying that the following sequences
\begin{equation*}
	\xymatrix{F(U) \ar[r]^-e & \prod F(U_i) \ar@<-.5ex>[r]_-{f^+} \ar@<.5ex>[r]^-{f^-} & \prod_{i<j} F(U_i\cap U_j)}
\end{equation*} 
\begin{equation*}
\xymatrix{\coprod_{i<j}\hF(U_i\cap U_j) \ar@<-.5ex>[r]_-{g^-} \ar@<.5ex>[r]^-{g^+} & \coprod \hF(U_i) \ar[r]^-u & \hF(U)}
\end{equation*}
are an equalizer and a co-equalizer respectively. 

\begin{exr}
Prove that the limit or colimit over the nerve of a cover can be determined after considering only the elements of the cover and their pairwise intersections. Do this by observing that the limit or colimit over the 1-skeleton\footnote{In higher homotopy analogs of sheaves and cosheaves one works over the whole \v{C}ech tower of a cover~\cite{dugger-hypercovers,douglas2007sheaves,gwilliam-stacks}.} of the nerve defines a cone or cocone over the whole nerve and employing universal properties. Then apply the equation from Theorem~\ref{thm:limits_equal} and its dual version to prove that the re-written axioms of Section~\ref{subsubsec:equalizers} and Section~\ref{subsubsec:exactness} are correct.
\end{exr}

To describe the maps explicitly requires some work. First, we choose an ordering of the indexing set of the cover $\covU=\{U_i\}_{i\in\Lambda}$. To specify a map to a product it suffices to specify maps to each factor of the product. Similarly, maps from a coproduct are specified by maps from each factor. This is summarized by the identities 
	\[
	\Hom(X,\prod_i Y_i)\cong \prod_i\Hom(X,Y_i) \qquad \mathrm{and} \qquad \Hom(\coprod_i X_i,Y)\cong\prod_i \Hom_i(X_i,Y).
	\]
To define the maps $e$ and $u$ we declare $e_i:=\rho_{U_i,U}$ and $u_i:=r_{U,U_i}$. For the maps $f^{\pm}$ and $g^{\pm}$ we define for each pair $i<j$ the maps
	\[
	f^+_{ij}:=\rho_{ij,j}\circ\pi_j \quad f^-_{ij}:=\rho_{ij,i}\circ\pi_i \qquad g^+_{ij}:=r_{j,ij}\circ\iota_{ij} \quad g^-_{ij}:=r_{i,ij}\circ\iota_{ij}
	\]
where $\pi_i:\prod F(U_i)\to F(U_i)$ is the natural projection and $\iota_{ij}:\hF(U_{ij})\to \coprod \hF(U_{ij})$ is the natural inclusion.
	
The reader might find it helpful to think of the maps in between the products as being represented by matrices. In the case of a cover with three elements $\covU=\{U_1,U_2,U_3\}$ all of whose pairwise intersections are non-empty, we can write
	\[
	f^+=\begin{bmatrix} \ast & \rho_{12,2} & \ast  \\  \ast & \ast & \rho_{13,3} \\ \ast & \ast & \rho_{23,3}\end{bmatrix} \quad f^-=\begin{bmatrix} \rho_{12,1} &\ast  &\ast  \\ \rho_{13,1} & \ast  & \ast \\ \ast & \rho_{23,2} &\ast \end{bmatrix}.
	\]
	The equalizer condition now reads that $f^+(s_1,s_2,s_3)=f^-(s_1,s_2,s_3)$, i.e. 
	\[
	(\rho_{12,2}(s_2),\rho_{13,3}(s_3),\rho_{23,3}(s_3))=(\rho_{12,1}(s_1),\rho_{13,1}(s_1),\rho_{23,2}(s_2)).
	\]
	
\subsection{Rephrased as Exactness}\label{subsubsec:exactness}
\index{sheaf!axiom as an exact sequence}\index{cosheaf!axiom as an exact sequence}

If $\dat=\Vect$, then we can add and subtract maps and look for kernels and cokernels instead of equalizers and co-equalizers. The sheaf and cosheaf axioms then reduce to linear algebra. The modified axioms now read as
\[
	\xymatrix{0 \ar[r] & F(U) \ar[r] & \prod F(U_i) \ar[r]^-{d^0} & \prod_{i<j} F(U_i\cap U_j)}
\]
\[
	 \xymatrix{\bigoplus_{i<j}\hF(U_i\cap U_j) \ar[r]^-{\partial_1} & \bigoplus \hF(U_i) \ar[r] & \hF(U) \ar[r] & 0}
\]
where $d^0$ is the matrix whose rows are parametrized by pairs $i<j$ and whose columns are parametrized by $k$ with entries given by $d^0_{ij,k}=[k:ij]\rho_{ij,k}$ where 
	\[
		[k:ij]=\left\{\begin{array}{ll} 0 & \mathrm{if} \,\, k\neq i\neq j \\
		1 & \mathrm{if} \,\, k=j\\
		-1 & \mathrm{if} \,\, k=i\end{array}\right.
	\]
The matrix $\partial_1$ is similarly defined except that the rows are indexed by $k$ and columns are indexed by pairs $i<j$ with entries $(\partial_1)_{k,ij}=[k:ij]r_{k,ij}$. Thus the sheaf axiom says that $F(U)\cong \ker(d^0)$ and the cosheaf axiom says that $\hF(U)\cong \coker (\partial_1)$.
	
In our example of a three set cover $\covU=\{U_1,U_2,U_3\}$ all of whose pairwise intersections are non-empty, the definition of $d^0$ corresponds to taking $f^+ - f^-$, i.e.
	\[
		d^0=f^+ - f^- =\begin{bmatrix} -\rho_{12,1} & \rho_{12,2} & 0 \\ -\rho_{13,1} & 0 & \rho_{13,3} \\ 0 & -\rho_{23,2} & \rho_{23,3}\end{bmatrix}
	\]
where each of the $\rho_{ij,k}$'s need to be filled in with some matrix representative of that linear map. The kernel is then identified with $F[\covU]$.

\section{\v{C}ech Homology and Cosheaves}\label{subsec:cech}
\index{Cech@\v{C}ech homology}

In Section~\ref{subsubsec:equalizers} we rephrased the limits and colimits coming from covers as equalizers and coequalizers. For the data category $\dat=\Vect$ we showed how to reinterpret this as an exact sequence. This perspective is indicative of a deeper and more computational idea, namely that of homology. We now show how to associate to any pre-cosheaf\footnote{Or pre-sheaf, but we'll leave it to the reader to dualize.} of vector spaces $\hF$ and an open cover $\covU=\{U_i\}_{i\in\Lambda}$ a complex of vector spaces whose zeroth homology computes $\hF[\covU]$. This allows us to compute the homology of data.

\begin{defn}[\v{C}ech Homology]
	Given a pre-cosheaf of vector spaces $\hF$ and an open cover $\covU=\{U_i\}_{i\in\Lambda}$, we define the \textbf{\v{C}ech homology on $\covU$} to be the homology of the complex
\[
	(\check{C}_{\bullet}(\covU;\hF),\partial_{\bullet}) \qquad \mathrm{where} \qquad \check{C}_p(\covU;\hF):=\bigoplus_{|I|=p+1} \hF(U_I) \qquad \mathrm{for} \qquad I\in N(\covU).
\]
By choosing an ordering on the index set $\Lambda$, we define the differential by extending the formula defined on elements $s_I\in\hF(U_I)$ by linearity, i.e.
\[
	\partial_p:C_p(\covU;\hF)\to C_{p-1}(\covU;\hF) \qquad \partial_p(s_I):=\sum_{k=0}^p (-1)^k r_{U^{(k)}_I,U_I}(s_I),
\]
where the symbol $U^{(k)}_I=U_{i_0}\cap \ldots \cap U_{i_{k-1}}\cap U_{i_{k+1}} \cap \ldots  U_{i_p}$ indicates the intersection that omits the $k$th open set.
Thus we can define by the usual formula the $p$th \textbf{\v{C}ech homology group}
\[
	\check{H}_p(\covU;\hF):=\frac{\ker \partial_p}{\im \partial_{p+1}} \qquad \mathrm{i.e.} \qquad H_p(\check{C}_{\bullet}(\covU;\hF)).
\]
\end{defn}

To guarantee that \v{C}ech homology is well-defined we verify the following lemma:
\begin{lem}\index{Cech@\v{C}ech homology!$\partial^2=0$}
	The differential $\partial$ in the \v{C}ech complex for a cover $\covU$ and a pre-cosheaf $\hF$ of vector spaces satisfies $\partial_{p}\circ\partial_{p+1}=0$.
\end{lem}
\begin{proof}
The combinatorial nature of the nerve of a cover guarantees that $\partial^2=0$. Specifically, there are two ways of going between incident simplices of dimension differing by two. Thus, we get the following diagram of open sets and data:
\[
	\xymatrix{ & U^{(j,k)}_I & \\
	U^{(j)}_I \ar[ur] & & U^{(k)}_I \ar[ul] \\
	& U_I \ar[ul] \ar[ur] \ar[uu] &}
	\qquad \qquad
	\xymatrix{ & \hF(U^{(j,k)}_I) & \\
	\hF(U^{(j)}_I) \ar[ur] & & \hF(U^{(k)}_I) \ar[ul] \\
	& \hF(U_I) \ar[ul] \ar[ur] \ar[uu] &}
\]

Let's follow a typical element $s_I\in\hF(U_I)$ through the diagram on the right upon applying the formula $\partial\circ\partial$. First note that the fact that $\hF$ is a pre-cosheaf implies that the square commutes, i.e. 
\[
r_{U^{(j,k)}_I,U^{(j)}_I}\circ r_{U^{(j)}_I,U_I}(s_I)=r_{U^{(j,k)}_I,U^{(k)}_I}\circ r_{U^{(k)}_I,U_I}(s_I)=r_{U^{(j,k)}_I,U_I}(s_I).
\] 
The first application of $\partial$ yields $(-1)^j r_{U^{(j)}_I,U_I}(s_I)$ and $(-1)^k r_{U^{(k)}_I,U_I}(s_I)$ as just two components in the formula for $\partial(s_I)$. Assuming $j<k$ and applying the definition of the boundary map to elements in $\hF(U^{(j)}_I)$ implies that we must actually delete the $k-1$st entry of $I-\{j\}$ since removing $j$ has caused everything above $j$ to shift down in the ordered list. Thus the image of $\partial^2(s_I)$ in $\hF(U^{(j,k)}_I)$ is
\[
	(-1)^{k-1}(-1)^jr_{U^{(j,k)}_I,U_I}(s_I) + (-1)^{k}(-1)^jr_{U^{(j,k)}_I,U_I}(s_I)=0.
\]
\end{proof}

\begin{ex}\index{Cech@\v{C}ech homology!example of}
	Consider the covers in Figure~\ref{fig:nerves}. The pre-cosheaf we described there assigned to each connected component of an open set a copy of the field $k$. First we consider the cover on the left of Figure~\ref{fig:nerves} with three open sets. We label the three vertices of the nerve, starting with the bottom left one and working counter-clockwise, $x,y$ and $z$ respectively. The \v{C}ech complex takes the form
	\[
		\xymatrix{k_{xyz} \ar[r]^-{\partial_2} & k_{xy}\oplus k_{xz}\oplus k_{yz} \ar[r]^{\partial_1} & k_x\oplus k_y \oplus k_z \ar[r] & 0 }
	\]
	where, using the lexicographic ordering for a basis, the matrix representatives for $\partial_2$ and $\partial_1$ take the following form:
	\[
		\partial_2=\begin{bmatrix} 1 \\ -1 \\ 1\end{bmatrix} \qquad \qquad  \partial_1=\begin{bmatrix} -1 & -1 & 0 \\ 1 & 0 & -1 \\ 0 & 1 & 1\end{bmatrix}
	\]
	One can easily verify that the $\ker \partial_1=\im \partial_2$ and consequently $\check{H}_1=0$. Furthermore, $\check{H}_0\cong k$, which happens to reflect that the union has one connected component.
	Similarly, one can consider the cover at the right of Figure~\ref{fig:nerves}. The \v{C}ech complex for this cover and the same pre-cosheaf is as follows:
	\[
		\xymatrix{k^3 \ar[r]^{\partial_1} & k^2 \ar[r] & 0 } \qquad \mathrm{where} \qquad \partial_1=\begin{bmatrix} -1 & -1 & -1 \\ 1& 1 & 1\end{bmatrix}
	\]
	Clearly $\check{H}_0\cong k$, whose dimension agrees with the number of connected components of the union, but also $\check{H}_1\cong k^2$, which witnesses the presence of two holes in the union.
\end{ex}

One can dually define \v{C}ech cohomology with coefficients valued in a pre-sheaf $F$. The discussion of Section~\ref{subsubsec:exactness}, along with the examples just presented, can be interpreted as saying a pre-sheaf $F$ or pre-cosheaf $\hF$ is a sheaf or cosheaf if and only if
	\[
	F(U)\cong \check{H}^0(\covU;F) \qquad \mathrm{or} \qquad \check{H}_0(\covU;\hF)\cong \hF(U).
	\]
for any choice of cover $\covU$ of $U$.
	
We would like to use this isomorphism to supply examples of sheaves and cosheaves from standard machinery in algebraic topology. Suppose one has an independent notion of homology, such as singular homology, and one can show it is isomorphic to \v{C}ech homology for suitably fine covers (see Section~\ref{subsec:refinement} to see why fineness matters) on nice spaces, then one can also define a cosheaf using those values. To make this rigorous, and to also provide a useful criterion for proving when a pre-cosheaf is a cosheaf, we recall a theorem:

\begin{thm}\label{thm:MV_cosheaf}\index{cosheaf!Mayer-Vietoris}\index{Mayer-Vietoris}
	Suppose $\hF$ is a pre-cosheaf, then $\hF$ is a cosheaf if and only if the following the following two properties hold
	\begin{itemize}
		\item For all open sets $U$ and $V$ the following sequence is exact
		\[
			\hF(U\cap V) \to \hF(U)\oplus \hF(V) \to \hF(U\cup V) \to 0.
		\]
		The first morphism is $(-r_{U,U\cap V},r_{V,U\cap V})$ and the second is $r_{U\cup V,U}+r_{U\cup V,V}$.
		\item If  $\{U_{\alpha}\}$ is directed upwards by inclusion, i.e. for for every pair $U_{\alpha}$ and $U_{\beta}$ there exists $U_{\gamma}$ containing both, then the canonical map
		\[
			\varinjlim_{\alpha} \hF(U_{\alpha}) \to \hF(\cup U_{\alpha})
		\]
		is an isomorphism. 
	\end{itemize}
	Dually, turning arrows around and using inverse limits gives a useful criterion for determining when a pre-sheaf is a sheaf.
\end{thm}
\begin{proof}
	Using induction one can prove that the cosheaf property for two sets implies the cosheaf property for finitely many sets (see~\cite[p. 418]{Bredon} for a proof). We now show that this implies the cosheaf axiom for arbitrary covers. Suppose $\{U_{\alpha}\}_{\alpha\in\Lambda}$ is a cover indexed by a potentially large, but ordered set $\Lambda$. For each finite subset $I\subset \Lambda$ we know that
	\[
	\bigoplus_{\alpha<\beta\in I} \hF(U_{\alpha,\beta}) \to \bigoplus_{\alpha\in I} \hF(U_{\alpha}) \to \hF(\bigcup_{\alpha\in I} U_{\alpha}) \to 0
	\]
	is exact. We know that the collection of finite subsets $I$ forms a directed system and that in $\Vect$ direct limits preserve exactness. As such we have that
	\[
	\varinjlim_I \bigoplus_{\alpha<\beta\in I} \hF(U_{\alpha,\beta}) \to \varinjlim_I \bigoplus_{\alpha\in I} \hF(U_{\alpha}) \to \varinjlim_I \hF(\bigcup_{\alpha\in I} U_{\alpha}) \to 0
	\]
	is exact as well, but by using the second property and the fact that the direct limit of the $I$'s is $\Lambda$ we have
	\[
		\bigoplus_{\alpha<\beta\in \Lambda} \hF(U_{\alpha,\beta}) \to \bigoplus_{\alpha\in \Lambda} \hF(U_{\alpha}) \to \hF(\bigcup_{\alpha\in \Lambda} U_{\alpha}) \to 0
	\]
	is exact. This proves the reverse direction. The other direction is clear.
\end{proof}

This theorem then provides us with a useful example of a cosheaf that we have implicitly used to generate examples. We now make this example explicit.

\begin{ex}\index{cosheaf!given by $H_0$}
	The assignment to an open set $U$ the $0$th singular homology of $U$
	\[
		U \rightsquigarrow H_0(U;k)
	\]
	is a cosheaf. This follows from the fact that the singular chain complex (see later for a definition) $C_{\bullet}(-;k)$ can be defined for any subset $U$ of $X$ and homology commutes with direct limits, thus the second property of the theorem holds. The first property in the theorem follows from exactness at the last two spots in the Mayer-Vietoris sequence:
	\[
	\xymatrix{
		            & H_1(U\cap V;k)\ar@{->}[r] & H_1(U;k)\oplus H_1(V;k) \ar@{->}[r]
		                   & H_1(U\cup V;k)  \ar `r[d] `_l[lll] ^{} `^d[dlll] `^r[dll] [dll] \\
		             & H_{0}(U\cap V;F)\ar@{->}[r] & H_0(U;k) \oplus H_0(V;k) \ar@{->}[r]
		                   & H_0(U\cup V;k) \ar[r] & 0 \\
		}
	\]
\end{ex}

The moral from this example is that, in essence,
\begin{quote}
	\begin{center}
	\emph{Any functor that satisfies Mayer-Vietoris is a cosheaf.}
\end{center}
\end{quote}

\section{Refinement of Covers}\label{subsec:refinement}
\index{cover!refinement of}

We have defined the sheaf and cosheaf axioms for a cover $\covU$. The coarsest possible cover of an open set $U$ is the cover with one element $\{U\}$. Thus, one way of interpreting the sheaf and cosheaf axiom is that $F[\covU]$ and $\hF[\covU]$ are independent of the cover chosen. A logical question to ask is if the axiom holds for some cover, but not all, then for what other covers does the axiom hold? To answer this question, we review some relevant concepts.

\begin{defn}[Refinement of Covers]
	Suppose $\covU$ and $\covU'$ are covers of $U$, then we say that $\covU'$ \textbf{refines} $\covU$ if for every $U_i'\in\covU'$ there is a $U_j\in\covU$ and an inclusion $U_i'\to U_j$. Note that every cover refines the trivial cover $\{U\}$.
\end{defn}

\begin{defn}\index{category!of covers}
The refinement relation endows the collection of covers of $U$ with the structure of a category $\Cov(U)$, whose objects are covers $\covU$ with a unique morphism $\covU'\to\covU$ if the former refines the latter.
\end{defn}

Note that if $U_{i_1}'\to U_{j_1}$ and $U'_{i_2}\to U_{j_2}$, then $U'_{i_1}\cap U_{i_2}'\to U_{j_1}\cap U_{j_2}$. So a refinement induces a functor between nerves, but it depends on which inclusions were chosen.
\[
\xymatrix{ & \Open(X) & \\ N(\covU') \ar[ru] \ar[rr] & & N(\covU) \ar[lu]} \qquad \xymatrix{ & \Open(X)^{op} & \\ N(\covU')^{op} \ar[ru] & & N(\covU)^{op} \ar[lu] \ar[ll]}
\]

The next lemma shows that these choices don't matter on the level of limits and colimits for pre-sheaves and pre-cosheaves.

\begin{lem}\index{presheaf!behavior under refinement}\index{pre-cosheaf!behavior under refinement}
	Let $\hF$ and $F$ be a pre-cosheaf and a pre-sheaf respectively. Suppose $\covU'$ refines another cover $\covU$ of an open set $U$. Then there are well-defined maps
	\[
	\hF[\covU']\to \hF[\covU] \qquad \mathrm{and} \qquad F[\covU]\to F[\covU'],
	\]
	i.e. we get functors $\hF:\Cov(U)\to\dat$ and $F:\Cov(U)^{op}\to\dat$.
\end{lem}
\begin{proof}
We'll detail the proof for a pre-cosheaf $\hF$ since the case for pre-sheaves can be found in the literature or obtained here via dualizing appropriately. A refinement $\covU'\to\covU$ defines a natural transformation $\hF\circ \iota_{\covU'}\Rightarrow \hF\circ \iota_{\covU}$. The colimit defines a natural transformation from $\hF\circ \iota_{\covU}$ to the constant diagram $\hF[\covU]$. Since the composition of natural transformations is a natural transformation, this induces a cocone $\hF\circ\iota_{\covU'}\Rightarrow \hF[\covU]$ which, by the universal property of the colimit, defines a unique induced map there, i.e.
\[
\hF\circ \iota_{\covU'}\Rightarrow \hF\circ \iota_{\covU}\Rightarrow \hF[\covU] \quad \mathrm{implies} \quad \exists!\, \hF[\covU']\to\hF[\covU].
\]
However, if in choosing the inclusions for the refinement we had made a different set of choices, $U_i'\to U_k$ rather than $U_j$, then a priori we might have expected different maps $\hF[\covU']\to\hF[\covU]$. Let us show this choice does not matter. If there is a choice, then we can take $U'_i\to U_j\cap U_k$ as a common refinement. As a consequence of $\hF$ being a functor from the open set category, the different maps to the colimit must agree, as they factor through whatever is assigned on the intersection, i.e. the following diagram commutes:
\[
\xymatrix{& & \hF(U_j) \ar[d] \\
\hF(U_i')  \ar[r] \ar[rru] \ar[rrd]& \hF(U_j\cap U_k) \ar[ru] \ar[rd] \ar[r] & \hF[\covU] \\
& & \hF(U_k) \ar[u]} 
\]
\end{proof}

\begin{cor}\label{cor:refined_axiom}\index{sheaf!behavior under refinement}\index{cosheaf!behavior under refinement}
	If $\hF$ is a cosheaf or $F$ is a sheaf for the cover $\covU'$, then it is a cosheaf or sheaf for every cover it refines.
\end{cor}
\begin{proof}
	Suppose we have a series of refinements
	\[
	\covU'\to\covU\to \{U\}.
	\] 
	To say that $\hF$ or $F$ is cosheaf or sheaf for $\covU'$ is to say that the following induced maps are isomorphisms:
	\[
	\xymatrix{\hF[\covU'] \ar[r] \ar@/^2pc/[rr]^{\cong} & \hF[\covU] \ar[r] & \hF(U)} \qquad \xymatrix{F(U) \ar[r] \ar@/^2pc/[rr]^{\cong} & \hF[\covU] \ar[r] & \hF[\covU']}
	\]
	However, by functoriality, the factored maps must themselves be isomorphisms, i.e.
	\[
	\xymatrix{\hF[\covU]\ar[r]^{\cong} & \hF(U)} \qquad \xymatrix{F(U) \ar[r]^{\cong} & F[\covU]}.
	\]
\end{proof}
We will make use of this corollary as we begin to consider sheaves and cosheaves on spaces where there is a finest cover. Checking the sheaf or cosheaf axiom there then guarantees it for all covers.

\section{Generalities on Sheaves and Cosheaves}
\label{subsec:generalities}

Sheaves have proved to be highly successful tools in pure mathematics over the past 60-70 years. This is largely because sheaves provide precise mechanisms for determining global unknowns from local knowns. These mechanisms are greatly enhanced by considering the operations such as $\Hom$ and $\otimes$ on sheaves, as well as the pushforward and pullback of a sheaf along a map, which we define in Section~\ref{subsec:six_ops}. Some of these operations are only available after applying a certain repair to turn a presheaf into a sheaf, also known as \textbf{sheafification}, which we define in Section~\ref{subsec:sheafification} after some preliminary discussion of examples.

One would like to know if a similar story is true for cosheaves. After all, any functor $\hF:\Open(X)\to\dat$ is exactly equivalent to specifying a functor $F:\Open(X)^{op}\to\dat^{op}$. However, certain asymmetries prevent such an observation from being as useful as one might hope. These asymmetries are outlined in Section~\ref{subsec:failure_to_commute} and obstruct the use of Grothendieck's version of sheafification to define an analogous \textbf{cosheafification}. However, we prove that such a device must exist in Section~\ref{subsec:cosheafification} without knowing a particular construction.

\subsection{Pre-sheaves and their Associated Sheaves}\label{subsec:sheafification}

Sheaves are fundamentally local structures. Informally stated, a pre-sheaf can fail to be a sheaf in two independent ways:
\begin{itemize}
\item \emph{Non-Local}: If a pre-sheaf has a section $s\in F(U)$ that cannot be constructed from sections over smaller open sets in $U$ --- a cover of $U$, for example --- then $F$ fails to be a sheaf.
\item \emph{Inconsistent}: If a pre-sheaf has a pair of sections $s\neq t\in F(U)$ such that when restricted to every smaller open set they define the same section, then $F$ fails to be a sheaf.
\end{itemize}

Let's illustrate both of these failures with two examples.
\begin{ex}[Non-local]\label{ex:nonlocalpresheaf}\index{presheaf!non-local example}
Let $F$ be a presheaf of vector spaces over the real line $\R$, defined as follows:
\begin{equation*}
 F(U)=
\begin{cases}
 k & \text{if } (-1,1)\subset U \\
0 & \text{o.w.}
\end{cases}
\end{equation*}
In particular, $F$ assigns the zero vector space to every open ball $B_r(x)$ centered at $x\in\R$ with $r\leq 1/2$. This collection of balls covers the real line thus if $F$ were a sheaf, then $F(\R)=0$, but it is the vector space $k$ instead.

An incarnation of this example, depicted in Figure~\ref{fig:nonlocal_presheaf}, is the pre-sheaf that assigns to every open set $U$, the first cohomology of the inverse image of $U$ under a map $f:S^1\to\RR$, i.e. 
\[
U\rightsquigarrow H_1(f^{-1}(U);k).
\] 
\begin{figure}[ht]
\centering
\includegraphics[width=.7\textwidth]{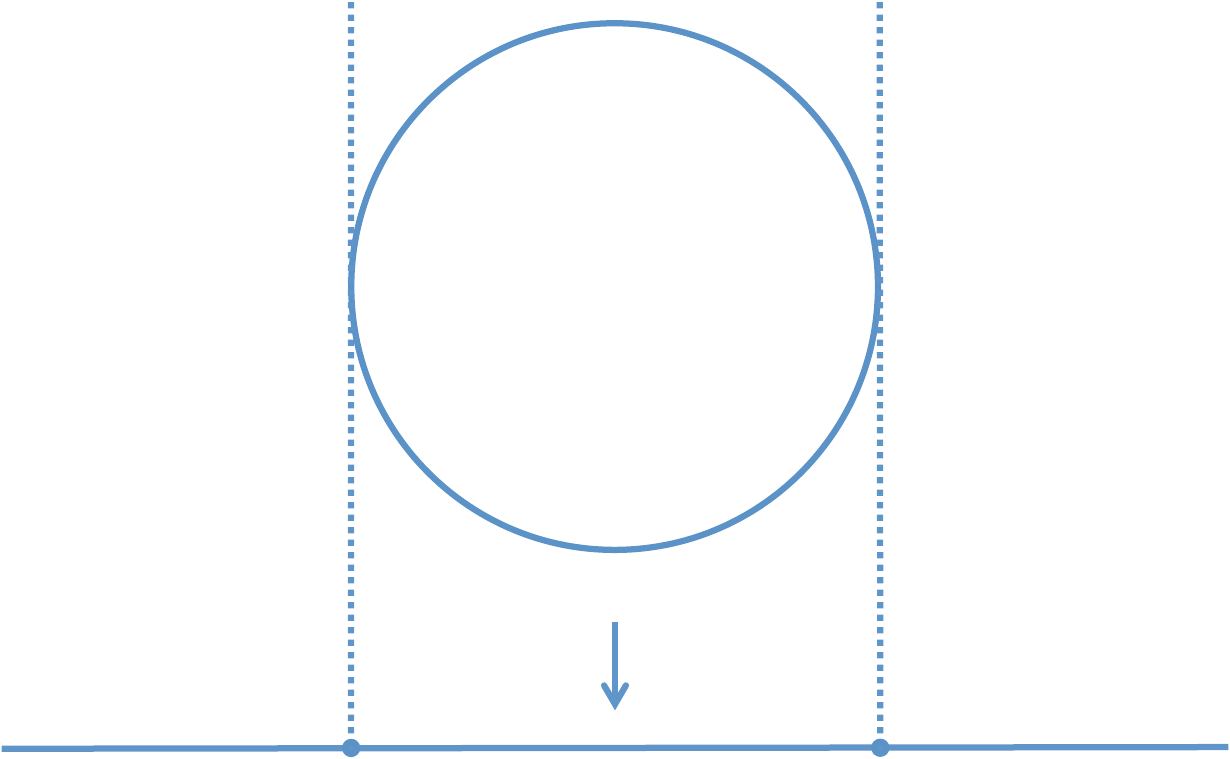}
\caption{Cohomology Pre-sheaf is Non-Local}
\label{fig:nonlocal_presheaf}
\end{figure}

\end{ex}

\begin{ex}[Inconsistent]\label{ex:inconsistentpresheaf}\index{presheaf!inconsistent example}
Let $F$ be a presheaf of sets over the real line $\R$, defined as follows:
\begin{equation*}
 F(U)=
\begin{cases}
 \{a,b\} & \text{if } U=X \\
\{\star\} & \text{o.w.}
\end{cases}
\end{equation*}
This presheaf is like two friends that agree on every possible political issue, but still belong to two different political parties.
\end{ex}

Fortunately, there is a general method of repair that can make any presheaf $F$ into a sheaf $\widetilde{F}$. Although this method modifies the value of $F$ on open sets, it leaves at least one feature of the presheaf unchanged. This is the stalk of the presheaf.

\begin{defn}\index{stalk}\index{presheaf!stalk}\index{sheaf!stalk}
Let $F$ be a pre-sheaf on a topological space $X$ and $x\in X$ a point. The \textbf{stalk} of $F$ at $x$ is the direct limit of $F$ over open sets $U$ containing $x$:
\[
F_x:=\varinjlim_{U\ni x} F(U)
\]
The stalk is the ``local value'' of a presheaf at $x$. Notice that every element $t\in F(U)$ with $x\in U$ has an associated value $t_x$, which is the image of $t$ in the direct limit.
\end{defn}
\begin{rmk}
 Since $F$ only assigns data to open sets, one often uses the direct limit construction to assign data to arbitrary sets of $X$; the stalk is just a special example of this principle.
 \end{rmk}

Now we introduce the procedure for turning an arbitrary presheaf of sets, vector spaces or groups into a sheaf of the same type.

\subsubsection{Sheafification}

We begin our introduction to the sheaf associated to a presheaf with a careful introduction to products and disjoint unions, following~\cite{deville2012dynamics}.

\begin{defn}[Disjoint Unions and Products in $\Set$]\index{product}\index{coproduct}\index{disjoint union}
Suppose $\{X_s\}_{s\in S}$ is a family of sets indexed by $S$. The \textbf{disjoint union} is a union that tracks the indexing set:
\[
\bigsqcup_{s\in S} X_s := \bigcup_{s\in S} X_s\times \{s\}
\]
The \textbf{product} can be written as a set of maps:
\[
\prod_{s\in S} X_s := \{ f: S \to \bigsqcup_{s\in S} X_s \,|\, f(s)\in X_s\, \forall s\in S\}
\]
The projection maps 
\[
\pi_{s'}:\prod_{s\in S} X_s \to X_{s'}
\]
are defined by evaluation $\pi_{s'}(f)=f(s')$.
\end{defn}

\begin{defn}[Sheafification]\label{defn:sheafification}\index{sheaf!sheafification}
Let $F:\Open(X)^{op}\to\Set$ be a presheaf of sets and let $F_x$ denote the stalk of $F$ at $x$. Now for each open set $U$, form the product $\prod_{x\in U}F_x$. The \textbf{sheafification} $\widetilde{F}$ of $F$ assigns to every open set $U$ the functions in $\prod_{x\in U}F_x$ that ``locally extend,'' i.e.
\[
\widetilde{F}(U):=\{s\in \prod_{x\in U}F_x\,|\, \forall x\in U\, s(x)\in F_x\,, \exists V\ni x\, V\subset U\,t\in F(V)\,\mathrm{s.t.}\, t_y=s(y)\,\forall y\in V\}
\]
There is a natural transformation $\theta:F\to \widetilde{F}$ that takes every element $s\in F(U)$ to the map $s:x\in U \mapsto s_x\in F_x$. In particular, $\theta_x:F_x\to\widetilde{F}_x$ is an isomorphism.
\end{defn}

One can summarize sheafification more elegantly in the language of categories. Since every sheaf is also a pre-sheaf, we have an inclusion functor
\[
\iota:\Shv(X;\Set)\hookrightarrow \Fun(\Open(X)^{op},\Set)=:\Preshv(X;\Set)
\]
that has a left adjoint, i.e. there is a universal natural transformation $\theta:\id_{\Preshv}\to \iota\circ \widetilde{(-)}$, see Section~\ref{sec:adjunctions} for a reminder. Such a subcategory is called \textbf{reflective}. This guarantees, for example, that if $F$ is an arbitrary pre-sheaf and $G$ is a sheaf regarded as a pre-sheaf $G=\iota(G)$, then we have the following universal property:\index{sheaf!sheafification!universal property as a left adjoint}
\[
\xymatrix{ & \widetilde{F} \ar@{.>}[d]^{\exists !} \\ F \ar[ur]^{\theta_F} \ar[r]^{\varphi} & \iota(G)}
\]
Pulling back along $\theta_F$ induces the natural isomorphism of $\Hom$-sets:
\[
\Hom_{\Shv}(\widetilde{F},G)\cong \Hom_{\Preshv}(F,\iota(G))
\]

\subsection{Grothendieck's Operations}\label{subsec:six_ops}
\index{Grothendieck!Six Functor Formalism}

What makes sheaf theory such a powerful machine is that there are many natural operations on sheaves and well understood adjunctions between these operations. However, many of these operations only exist with the aid of sheafification. In particular, there are the following six operations, grouped into three adjoint pairs, the third of which exists only in a suitable enlargement of the category of sheaves.
\[
(f^*,f_*) \qquad (\otimes,\SHom) \qquad (f_!,f^!)
\]
Here we will consider only four out of the six in order to forego this extra difficulty of ``enlarging'' the category of sheaves.

\begin{defn}[Pushforward Sheaf]\index{sheaf!pushforward}
Let $f:Y\to X$ be a continuous map and let $G$ be a sheaf on $Y$. The \textbf{pushforward sheaf} is defined by the formula:
\[
f_*G(U):=G(f^{-1}(U))
\]
\end{defn}

There should be an inverse operation that takes a sheaf $F$ on $X$ and pulls back along $f:Y\to X$. After all, if $i:W\hookrightarrow X$ is the inclusion of an open set, then a natural candidate for the pullback sheaf $i^*F$ would be the restriction of the domain of definition of $F$ to only those open sets contained in $W$.
\[
F|_W(U)=F(U)
\]
However, if $f:Y\to X$ is not an open map, then there is no hope for an easy definition. Sheafification, however, comes to the rescue.

\begin{defn}[Pullback Sheaf]\index{sheaf!pullback}\index{pullback!of sheaves}
Let $f:Y\to X$ be a continuous map and $F$ a sheaf on $X$. The \textbf{pullback sheaf}, written $f^*F$ is the sheafification of the pre-sheaf
\[
U\rightsquigarrow \varinjlim_{V\supset f(U)} F(V)
\]
\end{defn}

\begin{ex}[The Stalk]\index{stalk}\index{sheaf!stalk}
Let $i:\{x\}\hookrightarrow X$ be the inclusion of a point into a space with a sheaf $F$ defined on it. The sheaf $i^*F\cong F_x$ is the stalk at $x$.
\end{ex}
\begin{exr}
Verify for $i:W\hookrightarrow X$ that $i^*F=F|_W$
\end{exr}

The next pair of interest is the middle pair.
\begin{defn}[Sheaf Hom]\index{sheaf!Hom@$\SHom$}
Suppose $F$ and $G$ are sheaves of abelian groups over a single space $X$. The \textbf{sheaf hom} $\SHom(F,G)$ assigns to every open set
\[
\SHom(F,G)(U):=\Hom_{\Shv(U)}(F|_U,G|_U)
\]
\end{defn}

For the algebraically minded, there should be a knee-jerk response for an associated tensor sheaf, however the na\"ive assignment needs to be sheafified.

\begin{defn}\index{sheaf!tensor@$\otimes$}
Suppose $F$ and $G$ are sheaves of abelian groups over a single space $X$. The \textbf{tensor product of sheaves} $F\otimes G$ is defined to be the sheafification of the assignment
\[
U\rightsquigarrow F(U)\otimes G(U)
\]
\end{defn}

The reader is encouraged to work through the following exercise, borrowed from~\cite{achar-lsu} with a few extra hints.

\begin{exr}
Let $Q$ be the sheaf of sections of the map $f:S^1\to S^1$ defined via complex coordinates as $f(z)=z^2$, i.e.
\[
Q(U):=\{s:U\to S^1\,|\, f\circ s(z)=z\}.
\]
Check that this sheaf has no global sections. Now let $Q_{k}$ be the sheaf which assigns to each open set $U$ the $k$ vector space freely generated by the set $Q(U)$. Show by taking a carefully chosen cover of $S^1$ that
\[
F:U\rightsquigarrow Q_{k}(U)\otimes Q_{k}(U)
\]
is not a sheaf. Observe that we have a natural method for tensoring elements of $Q_{k}(U)$ together via pointwise multiplication. Any element $s\in Q_{k}(U)$ satisfies $(s\otimes s)(z)=s(z)^2=z$, but there are interesting cross-multiple terms.
\end{exr}

Although working out the above exercise is rewarding, the category theorist knows that since tensor products are colimit constructions and the sheaf axiom involves limits, one should instantly be suspicious of such a construction defining a sheaf. However, one can construct cosheaf-theoretic analogs of the above functors and there the tensor cosheaf is naturally a cosheaf, but cosheaf Hom needs to be ``cosheafified,'' if such a thing exists.

\subsection{Failures to Commute}\label{subsec:failure_to_commute}

Unfortunately, the universe appears to have a sort of handedness that makes certain constructions for sheaves natural, but not so for cosheaves. This is because most data categories $\dat$, such as $\Set,\Vect$ or $\Ab$, are not equivalent to their opposite categories. Thus the topological simplification of reducing cosheaves to sheaves comes at the cost of making the algebraic thinking more difficult. In particular, certain properties of $\dat=\Set,\Vect$ or $\Ab$ are used in the development of sheaf theory, which do not necessarily hold in $\dat^{op}$. The centerpiece of this discussion will be understanding that filtered colimits commute with finite limits in $\dat$, but cofiltered limits do not necessarily commute with finite colimits in $\dat$. Let us now relay the necessary definitions.

\begin{defn}\index{category!co@(co)filtered}
	A non-empty category $\cat$ is called \textbf{filtered} if the following two properties are satisfied: 
	\begin{itemize}
		\item For every pair of objects $x, y\in \cat$ there is a third object $z\in \cat$ with $x\to z$ and $y\to z$.
		\item For every pair of parallel morphisms $f,g:x\to y$ there is a third object and morphism $h:y\to z$ such that $h\circ f=h\circ g$.
	\end{itemize}
	A category $\cat$ is called \textbf{cofiltered} if $\cat^{op}$ is filtered. Sometimes, when the category is especially simple, we will simply call a cofiltered category filtered.
\end{defn}

\begin{ex}\index{category!of covers}
In Section~\ref{subsec:refinement} we considered the category of covers $\Cov(U)$ of an open set $U$. By noting that any two covers have a common refinement, one sees that this is an example of a cofiltered (or \emph{cofiltrant}) category.
\end{ex}

\begin{defn}\index{limit!cofiltered}\index{colimit!filtered}
	Suppose $I$ is a filtered indexing category with $F:I\to\dat$ and $G:I^{op}\to\dat$ diagrams in some category. We will call the colimit of $F$ a \textbf{filtered colimit} and the limit of $G$ a \textbf{cofiltered limit}.
\end{defn}

Now we give an example already introduced in the context of sheafification.
\begin{ex}[(Co)Stalks]\index{costalk}\index{stalk}
	Suppose $X$ is a topological space and $x$ is a point in $X$. The set of open sets containing $x$ defines a cofiltered subcategory of $\Open(X)$ or a filtered subcategory of $\Open(X)^{op}$. Consequently for a pre-sheaf $F$ or a pre-cosheaf $\hF$, the following
	\[
	F_x:=\varinjlim_{U\ni x} F(U) \qquad \varprojlim_{U\ni x}\hF(U)=:\hF_x
	\]
	define a filtered colimit and cofiltered limit, respectively. These are called the \textbf{stalk at x} of a pre-sheaf $F$ and the \textbf{costalk at x} of a pre-cosheaf $\hF$. Of course, to make such a statement meaningful, one needs to assume the data category $\dat$ has the relevant limits and colimits.
\end{ex}

The following theorem illustrates one of the fundamental differences between sheaves and cosheaves. It is expressed through the following algebraic fact, which the reader might like to compare with the Fubini theorem.
\begin{thm}\label{thm:filteredcolimits}\index{colimit!filtered!commutes with finite limits}
	Let $I$ be a filtered indexing category and $J$ a finite category. Then any functor $\alpha:I\times J\to\dat$ where $\dat=\Set,\Vect,$ or $\Ab$, has the property that the natural map
	\[
		\varinjlim_I\varprojlim_J \alpha(i,j) \to \varprojlim_J\varinjlim_I \alpha(i,j)
	\]
	is an isomorphism. We say for short that ``filtered colimits and finite limits commute'' in these categories.
\end{thm}
\begin{proof}
For a proof of this statement, we refer the reader to theorem 3.1.6 in~\cite{KS-CS}. First note that the product category $I\times J$ is just a product in $\Cat$ - the category of all categories. We can describe the objects of $I\times J$ as pairs of objects, one from each category, and the morphisms as tuples of morphisms, one for each object. One can take $\alpha$ and then define a new functor $\varprojlim_J \alpha:I\to\dat$ gotten by assigning to each object $i\in\obj(I)$ the limit over $J$ of $\alpha(i,-):J\to\dat$. Taking the colimit over $I$ defines the first expression. Similar reasoning defines the second.
\end{proof}

\index{exact!sequence of presheaves}\index{exact!sequence of sheaves}\index{sheaves!exact sequence of}\index{presheaf!exact sequence}
For an application of this theorem we introduce some ideas to the world of pre-sheaves valued in $\Vect$. Recall that a complex of vector spaces 
\[
\cdots \to V_1\to V_2 \to V_3 \to \cdots
\] 
is exact at a term in a sequence if the image of the incoming map coincides with the kernel of the outgoing map. A sequence of pre-sheaves\footnote{The corresponding statement for sheaves is \emph{not} true.} is exact if and only if for each open set $U$ the associated sequence of vector spaces is exact, i.e.
\[
	0 \to E \to F \to G \to 0 \qquad \mathrm{iff} \qquad 0 \to E(U) \to F(U) \to G(U) \to 0
\]
Theorem~\ref{thm:filteredcolimits} then implies that for any point $x\in X$ the induced sequence of stalks
\[
	0 \to E_x \to F_x \to G_x \to 0
\]
is exact. Intuitively this is because we can view $E$ as a kernel of the pre-sheaf map $F\to G$ and as already demonstrated, kernels are examples of finite limits. Thus taking the kernel of the stalk map $F_x\to G_x$ is the same as taking the stalk of the kernel of $F\to G$.

\begin{prop}\label{prop:dont_commute}\index{limit!cofiltered!do not commute with finite colimits}
	For $\dat=\Set,\Vect$ or $\Ab$ it is not true that cofiltered limits and finite colimits commute. Consequently, if $A,B,C:\NN^{op}\to\Ab$ (or $\Vect$) are functors from the category of natural numbers equipped with the opposite ordering, with natural transformations $A\to B \to C$ such that
	\[
		\xymatrix{ 0 \ar[r] & A_i \ar[r] & B_i \ar[r] & C_i \ar[r] & 0}
	\]
	is exact for every $i$, then it is \emph{not always} the case that the induced sequence on limits is exact.
	\[
		\xymatrix{ 0 \ar[r] & \varprojlim A \ar[r] & \varprojlim B \ar[r] & \varprojlim C \ar[r] & 0}
	\]
\end{prop}
\begin{proof}
	We borrow an example from Jason McCarthy's notes~\cite{limits-jason}. Consider the following system of short exact sequences of groups:
	\[
		\xymatrix{\vdots \ar[d] & \vdots \ar[d] & \vdots \ar[d] & \vdots \ar[d] & \vdots \ar[d] \\
		0 \ar[d] \ar[r] & \ZZ \ar[d]_{n+1} \ar[r]^n & \ZZ \ar[d]_{n+1} \ar[r] & \ZZ/n \ar[r] \ar[d]_{\id} & 0 \ar[d] \\
		0 \ar[d] \ar[r] & \ZZ \ar[d]_{n+1} \ar[r]^n & \ZZ \ar[d]_{n+1} \ar[r] & \ZZ/n \ar[r] \ar[d]_{\id} & 0 \ar[d] \\
		0  \ar[r] & \ZZ  \ar[r]^n & \ZZ  \ar[r] & \ZZ/n \ar[r] & 0}
	\]
	The inverse limit of the first (and second) column with non-zero entries must be zero. To see why, note that the inverse limit can be described as
	\[
		\varprojlim_{\NN^{op}} \ZZ_i = \{(x_i)\in \prod_i \ZZ_i\, |\, x_i=(n+1)^{j-i}x_j \,\forall j\geq i\},
	\]
	where we have viewed the indexing category as the natural numbers with the opposite ordering. Any non-zero element of the limit would have some non-zero factor $x_i$ and consequently all other factors would be non-zero (since the map $a\to (n+1)a$ is injective). In particular, all higher $x_j$ must be equal to $x_i/(n+1)^{j-i}$, but letting $j$ be suitably large would imply that $x_j$ must be less than one --- an impossibility. Thus the induced map of inverse limits is
	\[
		\xymatrix{0 \ar[r] & 0 \ar[r] & 0 \ar[r] & \ZZ/n\ZZ \ar[r] & 0}
	\]
	which is not exact. If the reader prefers an example in the category of vector spaces, one should see Schapira's example 4.2.5 in his notes~\cite{schapira-AlTo}.
\end{proof}

Thus the statement that short exact sequences of pre-cosheaves induces a short exact sequence on costalks cannot be guaranteed. There is a subtle work-around that says under suitable hypotheses\footnote{i.e. the Mittag-Leffler condition. See~\cite[Sec. 1.12]{KS} or~\cite[pp. 211-214]{alpine} for more details.} exactness can be guaranteed. This holds for categories like $\vect$, the category of finite-dimensional vector spaces, and $\mathbf{ab}$, the category of finite abelian groups, because this is where the descending chain condition holds~\cite{atiyah-mac}.

This last comment about $\vect$ provides justification for performing some dualization to obtain results about cosheaves from sheaves. After all, for finite-dimensional vector spaces it is true that
\[
	\Hom_{\vect}(-,k):\vect^{op}\to\vect
\]
establishes an equivalence of categories.\footnote{This does \emph{not} extend to an equivalence between $\Vect$ and its opposite category. In fact, $\Vect^{op}$ is equivalent to the category $\pro-\vect$, cf.~\cite[Rmk.6.2]{isaksen-pro}.\index{category!opposite!of Vect@$\Vect$}} However, issues of stalks versus costalks is not the primary obstacle that the asymmetry of Theorem~\ref{thm:filteredcolimits} presents. That obstacle has to do with a process known as \textbf{sheafification}, which provides a universal tool for turning any pre-sheaf into a sheaf. For most texts on sheaf theory it is presented before almost any other theory is developed.

\index{sheaf!sheafification}
The most general sheafification process outlined by Grothendieck takes a pre-sheaf $F$ and defines a new pre-sheaf $F^+$ that assigns to each open set $U$ the filtered colimit of $F:\Cov(U)^{op}\to\dat$, see~\cite[Sec. 17.4]{KS-CS} for a modern exposition. Applying this construction twice defines a sheaf. However, in order to guarantee that this $F^{++}$ is a sheaf one uses the properties of Theorem~\ref{thm:filteredcolimits}. This now gets us to the more fundamental reason why the study of cosheaves may be so obscure: The non-exactness of $\varprojlim$ thwarts the Grothendieck prescription for cosheafification. \emph{For pre-cosheaves valued in $\Set, \Ab$ or $\Vect$ there is simply no hope in using the standard, most general, cosheafification.}
\index{cosheaf!cosheafification!approaches}

There are a very small handful of approaches that have been used to circumvent this problem:
\begin{enumerate}
 \item \textbf{\v{C}ech Homology and Smoothness}: One approach developed by Bredon~\cite{bredon-ch,Bredon} is to define an equivalence relation on pre-cosheaves, more nuanced than isomorphism, which is constructed through zig-zag diagrams of \emph{local} isomorphisms. Bredon develops an operation which uses \v{C}ech homology to take in one pre-cosheaf and produce another. In the event that the starting pre-cosheaf was equivalent to a cosheaf (Bredon calls such a pre-cosheaf \emph{smooth}), he proves that his construction yields a cosheaf.

 \item \textbf{Pro-Objects}: Another notable approach is to use pro-objects, i.e. functors $P:I^{op}\to\cat$ where $I$ is filtrant. This theory is engineered in such a way that all the desired algebraic properties exist. This approach was perhaps first used by Jean-Pierre Schneiders~\cite{schneiders-ch} to develop a rich theory of cosheaves. The problem with pro-objects is its conceptual and algebraic difficulty. For the visually minded, cosheaves of pro-objects are infinite diagrams of infinite diagrams, which obscure the many natural examples of pre-cosheaves and cosheaves that one might want to capture. More recent work~\cite{sugiki,prasolov}, has also used this setup for cosheaves.

 \item \textbf{Topology:} Here, one eschews full generality and works only with certain cosheaves known as \textbf{constructible cosheaves}, which can be thought of as cosheaves on particular finite spaces. Cosheafification in this setting exists and is natural. Often one does not even think about needing to cosheafify, because the diagrams are modeled on the points of the space. This school of thought, motivated by the vision and unpublished ideas of Bob MacPherson, has some recent trace in the literature, see~\cite{woolf, dt-lag}.
\end{enumerate}

For the most part, we choose to sidestep the issues of sheafification and cosheafification by focusing on the third approach. We believe that this provides a better way of learning sheaf theory as it removes the ever-present phrase ``let \emph{blank} be the sheafification of \emph{blank}'' and focuses on the more important technical machinery first. However, we provide a proof that cosheafification does exist in the next section.

\subsection{The Existence of Cosheafification}
\label{subsec:cosheafification}
\index{cosheaf!cosheafification}

Grothendieck gave us a general construction of the sheafification functor $\widetilde{(-)}$ that works in more general data categories $\dat$, which comes from applying a certain functor $(-)^{+}$ twice. The requirement on the category $\dat$ includes, among other things, that filtered colimits and finite limits commute. If we were to regard a pre-cosheaf $\hF:\Open(X)\to\dat$ as a pre-sheaf valued in $\dat^{op}$, then the condition that filtered colimits and finite limits commute in $\dat^{op}$ would become the condition that cofiltered limits and finite colimits commute in $\dat$, which is patently false when $\dat=\Set,\Vect$ or $\Ab$, as the example in Proposition~\ref{prop:dont_commute} showed.

Consequently, we do not have a clear answer to the question: Does the inclusion functor $\iota$ have a \emph{right} adjoint $(-)_{\#}$?\index{cosheaf!cosheafification!universal property as a right adjoint}
\[
\iota:\Coshv(X;\dat) \hookrightarrow \Fun(\Open(X),\dat)=:\Precoshv(X;\dat)
\]
so that the dual universal property is satisfied, i.e. if $\hF$ is a pre-cosheaf and $\hG$ is a cosheaf with a morphism $\hG\to\hF$, then there is a unique way of completing the diagram.
\[
\xymatrix{ & \hF_{\#} \ar[d] \\
\hG \ar@{.>}[ur] \ar[r] & \hF}
\]
\index{cosheaf!cosheafification!exists for Set@$\Set$}
In the case where $\dat=\Set$, Jon Woolf's paper~\cite{woolf} contains a construction of cosheafification. Unfortunately, this cannot be adapted to categories like $\Vect$ or $\Ab$. For a high-level reason why, Mac Lane and Moerdijk explain on page 95 of~\cite{mm-topos} that a sheaf of abelian groups can be identified with an abelian group object in the category of sheaves. Since sheafification preserves finite products, sheafification of pre-sheaves of sets lifts to a functor between abelian group objects. Moreover, since the forgetful functor $for:\Ab\to\Set$ preserves limits (but not colimits), any sheaf of groups defines a sheaf of sets. Trying to repeat this last line of reasoning for cosheaves of groups fails, i.e. a cosheaf of groups does not, by forgetting, define a cosheaf of sets.

Our approach is to verify abstractly whether cosheafification exists without constructing it. Of course, one would like to use Freyd's Adjoint Functor Theorem~\ref{thm:freyd}, but we will use a different theorem~\cite[Thm 6.28]{adamek-lpac} that is easier to check.

\begin{thm}\index{Vopenka's principle}
Assuming Vopenka's principle (a large cardinal axiom), every full subcategory $\bat$ of a locally presentable category $\cat$, where $\bat$ is closed under colimits, is coreflective, i.e. the inclusion functor $\iota:\bat\hookrightarrow\cat$ has a right adjoint (a cofree functor).
\end{thm}

We will leave Vopenka's principle as a black box and assume it, even though many category theorists cringe at its very name. We prove that cosheafification exists by verifying the hypotheses of the above theorem for our case of interest.

\begin{cor}
The category of cosheaves of vector spaces is a coreflective subcategory of $\Fun(\Open(X),\Vect)$, i.e. cosheafification exists.
\end{cor}
\begin{proof}
It is clear that the category of cosheaves is closed under colimits, since we can define the colimit to be the pre-cosheaf, which open-by-open assigns the colimit of vector spaces over that open set. This pre-cosheaf is a cosheaf, because for a fixed cover, each vector space in the diagram is expressed as a colimit and colimits commute with colimits.

It remains to be seen that the category of pre-cosheaves is locally presentable. this means that the category is locally small, has small colimits, has a small set of objects $\covS$ that generates $\Precoshv(X;\Vect)$ in the sense that every pre-cosheaf is a colimit of objects in $\covS$, and every object is small. The first two statements are easily addressed. $\Open(X)$ is a small category and $\Vect$ is locally small, so the functor category is locally small. Colimits of pre-cosheaves are defined open-by-open. Since $\Vect$ is cocomplete, pre-cosheaves valued in $\Vect$ is also cocomplete. Now we address the existence of a generating set. 

Define, for each open set $U\in \Open(X)$, the following pre-cosheaf:
\[
\hat{h}_U(V)=\left\{\begin{array}{ll} k & if\, U\subset V \\ 0 & o.w.\end{array}\right.
\]
We'd like to say that every pre-cosheaf is a colimit of pre-cosheaves of the above form. The corresponding statement for pre-sheaves is proved in pages 41-42 of~\cite{mm-topos}. We will go ahead and repeat the argument here. Note that if $\hG$ is an arbitrary pre-cosheaf, then
\[
\Hom_{\Precoshv}(\hat{h}_U,\hG)\cong \hG(U).
\]
Observe that if $U\subset U'$, then we get a map of pre-cosheaves (a natural transformation) $\hat{h}_{U'}\to \hat{h}_{U}$. This in turn induces a map
\[
\Hom_{\Precoshv}(\hat{h}_U,\hG)\to \Hom_{\Precoshv}(\hat{h}_{U'},\hG)
\]
which coincides with the internal extension map of $\hG$, that is $r_{U',U}:\hG(U)\to \hG(U')$. In other words, the functor
\[
R:\Precoshv(X;\Vect)\to\Precoshv(X;\Vect)
\] 
\[
G \rightsquigarrow (U\mapsto \Hom_{\Precoshv}(\hat{h}_U,\hG)\cong \hG(U))
\]
is isomorphic to the identity functor. Since adjoints are unique up to isomorphism, then we can conclude that its (left) adjoint must also be isomorphic to the identity functor. However, we will construct explicitly the adjoint, which, combined with the fact that it must be the identity functor, exhibits $\hG$ as the colimit of pre-cosheaves of the form $\hat{h}_U$.

For each pre-cosheaf $\hG$, define the following \textbf{category of elements} $J_{\hG}$.\footnote{It should be noted that in~\cite{mm-topos}, the category of elements is written $\int \hG$.} The objects of $J_{\hG}$ are pairs
\[
(U,x) \qquad \mathrm{where} \qquad U\in\Open(X) \quad x\in \hG(U).
\]
A morphism $(U,x)\to (U',x')$ is defined if $U\subset U'$ and $x'=r_{U',U}(x)$. Clearly, there is a projection functor $\pi_{\hG}:J_{\hG}\to\Open(X)$ and by formality, there is a dual functor $\pi^{op}_{\hG}:J_{\hG}^{op}\to\Open(X)^{op}$.

Denote by $\mathcal{Y}$ the functor
\[
\mathcal{Y}:\Open(X)^{op}\to\Precoshv(X;\Vect) \qquad \qquad U \rightsquigarrow \hat{h}_U.
\]
We claim that the left adjoint $L$ to the functor $R$ considered above can be constructed object-wise as follows: for each pre-cosheaf $\hG$ define
\[
L(\hG):=\varinjlim \mathcal{Y}\circ \pi_{\hG}^{op}.
\]
We claim that $\hG$ is the colimit. This is already given from the fact that $L$ must be isomorphic to the identity functor, but let's at least check how $\hG$ is a co-cone, to make the statement more plausible. For each object $(U,x)$ in $J_{\hG}$ the map to $\hG$ is defined by
\[
x\in \hG(U)\cong \Hom_{\Precoshv}(\hat{h}_U,\hG)\ni \psi_{U,x}
\]
where $\psi_{U,x}$ is the natural transformation that sends $1\in\hat{h}_U(U)$ to $x\in \hG(U)$ and then for any larger open set $U\subset U'$ sends $1\in\hat{h}_U(U')$ to $r_{U',U}(x)$.
Observe that if $U\subset U'$ and $(U',x')\to (U,x)$ is a morphism in $J_{\hG}^{op}$, so $r_{U',U}(x)=x'$, then we have the following commutative diagram:
\[
\xymatrix{(U',x')\ar[d] \ar@{~>}[r] & \hat{h}_{U'}\ar[d] \ar[rd]^{\psi_{U,x'}} & \\
(U,x) \ar@{~>}[r] & \hat{h}_{U} \ar[r]_{\psi_{U,x}} & \hG}
\]
At the risk of demonstrating the obvious, the above diagram commutes if the the following diagram commutes for an arbitrary triple of open sets $V\subset V'\subset V''$. We will check it for the interesting boundary case $U\subset U'\subset U''$.
\[
\xymatrix{\hat{h}_{U'}(U'')=k \ar[r]^1 & \hat{h}_{U}(U'')=k \ar[r]^{r_{U'',U'}(x)} & \hG(U'') \\
\hat{h}_{U'}(U')=k \ar[r]^1 \ar[u]^1 & \hat{h}_U(U')=k \ar[u]^1 \ar[r]^{r_{U,U'}(x)=x'} & \hG(U') \ar[u]_{r_{U'',U'}} \\
\hat{h}_{u'}(U)=0 \ar[r] \ar[u] & \hat{h}_U(U)=k \ar[r]^{x} \ar[u]^1 & \hG(U) \ar[u]_{r_{U',U}} }
\]
This completes the plausibility check. We use the observation that $L$ is isomorphic to the identity functor to conclude that $\hG$ is actually the colimit. The conclusion is that
\[
\hG\cong \varinjlim \mathcal{Y}\circ\pi_{\hG}^{op}
\]
i.e. $\hG$ is expressible as a small colimit of pre-cosheaves of the form $\hat{h}_U$ where the size of the indexing set is bounded by the product of the cardinality of $\Open(X)$ and the maximum cardinality of $\hG(U)$ over varying $U$.

Now it remains to check the smallness of objects in $\Precoshv(X;\Vect)$. An object $\hG$ is small if there exists a regular cardinal $\kappa$ such that $\Hom(\hG,-)$ commutes with directed colimits of diagrams indexed by categories of cardinality at most $\kappa$.

Observe that for one of our pre-cosheaves $\hat{h}_U$ is compact since if $(\hF_i)$ is a direct system of pre-cosheaves, then
\[
\varinjlim\Hom(\hat{h}_U,\hF_i)\cong \varinjlim \hF_i \cong \Hom(\hat{h}_U,\varinjlim \hF_i).
\]
As already shown, for every pre-cosheaf $\hG$ there exists a diagram whose cardinality is the cardinality of $J_{\hG}$, which we will call $\kappa_J$. Now, we know how to express $\hG$ as a colimit of $\hat{h}_U$'s. Thus,
\begin{eqnarray*}
\Hom(\hG,\varinjlim_i \hF_i) &\cong & \Hom(\varinjlim_U \hat{h}_U,\varinjlim_i \hF_i)\\
&\cong &\varprojlim_U \Hom(\hat{h}_U,\varprojlim_i \hF_i) \\
&\cong &\varprojlim_U\varinjlim_i \Hom(\hat{h}_U,\hF_i) \\
&\cong &\varprojlim_i\varinjlim_U \Hom(\hat{h}_U,\hF_i) \\
&\cong &\varinjlim_i\Hom(\varinjlim_U \hat{h}_U,\hF_i) \\
&\cong &\varinjlim_i\Hom(\hG,\hF_i)
\end{eqnarray*}

The third line follows from compactness of $\hat{h}_U$. The fourth line follows from the fact that in $\Set$, $\kappa_J$ small colimits commute with $\kappa_J$ -filtered colimits. This completes the proof.

\end{proof}

%
%
\chapter{Preliminary Examples}
\label{sec:examples}

\begin{quote}
{\em``The content of a mathematical theory is never larger than the set of examples that are thoroughly understood.''}
\begin{flushright} --- Vladmir Arnol'd~\cite{arnol2004lectures} \end{flushright}
\end{quote}

Theories should be motivated by examples. In this chapter we develop the common themes these examples share. Broadly speaking, all sheaves are realized via local sections associated to a particular map. This principle is rigorously embodied by the \emph{\'etal\'e perspective}. Similarly, all cosheaves of sets are realized by connected components of the fiber of a map, embodied by the \emph{display perspective}, which is a generalization of the Reeb graph construction outlined in Definition~\ref{defn:reeb_graph}.

\section{Sheaves Model Sections}

Recall that if $f:Y\to X$ is a continuous map then a \textbf{section} is a continuous map $g:X\to Y$ such that $f(g(x))=x$ for all $x$. This definition has the property that $f$ is surjective. Sometimes a map admits a locally-defined section over a subset $U\subset X$, but not a global one. There is a sheaf that tracks this data.

\begin{defn}[Sheaf of Sections of a Map]\index{sheaf!of sections}
	Suppose $\pi:E\to X$ is a continuous map. Then we can associate a \textbf{sheaf of sections} to this map as follows:
	\[
	U\rightsquigarrow F(U):=\{s:U\to \pi^{-1}(U)\,\mathrm{continuous}\, | \pi(s(x))=x \}.
	\]
\end{defn}

Clearly, if $F(X)\neq \emptyset$, then we can answer positively the question ``Does $\pi:E\to X$ have a section?''
\[
\xymatrix{E\ar[d]_{\pi} \\ X \ar@/_/@{.>}[u]_{?}}
\]

To see why this is a pre-sheaf valued in $\dat=\Set$ note that what is assigned to an open set $U$ is a set of maps. A map whose domain of definition is $U$ can always be restricted to a smaller open subset $V\subset U$ to define a map on $V$. This process of restricting the domain of definition we write as $\rho_{V,U}(s):=s|_V$, which is what makes this assignment a pre-sheaf. 

Let us prove this defines a sheaf. Suppose $U\subset X$ and $\covU=\{U_i\}_{i\in\Lambda}$ is an arbitrary open cover of $U$. We must prove that the map
\[
F(U)\to F[\covU]:=\varprojlim(N(\covU)^{op}\to\Open(X)^{op}\to \Set)
\]
is an isomorphism.
Recall that the limit can be described in terms of products and equalizers. As such, every element of the limit is described by a collection of continuous sections $s_i:U_i\to \pi^{-1}(U_i)$, one for each element of the cover, such that on intersections $\rho_{ij,i}(s_i)=\rho_{ij,j}(s_j)$.\footnote{Here we have adopted the shorthand of referring to open sets via elements of the nerve.} The natural map from $F(U)$ to $F[\covU]$ simply takes a section $s\in F(U)$ to the collection of restricted sections $\{s_i:=s|_{U_i}\}$. If two sections over $U$ differ at a point $x$, then they will define different sections over $U_i\ni x$, thus the natural map is injective. To check surjectivity, note that an element in the limit defines a section over $U$ by setting $s(x)=s_i(x)$ if $x\in U_i$ and this will be continuous by the pasting lemma described at the beginning of Chapter~\ref{sec:abstract_sheaves}.

\begin{figure}[ht]
\centering
\includegraphics[width=.7\textwidth]{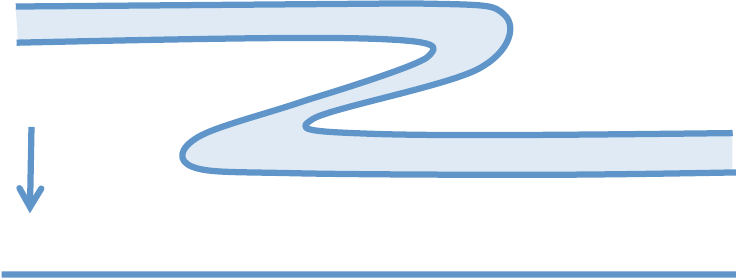}
\caption{Is There a Section?}
\label{fig:evade_sec}
\end{figure}

\begin{ex}
For a simple example, consider the projection onto the first coordinate $\pi_t:[0,1]\times [0,1]\to [0,1]$, which we regard as taking a time-space coordinate $(t,x)$ to its time coordinate $t$. There are lots of sections of this map. The map that assigns to each time $t$ a fixed position $p\in[0,1]$ defines a section, so there are uncountably many sections.

Now consider a different map that comes from restricting the time projection map to a subset $E\subseteq [0,1]\times[0,1]$, i.e. $\pi:=\pi_t|_E:E\to[0,1]$ is the restricted map. A drawing can be found in Figure \ref{fig:evade_sec} where $E$ is the region bound between the two curves. Does it have any global sections, i.e. is $F(X)\neq\emptyset$?
\end{ex}

The answer is clearly no. The example in Figure \ref{fig:evade_sec} illustrates a concept central to sheaf theory. Although about each point in time $t$ there is some $\epsilon>0$ such that on the open set $(t-\epsilon,t+\epsilon)$ a continuous section can be defined, there is no globally defined section. Thus local sections (local solutions) exist, but they do not always glue together to define a global section (global solution). This is why we say
\begin{quote}
	\begin{center}
	\emph{Sheaves mediate the passage from local to global.}
\end{center}
\end{quote}

\begin{ex}[Square Map]
	Suppose $f:\CC\to\CC$ is the map sending a complex number $z$ to $z^2$. For a  point $w=re^{i\theta}$ with $r\neq 0$ there are two points in the fiber: $z=\sqrt{r}e^{i\theta/2}$ and $z'=\sqrt{r}e^{i\theta/2 +\pi}$. Consequently, for a small connected neighborhood about $w$ there are two corresponding continuous sections. There is no global section because the square root map is necessarily multi-valued when considered all the whole complex plane.
\end{ex}

Lists of similar examples abound in geometry and topology, most of which are concerned with the following mathematical structure.

\begin{defn}\label{defn:fiber_bundle}\index{fiber bundle}\index{covering space}
A \textbf{fiber bundle} over $X$ consists of total space $E$ equipped with a continuous surjective map $\pi:E\to X$ satisfying the property that for each point $x\in X$ there exists an open neighborhood $U$ such that the following diagram commutes.
\[
\xymatrix{\pi^{-1}(U) \ar[rr]^{h_U} \ar[rd]_{\pi} & & U\times F \ar[ld]^{\pi_U} \\ & U & } 
\]
Here $F$ is the fiber space, $h_U$ is a homeomorphism and $\pi_U$ is projection onto the first factor. If $F$ is a discrete space then we usually write $\tilde{X}$ instead of $E$ and say that $\pi:\tilde{X}\to X$ is a \textbf{covering space}. If each fiber $\pi^{-1}(x)$ is endowed with the structure of a group, i.e. $F=G$ with the discrete topology, so that $h_U$ induces a group isomorphism between $\pi^{-1}(x)$ and $G$, then $E$ is called a \textbf{bundle of groups}. Analogous definitions hold for fiber a ring or a module.
\end{defn}

\index{M\"obius bundle}
The map $\pi:M\to S^1$ where $M:=S^1\times \RR/\sim$ with $(x,y)\sim (x+2\pi,-y)$ is an example of a fiber bundle over $S^1$. Restricting the domain of $\pi$ to the subspace $S^1\times [-1,1]$ allows one to think of this map as projecting the M\"obius bundle to its core circle. The projection $\pi$ has a section that embeds $S^1$ as the zero section, but there are no sections which avoid $S^1\times\{0\}$. The ``hairy ball'' theorem is the analogous statement except for the tangent bundle to the two sphere $S^2$. Sheaf theory is the \emph{lingua franca} for bundle theory and category theory. Thus even the most trivial example of a product bundle, $E=X\times k\to X$ where $k$ is a field, is of interest.

\begin{defn}[Constant Sheaf]\label{defn:const_sheaf}\index{sheaf!constant sheaf}\index{sheaf!locally constant sheaf}
	Suppose $A$ is an $R$-module equipped with the discrete topology and $E=X\times A\to X$ is the product bundle. The sheaf of sections of this map is called the \textbf{constant sheaf} $A_X$. If $A=R$ is a field $k$ or the ring $\ZZ$ we will just say the constant sheaf and write $k_X$ or $\ZZ_X$. We will almost always work over a field $k$.
	
	If $\pi:E\to X$ is not necessarily the product bundle, but has fiber $A$, then we call the sheaf of sections of $\pi$ the \textbf{locally constant sheaf} $F$ with value $A$.  
\end{defn}

For locally connected spaces the constant sheaf $k_X$ assigns to any open set $U$ the product of the field $k$ for as many connected components as $U$ has. To see this, let us be more precise about how the algebraic structure of the module/vector space interacts with the topological structure of a fiber bundle~\cite[Sec. 7.2.1]{macpherson-ih-notes}.

\begin{defn}[Local Systems]\label{defn:local_systems}\index{local system}
Let $k$ be a field viewed as a topological space with the discrete topology. 
A \textbf{local system} is a covering space $\pi:L\to X$ equipped with the structure of a $k$-vector space on each fiber. Specifically, an \textbf{$n$-dimensional local system} on a topological space $X$ is a topological space $L$, a map of spaces $\pi:L\to X$, and, for each point $p\in X$, a $k$-vector space structure on $\pi^{-1}(p)$ with the following property: Every point $p\in X$ has a neighborhood $U$ such that there is a homeomorphism $h_U:\pi^{-1}(U)\to U\times k^n$ such that $\pi_U\circ h_U=\pi$ and for each $x\in U$ the vector space structure on $\pi^{-1}(x)$ is induced by $h$ from the one on $\{x\}\times k^n$.
\end{defn}

Our definition of the locally constant sheaf $F$ in Definition~\ref{defn:const_sheaf} is more accurately defined as an $n$-dimensional local system $L$, at least when $A$ is a $k$-vector space. Consider one of the distinguished open neighborhoods $U$ of a point $x\in X$ provided by the definition. Here we have a commutative diagram:
\[
\xymatrix{\pi^{-1}(U) \ar[rr]^{h_U} \ar[rd]_{\pi} & & U\times k^n \ar[ld]^{\pi_U} \\ & U & } 
\]
For a locally connected space we can assume $U$ is connected by replacing $U$ with whatever connected component of $U$ contains $x$. Consider a section of $\pi:L\to X$ over $U$. Since $k^n$ has the discrete topology any section $s$ over $U$ has to be constant since the image of a connected set is always connected. This implies that every section $s$ has the form $h_U\circ s(y)=(y,(v_1,\ldots,v_n))$ for every $y\in U$ and for some fixed vector $\bar{v} \in k^n$. Consequently, for the distinguished neighborhood $U$, the set of sections
\[
F(U)=A_X(U)\cong k^n\cong A.
\]
Moreover, by the local system condition, we can form any linear combination of sections $s_1,s_2\in F(U)$ to obtain a third section $\alpha s_1+\beta s_2\in F(U)$. This implies that the locally constant sheaf is actually a sheaf valued in $\Vect$ --- the category of vector spaces. Arguing in the same way for each connected component tells us that over a union $U'$ of disjoint, connected, distinguished neighborhoods, a locally constant sheaf has the value
\[
F(U')\cong A^{\pi_0(U')} \cong H^0(U';A)\cong H^0(U';k)^n.
\]
This illustrates that a locally constant sheaf can be thought of as taking $H^0$ of a space with ``twisted'' coefficients.

\section{Local Systems: A Bridge Between Sheaves and Cosheaves}
\label{subsec:local_systems}

Formulating locally constant sheaves as a twisted $H^0$ presents an obvious dualization in terms of $H_0$. Moreover this duality reaches higher by considering higher cohomology and homology as well. To understand this, we will need to understand local systems better. 

\begin{defn}\index{local system}\index{category! of local systems}
The collection of local systems over $X$ forms a category $\LocSys(X)$. A morphism of between two local systems $\pi:L\to X$ and $\pi':L'\to X$ is a map $\varphi:L\to L'$ such that $\pi'\circ \varphi =\pi$ and the restricted map $\varphi_x:\pi^{-1}(x)\to \pi^{'-1}(x)$ is a linear transformation.
\end{defn}

The following theorem is classical and allows us to use two definition of local systems interchangeably.

\begin{thm}\index{local system!representation of $\pijuan(X)$}\index{fundamental groupoid $\pijuan(X)$!local systems}\index{representation!of fund@$\pijuan(X)$}
If $X$ is a locally connected and locally simply connected space, then the category of local systems is equivalent to the category of representations of the fundamental groupoid of $X$, i.e. 
\[
\LocSys(X)\simeq \Rep(\pijuan(X)).
\]
Recall that the objects of $\Rep(\pijuan(X))$ are functors $\cL:\pijuan(X)\to \Vect$.
\end{thm}
\begin{rmk}
For a connected space $X$ fixing a base point $x_0$ provides a skeletal subcategory $\pijuan(X;x_0)\hookrightarrow \pijuan(X)$. Precomposing $\cL$ with this inclusion defines a representation of the fundamental group $\pijuan(X;x_0)$.
\end{rmk}
\begin{proof}[Proof (Idea)]
The functor that realizes this equivalence is very easy to describe. Given a local system $L$ one defines a representation $\cL:\pijuan(X)\to\Vect$ by assigning to points $x\in X$, the vector space $\pi^{-1}(x)=:\cL(x)$. Now suppose $\gamma:[0,1]\to X$ is a path connecting $x$ to $y$. Since $\pi:L\to X$ is a covering space, there is a unique lift $\tilde{\gamma}$ connecting any element of $\pi^{-1}(x)$ to an element in $\pi^{-1}(y)$. These lifts piece together to define a \textbf{monodromy map} $\mu_{\gamma}:\cL(x)\to\cL(y)$. Since local systems are fiber bundles, a homotopy of paths pulls back to a trivial bundle, which shows that the map $\mu_{\gamma}$ is invariant under homotopy classes of paths rel endpoints.

Moreover, one can construct a space associated to a functor $\cL:\pijuan(X)\to\Vect$ by considering the product $L:=\prod_{x\in X}\cL(x)$ and topologizing suitably. For example, one could consider a basis of open sets around a point $v\in \cL(x)$ given by the collection of elements $\{w\in L\,|\,\exists\gamma, s.t. \mu_{\gamma}(v)=w\}$. This construction mirrors the usual construction of a classifying space given in Hatcher~\cite[Sec. 1.3]{hatcher} or Munkres~\cite[Ch. 13]{munkres}.
\end{proof}

\begin{figure}
\centering
\includegraphics[width=.6\textwidth]{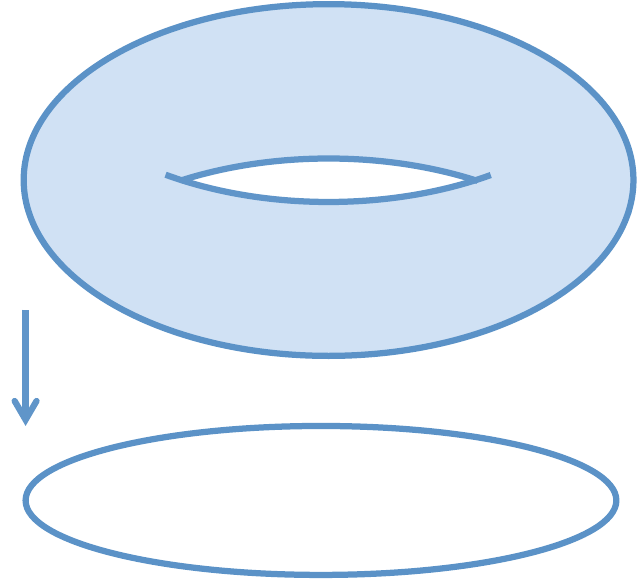}
\caption{Trivial Circle Bundle over the Circle}
\label{fig:circle_bundle}
\end{figure}

This equivalence allows us to define plenty of examples of local systems coming from fiber bundles.
\begin{prop}[Fiber Bundles Give Local Systems]\label{prop:bundles_locsys}\index{fiber byndle!gives local systems}\index{local system!from a fiber bundle}
	Suppose $\pi:E\to X$ is a fiber bundle, then for each $i$ the homology of the fiber $H_i(\pi^{-1}(x);k)$ defines a local system. Dually, for each $i$ the cohomology of the fiber $H^i(\pi^{-1}(x);k)$ defines a representation of the fundamental groupoid and consequently a local system.
\end{prop}
\begin{proof}
	This is easily seen because any path $\gamma:[0,1]\to X$ determines a pullback bundle $\gamma^*E$ which is trivial, so there is an isomorphism $H_i(\pi^{-1}(\gamma(0));k)\to H_i(\pi^{-1}(\gamma(1));k)$. Moreover, any homotopy of paths $H:[0,1]^2\to X$ determines a trivial pullback bundle $H^*E$.
\end{proof}

\begin{figure}
\centering
\includegraphics[width=\textwidth]{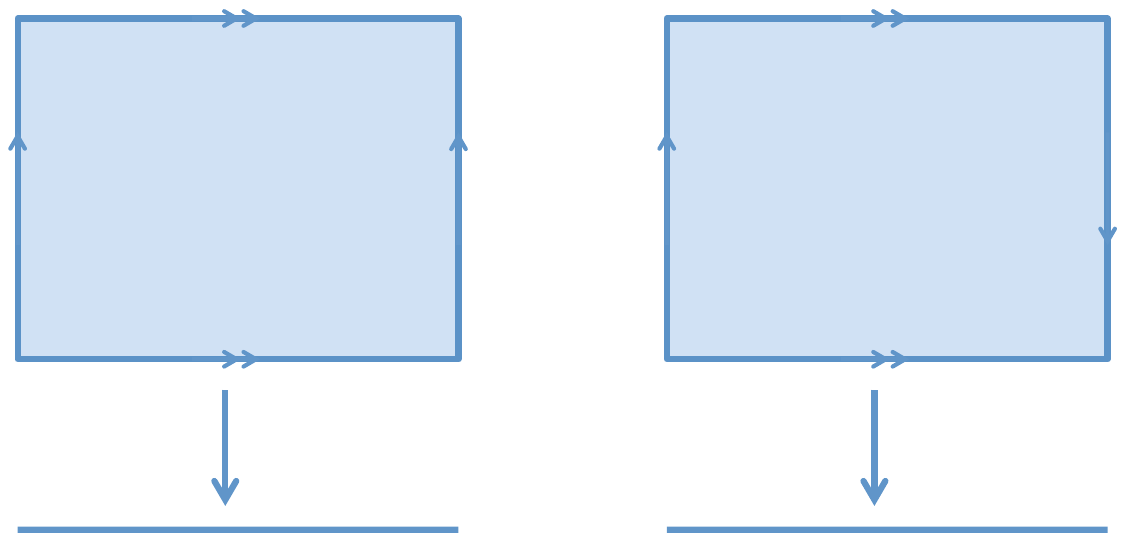}
\caption{Identification Spaces for the Torus and Klein Bottle}
\label{fig:torus_klein_id}
\end{figure}

\index{local system! Klein bottle and torus examples}
Let us now consider two examples. Firstly, in Figure~\ref{fig:circle_bundle} we drew a map from the torus to a circle. To define this map one considers the torus as the space $S^1\times S^1$ and defines the map to be projection onto the first factor. Secondly, consider the analogous projection map from the Klein bottle to the circle. An identification space model is depicted for both of these maps are drawn in Figure~\ref{fig:torus_klein_id} with the left hand side being the torus with its map and the right hand side being the Klein bottle with its map.

\begin{figure}
\centering
\includegraphics[width=\textwidth]{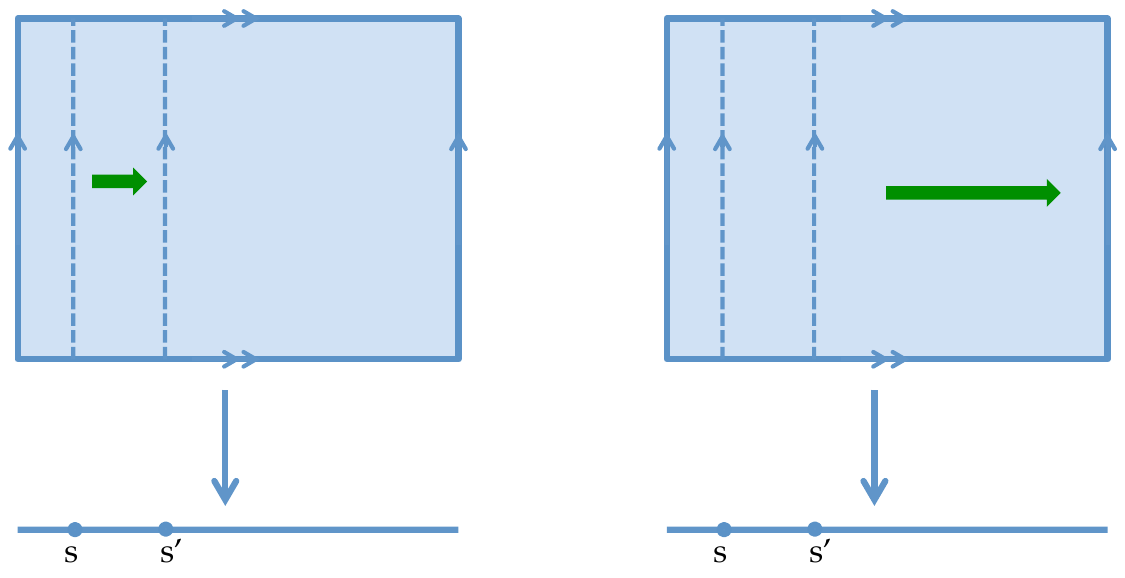}
\caption{Trivial Action with the Torus Map}
\label{fig:torus_action}
\end{figure}

Consider the local system gotten by taking $H_1(-;k)$ of the fiber $\pi_T:T\to S^1$. Between any two points $s$ and $s'$ there are two homotopy classes of paths connecting them: one that in the identification space model proceeds directly from $s$ to $s'$ and one that wraps around the circle using the implied identification. Choosing a basis for the vector space $H_1(\pi_T^{-1}(-);k)$ involves choosing a cycle along with an orientation. If one considers the monodromy map associated to either path, one can see that in the bases indicated for $H_1(\pi_T^{-1}(s);k)$ and $H_1(\pi_T^{-1}(s');k)$ in Figure~\ref{fig:torus_action} both monodromies are trivial (i.e. the identity map $k\to k$) as indicated by the green arrows.

\begin{figure}
\centering
\includegraphics[width=\textwidth]{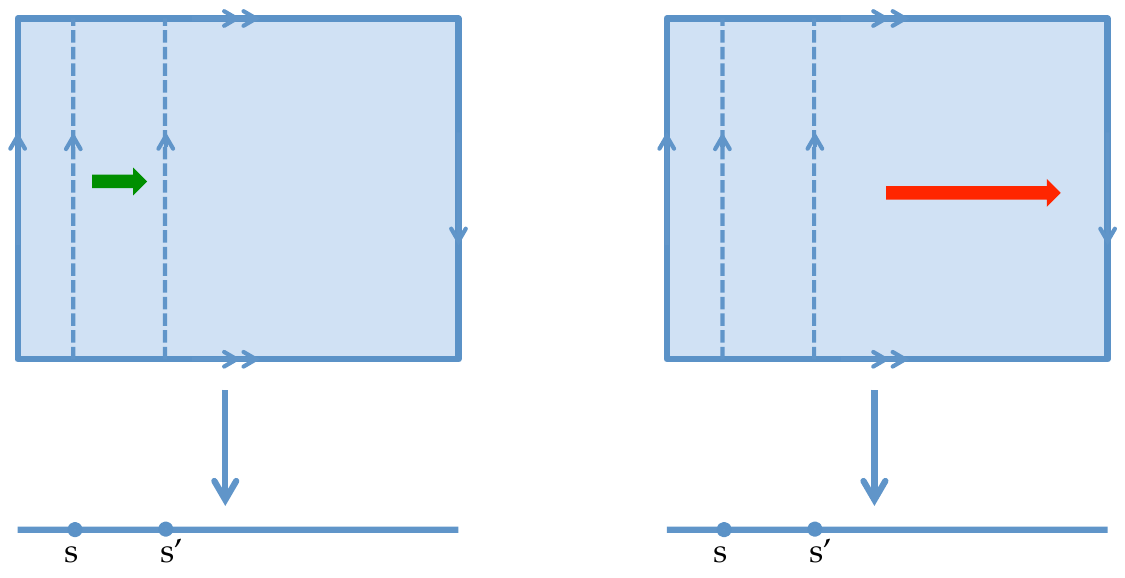}
\caption{Non-trivial Action from the Klein Bottle}
\label{fig:klein_action}
\end{figure}

Now consider the local system gotten by taking $H_1(-;k)$ of the fiber for $\pi_K:K\to S^1$. Choosing the same bases as before the monodromy associated to the longer path that wraps around the identification space is non-trivial
\[
\xymatrix{ H_1(\pi_K^{-1}(s');k) = k \ar[r]^{-1} & k = H_1(\pi_K^{-1}(s);k)}
\]
as indicated by the red arrows in Figure~\ref{fig:klein_action}.

We now will show that these examples actually provide examples of locally constant cosheaves, or sheaves if cohomology is taken. First, we will need some alternative definitions.

\begin{defn}\index{sheaf!locally constant sheaf}\index{cosheaf!locally constant cosheaf}
	Let $X$ be a locally connected space. A sheaf $A_X$ or cosheaf $\hat{A}_X$ on $X$ valued in $\Vect$ is \textbf{constant} with value $A$ if for every open set $U$ they make the following assignments:
	\[
	A_X: U \rightsquigarrow A^{\times \pi_0(U)} \qquad \hat{A}_X: U \rightsquigarrow A^{\oplus \pi_0(U)}.
	\]
	A sheaf $F$ or cosheaf $\cF$ is \textbf{locally constant} if for each point $x$ there is an open neighborhood $U$ such that $F$ or $\cF$ is constant, i.e. there is a vector space $A$ such that $F|_U\cong A_X$ or $\cF|_U\cong \hat{A}_X$.
\end{defn}

As a consequence of this definition and the topological assumptions on $X$, a locally constant sheaf or cosheaf possesses for each point $x\in X$ a collection of connected neighborhoods containing $x$ all of which take identical values. As a consequence $F(U)\to F_x$ or $\cF_x \to \cF(U)$ is an isomorphism. Moreover, for any other point $x'$ contained in $U$, the stalk or costalk at $x'$ can be chosen to be isomorphic to $F(U)$ or $\cF(U)$ respectively. By chaining together these sorts of isomorphisms, one can show the following theorem:

\begin{thm}\label{thm:lc_sheaf_rep}\index{sheaf!locally constant sheaf!representation of $\pijuan(X)$}\index{cosheaf!locally constant cosheaf!representation of $\pijuan(X)$}\index{fundamental groupoid $\pijuan(X)$!and locally constant (co)sheaves}\index{representation!of fund@$\pijuan(X)$}
	Suppose $X$ is a locally path connected, locally simply-connected paracompact Hausdorff space. A locally constant sheaf determines a a local system, where a local system is defined to be a representation of the fundamental groupoid of $X$, i.e.
	\[
		\Loc: \pi_1(X) \to \Vect.
	\]
	Similarly, any locally constant cosheaf valued in $\Vect$ determines a local system.
\end{thm}
\begin{proof}
	By taking stalks or costalks we can define the functor $\Loc$ on objects $x\in X$ to be $F_x$ or $\cF_x$, respectively. Since the theorem is well known (see~\cite{achar-lsu} for a proof, which we follow here) for sheaves we present the cosheaf-theoretic proof instead.
	
	Call a subset $K$ of $X$ \emph{fine}\footnote{In~\cite{achar-lsu} they use the word ``good.''} for a cosheaf $\cF$ if it is connected and is contained in a connected open set $V$ such that $\cF|_V$ is a constant cosheaf. For any set of points $\{x_i\}$ in a fine set $K$ we have a collection of isomorphisms
	\[
	\xymatrix{ & \cF(V) & \\
	\cF_{x_i} \ar[ur]^{\pi_i} & \cF_{x_j} \ar[u]^{\pi_j} & \cF_{x_k} \ar[lu]_{\pi_k} }
	\]
	that when composed together allows us to define an invertible map from $\cF_{x_i}\to \cF_{x_k}$ via $\pi_{k}^{-1} \pi_i$. Of course this map agrees with the composition of the analogously defined map $$\cF_{x_i}\to \cF_{x_j} \to \cF_{x_k}$$ because $\pi_k^{-1}\pi_j\pi_j^{-1}\pi_i=\pi_{k}^{-1} \pi_i.$
	
	Now we claim that given a path $\gamma:[0,1]\to X$ there exists a sequence of points $0=a_0< a_1 < \cdots < a_n=1$ so that for all $i$ the set $\gamma([a_i,a_{i+1}])$ is fine. This is the case because every point $\gamma(t)$ possesses a fine neighborhood and by continuity there are open intervals $V_t$ of $t$ such that $\gamma(V_t)$ is fine. If we choose intervals $[a,a']$ contained in each $V_t$, the interiors of these intervals will form an open cover of $[0,1]$. By compactness, finitely many of these intervals will do. Choosing such a finite list, merging and ordering the endpoints, gives the requested sequence.
	
	From the sequence we can define the map $\rho(\gamma):\cF_{\gamma(0)}\to\cF_{\gamma(1)}$ to be the composite
	\[
		\cF_{\gamma(a_0)} \to \cF_{\gamma(a_1)} \to \cdots \to \cF_{\gamma(a_n)}.
	\]
	This map is well defined by virtue of the fact that it is invariant under the addition of extra points $a'$ to the sequence above. Consequently, if any different sequence was chosen we could have merged it with this one and deduced that these maps were the same.
	
	A similar argument can be used to show that for homotopies $H:[0,1]\times [0,1] \to X$ there are sequences $\{a_i\}_{i=1}^n$ and $\{b_j\}_{j=1}^m$ so that the sets $H([a_i,a_{i+1}]\times [b_{j},b_{j+1}])$ are fine. Using the same concatenation of isomorphisms proves that if $\gamma$ and $\gamma'$ are homotopic relative endpoints, then the above defined maps $\cF_{\gamma(0)}\to \cF_{\gamma(1)}$ and $\cF_{\gamma'(0)}\to \cF_{\gamma'(1)}$ are the same.
\end{proof}

Moreover, one can show that representations of the fundamental groupoid give rise to locally constant sheaves and cosheaves. This will require a slightly more sophisticated version of van Kampen's theorem found in~\cite{brown-vkt, may-cat}.

\begin{prop}[van Kampen's Theorem]\index{van Kampen theorem}
	Suppose $X$ is a locally connected topological space, and suppose $\covU=\{U_i\}$ is a cover of $X$ by path-connected open subsets, then the van Kampen theorem states that
	\[
	\pijuan(X)\cong \varinjlim_{I\in N(\covU)} \pijuan(U_I),
	\]
	i.e. the functor $\pi_1:\Open(X)\to\Grpd$ is a cosheaf for the cover $\covU$. However, since any cover is refined by its connected components, which are open by assuming local connectivity, the arguments of Section~\ref{subsec:refinement} imply that the fundamental groupoid is a cosheaf. 
\end{prop}

\begin{thm}\label{thm:rep_lc_sheaf}\index{fundamental groupoid $\pijuan(X)$!determines locally constant (co)sheaves}
	A representation of a fundamental groupoid $\Loc:\pijuan(X)\to\Vect$ determines a locally constant sheaf and a locally constant cosheaf respectively.
\end{thm}
\begin{proof}
	Again the sheaf theoretic version of this statement is well known (see~\cite{achar-lsu}), so we carry out the cosheaf version. Assume we have a local system $\Loc$, then we define the associated cosheaf to be
	\[
	\widehat{L}:U \rightsquigarrow H_0(U;\Loc):=\varinjlim_{x\in U} \Loc|_U,
	\]
	which is a cosheaf on account of the fact that colimits commute with colimits. The fact that $\widehat{L}$ is locally constant comes from the fact that each point $x$ has a simply connected neighborhood $U$ for which the local system $H_0(U;\Loc)\cong \Loc(x)$ for any $x\in U$.
	
	Although it is not pointed out in the literature, the classical proof for sheaves follows by making the exact dual assignment.
	\[
		L: U \rightsquigarrow H^0(U;\Loc):=\varprojlim_{x\in U} \Loc|_U
	\]
\end{proof}
\begin{rmk}[Alternative Proof]
	The introduction of an apparently superfluous $H_0(-;\Loc)$ is an invocation of the principle that $H_0$ is a cosheaf. This principle, expressed in Theorem \ref{thm:MV_cosheaf}, actually states that ``$H_0$ for any homology theory that satisfies Mayer-Vietoris is a cosheaf.'' This is true again for this case, but it requires that the reader know that local systems allow us to define a homology theory with ``twisted coefficients.'' This theory, which uses singular chains with coefficients determined by $\Loc$, satisfies the Eilenberg-Steenrod axioms~\cite[Ch.~6]{whitehead} and thus Mayer-Vietoris~\cite[Ch.~4.6]{spanier}. To complete our alternative proof of the above theorem, we check one more hypothesis of Theorem \ref{thm:MV_cosheaf}. We observe that for an upward increasing sequence of open sets $\{U_i\}$ we have a directed system of chain complexes of twisted singular chains whose right most end point looks like the (non-exact) sequence
	\[
		C_1(U_i;\Loc) \to C_0(U_i;\Loc) \to 0.
	\]
	Taking $H_0$ involves only taking a cokernel, which commutes with direct limits. This proves that $H_0(-;\Loc)$ is a cosheaf.
\end{rmk}

This establishes the following corollary, justifying that local systems do indeed form a bridge between locally constant sheaves and cosheaves.

\begin{cor}
With the above topological assumptions on $X$, the category of locally constant sheaves and locally constant cosheaves are equivalent.
\end{cor}

We leave the details to the reader, who should note that any local system $\cL$ defines a sheaf by taking $H^0(-;\cL)$ and a cosheaf by taking $H_0(-;\cL)$.

\section{Cosheaves Model Topology}

The omnipresence of sheaves in geometry and topology should come with no surprise to many researchers in the algebraic cousins of these fields. Remarkably, cosheaves are just as abundant, but this fact is less well appreciated. This might stem from a desire to avoid excessive terminology as very classical constructions in topology might be called cosheaves, but we will briefly reverse this wisdom to provide ourselves with lots of examples.

Perhaps the closest parallel to the sheaf of sections is the cosheaf of pre-images, but the presence of topology makes it a richer object of study.

\begin{defn}[Cosheaf of Pre-images]\index{cosheaf!of pre-images}
	Suppose $f:Y\to X$ is a continuous map. We can define the pre-cosheaf of topological spaces $\hF:\Open(X)\to\Top$ by assigning to an open subset the pre-image $f^{-1}(U)\in\Open(Y)$ endowed with the subspace topology, i.e.
	\[
	U \rightsquigarrow f^{-1}(U).
	\]
	Since colimits in the open set category are just unions and $f^{-1}(\cup_i U_i)=\cup_i f^{-1}(U_i)$, this defines a cosheaf.
\end{defn}

\begin{ex}[Feature Function]\index{cosheaf!of features}
	Suppose we have a topological space $X$, populated with features of interest, expressed as a function $P:\{1,\ldots, n\}\to X$. We get a cosheaf of sets via $\hF(U)=P^{-1}(U)$. A slightly different cosheaf is gotten by letting $\hG(U)=U\cap\im(P)$, which cannot distinguish points with identical images.
\end{ex}

In the case $n=1$ we can linearize this last example to define an example analogous to an example commonly encountered when studying sheaves.

\begin{figure}[ht]
\centering
\includegraphics[width=.6\textwidth]{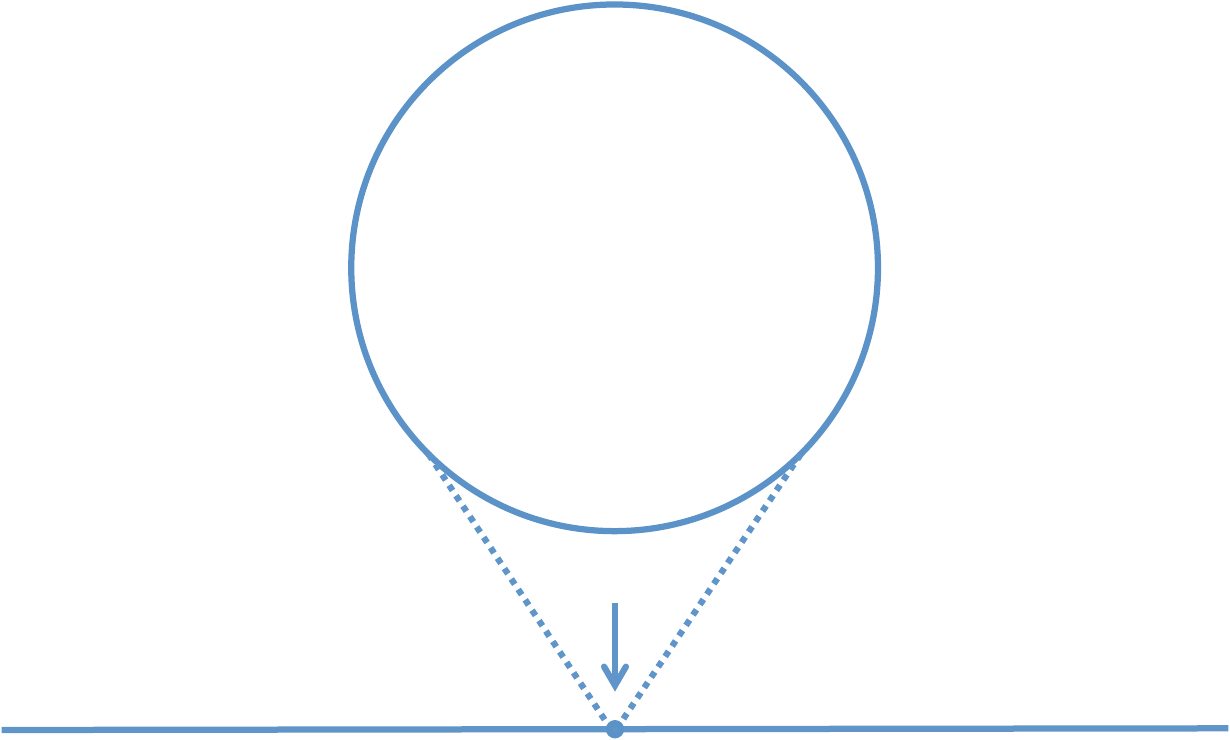}
\caption{Topological Model for Skyscraper Cosheaf}
\label{fig:skyscraper}
\end{figure}

\begin{defn}[Skyscraper Cosheaf]\index{cosheaf!skyscraper}
	Suppose $x\in X$ and $V$ is an $k$-vector space. Let's define the \textbf{skyscraper cosheaf} at $x$ with value $V$ to be
	\[
		\skycshv_x^V(U)=\left\{ \begin{array}{ll} V & \textrm{if $x\in U$}\\
		0 & \textrm{otherwise.}\end{array} \right.
	\]
	When $V=k$, we drop the superscript for notational convenience.

	A topological incarnation for the skyscraper is depicted in Figure~\ref{fig:skyscraper}. Here one makes the assignment
	\[
	U \rightsquigarrow H_0(f^{-1}(U);k)
	\]
	where $f$ is the map that maps the circle to the point $x$, i.e. the constant map with value $x$.
\end{defn}

\index{pointless topology}\index{locale}
We could adopt the perspective of cosheaves of pre-images as an alternative to continuous functions. This has been suggested in the past by John von Neumann and his derisively-named \textbf{pointless topology}, where in place of topological spaces one uses the poset of open sets as a primary notion --- an example of a \textbf{locale} --- and one observes that every continuous map of spaces $f:Y\to X$ induces a functor between categories $f^{\circ}:\Open(X)\to\Open(Y)$. This perspective will be of use later as we introduce operations on sheaves and cosheaves.

The cosheaf of pre-images will provide us with lots of examples of cosheaves pertinent to topology. However, viewing the entire information of the fiber (pre-image) is often too much to consider. Instead, one can consider invariants of the fiber and get a sometimes simpler, but still content-rich cosheaf (pending certain properties of the invariant).

\begin{defn}[Cosheaf of Connected Components]\index{cosheaf!of connected components}
	Given a continuous map of spaces $f:Y\to X$, one can define a pre-cosheaf of the components of the pre-image (not path components) $\hF:\Open(X)\to\Set$. This is done via the assignment
	\[
	U \rightsquigarrow \pi_0(f^{-1}(U)).
	\]
	This is not always a cosheaf. However, if $Y$ happens to be locally connected, i.e. the connected components of an open set are open, then it is. Alternatively, one can observe that the functor $\pi_0:\Top_{lc}\to\Set$ is left adjoint to the discrete space functor and so it preserves colimits~\cite{woolf}.
\end{defn}

\begin{ex}[The Square Map, Again]
Consider the cosheaf of connected components associated to the map $f:S^1\to S^1$ defined in complex coordinates as $f(z)=z^2$. For every point $p\in S^1$ there are two connected components in the fiber over $p$. However, there is only one connected component over the whole of $S^1$. This illustrates a sort of ``twisted'' $H_0$ already alluded to.
\end{ex}

\begin{exr}
Work out the cosheaf of connected components associated to the map $\pi:E\to X$ found in Figure~\ref{fig:evade_sec}.
\end{exr}

The next example provides a derived version of the principle that $H_0$ is a cosheaf.

\begin{ex}[Singular p-chains]\index{cosheaf!of singular chains}
	Fix $X$ a topological space and an open subset $U$. A singular p-chain on $U$ is nothing more than a $R$-linear combination of maps of the form $\sigma:\Delta^p\to U$. Since we can always post-compose a $p$-chain on $U$ with an inclusion $U\hookrightarrow V$, this defines a pre-cosheaf
	\[
	C_p(U)=\{\sum_{\sigma} \lambda_{\sigma}\sigma | \lambda_{\sigma}\in R, \sigma:\Delta^p\to U\}.
	\]
	This is, however, not a cosheaf as defined. Try writing down a chain on a union of two open sets as a linear combination of chains on the two sets.
	 A chain needs be sub-divided into pieces coming from each open set, each piece being represented as a map from a fixed simplex. As such, if we define
	\[
	\hat{C}_p(U):=\varinjlim C_p(U)
	\]
	where the colimit is being performed over iterated subdivision, then we obtain a cosheaf~\cite{Bredon}.
\end{ex}

\begin{figure}[ht]
\centering
\includegraphics[width=.4\textwidth]{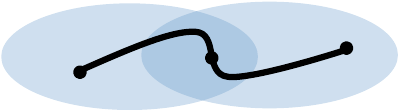}
\caption{Barycentric Subdivision of a Singular Chain}
\label{fig:sing}
\end{figure}

\begin{rmk}[Mayer-Vietoris and Cosheaves]\index{Mayer-Vietoris!as a cosheaf axiom}
	Another way of seeing that singular $p$-chains do not define a cosheaf is to recall that the proof of the \textbf{Mayer-Vietoris} theorem starts with the observation that the sequence
	\[
	0 \to C_p(U\cap V) \to C_p(U)\oplus C_p(V) \to C_p(U + V) \to 0
	\]
	is exact. Here the $C_p(U+V)$ is just notation for the cokernel of the previous map, thus the sequence is by definition exact. The elements of the cokernel  are linear combinations of singular chains strictly contained in either $U$ or $V$. One then uses barycentric subdivision to show that the \emph{complexes} $C_{\bullet}(U+V)$ and $C_{\bullet}(U \cup V)$ are chain homotopy equivalent. Letting $R=k$ be a field, this motivates defining a cosheaf valued in $\dat=K^b(\Vect_k)$ by assigning
	\[
	U \rightsquigarrow C_{\bullet}(U; k)
	\]
	and this will be a cosheaf.\footnote{The author has recently learned that Jacob Lurie calls this a \textbf{homotopy cosheaf}.} The category $K^b(\Vect_k)$ will be discussed later in the paper where it plays a more important role, but briefly stated it is the category whose objects are chain complexes of vector spaces of finite length and whose morphisms consist of equivalence classes of maps where we have identified those that are chain homotopic. This makes
	\[
		C_{\bullet}(U + V)\cong C_{\bullet}(U\cup V)
	\]
	thereby forcing the cosheaf axiom to hold. Of course the way this isomorphism is proven is via the use of barycentric subdivision, so we can avoid using cosheaves of chain complexes by working with the cosheaf $\hat{C}_p$ directly.
\end{rmk}

The cosheaves of singular chains serve a role precisely dual to the sheaves of co-chains commonly encountered in the literature. Consequently, homology is most naturally associated with cosheaf theory and cohomology is naturally associated with sheaf theory. However, there is a deeper duality between sheaves and cosheaves. When considering compactly-supported cohomology or closed (Borel-Moore) homology the natural habitats reverse. The kernel of this idea is present in the following example.

\begin{ex}[Compactly Supported Functions]\index{cosheaf!of compactly supported functions}
	Suppose $X$ is a locally compact Hausdorff space. Consider the following assignment:
	\[
	\Omega^0_c: U \rightsquigarrow \{f:U\to \RR \,|\, \mathrm{supp}(f) \; \mathrm{compact}\}
	\]
	Compactly supported functions defined locally can always be extended to larger open sets via extension by zero. If $X$ is a manifold, then we get more cosheaves of compactly supported differential $p$-forms $\Omega^p_c$ for $p\geq 0$.
\end{ex}

\section{Taming of the Sheaf... and Cosheaf}
\label{subsec:reeb}

As argued, the canonical example of a sheaf is the sheaf of sections of a map. This stands in contrast with the cosheaf of pre-images. However, a legitimate concern of both examples is its lack of computability. This concern is heightened given that the digital computer is becoming an increasingly common tool for modern mathematics.

A natural question might then be ``Can we store the sheaf of sections on a computer?'' Even in the example depicted in Figure \ref{fig:evade_sec}, it seems unlikely. On a small open set the sheaf of sections is in bijection with the set 
\[
\{f:(x-\epsilon,x+\epsilon)\to (a,b) \, |\, \mathrm{continuous}\},
\] 
which is uncountable. Moreover, for simple spaces like the closed unit interval with its Euclidean topology, there are uncountably many open sets that we need to assign data to.

To handle the first problem of ``too many sections'' in a somewhat ad hoc manner, we can conduct some pre-processing on the input data $\pi:E\to X$. As a motivating example, we can consider a construction normally defined when $X=\RR$.

\begin{defn}[Reeb Graph]\label{defn:reeb_graph}\index{Reeb graph}
	Suppose $Y$ is a topological space and $f:Y\to\RR$ is a continuous map. The \textbf{Reeb graph}~\cite{Reeb} is defined to be the quotient space $R(f):=Y/ \sim$ where $y\sim y'$ if and only if $y$ and $y'$ belong to the same connected component of the fiber $f^{-1}(t)$.
	\[
		\xymatrix{Y \ar[rr]^q \ar[rd]_{f} & & R(f) \ar[ld]^{\pi}\\
		& \RR &}
	\]
	Observe that $R(f)$ still possesses a map to $\RR$. There is clearly a direct generalization for arbitrary base spaces $X$.
\end{defn}

For an example of the Reeb graph, consider our zig-zag from Figure~\ref{fig:evade_sec}. Now let's work out what the sheaf of sections for $R(\pi)$ is and what the cosheaf of connected components is as well. Observe that we can probe the sheaf or the cosheaf on $[0,1]\subset\RR$ by asking what it assigns to open sets of the form $(x-\epsilon,x+\epsilon)$. Clearly it is constant except when the open set intersects a ``critical value.'' We express this observation by assigning values directly to cells in the visible decomposition of the codomain of the function. The data over incident edges and vertices are related, but the direction of that relation is dependent on whether we are considering a sheaf or a cosheaf. Making this observation rigorous has tremendous pay off because it allows us to avoid storing infinitely many open sets by instead working with finitely many cells.

\begin{figure}[ht]
\begin{minipage}[b]{0.47\textwidth}
\centering
\includegraphics[scale=.6]{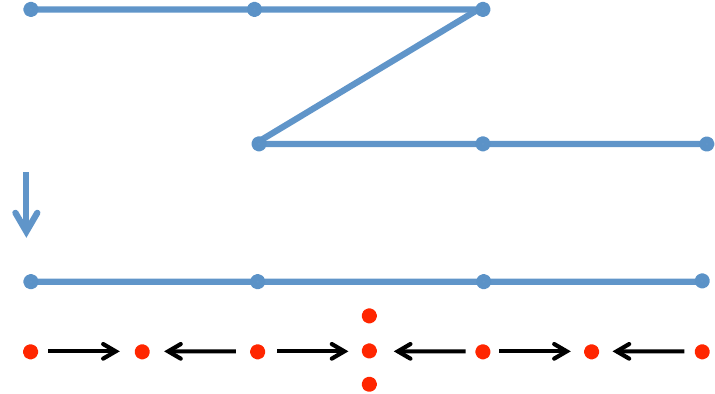}
\caption{Sheaf of Sections}
\label{fig:evade_sec_reeb}
\end{minipage}
\hfill
\begin{minipage}[b]{0.47\textwidth}
\centering
\includegraphics[scale=.6]{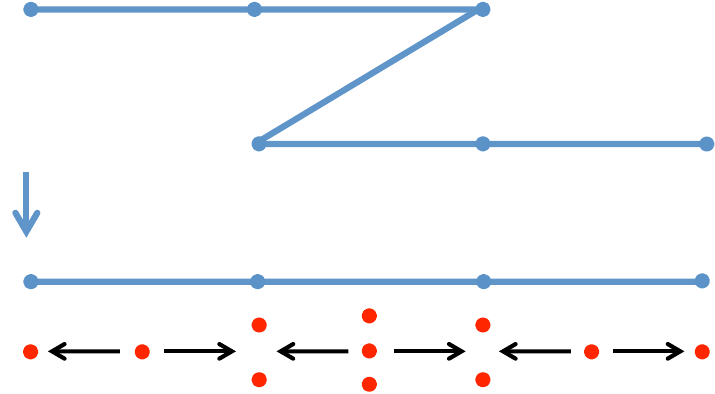}
\caption{Cosheaf of Components}
\label{fig:evade_comp}
\end{minipage}
\end{figure}

\ctparttext{In this part we emphasize that cellular sheaves and cosheaves are nothing more than linear algebra parametrized by a cell complex. The use of the term ``sheaf'' is justified by the Alexandrov topology, which makes functors out of posets into sheaves or cosheaves. Explicit proofs are presented since the primary reference of Shepard~\cite{shepard} is unpublished and not easily accessed. Cellular sheaf cohomology and cosheaf homology are presented computationally in Chapter~\ref{sec:homology} and then put on the firm foundation of derived categories in Chapter~\ref{sec:derived}. The novel contributions from this part, aside from working explicitly with cosheaves, are the introduction of ``barcodes'' to interpret cellular sheaf cohomology and cosheaf homology and exploiting the existence of enough projectives for cellular sheaves to define sheaf \emph{homology}.}
\part{Linear Algebra over Cell Complexes}
\label{part:easy_sheaves}

%
%

\chapter{Cellular Sheaves and Cosheaves}
\label{sec:cell_sheaves}

\begin{quote}
{\em``Sheaf theory is [where] you do topology horizontally and algebra vertically.''}
\begin{flushright} --- attributed to Maurice Auslander by~\cite{gray} \end{flushright}
\end{quote}

We can take it as an experimental observation from Figures \ref{fig:evade_sec_reeb} and \ref{fig:evade_comp} that in certain situations a sheaf or a cosheaf can be described as assigning data directly to the cells of a cell complex. Since cell complexes will be objects of primary importance to us, we review some definitions that may be non-standard.

\section{Cell Complexes and the Face-Relation Poset}

\begin{defn}[Regular Cell Complex~\cite{macpherson-seminar}]\index{regular cell complex}
	A {\bf regular cell complex} $X$ is a space equipped with a partition into pieces $\{X_{\sigma}\}_{\sigma\in P_X}$ such that the following properties are satisfied:
	\begin{enumerate}
		\item {\bf Locally Finite:} Each point $x\in X$ has an open neighborhood $U$ intersecting only finitely many $X_{\sigma}$.
		\item $X_{\sigma}$ is homeomorphic to $\RR^k$ for some $k$ (where $\RR^0$ is one point).
		\item {\bf Axiom of the Frontier}:\footnote{The \textbf{frontier} of a subspace $A$ is the complement of $A$ in its closure, i.e. $\fr(A):=\bar{A}-A$. In some forms this axiom reads: if $X_{\sigma}\neq X_{\tau}$ and $X_{\sigma}\cap\bar{X}_{\tau}\neq\emptyset$ then $X_{\sigma}$ is contained in the frontier of $X_{\tau}$.} If $\Bar{X}_{\tau}\cap X_{\sigma}$ is non-empty, then $X_{\sigma}\subseteq \Bar{X}_{\tau}$. When this occurs we say the pair are \textbf{incident} or that $X_{\sigma}$ is a face of $X_{\tau}$. The face relation makes the indexing set $P_X$ into a poset by declaring $\sigma\leq \tau$.
		\item The pair $X_{\sigma}\subset \Bar{X}_{\sigma}$ is homeomorphic to the pair $\mathrm{int}(B^{k})\subset B^k$, i.e. there is a homeomorphism from the closed ball $\varphi:B^k\to \Bar{X}_{\sigma}$ that sends the interior of the ball to $X_{\sigma}$.
	\end{enumerate}
\end{defn}

\begin{rmk}[Notation]
	Another common way of notating a cell complex is as a pair $(|X|,X)$ where $X$ is the set of cells and $|X|$ is the topological space being partitioned. To each cell $\sigma\in X$ there is a corresponding topological subspace $|\sigma|\subseteq |X|$. Our definition's notation says that $(X,P_X)$ is a cell complex. Our correspondence between cells and subspaces is $\sigma\rightsquigarrow X_{\sigma}$. However, we will have occasion to use both of these notations, and will sometimes use all three symbols $\sigma, |\sigma|$ and $X_{\sigma}$ to mean the same thing.
\end{rmk}

It is true that every regular cell complex can be further decomposed so that the resulting space is the homeomorphic image of a simplicial complex. However, for ease of computations we want to work with a class of spaces more general and natural than regular cell complexes. As such, we work with cell complexes, adopting the same convention as in Allen Shepard's thesis.

\begin{defn}[Cell Complex~\cite{shepard,macpherson-seminar}]\label{defn:face_poset}\index{cell complex}
 A \textbf{cell complex} is a space $X$ with a partition into pieces $\{X_{\sigma}\}$ that satisfies the first three axioms of a regular cell complex. Moreover, we require that when we take the one-point compactification of $X$, then the cells $\{X_{\sigma}\}\cup \{\infty\}$ are the cells of a regular cell complex structure on $X\cup\{\infty\}$.
\end{defn}

\begin{ex}
 The open interval $(0,1)$ decomposed with only one open cell is \emph{not} a cell complex. Its one-point compactification is the circle decomposed with one vertex $\{\infty\}$ and one edge $(0,1)$ whose attaching map is not an embedding, thus contradicting the fourth axiom.
\end{ex}

\begin{defn}[Cell category]\index{category!of a cell complex $\Cell(X)$}
 To a cell complex $(X,\{X_{\sigma}\}_{\sigma\in P_X})$ we can associate a category $\Cell(X;\{X_{\sigma}\})$, which is the indexing poset $P_X$ viewed as a category. This means that there is one object $\sigma$ for each $X_{\sigma}$ and a unique morphism $\sigma\to\tau$ for each incident pair $X_{\sigma}\subseteq \bar{X}_{\tau}$. When there is no risk of confusion, or a cell structure is specified at the beginning, then we will suppress the extra notation and just use $\Cell(X)$ or $X$.
\end{defn}

We now introduce diagrams indexed by the cell category. These were defined in Shepard's 1985 thesis~\cite[p.~6]{shepard}, but were known as ``stacks''\index{stacks} in the first published volume of Zeeman's 1954 thesis~\cite[p.~626]{dihom_1}. However, this term came to be used for an entirely different construction in the landmark paper of Pierre Deligne and David Mumford~\cite{DM-stacks}, which introduced the fundamental concept of algebraic and moduli stacks for algebraic geometry. Zeeman's usage of the term is now extinct, but his work anticipates MacPherson's cellular perverse sheaves (cf. Definition~\ref{defn:cellular_perverse_sheaves}) although MacPherson was unaware~\cite{macpherson} of the content of Zeeman's thesis~\cite{dihom_1,dihom_2,dihom_3}. To keep the presentation simple, we give Shepard's definition and its appropriate dualization.

\begin{defn}[Cellular Sheaves and Cosheaves]\label{defn:cell_sheaves}\index{cellular sheaves and cosheaves}\index{sheaf!cellular}\index{cosheaf!cellular}
	A \textbf{cellular sheaf} $F$ valued in $\dat$ on $X$ is a functor $F:\Cell(X)\to\dat$, i.e. it is 
	\begin{itemize} 
		\item an assignment to each cell $X_{\sigma}$ in $X$ an object $F(\sigma)$,
		\item and to every pair of incident cells $X_{\sigma}\subset \bar{X}_{\tau}$ a \textbf{restriction map}\footnote{Shepard calls these co-restriction maps since they point from faces to co-faces, but we will see they are restriction maps in the Alexandrov topology.} $\rho^F_{\sigma,\tau}:F(\sigma)\to F(\tau)$.
	\end{itemize}
	
	Dually, a \textbf{cellular cosheaf} $\hF$ valued in $\dat$ on $X$ is a functor $\hF:\Cell(X)^{op}\to\dat$, i.e. an assignment of an object $\hF(\sigma)$ for each cell, and an \textbf{extension map} $r_{\sigma,\tau}:\hF(\tau)\to\hF(\sigma)$ for every pair of incident cells $X_{\sigma}\subset\bar{X}_{\tau}$.
\end{defn}

Let us consider a few natural examples.

\begin{ex}[Bott's Torus]\index{cosheaf!cellular!height function on torus}
 The following example was first popularized by Raoul Bott in his book on Morse theory~\cite{macpherson}. Consider the height function on the torus, rotated by $90^{circ}$ so that the real line is underneath the torus, as shown in Figure~\ref{fig:bott_torus}.
\begin{figure}
  \centering
\includegraphics[width=.7\textwidth]{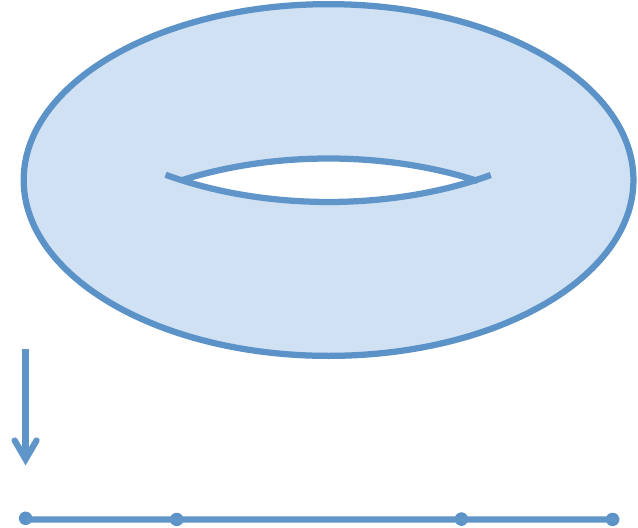}
\caption{Bott's Height Function on the Torus}
\label{fig:bott_torus}
 \end{figure}
 By taking the pre-image of the star of each cell, one obtains a diagram of spaces $\hF:X^{op}\to\Top$. By post-composing this diagram with $H_1(-;k)$, one obtains the cellular cosheaf indicated in Figure~\ref{fig:bott_torus_cosheaf}.
  \begin{figure}[ht]
  \centering
\includegraphics[width=.7\textwidth]{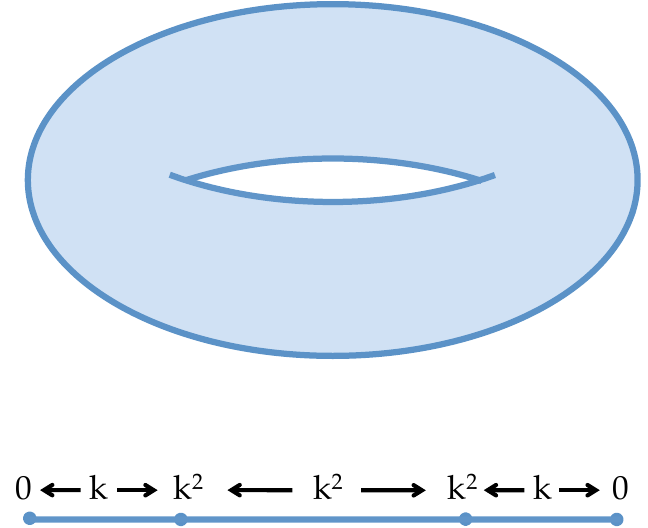}
\caption{Cellular Cosheaf by Taking $H_1$ of the Star pre-images}
\label{fig:bott_torus_cosheaf}
 \end{figure}
\end{ex}

\begin{ex}[Klein Bottle Revisited]\index{cosheaf!cellular!Klein bottle}
 As seen in Section~\ref{subsec:local_systems}, a Klein bottle can be viewed as a non-trivial $S^1$ bundle over the circle. We've already seen how the this leads to a representation of the fundamental groupoid $\pijuan(S^1)$. We can concoct a cellular cosheaf that describes this bundle in a different way. Let $s$ and $s'$ denote two vertices in a cell decomposition which includes two edges $a$ and $b$, as in Figure~\ref{fig:klein_action}. We can imagine calling the edge $a$ the short edge between $s$ and $s'$ and let $b$ be the long edge. To each cell $\hF$ assigns the homology of the fiber to that cell. The actions are encoded using maps between the cells. This gives us a diagram of vector spaces in the shape of a cell structure on $S^1$:
 \[
  \xymatrix{ & \ar[ld]_{-1} \hF(b) \ar[rd]^1 & \\
  \hF(s) & & \hF(s') \\
  & \hF(b) \ar[ul]^1 \ar[ur]_1 &}
 \]
\end{ex}

\begin{ex}
Let $Y=(0,1)$ be the open unit interval in $\RR$. Denote the inclusion of $Y$ into $\RR$ by $j$. To the constant sheaf on $Y$, written $k_Y$, we can associate two sheaves on $\RR$: $j_*k_Y$ and $j_! k_Y$. To describe these sheaves, we can think cellularly. For the first sheaf $j_*k_Y$ the cellular sheaf is simply the diagram of vector spaces
\[
0 \leftarrow k \rightarrow k \leftarrow k \rightarrow 0.
\]
For $j_!k_Y$ the cellular sheaf is the diagram of vector spaces
\[
0 \leftarrow 0 \rightarrow k \leftarrow 0 \rightarrow 0.
\]
The key difference being that the latter sheaf is zero on the endpoints $\{0\}$ and $\{1\}$. However, one can recover the classical, open set description of the sheaves $j_*k_Y$ and $j_! k_Y$ by considering any ordinary (Hausdorff) open set on the real line and then computing the limit of the diagram that lies over the open set. In Figure \ref{fig:star_shriek} we have drawn the two cellular sheaves of interest and the value of the limit over each open set.
	\begin{figure}
	\centering
	\includegraphics[width=\textwidth]{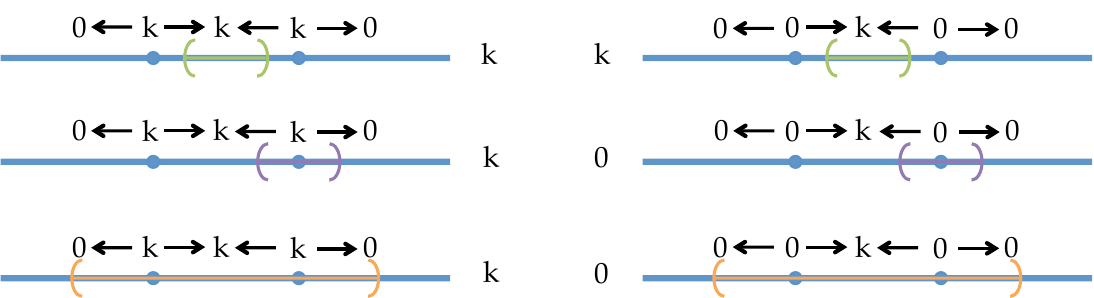}
	\caption{Cellular Description of $j_*k_Y$ and $j_!k_Y$ respectively}
	\label{fig:star_shriek}
	\end{figure}
Of course, one can dualize the discussion and consider cosheaves instead and use colimits to get functors from the open set category as defined in Chapter~\ref{sec:abstract_sheaves}. This perspective is developed more fully in Section~\ref{subsec:equivalence_cosheaves}.
\end{ex}

Since functors between categories assemble themselves into a category of their own, we get categories of cellular sheaves and cosheaves.
\begin{defn}\label{defn:cat_cell_sheaves}\index{category!of cellular sheaves $\Shv(X)$, cosheaves $\Coshv(X)$}
 We denote the category of cellular sheaves on $X$ by
\[
 \Shv(X;\dat):=\Fun(\Cell(X),\dat)
\]
and the category of cellular cosheaves by
\[
 \Coshv(X;\dat):=\Fun(\Cell(X)^{op},\dat).
\]
Morphisms are natural transformations of functors. If $\dat=\Vect$, then we will omit the notation after the semicolon and write $\Shv(X)$ and $\Coshv(X)$ instead.
\end{defn}

The notation deliberately coincides with the notation used for categories of sheaves and cosheaves on an arbitrary topological space, i.e. functors out of the open set category that satisfy the appropriate axiom. This conflict will be resolved in Section \ref{subsec:posets} in one way, and in Chapter \ref{subsec:strat} in an entirely different way.

\section{Partially Ordered Sets: Finite Spaces and Functors}
\label{subsec:posets}

Cellular sheaves and cosheaves earn their finiteness by assigning data directly to cells, rather than open sets. This turns out to not be entirely true; Cellular sheaves and cosheaves are simply operating on a different topology than the one we are accustomed to. Partially ordered sets can be endowed with a topology making cellular sheaves and cosheaves into actual sheaves and cosheaves on this topology.  

Here one can illuminate all of the general machinery of classical sheaf theory, but with a combinatorial finiteness that bends the theory to direct computation and understanding. Some of the explicit treatment of sheaves on posets is contained in the clear and concise work of Sefi Ladkani~\cite{sl-dereq}, but we streamline the discussion by using Kan extensions, which clarifies how cosheaves on a poset $X$ differ from sheaves on $X^{op}$.

\subsection{The Alexandrov Topology}
\label{subsubsec:alex}

In this section we introduce a class of non-Hausdorff spaces called \textbf{Alexandrov spaces}. The reader should note that although this topology is non-Hausdorff, it is highly relevant to concepts in algebraic topology. There is a remarkable theorem due to Michael McCord~\cite{mccord-finite} that states that every finite simplicial complex is weakly homotopy equivalence to an Alexandrov space. Thus, if one is interested in the topological properties of simplicial complexes, one should care about (non-Hausdorff) Alexandrov spaces. McCord even gives constructions of classical operations in algebraic topology, including suspension, in the Alexandrov setting. However, our ambitions for this section are far more limited. Let us begin with the necessary definitions.

\begin{defn}\index{preordered set}\index{poset}
	A \textbf{pre-order} consists of a set $P$ and a relation $\leq$ that is reflexive and transitive. A \textbf{poset} is a pre-order where the relation is also anti-symmetric, i.e. $x\leq y$ and $y\leq x$ implies $x=y$. A map $f$ of pre-orders is one that respects $\leq$. That is if $x\leq y$ then $f(x)\leq f(y)$. Pre-orders and order preserving maps form a category $\Preorder$. The collection of all posets form a subcategory of this category.
\end{defn}

Every pre-order can be equipped with a topology. However, it was first defined for finite posets by Pavel Alexandrov~\cite{alex-dr, alex-ct} and the general definition carries his name.

\begin{defn}\index{Alexandrov topology}
	On a pre-order $(P,\leq)$ define the \textbf{Alexandrov topology} to be the topology whose open sets are the sets that satisfy the following property:
	\[
		x\in U \qquad x\leq y \qquad \Rightarrow \qquad y\in U
	\]
	A basis is given by the sets of the form $U_x:=\{y\in P | x\leq y\}$ --- what we will call the \textbf{open star at x}. Similarly, we define the \textbf{closure of x} by $\bar{x}:=\{y\in P | y\leq x\}$. When $P$ is a finite poset, then a basis of closed sets is given by the $\bar{x}$'s.
\end{defn}

Any pre-order $P$ has an associated poset. This poset is gotten by defining an equivalence relation on $P$ via $x\sim y$ if and only if $x\leq y$ and $y\leq x$. One can check that this surjection is order-preserving. This construction defines a right adjoint to the inclusion of posets into pre-orders~\cite{woolf}.

\begin{rmk}[$P$ will mean a poset]
	Although spaces equipped with a pre-order are an interesting class of structures to consider, we will now work exclusively with posets. We do this to prevent closed loops from occurring in chains of related elements, as this would complicate our story.
\end{rmk}

\begin{ex}\label{ex:closed_gen}
	Consider $(\RR,\leq)$ with the usual partial order. The open sets are all those open or half open intervals such that the right-hand endpoint is $+\infty$. Observe that the closed set $(-\infty,0)$ cannot be written as an intersection of closed sets of the form $\bar{t}$. Thus the closures at $t$ do not form a basis.
\end{ex}

The dictionary between cellular complexes and Alexandrov spaces is easily described. First we introduce another definition.

\begin{defn}[Star]\index{star}
	Let $(X,\{X_{\sigma}\})_{\sigma\in P_X}$ be a cell complex. Every cell $X_{\sigma}$ has a \textbf{star}, which is a set that consists of all those cells $X_{\tau}$ such that $X_{\sigma}\leq X_{\tau}$. 
	\[
		\st(X_{\sigma}):=\{X_{\tau}\,|\,X_{\sigma}\leq X_{\tau} \}
	\]
	Since this definition only depends on the incidence relation of cells, we often drop the distinction between $X_{\sigma}$ and its label $\sigma$. Thus the star is also described as a subset of the poset $P_X$ consisting of those labels $\tau$ such that $\sigma\leq\tau$.
\end{defn}

The Alexandrov topology on the indexing poset $P_X$ of a cell complex allows us to define a continuous surjective map that comes from sending each cell $X_{\sigma}$ to its label $\sigma$.
This continuous surjective map gives an alternative way of describing how the Alexandrov topology arises. It is the quotient space where we identify two points $x$ and $y$ if and only if they belong to the same cell.
\[
	\xymatrix{X \ar[d]^{q} \\ P_X:=X/\sim}
\]
The inverse image of the star of $\sigma$ is an open union of cells, which is open. Thus this map is continuous and the topology that makes this map continuous is the Alexandrov topology.

\begin{figure}[ht]
	\centering
	\includegraphics[width=.5\textwidth]{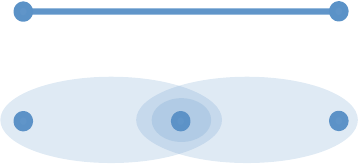}
	\caption{Alexandrov Space Associated to the Unit Interval}
	\label{fig:alex_interval}
\end{figure}

\begin{ex}[The Interval]\index{Alexandrov topology!unit interval example}
	Suppose $X=[0,1]$ is the unit interval given a cell complex structure with two vertices and one open interval. The face relation poset $P_X$ takes the following form:
	\[
		\xymatrix{ & \bullet & \\ \bullet \ar[ru] & & \bullet \ar[lu]}
	\]
	The Alexandrov topology has basic open sets corresponding to the star of each cell. The stars of the two vertices intersect each other. In Figure \ref{fig:alex_interval}, we have drawn the basic open sets.
\end{ex}

\subsection{Functors on Posets}
\label{subsubsec:fun_kan}
We want to understand how data modeled on posets can be treated as a sheaf or cosheaf on the Alexandrov topology. To do so we use the elegant, but sophisticated, approach of Kan extensions. To motivate this concept we will consider the relationship between a poset and its topology.

Observe that the correspondence between the relation internal to the poset $P$ and the containment relation for the open sets in the Alexandrov topology is order-reversing. Said more succinctly, we have an inclusion functor that is contravariant, i.e.
\[
\iota: P\to \Open(P)^{op} \qquad p \mapsto U_p.
\]
A natural question to ask is 
\begin{quote}\index{Kan extension!motivation}
	``Given a functor $F:P\to\dat$, is there a consistent way of extending $F$ to a functor $R:\Open(P)^{op}\to\dat$?''
\end{quote} 
One can hope to perform this extension since the image of the inclusion $\iota:P\to\Open(P)^{op}$ is a basis for the topology. Consequently, we can express arbitrary open sets as unions (colimits or limits in the opposite category) of basic open sets $\iota(p)=U_p$. A candidate extension would be to define 
\[
F(U):=\varprojlim_{U_p\subset U} F(p)
\] 
or as the colimit of $F$ over $U_p\subset U$. However, we should have some consistency. If one views $U_p=\{p'|p\leq p'\}$ as a subcategory of the category $P$, then it has an initial object $p$ and thus the limit of the diagram $F|_{U_p}$ is $F(p)$, i.e. 
\[
	\varprojlim_{p\leq p'} F(p') \cong F(p).
\]
This guides us to the following possible extension.
\[
	\xymatrix{P \ar[r]^F \ar[d]_{\iota} & \dat \\
	\Open(P)^{op} \ar[ur]_-{\varprojlim F(p)} & }
\]
This extension is nice for many reasons. By using limits to define data on larger open sets we have forced the sheaf axiom to hold, so this extension is in fact a sheaf. Moreover it illustrates through example a more general concept, which we now define.

\begin{rmk}[Caveat]
	We will make use of Kan extensions at a few points throughout the paper, but its immediate application is a theorem that says functors out of posets can be identified with sheaves. The proof of this theorem is described casually without the language of Kan extensions in~\cite{sl-dereq}, but adopting this language will be powerful and will make certain categorical properties transparent. 
\end{rmk}

\begin{defn}[Kan Extensions]\index{Kan extension}
	Suppose $\bat,\cat$ and $\dat$ are categories, $F:\bat\to\dat$ and $E:\bat\to\cat$ are functors, then the \textbf{right Kan extension of $F$ along $E$} written $R=\Ran_E F:\cat\to\dat$ is a functor and a natural transformation $\epsilon:RE\to F$ that is universal in the following sense. For every functor $H:\cat\to \dat$ with a natural transformation $\alpha:H\circ E\to F$ there exists a unique natural transformation $\sigma: H\to R$, i.e. $\Nat(H,R)\cong\Nat(H\circ E, F)$.
	\[
		\xymatrix{\bat \ar[r]^F \ar[d]_{E} & \dat \\
		\cat \ar[ur]_-{R=\Ran_E F} & }
	\]
	The \textbf{left Kan extension of F along E} written $L=\Lan_E F:\cat\to\dat$ is a functor with a natural transformation $\eta:F\to L\circ E$ that is universal as well. If $H:\cat\to\dat$ is a functor with a natural transformation $\omega: F\to H\circ E$, then there exists a unique $\tau:L\to H$, i.e. $\Nat(L,H)\cong \Nat(F,H\circ E)$.
	\[
		\xymatrix{\bat \ar[r]^F \ar[d]_{E} & \dat \\
		\cat \ar[ur]_-{L=\Lan_E F} & }
	\]
\end{defn}

\begin{rmk}[Existence of Kan Extensions]
	Kan extensions do not always exist, but we have already alluded to a situation where they do. If $\dat$ has enough limits and colimits, then we can give point-wise formulae for the left and right Kan extensions respectively:
	\[
		\Lan_E F (c):=\varinjlim_{E(b)\to c} F(b) \qquad \Ran_E F (c):=\varprojlim_{c\to E(b)} F(b)
	\]
\end{rmk}

One of the reasons that sheaves and cosheaves on Alexandrov spaces are so well-behaved is that every open set has a finest cover, so in particular, by Corollary \ref{cor:refined_axiom}, we only need to check the (co)sheaf axiom on this cover, and it will be guaranteed for all others. Furthermore, every point in an Alexandrov space has a smallest open neighborhood, and the (co)stalks are just the values on these minimal open sets. This is how we can use Kan extensions to create a dictionary between (co)sheaves on Alexandrov spaces and functors out of posets. 

\begin{thm}\label{thm:posets_sheaves}\index{sheaf!on a poset}\index{cosheaf!on a poset}
	Let $P$ be a poset and $\dat$ a category that is both complete and co-complete. Then the following categories are equivalent
	\[
		\Fun(P,\dat)\cong\Shv(P;\dat) \qquad \Fun(P^{op},\dat)\cong\Coshv(P;\dat)
	\]
\end{thm}
\begin{proof}
	We claim that taking the right Kan extension of $F:P\to\dat$ along the inclusion $\iota:P\to\Open(P)^{op}$ produces a sheaf. Suppose $U$ is an open set in the Alexandrov topology, i.e. one for which $p\in U$ and $p\leq p' \Rightarrow p'\in U$. It is true that every open set can be expressed as a union $U=\cup_{p\in U} U_p$ and thus the finest possible cover is $\{U_p\}_{p\in U}$. The right Kan extension then defines $F(U):=F[\{U_p\}_{p\in U}]$ so the sheaf axiom holds for that cover, but by Corollary \ref{cor:refined_axiom}, this means that $F$ is a sheaf. To go from a sheaf to the diagram, one simply takes stalks at every point. Since the smallest neighborhood containing $p$ is $U_p$, we get that $F_p=F(U_p)=F(p)$. 
	
	The dual argument for cosheaves is completely analogous: we take the left Kan extension of $\hF:P^{op}\to\dat$ along the inclusion $\iota:P^{op}\to \Open(P)$ to get a cosheaf. Taking costalks returns a diagram from a cosheaf.
\end{proof}

\begin{rmk}[Stalks and Costalks on Posets]
	To elaborate on the proof, let us compute some invariants. Recall that the stalk and costalk at a point $p\in P$ for a sheaf and cosheaf respectively is described via the use of filtered colimits and limits.
	\[
	F_p:=\varinjlim_{U\ni p} F(U) \qquad \mathrm{and} \qquad \hF_p:=\varprojlim_{U\ni p} \hF(U)
	\]
	In both cases when $P$ is a poset with the Alexandrov topology there is a smallest open set containing $p$, namely $U_p=\{q|p\leq q\}$, so $F_p\cong F(U_p)=F(p)$ and $\hF_p=\hF(U_p)=\hF(p)$.
\end{rmk}

\begin{defn}[Sections]\label{defn:sections}\index{sheaf!sections over a poset}\index{cosheaf!sections over a poset}
	Let $(P,\leq)$ be a poset and $F:P\to\dat$ a sheaf and $\hF:P^{op}\to\dat$. let $Z\subset P$ be any subset. We define the \textbf{sections over $Z$} to be
	\[
		\Gamma(Z;F):=\varprojlim F|_Z \qquad \mathrm{and} \qquad \varinjlim \hF|_Z=:\Gamma(Z;\hF).
	\]
	When $Z=P$, we call these \textbf{global sections}. Note that $\Gamma(Z;-)$ is context dependent: different definitions are used pending whether a sheaf or cosheaf is used.
\end{defn}

The above theorem provides the simplest explanation of why cellular sheaves and cosheaves deserve to be called sheaves and cosheaves. When Theorem \ref{thm:posets_sheaves} is specialized to the face relation poset $P_X$ of a cell complex, also called the cell category $P_X=\Cell(X)$ in Definition \ref{defn:face_poset}, we get that the category of sheaves in Definition \ref{defn:cat_cell_sheaves}. We summarize these observations in the following corollary. 

\begin{cor}
	Let $(X,P_X)$ be a cell complex. A \textbf{cellular sheaf on $X$} is a sheaf on $P_X$ equipped with the Alexandrov topology. Such a sheaf is uniquely determined by a functor $F:P_X\to\dat$. A \textbf{cellular cosheaf on $X$} is a cosheaf on $P_X$ with the Alexandrov topology. Such a cosheaf is uniquely determined by a functor $\hF:P_X^{op}\to\dat$.
\end{cor}

To close, we point out one of the symmetries that Alexandrov spaces possess.

\begin{clm}
	In the Alexandrov topology, arbitrary intersections of open sets are open and arbitrary unions of closed sets are closed. Thus, every Alexandrov space possesses a dual topology by exchanging open sets with closed sets.
\end{clm}
\index{cosheaf!on closed sets}\index{sheaf!on closed sets}
This observation would have pleased Leray. It demonstrates that one can also think of a functor $F:P\to\dat$ on a poset as either a sheaf or as a ``cosheaf on closed sets.'' What distinguishes these two though is whether we use limits or colimits to extend to larger sets. We will consider this perspective in greater detail in Section \ref{subsubsec:closed_push}.

%
%

\chapter{Functors Associated to Maps}
\label{sec:maps}

\begin{quote}
{\em``Qui s\`eme le foncteur r\'ecolte la structure.''}\footnote{``Who sows the functor reaps the structure.''}
\begin{flushright} --- Bourbaki \end{flushright}
\end{quote}

Since sheaves and cosheaves as defined here assign data to open sets, maps between spaces should only make reference to open sets. In Section~\ref{subsec:six_ops} we briefly introduced how to pushforward or pullback a sheaf along a map between spaces. In the case where our spaces are partially ordered sets endowed with the Alexandrov topology, it suffices to work directly with points since they are in bijection with a basis for the topology. However, playing these perspectives off of each other adds depth to the theory. In particular, by restricting our attention to these spaces, and using Kan extensions, we define the basic functors on (co)sheaves without making use of (co)sheafification. Pedagogically this is advantageous because the operation of sheafification tends to obfuscate the underlying ideas of sheaf theory. The lack of an explicit cosheafification process has historically been a stumbling block for the theory.

Recall that the definition of a continuous map $f:X\to Y$ says that the inverse image of an open set of $Y$ is an open set of $X$. This observation can be expressed by saying that we have a functor
\[
\mathring{f}:\Open(Y)\to\Open(X) \qquad U\subseteq Y \rightsquigarrow f^{-1}(U)\subseteq X.
\]
By formality, we also have a functor from the corresponding opposite categories
\[
\mathring{f}:\Open(Y)^{op}\to\Open(X)^{op}.
\]
We have purposely suppressed the $op$ superscript on $\mathring{f}$ for legibility. Now suppose we are given a pre-sheaf $G$ on $X$, then we get a naturally associated pre-sheaf on $Y$ by observing that the diagram
\[
\xymatrix{\Open(X)^{op} \ar[r]^-G & \dat \\
\Open(Y)^{op} \ar[u]^{\mathring{f}} \ar@{-->}[ru] & }
\]
has a natural completion given by pre-composition. 

\begin{defn}[Pushforward Sheaf and Cosheaf]\index{sheaf!pushforward}\index{cosheaf!pushforward}\index{direct image}
	Suppose $f:X\to Y$ is a continuous map of spaces, $G$ and $\hG$ a pre-sheaf and pre-cosheaf respectively, then define the \textbf{pushforward} or \textbf{direct image} pre-sheaf and pre-cosheaf via $f_*G(U):=G(f^{-1}(U))$ and $f_*\hG(U):=\hG(f^{-1}(U))$. Because $f^{-1}$ commutes with unions, we get that if $G$ or $\hG$ is a sheaf or cosheaf, then so is the pushforward. Moreover, this operation is functorial with respect to maps between (co)sheaves, so we get functors
	\[
	f_*:\Shv(X;\dat) \to \Shv(Y;\dat) \qquad f_*:\Coshv(X;\dat)\to\Coshv(Y;\dat).
	\]
\end{defn} 

There is also a pullback functor associated to a continuous map $f:X\to Y$, but it's construction is less obvious. Namely, if $G$ is a pre-sheaf on $Y$, then there is no clear way to define a pre-sheaf on $X$ because for an open set on $X$, $f(U)$ may not be open. The solution usually used is to take a system of approximations of $f(U)$ by open sets and to define the pullback sheaf as the limit of these approximations.
\[
	f^*G(U):=\varinjlim_{V\supset f(U)} F(V).
\]

Thinking categorically, the problem of ``approximation'' has been encountered before. Namely, how can we complete the following diagram?
\[
\xymatrix{\Open(Y)^{op} \ar[r]^-G \ar[d]_{\mathring{f}} & \dat \\
\Open(X)^{op} \ar@{-->}[ru]_? & }
\]
Again, by assuming that $\dat$ has sufficient colimits, we can could fill in the diagram by taking the left Kan extension of $G$ along $\mathring{f}$ and that will yield the candidate formula for the pullback just presented. Unfortunately, this definition for the pullback of a sheaf does not always define a sheaf. In sheaf theory over general spaces, this defect is circumvented by sheafification, as in Section~\ref{subsec:six_ops}. Fortunately, for the Alexandrov topology, this circumvention is unnecessary.

\section{Maps of Posets and Associated Functors}
\label{subsec:map_functors}

As already noted, sheaves and cosheaves on posets are easier to manipulate. The functors associated to a map of posets are explicitly defined without extra processing. Since posets can be made into a topological space the functors which exist for sheaves on general spaces can be studied here is a more tightly controlled laboratory. Moreover, since Alexandrov spaces have extra symmetries new functors not normally encountered exist here. 

\begin{defn}[Map of Posets]\index{poset!maps}
	Suppose $(X,\leq_X)$ and $(Y,\leq_Y)$ are posets. A map of posets is a map of sets $f:X\to Y$ that is order-preserving, i.e. if $x\leq_X x'$ then $f(x)\leq_Y f(x')$. Alternatively, since a poset can be viewed as a category, a map of posets is just a functor. When it is clear from context we will abbreviate $(X, \leq_X)$ by just $X$.
\end{defn}

\begin{rmk}[Notation and Cell Complexes]\label{rmk:cell_notation}
	In our effort to treat posets as spaces, we have used $X$ and $Y$ to denote partially ordered sets equipped with the Alexandrov topology. This might cause confusion since our canonical example of a poset will be the indexing poset of a cell complex $(X,P_X)$. Note that cell complexes consist of a pair of spaces, one is $X$, the Hausdorff space that is partitioned into pieces $X_{\sigma}$, the other is $P_X$, the poset of labels $\sigma$. From here on out we will work primarily with the poset $P_X$ as this is the combinatorial approximation to $X$. 
Thus, keeping in line with Shepard~\cite{shepard}, we change our notation from $(X,P_X)$ to $(|X|,X)$. Thus we have the following dictionary:
\begin{center}
\begin{tabular}{c |  c  c}
	 Dictionary & Old Notation & New Notation \\
	\hline
	Underlying Hausdorff Space & $X$ & $|X|$ \\
	Underlying Alexandrov Space & $P_X$ & $X$ \\
	Set of Points in a ``Cell'' & $X_{\sigma}$ & $|\sigma|$ \\
	``Cell'' viewed as a point & $\sigma$ & $\sigma$ \\
	Cellular Sheaf & $F:P_X\to\dat$ & $F:X\to\dat$
\end{tabular}
\end{center}
\end{rmk}

\subsection{Pullback or Inverse Image}
\label{subsubsec:pullback}

Recall that a sheaf on a poset $(Y,\leq_Y)$ is a functor $G:Y\to\dat$.  Similarly, a cosheaf is a functor $\hF:Y^{op}\to\dat$. For both structures the pull-back functor $f^*$ is easily described. It is the obvious pre-composition that completes the following diagram.
\[
\xymatrix{Y \ar[r]^-G & \dat \\
X \ar[u]^{f} \ar@{-->}[ru] & }
\]

\begin{defn}[Pullback for Poset Maps]\index{pullback!of a sheaf over a poset}
	Given a sheaf $G$ on $Y$ and a map of posets $f:X\to Y$, we can define the \textbf{pullback} or \textbf{inverse image} sheaf $f^*G$ on $X$ as follows:
	\begin{itemize}
		\item $f^*G(x)=G(f(x))$
		\item If $x\leq x'$, then let $\rho^{f^*G}_{x',x}=\rho^G_{f(x'),f(x)}$
		\item If $\eta: G\to H$ is a morphism in $\Shv(Y)$, i.e. a natural transformation of diagrams over $Y$, then $f^*\eta: f^*G\to f^*H$ is a morphism in $\Shv(X)$ defined by declaring $f^*\eta(x):f^*G(x)\to f^*H(x)$ to be equal to $\eta(f(x)):G(f(x))\to H(f(x))$.
	\end{itemize}
	
	The same definition and arguments go through for a cosheaf on $Y$ with suitable modification, i.e $r^{f^*G}_{x,x'}=r^G_{f(x),f(x')}$. Thus, we get functors
	\[
	f^*:\Shv(Y;\dat) \to \Shv(X;\dat) \qquad f^*:\Coshv(Y;\dat)\to\Coshv(X;\dat).
	\] 
\end{defn}

The definition of the pullback seems almost too good to be true, but one can check that the pre-sheaf description we outlined earlier agrees with this definition. Observe that if one applies that definition then
\[
	f^*F(U_x):=\varinjlim_{V \supset f(U_x)} F(V) \cong F(V_{f(x)})=F(f(x)),
\]
where we have used the fact that the smallest open set containing $f(U_x)=f(\{x' | x\leq x'\})$ is $V_{f(x)}=\{y|f(x)\leq y\}$.

\begin{ex}[Constant Sheaf and Cosheaf]\index{sheaf!constant sheaf!pullback from a point}
 Consider the constant map $p:X\to\star$. A sheaf $G$ on $\star$ consists of a single vector space $W$ and the identity morphism so we'll just call $G$ by the name $W$. We define the constant sheaf on $X$ with value $W$ to be $W_X:=p^*W$. One sees that it is a sheaf that assigns $W$ to every cell with all the restriction maps being the identity. Similarly, the constant cosheaf with value $W$ is $\hat{W}_X:=p^* W$.
\end{ex}

\subsection{Application: Subdivision}
\label{subsubsec:subdivide}
In the case where the poset is the face relation of a cell complex certain natural maps present themselves, such as subdivision.
	
\begin{defn}[\cite{shepard} 1.5, p.29]\index{subdivision}\index{cell complex!subdivision}
 A \textbf{subdivision} of a cell complex $X$ is a cell complex $X'$ with $|X'|=|X|$ and where every cell of $X$ is a union of cells of $X'$.
\end{defn}

Untangling the definition a bit we see that if $\sigma$ is a cell of $X$, then there is a collection of cells $\{\sigma_i'\}$ such that $\cup_i |\sigma_i'|=|\sigma|$. As such, we can define a surjective map of posets $s:X'\to X$ defined by making $s(\sigma')=\sigma$ if $|\sigma'|\subseteq |\sigma|$.

\begin{clm}
 Subdivision of a cell complex $X$ induces an order preserving map $s:X'\to X$ of the corresponding face-relation posets.
\end{clm}
\begin{proof}
 The ordering on $X'$ is given by the face relation. Suppose $\sigma'\leq \tau'$, then either $s(\sigma')=s(\tau')$ or not. If not, then $\sigma'$ and $\tau'$ belong to the subdivision of two cells $\sigma\leq \tau$.
\end{proof}

We are going to use this fact to define the subdivision of a sheaf in a cleaner manner than is found in~\cite{shepard}.

\begin{defn}\index{subdivision!of a cellular sheaf}\index{sheaf!cellular!subdivision}
 Suppose $F$ is a sheaf on $X$ and $s:X'\to X$ is a subdivision of $X$, then we define the subdivided sheaf $F':=s^*F$.
\end{defn}

\subsection{Pushforward or Direct Image}
\label{subsubsec:pushforward}

By adopting a point-theoretic picture rather than an open set-theoretic picture of sheaves and cosheaves over posets, we got an easy definition for the pullback functor. In the introduction we outlined a general definition for the pushforward functor $f_*$ on sheaves and cosheaves on an arbitrary topological space. Interestingly enough, although $f_*$ had a simple description using open sets, the point-level description requires thought.

\begin{defn}[Pushforward for Poset Maps]\index{sheaf!pushforward}\index{cosheaf!pushforward}\index{direct image}

	 Given a sheaf $F$ on $X$ and a map of posets $f:X\to Y$ we can define a sheaf on $Y$ as follows:
	\begin{itemize}
	 \item The \textbf{pushforward} of a sheaf is the right Kan extension of $F$ along $f$, i.e. $\Ran_f F$. \[f_*F(y)=\varprojlim_{f(x)\geq y} F(x)\]
	 \item Suppose $y\leq y'$, then $\{x|f(x)\geq y'\}\subseteq \{x|f(x)\geq y\}$. Any limit over the bigger set defines a cone over the smaller set by restriction, thus the universal property of limits guarantees the existence of a unique map $f_*F(y)\to f_*F(y')$ that we will call $\rho^{f_*F}_{y',y}$.
	 \item Suppose $\eta:F\to G$ is a map of sheaves, i.e. a natural transformation of diagrams over $X$. Then for any sub-poset $U$ of $X$, post-composing the limit over $U$ of $F$ with the arrows in the natural transformation defines a cone over $G$ restricted to $U$. By the universal property of limits there must be an induced map. $$\varprojlim_{x\in U} F \to \varprojlim_{x\in U} G$$
	\end{itemize}
	
	For cosheaves, the dual arguments go through with the slight modification that we use the left Kan extension along $f^{op}:X^{op}\to Y^{op}$.
	\[f_*\hF(y):= \varinjlim_{f(x)\geq y}\hF(x).\]
	Since both of these constructions are functorial, we have redefined two functors:
	\[
	f_*:\Shv(X;\dat) \to \Shv(Y;\dat) \qquad f_*:\Coshv(X;\dat)\to\Coshv(Y;\dat)
	\]
\end{defn}

\begin{ex}[Global Sections]\index{global sections!pushforward to a point}\index{pushforward to a point}
This functor is extremely useful as it gives us a way of defining the global sections of a sheaf or a cosheaf. For the constant map $p:X\to \star$ we offer the following definitions:
\[
	p_* F(\star)\cong F(X) = \Gamma(X;F) = H^0(X;F) \qquad p_*\hF(\star) \cong \hF(X) = \Gamma(X;\hF) = H_0(X;\hF)
\]
In Section \ref{sec:derived} we will use this definition as the prototype for defining ``higher'' pushforward or direct image functors.
\end{ex}

\subsection{$f_{\dd}$, Pushforwards and Closed Sets}
\label{subsubsec:closed_push}

One of the advantages of describing the standard functors of sheaf theory in the setting of posets is the presence of extra symmetries. Abstract definitions lend themselves to being dualized. In particular, in our point-theoretic definition of the pushforward we made use of Kan extensions, which come in two variants: left and right. In this section we consider the other variant and give a topological explanation for its origin.

\begin{rmk}[Caveat]\index{f@$f_{\dd}$ notation}
	The functor $f_{\dd}$ defined below is the left adjoint to $f^*$ (for cosheaves it will be the right adjoint). There seems to be a strong trend to call the left adjoint of $f^*$ by a different name: $f_!$. According to Joel Friedman (\cite{friedman} p. 22) the tradition goes back to Grothendieck in~\cite{sga4.1} SGA Expos\'e I, Proposition 5.1. The same notation is used by Ladkani~\cite{sl-dereq}, Lurie, Beilinson, Bernstein and others.

	This is unfortunate, since the notation $f_!$ is perhaps even more firmly established for the pushforward with compact supports functor used in classical sheaf theory. The reason seems to be that for general sheaves, there is no left adjoint to $f^*$, so it would be clear from context which was meant. However, for cellular sheaves, both functors exist and are useful.
\end{rmk}

\begin{defn}[Pushforward with Open Supports]\index{pushforward with open supports $f_{\dd}$}\index{sheaf!pushforward $f_{\dd}$}
 Given a sheaf $F$ on $X$ and a map of posets $f:X\to Y$ we can define a sheaf on $Y$ as follows:
\begin{itemize}
 \item The \textbf{pushforward with open supports} of a sheaf is the left Kan extension of $F$ along $f$, i.e. $\Lan_f F$. \[f_{\dd}F(y)=\varinjlim_{f(x)\leq y} F(x)\]
 \item If $y\leq y'$, then $\{x|f(x)\leq y\}\subseteq\{x|f(x)\leq y'\}$ and since any colimit over the bigger set defines a cocone over the smaller set by restriction, we get a unique map $\rho^{f_{\dd}F}_{y',y}:f_{\dd}F(y)\to f_{\dd}F(y')$.
 \item If we have a map of sheaves $\eta:F\to G$, then pre-composing the arrows for $\colim G$ with $\eta$ defines a co-cone over $F$. By universal properties we get an induced map \[\varinjlim_{x\in V} F\to \varinjlim_{x\in V} G.\]
\end{itemize}

Dually, for cosheaves we use the right Kan extension along $f^{op}$.
\[f_{\dd}\hF(y):= \varprojlim_{f(x)\leq y}\hF(x)\]
Both of these constructions are functorial and thus we have defined two functors:
\[
f_{\dd}:\Shv(X;\dat) \to \Shv(Y;\dat) \qquad f_{\dd}:\Coshv(X;\dat)\to\Coshv(Y;\dat)
\]
\end{defn}

This functor appears to be quite unusual, despite its naturality from the categorical perspective. To explain its topological origin, we revisit some of the original ideas of Alexandrov.

When Alexandrov first defined his topology he did two things differently:
\begin{enumerate}
	\item He only defined the topology for \emph{finite} posets.
	\item He defined the \emph{closed} sets to have the property that if $x\in V$ and $x'\leq x$, then $x'\in V$.
\end{enumerate}
\index{cosheaf!on closed sets}\index{sheaf!on closed sets}
Let us repeat the initial analysis of sheaves and diagrams indexed over posets, where we now put closed sets on equal footing with open sets. Observe that as before we have an inclusion functor:
\[
 j:(X,\leq) \to \Closed(X) \qquad x \rightsquigarrow \bar{x}:=\{x'| x'\leq x\}
\]
Consequently, we have a similar diagram for a functor $F:X\to\dat$ as before.
\[
	\xymatrix{X \ar[r]^F \ar[d]_{j} & \dat \\
	\Closed(X) \ar@{.>}[ur]_{?} & }
\]

If we choose the left Kan extension, we'd like to say the extended functor is a cosheaf on closed sets, i.e. use the definition of a cosheaf but replace open sets with closed sets. Unfortunately, this concept is not well defined for general topological spaces because the arbitrary union (colimit) of closed sets is not always closed. For Alexandrov spaces this property does hold and this illustrates one of the extra symmetries this theory possesses.

However, in order for the Kan extension to take a diagram and make it into a cosheaf, we need to know whether the image of the inclusion functor defines a basis for the closed sets. In Example \ref{ex:closed_gen} we showed that this is not always the case. The topology generated by the image of this functor is called the \textbf{specialization} topology and it suffers from certain technical deficiencies. In particular, order-preserving maps are not necessarily continuous in this topology, thus it fails to give a functorial theory. Fortunately, for finite posets these topologies agree and we can talk about cosheaves on closed sets without any trouble.

We now can give a topological explanation for the existence of the functor $f_{\dd}$. It is the functor analogous to ordinary pushforward where we have adopted closed sets as the indexing category for cosheaves and sheaves. If $f:X\to Y$ is a map of posets, then $f^c$ is the induced map between closed sets. The dagger pushforward is then the obvious completion of the diagram.
\[
\xymatrix{\Closed(X) \ar[r]^-F & \dat \\
\Closed(Y) \ar[u]^{f^c} \ar@{-->}[ru]_{f_{\dd}F} & }
\]

In Section \ref{subsec:sheaf_homology} this functor provides the foundation for defining \textbf{sheaf homology} and \textbf{cosheaf cohomology} --- theories that don't exist for general spaces.

\subsection{$f_!$: Pushforward with Compact Supports on Cell Complexes}
\label{subsubsec:shriek}

The three functors $f^*,f_*$ and $f_{\dagger}$ induced by a map of posets are well defined for any poset and any diagram. However, when Shepard wrote his thesis the only posets that he considered were posets coming from cell complexes. By working in this smaller class and imitating the theory of constructible sheaves, Shepard described another functor that is not defined for arbitrary Alexandrov spaces: the pushforward with compact supports $f_!$.

This fourth functor is meant to provide a cellular (constructible) analog of a functor naturally defined for sheaves on more general topological spaces and the name is borrowed from there. The reader must keep this in mind since every set in a finite Alexandrov space is compact. Thus, when we say ``pushforward with compact supports'' we mean a discrete model for the pushforward with compact supports functor defined for locally compact Hausdorff spaces.

Following Shepard, this functor $f_!$ will only be defined for cellular maps, which are stratified (or even definable) maps naturally adapted to cell complexes.

\begin{defn}[Cellular Map~\cite{shepard} pg. 32]\label{defn:cell_map}\index{cellular maps}
	Let $X$ and $Y$ be cell complexes. A \textbf{cellular map} $(|f|,f)$ consists of a map of posets $f:X\to Y$ and a continuous ``geometric'' map $|f|:|X|\to|Y|$ satisfying the following compatibility conditions:
	\begin{enumerate}
		\item For every $\sigma\in X$, $|f|(|\sigma|)$ is the cell $|f(\sigma)|$.
		\item The restricted map $|f|_{|\sigma|}:|\sigma|\to |f(\sigma)|$ is the projection $\RR^{n+k}\to\RR^n$ onto the first $n$ coordinates.
		\item Given $\sigma\in X$ and $y,z\in |f(\sigma)|$, $|f|^{-1}(y)\cap \bar{|\sigma|}$ is compact if and only if $|f|^{-1}(z)\cap \bar{|\sigma|}$ is.
	\end{enumerate}
\end{defn}

\begin{rmk}
The first and second conditions clearly restrict the types of maps of posets that can be considered. It appears that the third condition is redundant given the first two, but this is how it is recorded in Shepard's thesis.
\end{rmk}
 
\begin{ex}
Let $X=[0,1)$ be given the cell structure $x=0$ and $b=(0,1)$. Let $Y=[0,1)\times [0,1)$ be given the simplest possible cell structure. The underlying posets for these spaces are as follows:
\[
  \xymatrix{ & b\\x \ar[ur] &} \qquad \xymatrix{ & \sigma & \\ a \ar[ur] & & b \ar[ul]\\ & x \ar[ul] \ar[uu] \ar[ur] &}
\]
Here $x$ refers to the vertex, $a$ and $b$ the open edges, and $\sigma$ is the open face. Clearly, $f(x)=a$ and $f(b)=\sigma$ would be a map of these posets, but it is not a cellular map.
\end{ex}

The definition of $f_!$ uses kernels and other standard linear algebra operations. As such, we now assume $\dat=\Vect$ and suppress it from our notation.

\begin{defn}[Pushforward with Compact Supports]\index{pushforward with compact supports $f_{"!}$}\index{sheaf!pushforward with compact supports $f_{"!}$}
Given a sheaf $F$ on $X$ and a cellular map $f:X\to Y$, we can define the \textbf{pushforward with compact supports} sheaf on $Y$ as follows:
\begin{itemize}
 \item $f_!F(\tau)=\{s\in\Gamma(f^{-1}(\tau);F)\,|\,s(\sigma)=0\,\, \mathrm{if}\,\, |\bar{\sigma}|\cap f^{-1}(y)\,\,\mathrm{not}\,\,\mathrm{compact}\,\,\mathrm{for}\,\, y\in|\tau|\,\}$
\item Let $\gamma\leq \tau$ be cells in $Y$, and let $s\in f_!F(\gamma)$ and $t\in f_!F(\tau)$. We define $\rho^{f_!F}_{\tau,\gamma}(s)=t$ if for every $\sigma\in f^{-1}(\tau)$ and every $\lambda\in f^{-1}(\gamma)$ such that $\lambda\leq \sigma$ $\rho^F_{\sigma,\lambda}(s(\lambda))=t(\sigma)$. If there is no such $t\in f_!F(\tau)$ then we define $\rho^{f_!F}_{\tau,\gamma}(s)=0$.
\end{itemize}
The notation $\Gamma(-;F)$ for sections is explained in Definition \ref{defn:sections}. The verification that $f_!F$ is actually a sheaf and that it is functorial, is much more drawn out and is done in detail in~\cite{shepard} pp. 35-38. As such we have defined a functor
\[
	f_!:\Shv(X)\to\Shv(Y)
\]
\end{defn}

\begin{rmk}[Compact Supports for Cosheaves]\index{cosheaf!pushforward with compact supports}
	The definition for cosheaves cannot be written so simply because the vector space of ``compactly supported'' sections of a cosheaf, is a quotient of the space of all sections. The simplest definition would be, assuming $\hF:X^{op}\to\vect$, to take transposes and turn $\hF$ into a sheaf $F$ and apply the definition above. We will not make use of the cosheaf version of this functor.
\end{rmk}

\section{Calculated Examples}\label{subsec:calculated_examples}
In this section we compute explicit examples of the functors defined above. To avoid clutter, we consider only sheaves and leave it the reader to dualize and check the corresponding functors on cellular cosheaves. We further assume that $\dat=\Vect$ and leave it as implicit that all operations are to be performed in vector spaces. 

The notation $\Box$ will be a place holder for any one of the three symbols $*,\dd,!$.

\subsection{Projection to a point}\index{cellular maps!projection to a point}
We consider the constant map $p:X\to\star$. The output of $p_{\Box}F$ is a single vector space, namely $p_{\Box}F(\star)$. 

\begin{figure}[!ht]
\centering
\includegraphics[width=.7\textwidth]{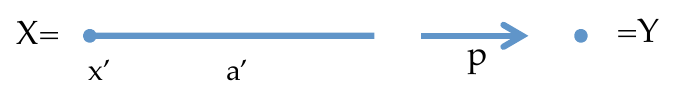}
\caption{Projection to a Point}
\label{fig:shv_fun_ex1}
\end{figure}

Without too much effort we compute the following:
\begin{itemize}
 \item $p_*F=\varprojlim\{F(x')\to F(a')\}\cong F(x')$
 \item $p_{\dd}=\varinjlim\{F(x')\to F(a')\}\cong F(a')$
\end{itemize}

For the pushforward with compact supports, we will be extra careful. Recall the definition states that $p_!F(\tau)=\{s\in\Gamma(p^{-1}(\tau);F)|s(\sigma)=0\, \mathrm{if}\, |\bar{\sigma}|\cap p^{-1}(y)\,\mathrm{not}\,\mathrm{compact}\,\mathrm{for}\, y\in|\tau|\}$.

In our example $y$ can be the only point $\star$ and $p^{-1}(\star)=X$. Thus we have only two cells to check whether their closures are compact or not. Clearly $\bar{x}'=x'$ is compact, but $\bar{a}'=X$ is not compact. The definition then says that we only allow sections whose value on $a'$ is zero.
\begin{itemize}
 \item $p_!F=\ker (\rho_{x',a'}:F(x')\to F(a'))$
\end{itemize}

\subsection{Inclusion into a Closed Interval}\index{cellular maps!inclusion of a point}
Here we encounter an open inclusion $j:X\to Y$. The first thing to note is that in this case, the value of $j_!F$ is not going to change since either $j^{-1}(y)=\{x\}$ or it is empty. Since points are closed and bounded, the compactness condition on $|\bar{\sigma}|\cap \{j^{-1}(y)\}$ is always satisfied.

\begin{figure}[!ht]
\centering
\includegraphics[width=\textwidth]{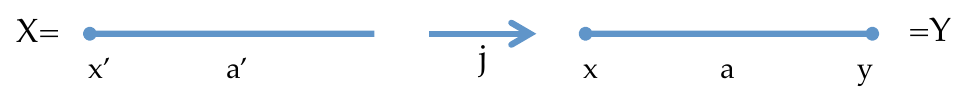}
\caption{Inclusion into a Closed Interval}
\label{fig:shv_fun_ex2}
\end{figure}

We see in this example that
\begin{itemize}
 \item $j_*F(x)=\varprojlim\{F(x')\to F(a')\}\cong F(x')$, $j_*F(a)=F(a')$, and less intuitively, $j_*F(y)\cong F(a')$.
 \item $j_{\dd}F(x)= F(x')$, $j_{\dd}F(a)\cong F(a')$, and $j_{\dd}F(y)=\varinjlim\{\emptyset\}=0$.
 \item $j_!F\cong j_{\dd}F$.
\end{itemize}

\subsection{Map to a Circle}\index{cellular maps!bijection with a circle}
Here is an example where the function is bijective and continuous (in both topologies), but not an embedding, i.e. the domain is not homeomorphic with its image.

\begin{figure}[h!]
\centering
\includegraphics[width=\textwidth]{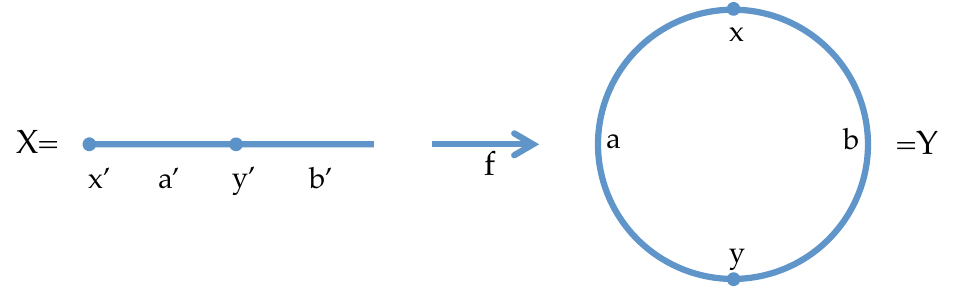}
\caption{Map to a Circle}
\label{fig:shv_fun_ex3}
\end{figure}

All three sheaves agree on the values and the restriction maps $f_{\Box}F(y)\cong F(y')\to F(a')\cong f_{\Box}F(a)$. We concentrate on the other two cells.
\begin{itemize}
 \item Here diagram we are taking the limit over is disconnected because the inverse image of the star of $x$ in $Y$ is disconnected. Consequently, $f_*F(x)=\varprojlim\{F(x')\to F(a') \quad F(b')\}\cong F(x')\oplus F(b')$ and $f_*F(b)=F(b')$.
 \item Here $f_{\dd}F(x)=F(x')$, but for similar reasons as before $f_{\dd}F(b)=\varinjlim\{F(x')\quad F(y')\to F(b')\}=F(x')\oplus F(b')$.
 \item As noted, $f$ is injective thus the value of $f_!F$ on any cell in the image is un-changed. However, we need to pay careful attention to how the restriction map is defined. The map $f$ is injective so the fiber over $a$ is $a'$ and over $x$ is $x'$, but $x'\nleq a'$ so the restriction map must be zero.
\end{itemize}

\section{The Push-Pull Adjunctions}\label{subsec:adjunctions}
\index{sheaf!push-pull adjunctions}\index{cosheaf!push-pull adjunctions}

Recall from Section~\ref{sec:adjunctions} that adjunctions allow us transform a complicated problem into an easy one. To derive these adjunctions, we can take two approaches: Use Freyd's adjoint functor theorem~\ref{thm:freyd}, or explicitly construct the adjunction. Since in our construction of the functors associated to a map, we made explicit use of limits and colimits, corresponding to the right and left Kan extensions respectively, and (co)limits commute with (co)limits, the following theorems are automatic. However, we check them explicitly for sheaves and leave the dual proof for the reader to fill out on their own.

\begin{thm}
The functors $f^*:\Shv(Y)\to \Shv(X)$ and $f_*:\Shv(X)\to\Shv(Y)$ form an adjoint pair $(f^*,f_*)$ and thus
\[
\Hom_{\Shv(X)}(f^*G,F)\cong \Hom_{\Shv(Y)}(G,f_*F).
\]
Dually, the functors for cosheaves satisfy the opposite adjunction $(f_*,f^*)$
\[
\Hom_{\Coshv(Y)}(f_*\hF,\hG) \cong \Hom_{\Coshv(X)}(\hF,f^*\hG).
\]
\end{thm}
\begin{proof}
	Recall that $f^*(f_*F)(x)=(f_*F)(f(x))$. Using the fact that $(f_*F)(f(x))=\varprojlim\{F(z)|f(z)\geq f(x)\}$, we get a map to $F(x)$ since $x\in f^{-1}(f(x))$ and this morphism is final for each $x$. This implies there is a natural transformation of functors $f^*\circ f_*\to\id$, which is universal (final). 
	
	Similarly, $f_*(f^*G)(y)=\varprojlim\{f^*G(x)=G(f(x))|f(x)\geq y\}$ and since $y\leq f(x)$ we can use the restriction map $\rho^G_{f(x),y}:G(y)\to G(f(x))$. The universal property of the limit guarantees a map $G(y)\to\varprojlim G(f(x))=f_*f^*G(y)$ and thus a natural transformation of functors $\id\to f_*f^*$.
\end{proof}

\begin{thm}
	The functors $f_{\dd}:\Shv(X)\to\Shv(Y)$ and $f^*:\Shv(Y)\to\Shv(X)$ form an adjoint pair $(f_{\dd},f^*)$ and thus
	\[
	\Hom_{\Shv(Y)}(f_{\dd}F,G)\cong\Hom_{\Shv(X)}(F,f^*G).
	\]
	Dually, the functors for cosheaves satisfy the opposite adjunction $(f^*,f_{\dd})$
	\[
	\Hom_{\Coshv(X)}(f^*\hG,\hF) \cong \Hom_{\Coshv(Y)}(\hG,f_{\dd}\hG).
	\]
\end{thm}
\begin{proof}
	$f_{\dd}(f^*G)(y)=\varinjlim \{G(f(x))|f(x)\leq y\}$ so again we can use the restriction maps to define maps to $G(y)$. The universal property of colimits gives a map $f_{\dd}f^*G(y)\to G(y)$ and thus a map of functors $f_{\dd}f^*\to\id$. Similar arguments give a map $\id\to f^*f_{\dd}$
\end{proof}

To conclude, we derive the first interesting consequence of an adjunction. In effect it reduces all the possible natural transformations between a certain pair of functors to a single vector space.
\begin{prop}
	If $F:X\to\Vect$ is a sheaf and $p:X\to\star$ is the constant map, then
	\[
		\Hom_{\Shv(X)}(p^*k,F)\cong \Hom_{\Vect}(k,p_*F)\cong F(X) = H^0(X;F).
	\]
\end{prop}
\begin{proof}
	The first isomorphism is the adjunction $(p^*,p_*)$. The second isomorphism is simply the observation that every linear map is determined by where it sends 1, i.e. $\Hom_{\Vect}(k,W)\cong W$.
\end{proof}

%
%

\chapter{Homology and Cohomology}
\label{sec:homology}

\begin{quote}
{\em``The de Rham complex may be viewed as a God-given set of differential equations, whose solutions are the closed forms.... A measure of the size of the space of `interesting' solutions is the definition of the de Rham cohomology.''}
\begin{flushright} --- Raoul Bott and Loring Tu~\cite[p. 15]{bott1982differential} \end{flushright}
\end{quote}

In Section~\ref{subsec:posets} and Chapter~\ref{sec:maps} we worked over arbitrary posets. We did this because it was natural and some applications may need this level of generality. In this section, we eschew this generality and restrict ourselves to posets arising as the face relation of a finite cell complex. This is beneficial not only because cell complexes are of great interest, but because sheaves and cosheaves over them have easily defined cohomology and homology theories.

We will start by describing a simple generalization of cellular cohomology and homology where we have augmented the coefficients by placing vector spaces over individual cells and linear maps between incident cells. This is a generalization in the sense that if one restricts to the case where every cell is assigned the one-dimensional vector space $k$ and all the incident linear maps are the identity, we recover classical cellular (co)homology. However interesting this special case may be, it misses a theory general enough to compute homological invariants of data varying over a cell complex.

The theory presented is combinatorial and computable. One needs only a good working knowledge of linear algebra to be able to use it. However, one can compute cellular sheaf cohomology without understanding it. To clarify the meaning of these computations we adopt a representation-theoretic perspective. This allows us to break up sheaves and cosheaves into the basic building blocks of indecomposable representations of the cell category. Thus, borrowing terminology from the persistent homology community, we use ``generalized barcodes'' to see the topology of data in a wider world of applications. These ideas are be put into practice in Chapters \ref{sec:barcodes}, \ref{sec:nc_coding}, and \ref{sec:sensors}, where many examples are considered.

\section{Chain Complexes and Homology}

\begin{defn}\index{graded vector space}
A \textbf{$\Z$-graded vector space} $V^{\ast}$ is a collection of vector spaces $\{V^i\}_{i\in\ZZ}$ with one for each integer. A \textbf{graded map} is a collection of linear maps $f:V^i\to W^i$. The \textbf{category of $\Z$-graded vector spaces}, $\grVect$, has graded vector spaces $V^{\ast}$ for objects and graded maps for morphisms.
\end{defn}

A (co)chain complex is a graded vector space with extra structure.

\begin{defn}\index{Cochain@(co)chain complex}
A \textbf{cochain complex} consists of a collection of vector spaces called \textbf{cochain groups} $\{V^i\}_{i\in\Z}$ and a collection of linear maps called \textbf{differentials} $d^i:V^i\to V^{i+1}$ that satisfy $d^{i+1}\circ d^i=0$ for every $i\in\Z$. We denote a cochain complex by $(V^{\bullet},d^{\bullet})$. Alternatively said, a cochain complex is a graded vector space equipped with a degree one increasing map that when composed twice gives the zero map.

A \textbf{chain complex} is a cochain complex with different notation. The \textbf{chain groups} $\{V_i\}_{i\in\Z}$ and \textbf{boundary maps} $\partial_i:V_i\to V_{i-1}$ are decorated with subscripts; this is the only difference. The maps satisfy $\partial_{i-1}\circ\partial_i=0$ for every $i\in\Z$. We denote a chain complex by $(V_{\bullet},\partial_{\bullet})$.
\end{defn}

\begin{rmk}
Since chain complexes and cochain complexes are the same thing, merely dressed up in different notation, we will usually just say ``Let $(V^{\bullet},d^{\bullet})$ be a chain complex'' and let the mathematical notation be precise. As an aside, one can also say that a chain complex is \emph{homologically indexed} if it is written as $(V_{\bullet},\partial_{\bullet})$ or \emph{cohomologically indexed} if it is written as $(V^{\bullet},d^{\bullet})$.
\end{rmk}

\begin{rem}
Sometimes we drop the subscript or superscript $\bullet$ and write $(V,\partial)$ or $(V,d)$ to refer to a chain complex or cochain complex. Dropping the superscript can lead to overloaded notation. For example, the expression $d^2=0$ is a synonym for ``$d$ is a differential,'' i.e. $d^{i+1}\circ d^i=0$ for ever $i\in\Z$, but it could also mean that the map $V^2\to V^3$ is zero. This is one of the perils of cohomological indexing for chain complexes, but the ambiguity is resolved by scoping the context. If we are speaking at the high-level of viewing a chain complex as a different sort of structure, then the former interpretation is intended. If we are talking about the particulars of a given chain complex, then the latter is meant. 
\end{rem}

\begin{defn}\index{category!C@$\Ch^{\bullet}(\Vect)$ of chain complexes}
The \textbf{category of chain complexes}, $\Ch^{\bullet}(\Vect)$, has chain complexes for objects, and chain maps $f^{\bullet}:(V,d_V)\to (W,d_W)$ for morphisms, i.e. a collection of maps $f^i:V^i\to W^i$ such that $f^{i+1}\circ d^i_V=d^i_W\circ f^i$.
\end{defn}

There is a natural functor $\iota: \grVect\to\Ch(\Vect)$ that treats a $\Z$-graded vector space as a chain complex with zero differentials, i.e.
\[
\{V^i\}_{i\in\Z} \rightsquigarrow (V^{\bullet},0)
\]
Taking cohomology of a chain complex defines a functor going the other way.

\begin{defn}\index{cohomology}
\textbf{Cohomology} is a functor $H^{\ast}:\Ch(\Vect)\to\grVect$, which takes a chain complex $(V,d)$ and places the quotient vector space 
\[
H^i(V,d):=\ker(d^i)/\im(d^{i-1})
\] 
in degree $i$. Without too much work one can show that a chain map $f^{\bullet}$ induces maps of the associated cohomology spaces $H^i(f):H^i(V,d_v)\to H^i(W,d_w)$, making $H^{\ast}$ into a functor.
\end{defn}

\subsection{The Combinatorics of Cell Complexes and Homology}

The motivation for chain complexes and homology comes from computing invariants of topological spaces. As already indicated, posets can be regarded as topological spaces, but not every poset has the nice structure that the face-relation poset of a cell complex has. This nice structure is what determines whether certain sequence of vector spaces and maps defines a chain complex.

\begin{defn}
We write $\sigma \leq_i \tau$ if the difference in dimension of the cells is $i$.
\end{defn}

\begin{lem}\label{lem:2cells}
 If $\sigma\leq_2\tau$, then there are exactly two cells $\lambda_1,\lambda_2$ where $\sigma\leq_1\lambda_i\leq_1\tau$.
\end{lem}

We want to invent a sign condition that distinguishes these two different sequences of incidence relations.
\[
	\xymatrix{ & \tau & \\ \lambda_1 \ar[ru] & & \lambda_2 \ar[lu]\\ & \ar[lu] \sigma \ar[ru] &}
\]

\begin{defn}[Signed Incidence Relation]\index{incidence relation}
 A \textbf{signed incidence relation} is an assignment to any pair of cells $\sigma,\tau\in X$ a number $[\sigma:\tau]\in\{0,\pm 1\}$ such that 
 \begin{itemize}
 \item if $[\sigma:\tau]\neq 0$, then $\sigma\leq_1\leq _1 \tau$, and
 \item if $\gamma$ and $\tau$ are any pair of cells, the sum $\sum_{\sigma}[\gamma:\sigma][\sigma:\tau]=0$.
 \end{itemize}
\end{defn}

One way to get a signed incidence relation is to choose a \textbf{local orientation}\index{local orientation} (via the homeomorphism of each cell $|\sigma|$ with $\RR^k$) for each cell without regard to global consistency. Then for every pair of incident cells $\sigma\leq\tau$ we have a number $[\sigma:\tau]=\pm1$ given by $+1$ if the orientations agree and $-1$ otherwise.

Another way is motivated by working with regular cell complexes, where we can subdivide so that we have a simplicial complex. We can refer to any cell by a list of its vertices. If we order the set of vertices, then we have a procedure for orienting the cells. A local orientation of a cell $\sigma\in X$ consists of divvying up the set of ordered lists representing $\sigma$ into classes each of which are invariant under even permutations. We can then pick the class with the list of vertices in increasing order as ``the'' orientation. Either method enables us to define a chain complex associated to a cell complex.

\begin{prop}[Cellular Cohomology]\index{cohomology!of a cell complex}
Let $X$ be a cell complex equipped with a sign relation. Let $C^n(X;k)$ be the vector space spanned by the $n$-dimensional calls of $X$. We define a map $\delta:C^n\to C^{n+1}$ on the basis by defining $\delta(\sigma)=\sum_{\tau} [\sigma:\tau]\tau$. Clearly $\delta^2=0$.
\end{prop}

\section{Computational Sheaf Cohomology and Cosheaf Homology}
\label{subsec:comp_homology}

We now provide formulae for computing cellular sheaf cohomology and cellular cosheaf homology that is completely analogous to cellular cohomology.

\subsection{Cellular Sheaf Cohomology}

\begin{defn}[\cite{dihom_1, shepard}]\index{cohomology!of a cellular sheaf!compactly-supported}\index{sheaf!cellular!compactly supported cohomology of}
 Given a cellular sheaf $F:X\to\Vect$ we define its \textbf{compactly supported $k$ co-chains} to be the product\footnote{Here we implicitly assume that $X$ has finitely many cells in a given dimension so products and direct sums agree.} of the vector spaces residing over all the $k$-dimensional cells.
\[
 C^k_c(X;F)=\bigoplus_{\sigma^k}F(\sigma^k)
\]
 These vector spaces are graded components in a complex of vector spaces $C^{\bullet}_c(X;F)$. The differentials are defined by
\[
 \delta^k_c=\sum_{\sigma\leq\tau} [\sigma^k:\tau^{k+1}]\rho_{\tau,\sigma}.
\]
The cohomology of this complex
\[
 \xymatrix{0 \ar[r] & \oplus F(\mathrm{vertices}) \ar[r]^-{\delta_c^0} & \oplus F(\mathrm{edges}) \ar[r]^-{\delta_c^1} & \oplus F(\mathrm{faces}) \ar[r] & \cdots & = C^{\bullet}_c(X;F)}
\]
is defined to be the \textbf{compactly supported cohomology} of $F$, i.e. $H^k_c(X;F)=\ker\delta^k_c/\im \delta^{k-1}_c$.
\end{defn}

\begin{lem}
 $(C^{\bullet}_c(X;F),\delta^{\bullet}_c)$ is a chain complex.
\end{lem}
\begin{proof}
To see why the chain complex condition $\delta_c^{k+1}\delta_c^{k}=0$ is assured, Lemma \ref{lem:2cells} is crucial. This is the very same lemma that proves that ordinary cellular homology is computed via a chain complex. One must now observe that varying data over the cells does not change the result.
\begin{eqnarray*}
 \delta_c\delta_c &=& \sum_{\sigma\leq_1\tau} [\sigma:\tau]\rho_{\tau,\sigma}(\delta_c) \\
 &=& \sum_{\sigma\leq_1\tau} [\sigma:\tau]\rho_{\tau,\sigma}(\sum_{\gamma\leq_1 \sigma} [\gamma:\sigma]\rho_{\sigma,\gamma}) \\
 &=& \sum_{\gamma\leq_1\sigma\leq_1\tau} [\gamma:\tau] \rho_{\tau,\sigma}\rho_{\sigma,\gamma} \\
 &=& \sum_{\gamma\leq_1\sigma\leq_1\tau} [\gamma:\tau] \rho_{\tau,\gamma} \\
 &=& \sum_{\gamma\leq_1\sigma\leq_1\tau} ([\gamma:\sigma_1][\sigma_1:\tau]+[\gamma:\sigma_2][\sigma_2:\tau]) \rho_{\tau,\gamma} \\
 &=& 0
\end{eqnarray*}
\end{proof}

To define the arbitrarily-supported cochain complex associated to a cellular sheaf $F$ on $X$, we simply remove all the cells from $X$ without compact closures and apply the same formula.

\begin{defn}[Ordinary Cohomology]\index{cohomology!of a cellular sheaf}\index{sheaf!cellular!cohomology of}
    Let $X$ be a cell complex and $F:X\to\Vect$ a cellular sheaf. Let $j:X'\to X$ be the subcomplex consisting of cells that do not have vertices in the one-point compactification of $X$. Define the ordinary cochains and cohomology by
\[
	C^{\bullet}(X;F)=C^{\bullet}_c(X';j^*F) \qquad H^i(X;F):=H^i_c(X';j^*F)
\]
\end{defn}

The situation may seem a bit unusual. The naturally defined chain complex computes a more restrictive type of cohomology. To get the standard cohomology, one needs to remove non-compact cells. When we define cohomology via the derived perspective of Section \ref{sec:derived}, this quirk of linear algebra disappears. Ordinary cohomology will fall out naturally using limits and injective resolutions, and compactly-supported sheaf cohomology will require some finesse.

\begin{ex}[Compactly Supported vs. Ordinary Cohomology]
To see why the na\"ive chain complex computes compactly supported cohomology, consider the example of the half-open interval $X=[0,1)$ decomposed as $x=\{0\}$ and $a=(0,1)$. Now consider the constant sheaf $k_X$. To compute compactly supported cohomology, we must first pick a local orientation of our space. By choosing the orientation that points to the right, we get that $[x:a]=-1$. The cohomology of our sheaf is computed via the complex
\[
  \xymatrix{k \ar[r]^{-1} & k},
\]
which yields $H^0_c=H^1_c=0$. If we follow the prescription for computing ordinary cellular sheaf cohomology, then we must remove the vector space sitting over $a$ in our computation. The resulting complex is simply the vector space $k$ placed in degree 0, so $H^0(X;k_X)=k$ and is zero in higher degrees.
\end{ex}

\begin{figure}[ht]
\centering
\includegraphics[width=.3\textwidth]{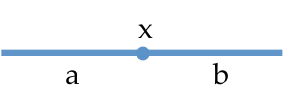}
\caption{Minimal Cell Structure on an Open Interval}
\label{fig:open_interval}
\end{figure}

\begin{ex}[Open Interval]
 If we pretended for a moment that the pure stratum $Y=(0,1)$ is a cell complex\footnote{Recall that we require a cell complex to have a one point compactification that is a regular cell complex.} with no other cells, then computing the compactly supported cohomology of the constant sheaf would yield a vector space in degree one and nowhere else, hence $H^1_c(Y;k_Y)=k$.

 To make this example a legitimate example, as in Figure \ref{fig:open_interval}, we place a vertex at $x=1/2$. We call our new cells $a=(0,1/2)$ and $b=(1/2,1)$. If we orient our 1-cells to point to the right, then $[x:a]=1$ and $[x:b]=-1$. Using the lexicographic ordering on our cells to get a basis for $C^1_c(Y;k_Y)$ we can compute explicitly the compactly supported cohomology.
 \[
  \delta^0_c=\begin{bmatrix}1 \\ -1\end{bmatrix} : k_x \to k_a\oplus k_b \qquad \Rightarrow \qquad H^0_c=0 \quad H^1_c=k
 \]
\end{ex}

\subsection{Cellular Cosheaf Homology}

For cellular cosheaves the exact dual construction works, but the terminology is slightly different.

\begin{defn}[Borel-Moore Cosheaf Homology]\index{homology!of a cellular cosheaf!Borel-Moore}\index{cosheaf!cellular!Borel-Moore homology of}
 Let $X$ be a cell complex and let $\hF:X^{op}\to\Vect$ be a cellular cosheaf. Define the \textbf{Borel-Moore homology of} $H_{\bullet}^{BM}(X;\hF)$ to be the homology of the following complex:
\[
 \xymatrix{C^{BM}_{\bullet}(X;\hF) = & \cdots \ar[r] & \oplus \hF(\mathrm{faces}) \ar[r]^-{[e:\sigma]r_{e,\sigma}} & \oplus \hF(\mathrm{edges}) \ar[r]^-{[v:e]r_{v,e}} & \oplus \hF(\mathrm{vertices}) \ar[r] & 0 }
\]
\end{defn}

\begin{defn}[Ordinary Cosheaf Homology]\index{homology!of a cellular cosheaf}\index{cosheaf!cellular!homology of}
	Let $X$ be a cell complex and let $\hF:X^{op}\to\Vect$ be a cellular cosheaf. By discarding all the cells without compact closure, we obtain the maximal compact subcomplex $X'$. If we write $j:X'\hookrightarrow X$ for the inclusion, then we can define the ordinary chain complex to be
	\[
	C_{\bullet}(X;\hF)=C^{BM}_{\bullet}(X';j^*\hF).
	\]
	Applying the definition above gives the ordinary \textbf{cosheaf homology} $H_{\bullet}(X;\hF)$ of a co-sheaf.
\end{defn}

All of the examples of cellular sheaf cohomology dualize to give interesting examples of cellular cosheaf homology. Let us define the functor that performs this operation.
\begin{defn}[Linear Duality]\index{duality!linear}\index{linear duality}
	Let $\widehat{V}: \Shv(X;\vect_k)^{op} \to \Coshv(X;\vect_k)$ be the contravariant equivalence from sheaves to cosheaves, both valued in \emph{finite} dimensional vector spaces, defined as follows:
	\[
		\xymatrix{F(\tau) \ar@{~>}[r] & F(\tau)^*  \ar[d]^{\rho^*_{\tau,\sigma}} & \widehat{V}(F)(\tau) \ar[d]^{r_{\sigma,\tau}} \\
		F(\sigma) \ar[u]^{\rho_{\tau,\sigma}} \ar@{~>}[r] & F(\sigma)^* & \widehat{V}(F)(\sigma)}
	\]
\end{defn}

\begin{lem}
 Taking linear duals preserves cohomology, i.e. $H_i(X;\widehat{V}(F))\cong H^i(X;F)$ and $H_i^{BM}(X;\widehat{V}(F))\cong H^i_c(X;F)$.
\end{lem}

\section{Explaining Homology and Cohomology via Indecomposables}
\label{subsec:rep_thy}

Sheaf cohomology is notoriously difficult to interpret. Every time a successful interpretation is discovered, a cornerstone of a theory waiting to be fleshed out is put into place. For example, the Cousin problems of complex analysis asks whether a meromorphic function with a given divisor (zeros and poles) exists or not. When Cartan and Serre interpreted this problem in terms of sheaf theory, sheaf cohomology groups gave a complete classification and obstruction theory; see~\cite{gray} p. 17. The narrative that falls out of those historical successes is that sheaf cohomology gives calculable obstructions to finding solutions.

However, when the above interpretation fails, we need to compute examples and extract new interpretations. When computing sheaf cohomology, one encounters a plethora of choices that obfuscate the natural meaning of the vector spaces $H^i(X;F)$: picking ordered bases for each $F(\sigma)$, choosing local orientations, computing kernels and quotients, taking representative elements of cohomology or homology, etc. Each of these lead one farther from a workable interpretation of the topology of data.

The experience of the author in computing examples of sheaf cohomology has led him to believe that the best way of circumventing these issues is to borrow an idea from the representation theory of quivers. Specifically, if one knows the direct sum decomposition of a sheaf into indecomposable sheaves, then one gets a distinguished basis for sheaf cohomology. These indecomposables allow one to see how data travels through a space.

\subsection{Persistence Modules and Barcodes}
\label{subsubsec:barcodes}

To begin the introduction of representation theory gently, we will describe a convenient visualization technique called a \textbf{barcode}, which was first described by Carlsson, Zomorodian, Collins and Guibas~\cite{carlsson2004persistence}. The motivation for those authors was to provide a simple shape descriptor for data that could be used by scientists not trained in representation theory, but we will adapt it for understanding sheaves and cosheaves. 

To begin, let us recast a chain complex as a special instance of the following structure:

\begin{defn}\label{defn:persistence_module}\index{persistence module}\index{module!persistence}
A \textbf{persistence module} consists of a collection of vector spaces $\{V^i\}_{i\in\Z}$, one for each integer, and a collection of linear maps $\rho^i_V:V^i\to V^{i+1}$. If $i\leq j$, then we define $\rho^V_{i,j}:=\rho^{j-1}_V\circ\cdots\circ \rho^i_V$ to be the uniquely determined map from $\rho^V_{i,j}:V^i\to V^j$. We denote a persistence module by $(V,\rho_V)$, but we may suppress the $V$ in $\rho_V$ or even drop the $\rho_V$ all together.
\end{defn}

Observe that one can add two persistence modules to create a third persistence module, i.e. if $(V^{\bullet},\rho_V)$ and $(W,\rho_W)$ are two persistence modules, then one obtains a third persistence module $(U,\rho_U)$ by defining $U^i:=V^i\oplus W^i$ and $\rho^i_U:=\rho^i_V\oplus \rho^i_W$. We denote the sum by $(V\oplus W,\rho_V\oplus\rho_W)$ or more simply by $V\oplus W$.

There is a fundamental structure theorem for persistence modules, due to Crawley-Boevey~\cite{crawley2012decomposition}, that explains how any persistence module can be written as a direct sum of simpler persistence modules. We now introduce these simpler persistence modules. 
\begin{defn}\label{defn:interval}\index{interval module}\index{module!interval}
Recall that an \textbf{interval} in $(\Z,\leq)$ is a subset $I\subset \Z$ having the property that if $i,k\in I$ and if there is a $j\in I$ such that $i\leq j\leq k$, then $j\in I$. An \textbf{interval module} $k_I$ assigns to each element $i\in I$ the vector space $k$ and assigns the zero vector space to elements in $\Z\setminus I$. All maps $\rho_{i,j}$ are the zero map, unless $i,j\in I$ and $i\leq j$, in which case $\rho_{i,j}$ is the identity map. 
\end{defn}

Since interval modules are completely determined by the interval where they assign non-zero vector spaces, we can draw a \textbf{bar} to represent an interval module. The following structure theorem shows that any persistence module can be represented by a collection of bars, called a \textbf{barcode}.\index{barcode}

\begin{thm}[Decomposition for Pointwise-Finite Persistence Modules]\label{thm:crawleyboevey}\index{persistence module!decomposition theorem}\index{decomposition!interval/barcode}

If $(V,\rho_V)$ is a persistence module for which every vector space $V^i$ is finite-dimensional, then the module is isomorphic to a direct sum of interval modules, i.e.
\[
V\cong \bigoplus_{I\in D} k_I.
\]
Here $D$ is a multi-set of intervals. A multi-set is a set allowing repetitions, i.e. a set equipped with a function $\mu$ indicating the multiplicity of each given element.

\end{thm}

\begin{rem}\index{barcode!length}

We will refer to the length of the bar as $\ell(I)=j-i$

\end{rem}

This theorem summarizes a great deal of elementary linear algebra and quiver representation theory. For linear algebra, it has the fundamental theorem of linear algebra as a consequence~\cite{strang1993fundamental}, i.e. any map of vector spaces $T:V\to W$ has a matrix representation that is diagonal with $0$ and $1$ entries, the number of 1s corresponding to the rank of the matrix, cf.~\cite{artin1991algebra} Chapter 4, Proposition 2.9. Said differently, there are vector space isomorphisms making the following diagram commute:

\[
    \xymatrix{ V \ar[r]^T \ar[d]_{\varphi}^{\cong} & W \ar[d]^{\psi}_{\cong} \\
    \im(T) \oplus \ker(T) \ar[r]_{\id\oplus 0} & \im(T)\oplus \cok(T)}
\]
Here $\im(T)$, $\ker(T)$, and $\cok(T)$ refer to the image, kernel and cokernel of $T$ respectively. Although the image of $T$ is properly a subspace of $W$, the first isomorphism theorem identifies it with $V$ modulo the kernel.

\begin{figure}
    \centering
    \includegraphics[width=.8\textwidth]{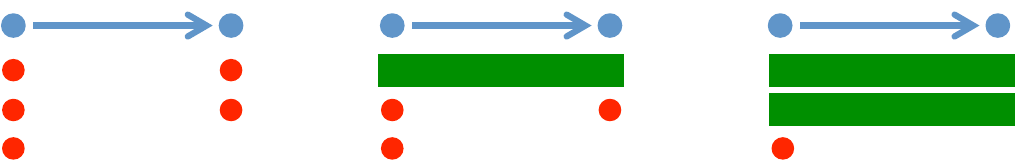}
    \caption{Barcodes associated to $T:\R^3\to\R^2$ for $\rank(T)=0,1,2$}
    \label{fig:rank012bars}
\end{figure}

\begin{ex}[Barcodes for Linear Algebra]\index{barcode!fundamental theorem of linear algebra}

Consider any linear map $T:\R^3\to\R^2$ as a persistence module by extending by zero vector spaces and maps. There are three isomorphism classes of such persistence modules determined by the rank of $T$. The associated barcodes are depicted in Figure~\ref{fig:rank012bars}.

\end{ex}

\begin{figure}[ht]
    \centering
    \includegraphics[width=\textwidth]{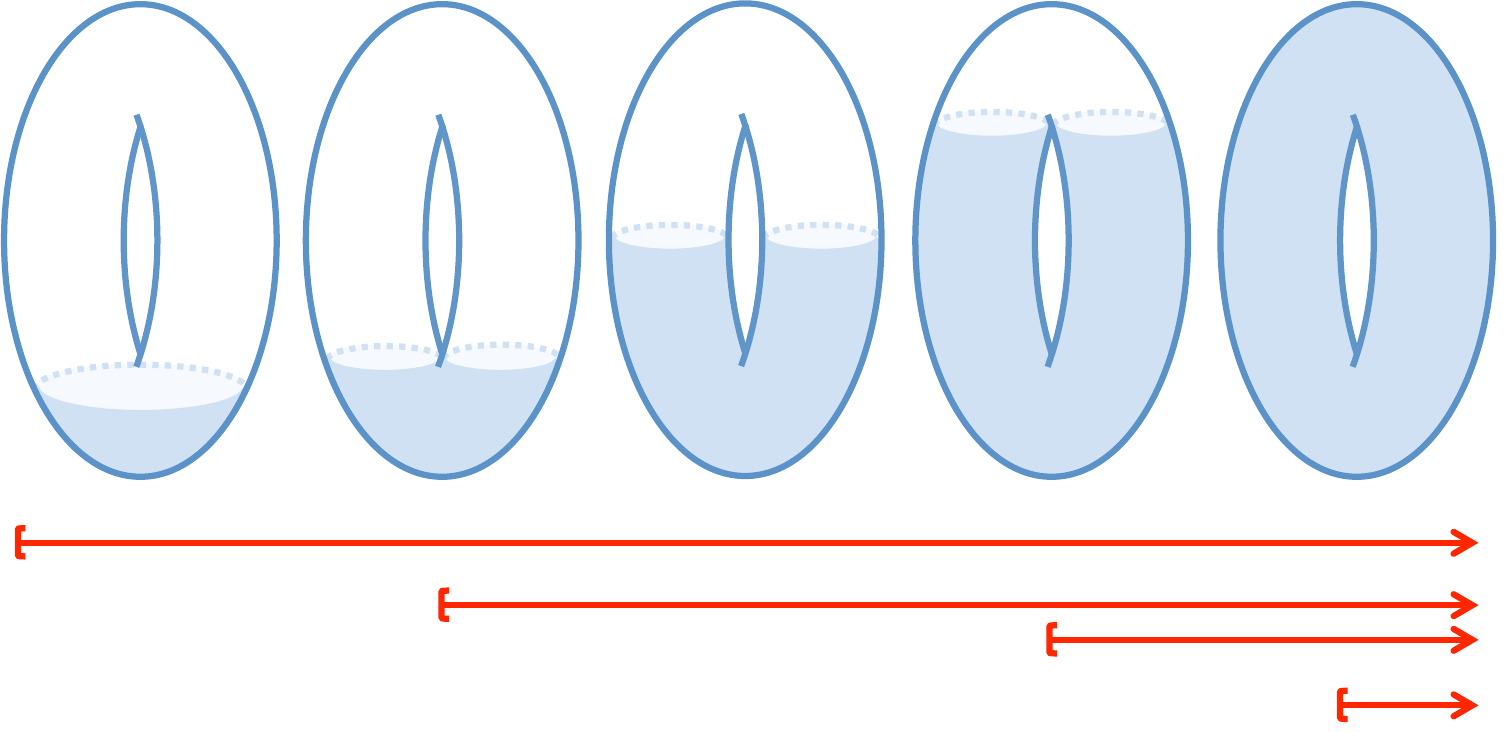}
    \caption{Barcodes for the filtration of a Torus}
    \label{fig:persistence_torus}
\end{figure}

\begin{ex}[Barcodes for Filtrations]\index{barcode!filtration of Bott's torus}

Barcodes with longer bars appear in filtrations of topological spaces. For example, consider the standard height function on the torus $h:X\to\R$. By choosing a discrete set of points $\{t_0 < t_1,\ldots \}$ to sample the function $h$ at, we get a sequence of spaces $\{X_{\leq t_i}=h^{-1}(-\infty,t_i]\}$, which after taking homology in some degree $i\geq 0$ defines a persistence module. Such an example is depicted in Figure \ref{fig:persistence_torus}.

\[
X_{\leq t_0} \hookrightarrow X_{\leq t_1} \hookrightarrow \cdots \qquad \rightsquigarrow \qquad H_i(X_{\leq t_0};k) \to H_i(X_{\leq t_1};k) \to \cdots
\]

\end{ex}

Now we reach an example that will be useful when we undertake the derived category of chain complexes.

\begin{figure}[ht]
    \centering
    \includegraphics[width=.8\textwidth]{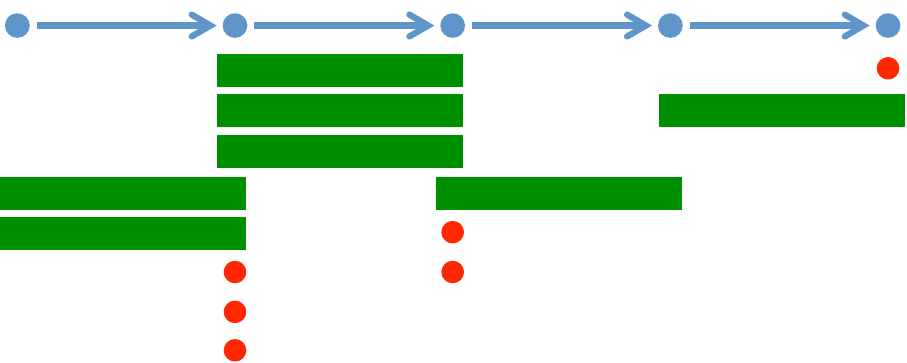}
    \caption{Barcodes for a Chain Complex}
    \label{fig:chaincomplex_bcs}
\end{figure}

\begin{ex}[Barcodes for Chain Complexes]\label{ex:chaincomplex_bc}\index{barcode!associated to chain complex}
As already remarked, a chain complex of vector spaces is a special example of a persistence module and, consequently, has a barcode decomposition. With a moment's reflection one can see that any chain complex can be written as the direct sum of two types of modules: 
the length zero interval modules
\[
S_i: \qquad \cdots \to 0 \to k \to 0 \to \cdots
\]
and the length one interval modules.
\[
P_i: \qquad \cdots \to 0 \to k \to k \to 0 \to \cdots
\]
Figure~\ref{fig:chaincomplex_bcs} gives a visual depiction of such a barcode decomposition.
\end{ex}

\subsection{Representation Theory of Categories and the Abelian Structure}
\label{subsubsec:rep_abs}
For the purposes of this section, there is no real difference between cellular sheaves and cosheaves --- they are both representations of the cell category $\Cell(X)$. Recall, for any category $\cat$ the \textbf{category of representations} is defined to be the category of functors to $\Vect$:\index{category!of representations}
\[
	\Rep(\cat):=\Fun(\cat,\Vect)
\]
This category has the structure of an abelian category, which we explain in this section. In effect, this means we can do everything in $\Rep(\cat)$ that we can do in $\Vect$: take kernels and cokernels of maps between representations, talk about images of maps, add maps and so on. We will introduce these properties as we need them.

\begin{clm}\index{category!exact}
	For $\cat$ a category, $\Rep(\cat)$ is an \textbf{exact category}. This means we can talk about exact sequences. Specifically:
	\begin{itemize}
		\item There is a \textbf{zero representation} given by sending all objects and morphisms to the zero object and the zero morphism.
		\item Between any two representations $F$ and $G$ there is a \textbf{zero map}, which can be factored through the zero representation.
		\item Since any morphism $\eta:F\to G$ is a natural transformation occurring inside $\Vect$, there are associated \textbf{kernel and cokernel representations} denoted $\ker(\eta)$ and $\coker(\eta)$ defined by taking kernels and cokernels object wise:
		\[
		\xymatrix{\ker(\eta(c')) \ar[r] & F(c') \ar[r]^{\eta(c')} & G(c') \ar[r] & \coker(\eta(c')) \\
		\ker(\eta(c)) \ar[r] \ar@{.>}[u] & F(c) \ar[r]^{\eta(c)} \ar[u]_{F(f)} & G(c) \ar[r] \ar[u]_{G(f)} & \coker(\eta(c)) \ar@{.>}[u] }
		\]
		\item There is an \textbf{image representation} $\im(\eta)$ defined as the object-wise image.
	\end{itemize}
\end{clm}

As usual, we say that a sequence of representations $A\to B \to C$ is \textbf{exact}\index{exact!sequence of representations} at $B$ if the kernel of the outgoing morphism is equal to the image of the incoming morphism. A longer sequence is exact if it is exact at each place with an incoming and outgoing morphism.

We are going to do a brief sketch of some representation theory for categories, using the terminology introduced.

\begin{defn}\index{representation!subrepresentation}
	A \textbf{subrepresentation} $E$ consists of a choice of subspace $E(c)\to F(c)$ for each object that is invariant under all the linear maps $F(f)$. Restriction of $F(f)|_{E(c)}=:E(f)$ makes $E$ into a representation of its own right. Said more succinctly, $E\to F$ is a natural transformation of functors that is object wise an inclusion, i.e.
	\[
		0 \to E \to F
	\] 
is an exact sequence. Dually, we can say $G$ is a quotient representation by saying $F\to G \to 0$ is an exact sequence.
\end{defn}

\begin{defn}\index{representation!direct sum}
	Suppose $F:\cat\to\Vect$ and $G:\cat\to\Vect$ are two representations of a small category $\cat$, then we can define the \textbf{direct sum} of these two representations $H=F\oplus G$ by defining on objects $H(c):=F(c)\oplus G(c)$ and on morphisms $H(f)=F(f)\oplus G(f):H(c)\to H(c')$.
\end{defn}

The above definition further clarifies the structure of $\Rep(\cat)$.

\begin{clm}\index{category!additive}\index{category!abelian}
	For $\cat$ a category, $\Rep(\cat)$ is both an exact and an \textbf{additive category}. This latter definition requires the following:
	\begin{itemize}
		\item For any two representations $F$ and $G$ the set $\Hom_{\Rep(\cat)}(F,G)$ has the structure of an abelian group (with the zero map being the additive identity) making composition bilinear.
		\item The direct sum of two representations is a representation.
	\end{itemize}
	A category that is exact and additive is defined to be \textbf{abelian}. Thus $\Rep(\cat)$ is an abelian category. 
\end{clm}

\begin{rmk}
In any additive category, it can be shown that having finite direct sums (finite coproducts) implies the existence of finite direct products (finite products) and these are isomorphic. 
\end{rmk}

\begin{defn}[Indecomposable]\index{representation!indecomposable}\index{indecomposable representation}
	A representation $F:\cat\to\Vect$ is called \textbf{indecomposable} if whenever $F$ is written as a direct sum of representations one of the representations is the zero one; i.e. every direct sum decomposition is trivial.
	
	Said using sequences, a representation $F$ is indecomposable if whenever we have a short exact sequence of representations
	\[
	0\to E \to F \to G \to 0
	\] 
	with neither $E$ nor $G$ the zero representation, then $F\ncong E\oplus G$, i.e. the sequence does not split.
\end{defn}

\begin{exr}
Verify that the interval modules in Definition \ref{defn:interval} are indecomposable representations. What is the underlying category that these modules represent?
\end{exr}

\begin{rmk}[Splitting Lemma]\index{splitting lemma}
There is a general lemma called the \textbf{splitting lemma}, which provides equivalent ways of saying that $F$ is indecomposable. It states that writing $F$ as a direct sum is equivalent to either having a map back from $F$ to $E$, which precomposed with the inclusion $E\to F$ yields the identity, or having a map back from $G$ to $F$, which post-composed with the surjection is the identity on $G$.
\end{rmk}

\begin{defn}[Remak Decomposition]\index{decomposition!Remak}
	A direct sum decomposition of an object $F\in\Rep(\cat)$
	\[
		F\cong F_1\oplus\cdots\oplus F_n
	\]
	where each $F_i$ is indecomposable and non-zero is called a \textbf{Remak decomposition}.
\end{defn}

A fact that we would very much like to know is whether every representation admits a Remak decomposition. The structure theorem \ref{thm:crawleyboevey} provides an example where this is the case. Sir Michael Atiyah considered such a question in the very general setting of abelian categories~\cite{atiyah_ks}. He developed a bi-chain condition and proved that under this condition every non-zero object admitted a Remak decomposition. We use a stronger condition of finite-dimensionality that Atiyah showed implied his bi-chain condition.

\begin{thm}[Krull-Schmidt Theorem for Representations~\cite{atiyah_ks}]\index{Krull-Schmidt Theorem}
	Suppose $\aat$ is an abelian category, further satisfying
	\begin{enumerate}
		\item For every pair of objects $\Hom_{\aat}(A,B)$ is a finite dimensional vector space, and
		\item Conjugation is linear, i.e. for every pair of morphisms $\varphi:A\to B$ and $\psi:B'\to A'$ the following map is linear
		\[
			\Hom_{\aat}(B,B')\to \Hom_{\aat}(A,A') \qquad \eta\mapsto \psi\circ\eta\circ\varphi
		\]
	\end{enumerate}
		then the Krull-Schmidt theorem holds. This says that every non-zero object $A$ has a Remak decomposition and for any two such decompositions
		\[
			A\cong A_1\oplus\cdots\oplus A_n \qquad A\cong A'_1\oplus\cdots\oplus A'_m
		\]
		$n=m$ and after re-ordering $A_i\cong A_i'$.
\end{thm}

For $\aat=\Rep(\cat)$ the second condition is certainly satisfied. The first condition imposes significantly stronger conditions. First of all, we must restrict to the full subcategory of finite dimensional representations.
\[
	\Rep_f(\cat):=\Fun(\cat,\vect) \subset \Fun(\cat,\Vect)=:\Rep(\cat)
\]
Secondly, one must observe that for any two representations $F$ and $G$ the space of natural transformations is a subspace of a potentially infinite product of finite dimensional spaces.
\[
	\Hom(F,G)\subseteq \prod_{c\in\cat}\Hom_{\vect}(F(c),G(c)) 
\]
One severe restriction one can make to insure that Atiyah's first condition holds is to assume that the category $\cat$ has finitely many objects. This is not strictly necessary, but it does provide us with the following corollary:

\begin{cor}[Sheaves and Cosheaves on Finite Posets have Remak Decompositions]\label{cor:sheaves_remak}
	Suppose $(X,\leq)$ is a finite poset, then $\Shv(X)$ and $\Coshv(X)$ satisfy the Krull-Schmidt theorem.
\end{cor}

The example that we have in mind, of course, is the poset associated to a a cell complex $X$. In this situation, one can recognize a large set of examples of indecomposable representations.

\begin{lem}[Constant (Co)Sheaves are Indecomposable]\label{lem:const_indecomp}\index{indecomposable representation!constant (co)sheaves on connected}
	Suppose $X$ is a connected cell complex, then the constant sheaf $k_X$ and the constant cosheaf $\hat{k}_X$ are indecomposable.
\end{lem}
\begin{proof}
	We'll state the proof for sheaves and leave it to the reader to dualize for cosheaves. Suppose for contradiction that $k_X\cong F\oplus G$ where neither $F$ nor $G$ is the zero sheaf. Now as a consequence of neither $F$ nor $G$ being zero, and $k_X$ being one dimensional on each cell, there must be a pair of cells $\sigma$ and $\tau$ such that one is in the support of $F$ and the other is in the support of $G$. We argue that we can choose $\sigma$ and $\tau$ such that one is the face of the other. If not, then the support of $F$ (or $G$) would be closed under the following operations
	\[
		\sigma\subset \supp(F) \quad \mathrm{and} \quad \tau\subset \bar{\sigma} \quad \mathrm{or} \quad \sigma \subset \bar{\tau} \qquad \mathrm{then} \qquad \tau\subset \supp(F)
	\]
	which by connectedness of $X$ would imply that $\supp(F)=X$; a contradiction to the supposition that neither $F$ nor $G$ was the zero sheaf. (To see why $\supp(F)=X$, one can imagine drawing the Hasse diagram of the poset $X$ and realizing that connectedness means that the diagram is connected.) Thus we have such a pair $\sigma\subset \tau$ with one in the support of $F$ and the other in the support of $G$, but this also can not occur since the identity cannot be written as a sum of zero maps.
	\[
		k\to k \neq (k\to 0)\oplus (0\to k).
	\]
\end{proof}

\subsection{Quiver Representations and Gabriel's Theorem}
\label{subsubsec:quivers}

\begin{quote}
{\em ``These graphs arise in a multitude of classification problems in mathematics, such as classification of simple Lie algebras, singularities, platonic solids, reflection groups, etc. In fact, if we needed to make contact with an alien civilization and show them how sophisticated our civilization is, perhaps showing them Dynkin diagrams would be the best choice!''}~\cite{etingof_rep}
\end{quote}

There are natural examples of representations of categories where these ideas and their consequences have been studied. One such example is the category associated to a directed graph. 

\begin{defn}\index{quiver}
A \textbf{quiver} or \textbf{directed graph} is defined by a pair of sets consisting of ``edges'' $E$ and ``vertices'' $V$ along with a pair of functions $h,t:E\to V$, which we think of as denoting the head and tail of a directed edge respectively. Alternatively, a quiver can be topologically regarded as a one-dimensional cell complex equipped with a local orientation of its edges.
\end{defn}

One should be careful to note that a directed graph is not a category in and of itself, but there is a natural category associated to a directed graph, which we now define.

\begin{defn}\index{category!path/free}
To a quiver we can associate a category called the \textbf{free category} or \textbf{path category} written $\Free(X)$ . The objects are vertices and the morphisms are directed paths between vertices. Since paths are just concatenated edges, we think of the morphisms as being freely generated by the edges. We must consider simply sitting at a vertex as the identity directed path connecting the vertex to itself. 
\end{defn}

\begin{defn}\index{quiver!representation}\index{representation!of a quiver}
A \textbf{quiver representation} is thus nothing more than a functor $F:\Free(X)\to \Vect$. Because a general path is simply a sequence of edges, such a functor is equivalent to specifying a vector space for each vertex in $V$ and a linear map for each edge in $E$ that goes from the source to the target.
\end{defn}

Every finite dimensional quiver representation can be decomposed into a direct sum of indecomposable representations. However, this list can be very unwieldy. Gabriel's theorem provides an precise description of which quivers admit a finite list of indecomposable representations.

\begin{figure}
    \centering
    \includegraphics[width=.6\textwidth]{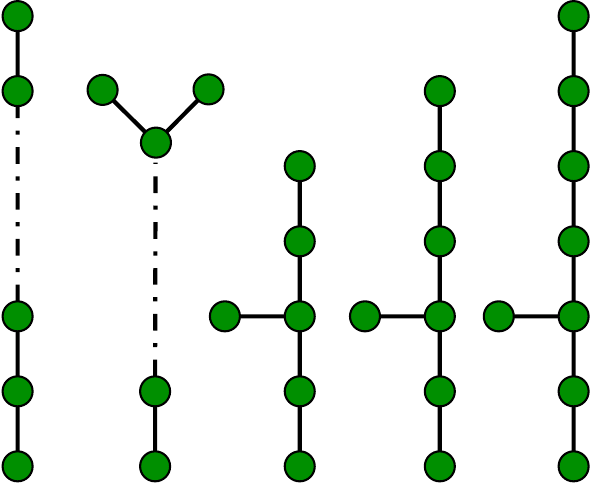}
    \caption{Simply Laced Dynkin Diagrams}
    \label{fig:sldd}
\end{figure}

\begin{thm}[Gabriel's Theorem~\cite{derksen}]\label{thm:gabriel}\index{Gabriel's theorem}
Let $Q$ denote a quiver. The category of representations $\Rep(Q)$ has finitely many indecomposables if and only if the underlying undirected graph of $Q$ is a union of Dynkin graphs of type $A_n$, $D_n$, $E_6$, $E_7$ or $E_8$. These are depicted in Figure \ref{fig:sldd}.
\end{thm}

\subsection{A Remark on Quivers and Perverse Sheaves}

From a quiver representation, one can always construct a cellular sheaf or cosheaf over a one-dimensional base space. This is done by turning every map
\[
	\xymatrix{F(s(e)) \ar[r]^{\rho_{t,s}} & F(t(e))}
\]
into one of the following diagrams:
\[
	\xymatrix{F(s(e)) \ar[r]^-{\rho_{t,s}} & F(t(e))=F(e) & \ar[l]_-{\id} F(t(e))} \qquad  \xymatrix{F(s(e)) & \ar[l]_-{\id} F(t(e))=F(e) \ar[r]^-{\rho_{t,s}} & F(t(e))}
\]
The former choice would make a quiver representation into a cellular sheaf, the latter into a cellular cosheaf.

There are dangers in trying to use quiver theory as a substitute for cellular sheaf or cosheaf theory. One might try to think of a poset as a certain type of quiver with vertices corresponding to elements and an edge between two elements if $s(e)\leq t(e)$. For example, consider the poset coming from the face relation of the cell complex $Y=[0,1)\times [0,1)$:
\[
	\xymatrix{ & \sigma & \\ a \ar[ur] & & b \ar[ul]\\ & x \ar[ul] \ar[uu] \ar[ur] &}
\]
A quiver representation produces a diagram of vector spaces
\[
\xymatrix{ & F(\sigma) & \\ F(a) \ar[ur] & & F(b) \ar[ul]\\ & F(x) \ar[ul] \ar[uu] \ar[ur] &}
\]
that does not commute. In contrast, if $F$ were a cellular sheaf, then the two triangles would commute. If we were to impose ``relations'' on the quiver representation by identifying different paths, then we could recover cellular sheaves. 

A relaxed version of this observation generates a combinatorial model for perverse sheaves, which was invented by Bob MacPherson~\cite{gmv} and explored by Maxim Vybornov~\cite{vybornov-mqa,vybornov-constr}.

\begin{defn}\index{perversity}
A \textbf{perversity} $p:\Z_{\geq 0} \to \Z$ is a function from the non-negative integers to the integers such that $p(0)=0$ and $p$ takes every interval $\{0,\ldots,k\}$ bijectively to an interval $\{a,\ldots,a+k\}$ where $a\in\Z_{\leq 0}$.
\end{defn}

\begin{defn}[Cellular Perverse Sheaves]\label{defn:cellular_perverse_sheaves}\index{cellular perverse sheaves}\index{sheaf!cellular perverse}
Let $X$ be a cell complex and $P_X$ its associated poset. Let $p:\Z_{\geq 0}\to\Z$ be a perversity. Define a quiver $Q_X$ whose vertices are the elements of $P_X$ and whose edges have the property that if $\tau\to\sigma$ is an edge, then $\sigma$ is incident to $\tau$ and $p(\dim\tau)=p(\dim\sigma)+1$. A \textbf{cellular perverse sheaf} assigns to each vertex of $Q_X$ a vector space $\cP(v)$ and to each edge from $\tau$ to $\sigma$ a linear map $\kappa_{\sigma,\tau}:\cP(\tau)\to\cP(\sigma)$. These maps satisfy the chain complex condition for any pair of vertices $\gamma,\tau$
\[
\sum_{\sigma} \kappa_{\gamma,\sigma}\circ \kappa_{\sigma,\tau}=0
\]
where $\sigma$ ranges over all vertices containing with an edge from $\tau$ and to $\lambda$.
\end{defn}

%
%

\chapter{The Derived Perspective}
\label{sec:derived}

\begin{quote}
{\em``The maker of a sentence launches out into the infinite and builds a road into Chaos and old Night, and is followed by those who hear him with something of wild, creative delight.''}
\begin{flushright} --- Ralph Waldo Emerson~\cite[p.59]{emerson-journals} \end{flushright}
\end{quote}

 \begin{figure}
    \centering
    \includegraphics[width=\textwidth]{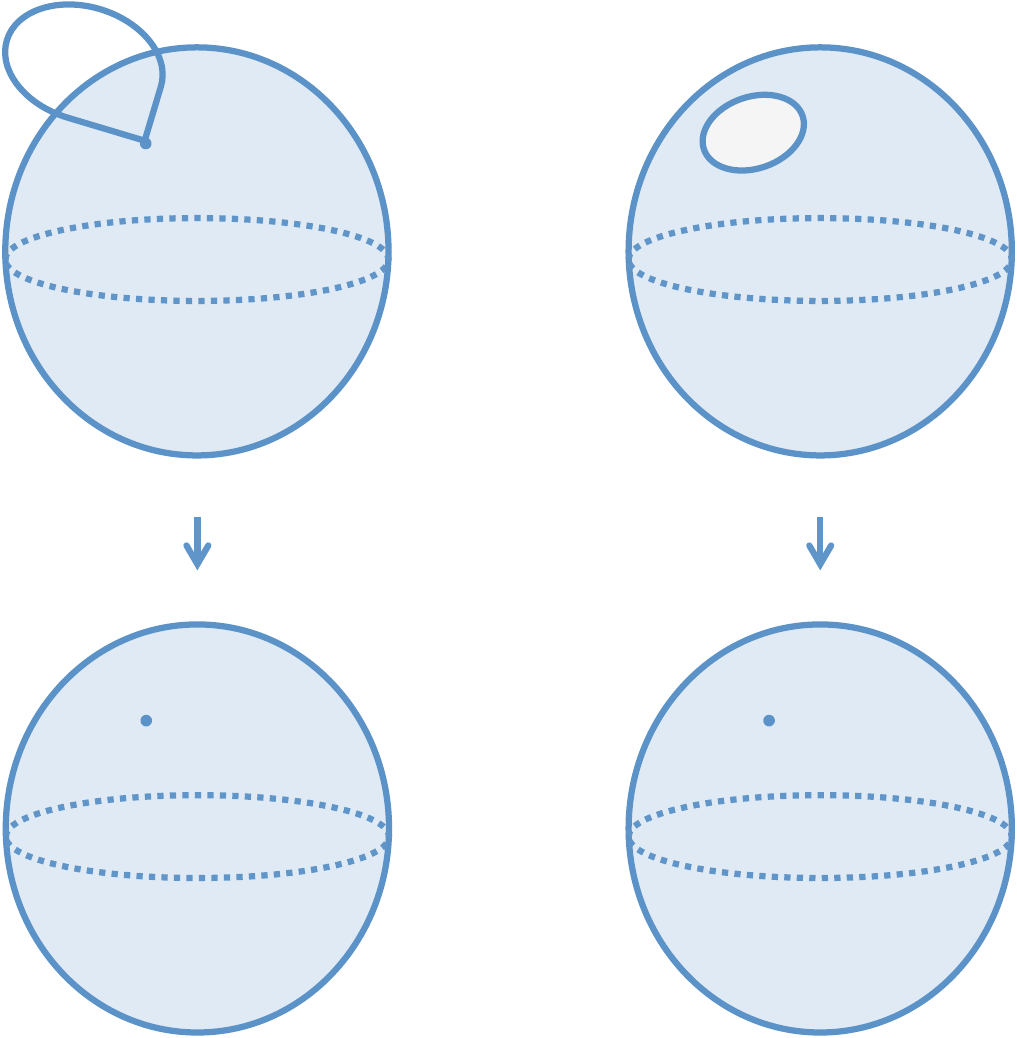}
    \caption{MacPherson's Motivating Example for the Derived Category}
    \label{fig:derived_category}
\end{figure}

The need for a derived perspective can be stated with one picture. In Figure~\ref{fig:derived_category} two maps are drawn to the two-sphere $S^2$. One is defined on the wedge sum $S^2\vee S^1$ and maps the $S^1$ to a point. The other is defined on the closed disk $\DD^2$ and maps the boundary circle to a point. If one is only allowed to look at the homology of the fiber for both of these maps, they will not be able to tell them apart. The derived category is the universal solution to this problem, as well as many others. 

In the derived category, one does not consider a sheaf in isolation, but rather one considers complexes of sheaves or, alternatively said, sheaves of complexes. In order to solve the problem presented by Figure~\ref{fig:derived_category}, one works with the sheaf of cochains on each fiber, along with their differentials. This transition is formally motivated via an analogy with Taylor series in Section \ref{subsec:taylor}. Injective and projective sheaves are introduced as the basic building blocks for the derived category, just as polynomials are the basis for Taylor series. Because the Alexandrov topology is so simple, we can described explicitly the elementary injective and projective (co)sheaves in Section \ref{subsubsec:elem_inj}. Injective and projective resolutions are then introduced in Section \ref{subsubsec:resolutions}.

Section \ref{subsec:htpy_chain} gives a high-level introduction to the homological algebra techniques necessary to understanding the derived category. The explicitness of cellular sheaves allows us to give concrete examples of what is usually taken on faith when first learning the subject. The notion that maps are unique up to homotopy and that sheaves can be ``quasi''-isomorphic without being isomorphic, are demonstrated in Examples \ref{ex:not_unique} and \ref{ex:quasi_iso}.

The derived definition of cosheaf homology is given in Section \ref{subsec:derived_sheaf_cohom} and the derived definition of sheaf cohomology can be dualized from there or looked up in Shepard's thesis~\cite{shepard}. These definitions should be regarded as the true definition of cosheaf homology and sheaf cohomology. The compactly supported variant, which we call Borel-Moore cosheaf homology, is defined in Section \ref{subsubsec:BM_cosheaf_homology}. The derived functor formalism allows us to resolve the question of invariance under subdivision in Section \ref{subsubsec:subdivision} with considerable ease.

Finally, we exploit the special features of the Alexandrov topology to develop two new theories: sheaf \emph{homology} and cosheaf \emph{cohomology}. Although these theories are invariant under subdivision in the domain of a map, they are not invariant under subdivision in the target of a map. These theories are sensitive to both the cell structure and the embedding. We compute some explicit examples of these theories in Section \ref{subsubsec:sheaf_homology_graphs}.

\section{Taylor Series for Sheaves}
\label{subsec:taylor}

When first learning about the derived perspective a helpful analogy might be the following. We can approximate suitably nice functions around a point via the use of Taylor series:
\[
	f(x)\simeq f(a) + f'(a)(x-a) + \frac{f''(a)}{2!}(x-a)^2+\cdots
\]  
The working physicist or engineer appreciates deeply how by only using a few terms, one can make serious headway into the analysis of integrals or other problems involving $f$.

In similar spirit one might start approximating or ``taking the Taylor series expansion'' of a topological space $X$ via its homotopy or homology groups:
\[
	\pi_0(X,x) \quad \pi_1(X,x) \quad \pi_2(X,x) \quad \cdots \qquad | \qquad \cdots \quad H_2(X) \quad H_1(X) \quad H_0(X)
\]
One should realize that both of these series expansions arise from more fundamental sequences:
\[
	X\to \Omega_x X \to \Omega_x^2 X \to \cdots\qquad | \qquad \cdots \to C_2(X;k)\to C_1(X;k) \to C_0(X;k)
\]
Here $\Omega_x X$ denotes the space of loops in $X$ based at $x$ (and iterated applications thereof) and $C_p(X;k)$ denotes the $p$-chains.

For a sheaf $F$ on a topological space $X$ one also has a similar process. Namely, there is an exact sequence called a \textbf{resolution}
\[
	0 \to F\to I^0 \to I^1 \to I^2 \to \cdots
\]
that when evaluated on an open set $U\subset X$ produces a sequence
\[
	0 \to F(U) \to I^0(U) \to I^1(U) \to I^2(U) \to \cdots
\]
that is exact at $F(U)$ and $I^0(U)$.\footnote{There is no guarantee for exactness at higher terms.} Like the physicist with their Taylor series, one can discard the original sheaf and work solely with the terms in the sequence $I^{\bullet}(U)$. This sequence is a chain complex with potentially interesting cohomology. For each $i$, these cohomologies piece together to provide a pre-sheaf description of the complex of sheaves $I^{\bullet}$ --- the derived replacement for $F$.
\[
	U \rightsquigarrow H^i(I^{\bullet}(U))=:H^i(U;F).
\]
If we specialize to the constant sheaf $F=k_X$, then we obtain another familiar series expansion of the space $X$: the cohomology. However, this series is much more general, as it encodes the cohomology of each open set in $X$. Consequently, even if one embeds $X$ into the contractible cone $CX$, the constant sheaf and its derived replacement will remember the topology on $X$.

However, just as the reason that Taylor series are amenable to analysis because polynomials have simple properties, for general sheaves we must develop an algebraic analogue of a polynomial, which are the injective sheaves.

\subsection{Elementary Injectives and Projectives}
\label{subsubsec:elem_inj}

In this section we consider the basic building blocks of the derived category, which are injective or projective objects. These objects are characterized by universal mapping properties. For cellular sheaves and cosheaves the injectives and projectives can be described explicitly.

\subsubsection{Injectives}
\begin{defn}\index{injective!objects}
 A representation of a small category $I:\cat\to\Vect$ is \textbf{injective} if, for any natural transformation $\eta:A\to I$ and any injection $\iota:A\hookrightarrow B$, there is an extension $\tilde{\eta}:B\to I$ such that $\eta=\tilde{\eta}\circ\iota$. Said using diagrams, they are characterized by the usual universal property:
\[
 \xymatrix{0 \ar[r] & A \ar[r]^{\iota} \ar[d]_{\eta} & B \ar@{.>}[dl]^{\exists \tilde{\eta}} \\ & I &}
\]
\end{defn}

\begin{exr}\label{exr:inj_props}\index{injective!properties}
Use the universal property of an injective representation to prove the following statements:
 \begin{itemize}
  \item Every short exact sequence \[0\to I\to A \to Q \to 0,\] where $I$ is injective, splits, i.e. if $\iota: I\hookrightarrow A$ is an inclusion, then $A\cong I\oplus \cok(\iota)$.
  \item $\prod A_i$ is injective if and only if each $A_i$ is injective.
 \end{itemize}
 We will make use of these properties in Lemma~\ref{lem:inj}.
\end{exr}

For our first example of an injective representation, we consider an injective cell sheaf. These sheaves are supported on the closures of cells.

\begin{defn}\index{sheaf!elementary injective cellular sheaf $[\sigma]^W$}\index{injective!cellular sheaf}
 An \textbf{elementary injective cell sheaf} on $X$ concentrated on $\sigma\in X$ with value $W\in\Vect$ is given by
\begin{equation*}
[\sigma]^W(\tau)=
\begin{cases}
W & \text{if } \tau\leq \sigma,\\
0 & \text{other wise.}
\end{cases}
\end{equation*}
where the only possible non-zero restriction maps are the identity.
\end{defn}

In order to prove that this sheaf is actually injective we introduce an alternative definition of injective sheaves and cosheaves defined on arbitrary posets. This definition makes use of the functors $f_*$ and $f_{\dd}$.

\begin{defn}\index{injective!sheaf on a poset}\index{sheaf!elementary injective on a poset}\index{injective!cosheaf on a poset}\index{cosheaf!elementary injective on a poset}
 Let $i_x:\star \to X$ be the map that assigns to the one element poset the value $x\in X$, i.e. $x=i_x(\star)$. Define the \textbf{elementary injective sheaf} on $x\in X$ with value $W\in\Vect$ to be  $[x]^W=(i_x)_*W$ and the corresponding elementary injective cosheaf to be $\{\hat{x}\}^W:=(i_x)_{\dd}\hat{W}$. 
\end{defn} 

One can see that for cosheaves, the elementary injectives are concentrated on the open stars of cells. To prove these objects are actually injective we make use of the adjunctions already defined.

\begin{lem}
	The sheaf $[x]^W=(i_x)_*W$ and cosheaf $\{\hat{x}\}^W:=(i_x)_{\dd}\hat{W}$ are injective.
\end{lem}
\begin{proof}
	The proof is immediate from the following adjunctions
	\[
		\Hom_{\Shv(X)}(A,(i_x)_*W)\cong \Hom_{\Vect} (A(x),W)
	\]
	\[ 
		\Hom_{\Coshv(X)}(\hat{A},(i_x)_{\dd}\hat{W})\cong \Hom_{\Vect}(\hat{A}(x),\hat{W})
	\]
	and the fact that in the category of vector spaces every object is injective.
\end{proof} 

There is one last lemma that tells us that the considering sheaves of the form $[x]^W$ suffices for understanding all injective sheaves over a cell complex. The proof is presented in~\cite[Thm.~1.3.2, p.19-20]{shepard}, but we will give the direction needed for our derived equivalence proof given in Theorem~\ref{thm:equivalence}.

\begin{lem}[All Injectives are Sums]\label{lem:inj}\index{injective!cellular sheaves are sums}
Let $X$ be the face-relation poset of a cell complex. A sheaf $I$ is injective if and only if it is isomorphic to a one of the form $\oplus_{\sigma}[\sigma]^{V_{\sigma}}$.
\end{lem}
\begin{proof}
One can easily check that the direct sum of elementary injective sheaves is injective, so this proves the `if' direction.

To prove that every injective is isomorphic to a direct sum of injectives --- the `only if' direction --- requires a little work. Assume for induction that every injective sheaf $I$ that is non-zero on at most $k\leq n-1$ cells is isomorphic to $\oplus_{\sigma}[\sigma]^{V_{\sigma}}$. Now consider an injective sheaf that is non-zero on exactly $n$ cells. Let $\sigma$ be a cell of maximal dimension where $I(\sigma)=:V\neq 0$. Since $I$ is zero on all higher cells incidence to $\sigma$, there is a non-zero map $\eta$ from the skyscraper sheaf $S^V_{\sigma}$ to $I$ with $\eta(\sigma)=\id_V$. There is also a non-zero map $\iota:S_{\sigma}^V\to[\sigma]^V$. This gives us a diagram
\[
 \xymatrix{0 \ar[r] & S_{\sigma}^V \ar[r]^{\iota} \ar[d]_{\eta} & [\sigma]^V \ar@{.>}[dl]^{\exists \tilde{\eta}} \\ & I &}
\]
and the implicated existence of a map $\tilde{\eta}:[\sigma]^V\to I$. If $\tau\leq \sigma$, then by the fact that $\tilde{\eta}$ is a sheaf map,
\[
\id_V=\tilde{\eta}(\sigma)\circ \rho^{[\sigma]}_{\sigma,\tau}=\rho^{I}_{\sigma,\tau}\circ \tilde{\eta}(\tau)
\]
the map $\tilde{\eta}$ is injective. By the second property of Exercise~\ref{exr:inj_props}, we can deduce that $I\cong [\sigma]^V\oplus \cok(\tilde{\eta})$. Since $\cok(\tilde{\eta})$ is zero wherever $I$ is and also zero on $\sigma$, it is non-zero on at most $n-1$ cells and the induction hypothesis applies. The zero sheaf is clearly equal to a direct sum of elementary injectives with the zero vector space, which checks the base case, completing the induction.
\end{proof}

\subsubsection{Projectives}

There is a dual universal object that is called projective.
\begin{defn}\index{projective!objects}
 A representation of a small category $P$ is \textbf{projective} if for any natural transformation $\epsilon:P\to A$ and any surjection $\pi:B\twoheadrightarrow A$, there is a map $\tilde{\epsilon}:P\to B$ such that $\epsilon=\pi\circ \tilde{\epsilon}$. Said using diagrams:
\[
 \xymatrix{ & P\ar[d]^{\epsilon} \ar@{.>}[dl]_{\exists \tilde{\epsilon}} & \\ B \ar[r]_{\pi} &  A \ar[r] & 0}
\]
\end{defn}

As before, we have some dual consequences:
\begin{itemize}\index{projective!properties}
 \item Every short exact sequence \[0\to A \to B \to P \to 0\] where $P$ is projective, splits.
 \item $\oplus B_i$ is projective if and only if each $B_i$ is projective.
\end{itemize}

Since the adjunctions will be our guide we make the following definitions.

\begin{defn}\index{projective!sheaf on a poset}\index{sheaf!elementary projective on a poset}\index{projective!cosheaf on a poset}\index{cosheaf!elementary projective on a poset}
 Let $i_x:\star \to X$ be the map that assigns to the one element poset the value $x\in X$, i.e. $x=i_x(\star)$. Define the elementary projective sheaf on $x\in X$ with value $W\in\Vect$ to be  $\{x\}^W=(i_x)_{\dd}W$ and the corresponding elementary projective cosheaf to be $[\hat{x}]^W:=(i_x)_{*}\hat{W}$. 
\end{defn}

We leave it to the reader to check that these objects are actually projective.

\subsubsection{Mapping Identities}
\index{injective!mapping identities for (co)sheaves}
\index{projective!mapping identities for (co)sheaves}

Before moving on to the derived definition of sheaf cohomology, we record some useful identities that should be evident from the definition and the adjunctions.

\begin{equation*}
 \Hom_{\Shv}([\tau]^U,[\sigma]^W)=
\begin{cases}
 \Hom_{\Vect}(U,W) & \text{if } \sigma \leq \tau, \\
0 & \text{o.w.}
\end{cases}
\end{equation*}

\begin{equation*}
 \Hom_{\Shv}(\{\tau\}^U,\{\sigma\}^W)=
\begin{cases}
 \Hom_{\Vect}(U,W) & \text{if } \sigma \leq \tau, \\
0 & \text{o.w.}
\end{cases}
\end{equation*}

\begin{equation*}
 \Hom_{\Coshv}([\hat{\sigma}]^W,[\hat{\tau}]^U)=
\begin{cases}
 \Hom_{\Vect}(W,U) & \text{if } \sigma \leq \tau, \\
0 & \text{o.w.}
\end{cases}
\end{equation*}

\begin{equation*}
 \Hom_{\Coshv}(\{\hat{\sigma}\}^W,\{\hat{\tau}\}^U)=
\begin{cases}
 \Hom_{\Vect}(W,U) & \text{if } \sigma \leq \tau, \\
0 & \text{o.w.}
\end{cases}
\end{equation*}

\subsection{Injective and Projective Resolutions}
\label{subsubsec:resolutions}
\index{injective!resolutions}\index{projective!resolutions}

As promised, we aim to prove every sheaf has a resolution by injective sheaves. This follows from the following claim, which we now prove. Although this theorem is true for general spaces, we work with Alexandrov spaces arising as posets as usual.

\begin{clm}
	Every sheaf $F:X\to\Vect$ on a poset $(X,\leq)$ possibly of infinite size, $F$ admits an inclusion into an injective sheaf. Dually, every cosheaf admits is surjected onto by a projective cosheaf.
	\[
		0 \to F \to I
	\]
\end{clm}
\begin{proof}
	We construct $I$ explicitly. It is given by
	\[
		0 \to F \to I:=\prod_x [x]^{F(x)} = \prod_x (i_x)_*F(x).
	\]
	The map to $I$ is defined easily using the standard adjunctions
	\[
		\iota\in \Hom(F,\prod_x(i_x)_*F(x))\cong \prod_x \Hom(F,\prod_x(i_x)_*F(x)) \cong \prod_x \Hom(F(x),F(x)) \ni \prod_x \id_{F(x)}.
	\]
	We encourage the reader to describe this map is explicitly, by seeing how a single $\id_{F(x)}$ traces through this adjunction, which we'll call $\iota_x\in \Hom(F,(i_x)_*F(x))$.
	
	Similarly, for a cosheaf $\hF:X^{op}\to\Vect$ on an Alexandrov space we could have built a projective surjection by taking
	\[
	 \hat{P}_0:=\bigoplus_{x} (i_x)_*\hF(x) \to \hF \to 0
	\]
	where the map $\pi_0:P_0\to \hF$ is gotten through the corresponding adjunction for cosheaves and using the contravariance of $\Hom$ in the first slot
	\[
		\Hom(\bigoplus_{x} (i_x)_*\hF(x),\hF)\cong \prod \Hom((i_x)_*\hF(x),\hF)\cong \prod \Hom(\hF(x),\hF(x)).
	\]
\end{proof}

\begin{cor}
	Every sheaf $F:X\to\Vect$ has an injective resolution. Dually, every cosheaf $\hF:X^{op}\to\Vect$ has a projective resolution.
\end{cor}
\begin{proof}
	Since cokernels exist in the category of sheaves by taking element-by-element quotients and by iteratively applying the claim, we obtain an injective resolution of $F$:
	\[
	 \xymatrix{F \ar[r]^{\iota^0} & I^0 \ar[rd]_{\pi_0} \ar[rr]^{\iota^1=j_1\pi_0} & & I^1 \ar[rr]^{\iota^2=j_2\pi_1} \ar[rd]_{\pi_1} & & I^2 &\cdots \\
		      & & \cok(\iota^0) \ar[ru]_{j_1} & & \cok(\iota^1) \ar[ru]_{j_2} & &\cdots }
	\]
	
	Iterating the analogous process for projective cosheaves, replacing kernels where one sees cokernels above, one obtains an exact sequence of cosheaves called the \textbf{projective resolution} of $\hF$:
	\[
	 \cdots \hat{P}_2 \to \hat{P}_1 \to \hat{P}_0 \to \hF \to 0.
	\]
\end{proof}

These exact sequences can be used to replace $F$ or $\hF$ in a suitable sense, defined by the derived category. Before moving onto that discussion, we note one interesting point.

\begin{prop}
 The length of injective resolution of any sheaf $F\in\Shv(X)$ is bounded by the length of longest chain in the poset. In particular for $X$ a cell complex, it is bounded by the dimension.
\end{prop}
\begin{proof}
 Pick a maximal ordered subset in $X$ and consider its top element, say $x'$, then $I^0(x')=F(x')$ since nothing is larger than $x'$. The cokernel sheaf of $\iota^0$ evaluated on $x'$ is then $\cok(\id:F(x')\to F(x'))=0$. So for any maximally ordered chain in $X$, $I^1$ is zero on the top-most element. Arguing inductively finishes the proof.
\end{proof}

\section{The Derived Category and Homotopy Theory of Chain Complexes}
\label{subsec:htpy_chain}

The purpose of the derived category is to replace the category of sheaves with a category of complexes where certain operations are more natural. We have already shown that one can replace a sheaf by its injective resolution and a cosheaf by its projective resolution. This will define our derived replacement on the level of objects, but we have not yet shown how a map of sheaves or cosheaves induces a map on the level of resolutions. 

If $\phi:\hF\to\hG$ is a map of cosheaves, then it can be checked from the universal properties of projective objects, that this induces a map of complexes
\[
	\xymatrix{ \cdots & \hat{P}_2 \ar[r] \ar@{.>}[d]^{\phi_2} & \hat{P}_1 \ar[r] \ar@{.>}[d]^{\phi_1} & \hat{P}_0 \ar[r] \ar@{.>}[d]^{\phi_0} & \hF \ar[d]^{\phi} \ar[r] & 0\\ \cdots & \hat{Q}_2 \ar[r] & \hat{Q}_1 \ar[r] & \hat{Q}_0 \ar[r] & \hG \ar[r] & 0}
\]
where all the squares in sight commute. For a hint on how to see this, consider the composite map $\hat{P}_0\to\hF\to\hG$ and let $\hG=A$ and $B=\hat{Q}_0$ in the definition of the universal property defining a projective object. This induces our first map $\hat{P}_0\to\hat{Q}_0$. To get the next, all important step, one must recognize that having maps from $\hat{P}_0\to \hat{Q}_0$ and $\hF\to\hG$ induces maps between the kernels of the map $\hat{P}_0\to\hF$ and $\hat{Q}_0\to\hG$. Since $\hat{Q}_1$ surjects onto the kernel of the latter map repeating the initial argument provides a map from $\hat{P}_1$ to $\hat{Q}_1$. This shows that the projective replacement of cosheaves is functorial.

Aside from functoriality, there is one more snag that needs to be mentioned: For a sheaf or a cosheaf it is possible that the choice of injective or projective resolution is not unique. If one really wants to use these as replacements for the original sheaf or cosheaf, there must be a strong relationship between these two complexes. This is best seen by specializing the functoriality discussion above to the case $\phi=\id$.
\[
	\xymatrix{ \cdots & \hat{P}_2 \ar[r] \ar@{.>}[d]^{\phi_2} & \hat{P}_1 \ar[r] \ar@{.>}[d]^{\phi_1} & \hat{P}_0 \ar[r] \ar@{.>}[d]^{\phi_0} & \hF \ar[d]^{\id} \ar[r] & 0\\
	\cdots & \hat{Q}_2 \ar[r] & \hat{Q}_1 \ar[r] & \hat{Q}_0 \ar[r] & \hF \ar[r] & 0}
\]
The resulting map of complexes need not be a term-by-term isomorphism with all squares in sight commuting, but rather a more general notion must be substituted, namely the definition of chain homotopy. Before giving that, let us give an example.

\begin{ex}[Non-Unique Projective Resolutions]
	Let us work again over our test space of the closed unit interval $X=[0,1]$ stratified as $x=0$, $y=1$ and $a=(0,1)$. The constant cosheaf $\hat{k}_X$ is then modeled as
	\[
		\xymatrix{ & \ar[ld]_1 k \ar[rd]^1 & \\ k & & k}
	\]
	Revisiting the definition of the elementary projective cosheaves, there is one obvious projective resolution because the constant cosheaf on this stratification of the unit interval is already projective, so we have the identity map
	\[
		\hat{P}_{\bullet}: \qquad [\hat{a}] \to \hat{k}_X.
	\]
	On the other hand, following blindly the prescription provided for computing the projective resolution of an arbitrary cellular cosheaf would have lead us to the following ``canonical'' resolution:
	\[
		\hat{Q}_{\bullet}: \qquad [\hat{x}]\oplus [\hat{y}]\to [\hat{a}]\oplus [\hat{x}]\oplus [\hat{y}] \to \hat{k}_X
	\]
\end{ex}

\begin{defn}[Chain Homotopy]\index{Cochain@(co)chain complex!homotopy}
	Suppose $(A^{\bullet},d_A)$ and $(B^{\bullet},d_B)$ are two (cohomological) chain complexes and $\phi^{\bullet}$ and $\psi^{\bullet}$ are two chain maps, then \textbf{a chain homotopy} $h^{\bullet}$ is a chain map $h^i:A^{i}\to B^{i-1}$
	\[
		\xymatrix{\cdots\ar[r] & A^{i-2} \ar[d] \ar[r] & \ar[ld] A^{i-1} \ar[r] \ar[d]^{\psi^{i-1}}_{\phi^{i-1}} & \ar[ld] A^i \ar[r] \ar[d]^{\psi^i}_{\phi^i} & \ar[ld] A^{i+1} \ar[d]^{\psi^{i+1}}_{\phi^{i+1}} \ar[r] & \cdots \\
		\cdots \ar[r] & B^{i-2} \ar[r] & B^{i-1} \ar[r] & B^i \ar[r] & B^{i+1} \ar[r] & \cdots}
	\]
	such that 
	\[ 
		\phi^i-\psi^i=h^{i+1}d^i_A+d^{i-1}_Bh^i.
	\]
	In which case we say that $\phi\sim\psi$ are \textbf{chain homotopic}. 
	
	Consider now two chain maps $\phi:A^{\bullet}\to B^{\bullet}$ and $\psi:B^{\bullet}\to A^{\bullet}$, such that
	\[
		\phi\circ\psi \sim \id \qquad \mathrm{and} \qquad \psi\circ\phi \sim \id
	\]
	then one says $A^{\bullet}$ and $B^{\bullet}$ are \textbf{chain homotopy equivalent}.
\end{defn}

\begin{ex}[Non-unique, but equivalent]\label{ex:not_unique}
	Consider again the case of the two different projective resolutions of the constant sheaf $\hat{k}_X$ on the closed unit interval. On the one hand the composite
	\[
		\xymatrix{ 0 \ar[d] \ar[r] & [\hat{a}] \ar[d] \\
		[\hat{x}]\oplus [\hat{y}] \ar[d] \ar[r] & [\hat{x}]\oplus [\hat{a}]\oplus [\hat{y}] \ar[d] \\
		0 \ar[r] & [\hat{a}] }
	\]
	is clearly the identity on $\hat{P}_{\bullet}$, but the composite
	\[
		\xymatrix{
		[\hat{x}]\oplus [\hat{y}] \ar[d] \ar[r] & [\hat{x}]\oplus [\hat{a}]\oplus [\hat{y}] \ar[d] \\
		 0 \ar[d] \ar[r] & [\hat{a}] \ar[d] \\
		[\hat{x}]\oplus [\hat{y}] \ar[r] & [\hat{x}]\oplus [\hat{a}]\oplus [\hat{y}] }
	\]
	cannot possibly be the identity because one map factors through zero. However, if we employ a self-homotopy of $\hat{Q}_{\bullet}$ by defining a homotopy for the only possible degree to be
	\[
		h^0:[\hat{a}]\oplus[\hat{x}]\oplus[\hat{y}] \to  [\hat{x}]\oplus[\hat{y}]
	\]
	which is zero on the $a$ component and the identity elsewhere. One can then check that this defines a homotopy between the identity map and the map indicated in the second composite.
\end{ex}

The conclusion from the example should be that although one can use different projective resolutions, the choice is irrelevant up to homotopy. The derived category should not be able to discriminate between them. As such, we make the following definitions.

\begin{defn}\index{category!C@$C^b(\aat)$ of chain complexes}\index{category!K@$K^b(\aat)$ homotopy category of chain complexes}\index{H@$K^b(\aat)$ homotopy category of chain complexes}
	Let $\aat$ be an abelian category, such as the category of sheaves or cosheaves. The \textbf{category of chain complexes} in $\aat$, written $C^b(\aat)$ has objects that are chain complexes and morphisms that are chain maps.
	
	The \textbf{homotopy category of complexes} $K^b(\aat)$ of an abelian category $\aat$ has the same objects as $C^b(\aat)$, but where we have identified chain homotopic maps.
\end{defn}

\begin{defn}\index{Derived@$D^b(\aat)$ derived category}\index{category!Derived@$D^b(\aat)$ derived category}
	For $\aat=\Shv(X)$ we define the \textbf{bounded derived category of sheaves} $D^b(\Shv(X))$ to be $K^b(\mathrm{Inj}-\Shv(X))$ the homotopy category that uses only complexes of injective sheaves. 
	
	Similarly, for $\aat=\Coshv(X)$, we define the \textbf{bounded derived category of cosheaves} $D^b(\Coshv(X))$ by $K^b(\mathrm{Proj}-\Coshv(X))$ where complexes of projective cosheaves are used instead.
\end{defn}

This definition, is an equivalent reformulation of another definition of the derived category. This other perspective is built on the foundational notion of a quasi-isomorphism, which is in turn built on the idea of a cohomology sheaf or homology cosheaf.

\begin{defn}
	Suppose we are given a complex of cellular sheaves
	\[
		(F^{\bullet},d^{\bullet}): \qquad \cdots \to F^{i-1} \to F^{i} \to F^{i+1} \to \cdots,
	\]
	i.e. for each cell $\sigma$ we have a complex of vector spaces. For each $i$ we can define the $i$th \textbf{cohomology sheaf} as the assignment
	\[
		\cohomsheaf^i(F^{\bullet}): \qquad \sigma \rightsquigarrow H^i(F^{\bullet}(\sigma))
	\]
	which is a cellular sheaf. The restriction maps being defined as the induced map on cohomology for the chain map $F^{\bullet}(\sigma)\to F^{\bullet}(\tau)$ for $\sigma\leq\tau$. 
	
	Considering all $i$ at once defines a functor from the category of complexes of sheaves and the category of graded sheaves (sheaves of graded vector spaces with level preserving restriction maps)
	\[
		\cohomsheaf^*: C^b(\Shv(X))\to \Shv(X;\mathrm{gr}\Vect) \qquad F^{\bullet} \rightsquigarrow \bigoplus_i \cohomsheaf^i(F^{\bullet}).
	\]
	There are completely dual notions of \textbf{homology cosheaves}, where we generally use homological indexing and notation $(\hF_{\bullet},\partial_{\bullet})$.
\end{defn}

\begin{defn}[Quasi-Isomorphisms]\index{quasi-isomorphism}
	A map of complexes of sheaves (or cosheaves) $\alpha^{\bullet}:F^{\bullet}\to G^{\bullet}$ such that the induced map
	\[
		\cohomsheaf(\alpha^{\bullet}):\cohomsheaf^i(F^{\bullet})\to \cohomsheaf^i(G^{\bullet})
	\]
	is an isomorphisms for every $i$, is called a \textbf{quasi-isomorphism}.
\end{defn}
	
The term ``quasi-isomorphism'' reflects the fact that if $\alpha^{\bullet}:F^{\bullet}\to G^{\bullet}$ is a quasi-isomorphism, then there does not always exist an inverse map $\beta^{\bullet}$ that gives the identity, or even chain homotopic to the identity, any map back may simply not exist.
	
\begin{ex}\label{ex:quasi_iso}
	Consider again the unit interval $X=[0,1]$ decomposed into two vertices $x$ and $y$ and an open interval $a$. Consider the stalk sheaf $S_a$ that assigns $k$ to $a$ and is zero everywhere else. It's injective resolution defines a chain map
	\[
	\xymatrix{ 0\ar[r] \ar[d] & S_a \ar[r] \ar[d] & 0 \ar[d] \\
	0  \ar[r] & [a] \ar[r] & [x]\oplus [y] }
	\]
	which is a quasi-isomorphism. However, there does not exist a map of sheaves $[a]\to S_a$.
\end{ex}

The slogan most commonly associated with the derived category is that one ``formally inverts the quasi-isomorphisms.'' This is formalized by the process of localizing categories. Namely, if $Q$ is a collection of morphisms in $\bat$ that is closed under certain operations, then we can consider the following universal problem: suppose $L:\bat\to\cat$ is a functor such that if $\alpha\in Q$, then $L(\alpha)$ is an isomorphism, then every such functor factors through the \textbf{category localized at $Q$}, written $\bat[Q^{-1}]$.
\[
	\xymatrix{\bat \ar[rr]^L \ar[rd] & & \cat \\
	& \bat[Q^{-1}] \ar@{.>}[ur]_{\exists} & }
\]
An alternative approach to the derived category of an abelian category $\aat$ is to define
\[
	D(\aat):= K(\aat)[Q^{-1}] \qquad Q=\{\mathrm{quasi-isomorphisms}\}
\]
where we have removed the boundedness hypothesis.

One then proves the following claim to re-obtain the definition we provided here
\begin{thm}[\cite{aluffi} Thm 6.7]\label{thm:proj_or_inj}
	Suppose $A$ is an abelian category with enough projectives, then $D^-(A)\cong K^-(P)$ where $P$ denotes projective objects of $A$. Similarly, if $A$ has enough injectives then $D^+(A)\cong K^+(I)$ where $I$ denotes injective objects of $A$.
\end{thm}

\section{The Derived Definition of Cosheaf Homology and Sheaf Cohomology}
\label{subsec:derived_sheaf_cohom}

We are now in a position to give the derived definition of cosheaf homology and show that it agrees with the computational formula provided earlier. This discussion can be dualized and readily found in the literature. The proof that the formula for sheaves computes the cohomology as defined by taking an injective resolution and applying $\Gamma(X;-)=p_*$ can be found in~\cite{shepard} pp. 28-29. Let's more or less repeat the proof for cellular cosheaves since it is nowhere in the literature.

\begin{defn}\index{Left@$Lf_*$ left derived pushforward}
	Given a cosheaf $\hF$ on $X$ we define the \textbf{left derived pushforward} along $f:X\to Y$ by taking a projective resolution and applying pushforward term by term:
	\[
		Lf_*\hF:=f_*P_{\bullet}.
	\]
	We define the $i$th derived functor by
	\[
		L_if_*\hF:= \cohomsheaf_i(f_*P_{\bullet}).
	\]
	In the special case where $f=p:X\to\star$ we write
	\[
		H_i(X;\hF):=L_ip_*\hF
	\]
	for the $i$th cosheaf homology group of $\hF$.
\end{defn}

We now aim to prove the following theorem.
\begin{thm}\index{cosheaf!cellular homology!Left@$Lf_*$ left derived pushforward}
 Let $p:X\to\star$ be the constant map and $\hat{F}$ a cellular cosheaf on $X$. Then the left derived functors of $p_{*}$ agree with the computational formula for homology, i.e. $L_i p_{*}\hat{F}=H_i(X;\hat{F})$.
\end{thm}
\begin{proof}
Begin with a projective resolution of $\hat{P}_{\bullet}\to\hat{F}$ and then take cellular chains of each cosheaf to obtain the following double complex:
\[
 \xymatrix{ & \vdots & \vdots & \vdots & \\
\cdots &  C_1(X;\hat{P}_1) \ar[r]\ar[d] & C_1(X;\hat{P}_0) \ar[r]\ar[d] & C_1(X;\hat{F}) \ar[r]\ar[d] & 0 \\
\cdots & C_0(X;\hat{P}_1)\ar[r]\ar[d] & C_0(X;\hat{P}_0)\ar[r]\ar[d] & C_0(X;\hat{F})\ar[r]\ar[d] & 0 \\
\cdots & \colim \hat{P}_1 \ar[r] \ar[d]& \colim \hat{P}_0\ar[r]\ar[d] & \colim \hat{F}\ar[r]\ar[d] & 0 \\
& 0 & 0 & 0 &}
\]

Now we make use of the following two observations, which dualize~\cite[Thm. 1.3.10, 1.4.1]{shepard}.
\begin{lem}
 For $\hat{P}$ a projective cosheaf 
$$H_p(C_{\bullet}^{BM}(X;\hat{P}))\cong H_p(C_{\bullet}(X;\hat{P}))\cong 0$$ for $p>0$.
\end{lem}
\begin{proof}
 Observe that we can assume that $\hat{P}$ is an elementary projective co-sheaf with value $\RR$, i.e. $[\hat{\sigma}]$, since $C_{\bullet}^{BM}(X;\oplus A_i)=\oplus C_{\bullet}^{BM}(X;A_i)$.

Everything follows from the following consequence of our definition of a cell complex: In the one-point compactification of $X$, the closure of any cell $\sigma\in X$, call it $|\bar{\bar{\sigma}}|$, has the homeomorphism type of a closed $k$-simplex.

$C_{\bullet}(X;[\hat{\sigma}])$ is the chain complex that computes the cellular homology of $Y=|\{\tau\leq \sigma|\bar{\tau}\,\mathrm{is}\,\mathrm{compact}\}|$, which is a closed $k$-simplex minus the star of a vertex. On the other hand, $C^{BM}_{\bullet}(X;[\hat{\sigma}])$ is equal to the chain complex calculating the cellular homology of $|\bar{\bar{\sigma}}|$ except in degree zero if $|\bar{\sigma}|$ is not compact. Notice that $H_1$ for both of these complexes is the same, as $|\bar{\sigma}|$ and $|\bar{\bar{\sigma}}|$ are simply connected. This proves the claim.
\end{proof}

\begin{lem}
 For any cellular cosheaf $\hat{F}$ on a cell complex $X$ we have that $$\colim \hat{F}\cong\cok(C_1(X;\hat{F})\to C_0(X;\hat{F})).$$
\end{lem}
\begin{proof}
 First let us prove that taking the coproduct of $\hat{F}$ over all the cells obtains a vector space that surjects onto the colimit. As part of the definition of $\colim \hat{F}$ is a choice of maps $\psi_{\sigma}:\hat{F}(\sigma)\to\colim \hat{F}$.  Let $\Psi=\oplus \psi_{\sigma}:\oplus \hat{F}(\sigma)\to \colim \hat{F}$, now consider the factorization of this map through the image:
\[
 \xymatrix{\oplus \hat{F}(\sigma) \ar[rr]^{\Psi} \ar[rd] & & \colim \hat{F} \\
& \im \Psi \ar[ru]^j &}
\]
Now we can use the $\im \Psi$ to define a new co-cone over the diagram $\hat{F}$ simply by pre-composing the factorized map with the inclusions $i_{\sigma}:\hat{F}(\sigma)\to\oplus \hat{F}(\sigma)$. Since the colimit is the initial object in the category of co-cones, there must be a map $u:\colim \hat{F}\to \im \Psi$ and thus $u\circ j=\id$ since there is only one map $\colim\hat{F}\to\colim\hat{F}$.

Now observe that $C_0(X;\hat{F})=\oplus \hat{F}(v_i)$ surjects onto the colimit of $\hat{F}$ by virtue of the fact that since every cell $\sigma\in X$ has at least one vertex as a face, the map $\Psi$ factors through $\oplus \hat{F}(v_i)$. Thus there is a surjection from $\Psi':C_0(X;\hat{F})\to \colim\hat{F}$. Notice that by universal properties of the cokernel of $\partial_0:C_1(X;\hat{F})\to C_0(X;\hat{F})$ it suffices to check that $\Psi'\circ\partial_0=0$. However, this is clear since every $e$ edge has two vertices $v_1$ and $v_2$ (we've discarded all those edges without compact closures), then we need only check the claim for each diagram of the form
\[
 \xymatrix{ & \hat{F}(e) \ar[ld]_{r_{e,v_1}} \ar[rd]^{r_{e,v_2}} & \\
\hat{F}(v_1) & & \hat{F}(v_2)}
\]
where it is clear that the colimit can be written as $\hat{F}(v_1)\oplus \hat{F}(v_2)$ modulo the equivalence relation $(r_{e,v_1}(w),0)\simeq (0,r_{e,v_2}(w))$, i.e. $\partial_0|_e(w)=(-r_{e,v_1}(w),r_{e,v_2}(w))\simeq(0,0)$.
\end{proof}

From these two theorems we can conclude that the columns away from the chain complex of $\hat{F}$ are exact and thus $\mathrm{Tot}_{\bullet}(C_i(X;\hat{P}_j))$ induces quasi-isomorphisms between $\colim \hat{P}_{\bullet}$ and $C_{\bullet}(X;\hat{F})$. We have thus established the theorem.
\end{proof}

\subsection{Borel-Moore Cosheaf Homology}
\label{subsubsec:BM_cosheaf_homology}
\index{Borel-Moore cosheaf homology}\index{cosheaf!cellular!Borel-Moore homology}
\begin{defn}\label{defn:BM_cosheaf_homology}
 Suppose $\hF$ is a cellular cosheaf. Define $\Gamma^{BM}(X;F)$ to be the colimit of the diagram extended over the one-point compactification of $X$ where we define $\hF(\infty)=0$. Alternatively said we look at the inclusion $j:X\to X\cup\{\infty\}$ and define
\[
	\Gamma^{BM}(X;\hF):=p_*j_{\dd}\hF.
\]
Another possible definition is to dualize a cellular cosheaf of finite-dimensional vector spaces to a cellular sheaf by post-composing $\hF:X\to\Vect$ with $\Hom_{\vect}(-,k)$, apply $p_!$ and then dualize back.
\end{defn}

\begin{rmk}[Functoriality]
	The definitions that involve the one-point compactification are deficient in the following way. A map of cell complexes $f:X\to Y$ does not necessarily extend to a map between the one-point compactifications. It is for this reason that for functoriality, the definition using $p_!$ is preferred.
\end{rmk}

Now we can prove that the formula provided calculates the Borel-Moore homology of a cosheaf $\hat{F}$ by establishing the following lemma:

\begin{lem}\label{lem:BM_cosheaf_homology}
 For any cellular cosheaf $\hat{F}$ on a cell complex $X$ we have that $\Gamma^{BM}(X;\hat{F})\cong \cok(C_1^{BM}(X;\hat{F})\to C^{BM}_0(X;\hat{F}))$.
\end{lem}
\begin{proof}
 The proof above goes through until the last argument. Now we have edges $e$ with only one vertex. However, by extending and zeroing out at infinity to get that the colimit of
\[
 \xymatrix{ & \hat{F}(e) \ar[ld]_{r_{e,v}} \ar[rd]^0 & \\
\hat{F}(v) & & \hat{F}(\infty)=0}
\]
is exactly equal to the co-equalizer of $r_{e,v}:\hat{F}(e)\to F(v)$ and the zero morphism, i.e. the cokernel.
\end{proof}

\subsection{Invariance under Subdivision}
\label{subsubsec:subdivision}
\index{Borel-Moore cosheaf homology!invariance under subdivision}\index{cosheaf!cellular!Borel-Moore homology under subdivision}

Now we take up the question of invariance under subdivision by applying the derived perspective. For convenience, we work with sheaves, but the reasoning can be dualized.

\begin{defn}
 Suppose $F$ is a sheaf on $X$ and $s:X'\to X$ is a subdivision of $X$, then we define the subdivided sheaf $F':=s^*F$.
\end{defn}

For an example, let $X$ be the unit interval $[0,1]$ stratified in the obvious way with $x=0$, $y=1$ and $a=(0,1)$. Now consider a sheaf $F$ on $X$. We will want to investigate what happens to this sheaf as we subdivide the space. In this example, the barycentric subdivision of $X$ produces a space $X'$ with a third vertex $\bar{a}$ and two edges $a_x$ and $a_y$. The obvious way of defining a subdivided sheaf is to define $F'(\bar{a})=F'(a_x)=F'(a_y)=F(a)$ where we use the identity map for the two new restriction maps. Observe that if $F$ is the elementary injective sheaf $[a]$, then $F'$ is \emph{not} an injective sheaf, yet nevertheless $F'$ and $F$ have isomorphic cohomology.

More generally we are concerned with the following diagram of spaces (posets)
\[
 \xymatrix{X' \ar[rr]^s \ar[rd]_{p_{X'}} & & X\ar[ld]^{p_X}\\
	      &\star&}
\]
and the induced functors on sheaves. For example, if we analyze the ordinary pushforward functor, then we would obtain the following result, which is a simplified proof of one found in~\cite[Thm. 1.5.2]{shepard}:

\begin{thm}\label{thm:subdivision}
 Suppose $F$ is a sheaf on $X$ and $X'$ is a subdivision of $X$, then
\[
 H^{\bullet}(X;F)\cong H^{\bullet}(X';F')
\]
\end{thm}
\begin{proof}
Observe that since $p_{X'}=p_X\circ s$, then $(p_{X'})_*=(p_X)_*\circ s_*$. Now recall $$(p_{X'})_*F'=(p_{X'})_*s^*F=(p_{X})_*\circ s_* s^* F.$$ 
The question then boils down to understanding the relationship between $s_*s^*F$ and $F$. Unraveling the definition reveals
\begin{eqnarray*}
 s_*s^*F(y) & = &\varprojlim \{s^*F(x) | s(x)\geq y\} \\
&=&\varprojlim \{F(s(x)) | s(x)\geq y\} \\
(\mathrm{surjectivity})&=&\varprojlim \{F(x) | x\geq y\} \\
(\mathrm{sheaf-axiom})&=& F(U_y) \\
&=& F(y)
\end{eqnarray*}
So we have that for the subdivision map $s_*s^*F\cong F$ and as a consequence
$$(p_{X'})F'\cong p_X F.$$
Now we can just take the associated right derived functors to obtain the result.
\end{proof}

\section{Sheaf Homology and Cosheaf Cohomology}
\label{subsec:sheaf_homology}
\index{sheaf!homology, cellular}\index{cosheaf!cohomology, cellular}

There is a surprising symmetry in the land of cellular sheaves and cosheaves, which is unique to the land of Alexandrov spaces and deserves to be explored. Contrary to the existence of enough injective sheaves, which for general sheaves is gotten as a consequence of the target category, e.g. $\Ab,\Vect$, etc., the existence of enough projective sheaves is driven by the underlying topology of the space.

\begin{prop}
 Suppose $X$ is a topological space with the property that there is a point $x\in X$ such that for every open neighborhood $U\ni x$ there is a strictly smaller open neighborhood $V\subset U$. Then the category of sheaves on $X$ does not have enough projectives.\footnote{The author would like to acknowledge the contributions of Valery Alexeev, David Treumann, and Jon Woolf on mathoverflow in regards to this question.}
\end{prop}
\begin{proof}
 Consider the map $i:x\hookrightarrow X$ and the sheaf $i_*k$. Suppose it has a projective resolution, i.e. a projective sheaf $P$ and a surjection $P\to i_*k$. Now let's examine this map evaluated on an open set $U\ni x$. By assumption there is another open set $V\subset U$ and we can put the constant sheaf extended by zero on $V$, denote the inclusion by $j:V\hookrightarrow X$. Note that we have the following diagram of sheaves
\[
 \xymatrix{j_!\tilde{k}_V \ar[r] & i_*k \ar[r] & 0\\ & P \ar@{.>}[lu] \ar[u] &}
\]
whose value on the open set $U$ is
\[
 \xymatrix{j_!\tilde{k}_V(U)=0 \ar[r] & i_*k(U)=k \ar[r] & 0\\ & P(U) \ar@{.>}[lu] \ar[u] &}
\]
so in particular the surjection must factor through zero --- a contradiction.
\end{proof}

Contrary to sheaves on manifolds and other Hausdorff spaces, cellular sheaves are can be viewed as sheaves on finite posets and as such do not suffer from the above argument. In fact, computing a projective resolution is as easy as computing injective resolutions. To see how this goes recall we need to find a projective sheaf that surjects onto our sheaf of interest. 
\[
 P^0:=\oplus_{\sigma\in X} \{\sigma\}^{F(\sigma)} \to F \to 0
\]
 serves nicely and by finding the kernel sheaf (which is easier to understand than cokernels!) and then iterating this process will obtain a projective resolution
\[
 \cdots P^{-3}\to P^{-2} \to P^{-1} \to P^0 \to F \to 0.
\]

This motivates the following definitions:
\begin{defn}
 Given a cellular sheaf, we can construct its projective resolution $P_{\bullet}\to F$, calculate colimits of $P_{\bullet}$ and take the cohomology of the resulting complex of vector spaces. Assuming $F$ was in degree zero, this will be concentrated in negative degree and we define the \textbf{homology of a cellular sheaf} $F$ to be $H_i(X;F):=H^{-i}(p_{\dd}P_{\bullet})$.

Similarly we define the \textbf{cohomology of a cellular cosheaf} $\hat{F}$ by taking its injective resolution $\hat{I}^{\bullet}$, and taking limits, i.e. $H^i(X;\hat{F})=H^i(\hat{p}_*\hat{I}^{\bullet})$.
\end{defn}

The reasons for it's apocryphal nature are many:
\begin{enumerate}
 \item Only for (co)sheaves over finite spaces are there enough projectives and enough injectives.
 \item Spaces for which there is not a fixed $n$ so that every cell $\sigma$ contains in its star a cell $\tau$ such that $\dim\tau=n$ cannot hope to have the same computational formula for (co)homology because we can't treat the colimit (in the case of a sheaf) as a quotient object of $\oplus_{\dim\tau=n} F(\tau)$ and dually for limits of cosheaves.
 \item This defect, which is measured by the difference of $H^n(X;F)$ and $\colim F$, is only the first in a series of obstructions that appear to detect whether $X$ is a cell structure on a manifold.
\end{enumerate}

The evidence for the last two observations is further solidified in view of the following theorem.

\begin{thm}\label{thm:mfld_sheaf_cosheaf}\index{duality!over a manifold}
 Suppose $F$ is a cellular sheaf on a triangulated closed $n$-manifold $X$, then $F$ defines a cellular cosheaf on the dual triangulation and moreover all the homologies and cohomologies of both agree.
\end{thm}
\begin{proof}
 This is a consequence of the following simple observation:
\[
 \xymatrix{F(\sigma^i) \ar[r]^{\rho_{\sigma,\tau}} \ar@{~>}[d] & F(\tau^{i+1}) \ar@{~>}[d] \\
\hF(\tilde{\sigma}^{n-i}) \ar[r]^{\rho_{\tilde{\sigma},\tilde{\tau}}} & \hF(\tilde{\tau}^{n-i-1})}
\]
So the same abstract diagram of vector spaces $F:X\to\Vect$ defines a diagram over $\hF:\tilde{X}^{op}\to\Vect$, i.e. a cosheaf on the dual cell structure. Since they are the same diagrams everything about them is the same.
\end{proof}

\subsection{Invariance under Subdivision}
\index{sheaf!homology, cellular!subdivision}\index{cosheaf!cohomology, cellular!subdivision}

One can ask whether this new invariant is invariant under subdivision. In this section we show that it is invariant for the domain of a map, but is not invariant under subdivision of the target. One can see this latter claim by an earlier example already considered with the pushforward with open supports to a circle.

\begin{thm}
  Suppose $F$ is a sheaf on $X$ and $X'$ is a subdivision of $X$, then
\[
 (p_{X'})_{\dd}F'\cong (p_{X})_{\dd}F
\]
and consequently
\[
 (Lp_{X'})_{\dd}F'\cong (Lp_X)_{\dd}F
\]
thus sheaf homology is invariant under subdivision. Similarly, the same result should hold for cosheaf cohomology.
\end{thm}
\begin{proof}
Getting right down to it we see
\begin{eqnarray*}
 s_{\dd}s^*F(y) &=& \colim \{s^*F(x) |s(x)\leq y\} \\
&=& \colim \{ F(s(x)) | s(x)\leq y\} \\
(\mathrm{surjectivity})&=&\colim \{F(x) | x\leq y\} \\
(\mathrm{check-directly})&=& F(y)
\end{eqnarray*}
and thus
$$(p_{X'})_{\dd}s^*F=(p_X)_{\dd}s_{\dd}s^*F \cong (p_X)_{\dd}F.$$
Taking the left derived functors gives the higher result.
\end{proof}

\ctparttext{This part constitutes a first application of cellular sheaves and cosheaves to problems in science and engineering. 

Chapter~\ref{sec:barcodes} begins with a short introduction to persistent homology, which we reformulate using sheaves and cosheaves. The advantage of this reformulation is the ability to distribute homology computations and aggregate efficiently, as noted in Section~\ref{subsec:local-to-global}, which is part of joint work~\cite{DMT_sheaves} where discrete Morse theory is adapted to compute cellular sheaf cohomology. A theorem that uses the apparatus of spectral sequences to connect level set and sub-level set persistence is proved in Section~\ref{subsec:level2sub}. A motivating example for multi-dimensional persistence and an introduction of ``generalized barcodes'' concludes the chapter.

Chapter~\ref{sec:nc_coding} reviews an application of sheaves to network coding introduced first in~\cite{GH-ncs}. However, here the language of cellular sheaves is employed and the barcode method is used to visualize the flow of data. A duality theorem for network coding sheaves is proved in two different ways.

Chapter~\ref{sec:sensors} casts various sensor network problems in the language of sheaves. Any attempt to use level set persistence to study the intruder problem is proven to be a ``no-go'' using the machinery of cosheaves. The main contribution of the chapter is the introduction of a linearized model for multi-modal sensing in Section~\ref{subsec:multi_modal}. It was this model that first motivated the author to take up the theory of indecomposables as a way of interpreting sheaf cohomology.
}

\part{Applications to Science and Engineering}
\label{part:applications}

%
%
\chapter{Topological Data Analysis}
\label{sec:barcodes}


\begin{quote}
{\em``We are enabled to divide into forms, following the objective articulation; we are not to hack off parts like a clumsy butcher.''}
\begin{flushright} --- Plato's Phaedrus 265e, translated by R. Hackforth \end{flushright}
\end{quote}

Suppose two patients enter an office with a recent diagnosis of cancer. If they are very lucky, they might have their respective tumors biopsied and sent off for genetic analysis. The genetic analysis consists of measuring the gene expression levels of over a thousand genes in each of the respective tumors. These two sets of expression levels are then compared with a few hundred other tumor samples, where varying therapies were used to varying degrees of success. How should the doctor determine which course of action to take? If we are to carve the universe at its joints, on which side do these patients lie?

A similar, but apparently less dramatic, situation occurs in topology. Given two topological spaces $X$ and $Y$, how do we discern whether $X$ and $Y$ are essentially the same or different? The entire apparatus of algebraic topology was developed to address this problem. It turns out that these two situations are not just formally similar. Topological data analysis stems from the observation that data of the above form can give rise to topological spaces, which can in turn be discriminated using classical constructions such as homology~\cite{lum2013extracting}.

\section{Point Clouds and Persistent Homology}
\label{subsec:pointclouds}\index{point cloud}\index{persistent homology}

To illustrate why the applicability of homological methods is not so far-fetched, consider the following toy problem, which serves as the standard entr\'ee into \textbf{persistent homology}.
\begin{figure}[ht]
	\centering
	\includegraphics[width=.7\textwidth]{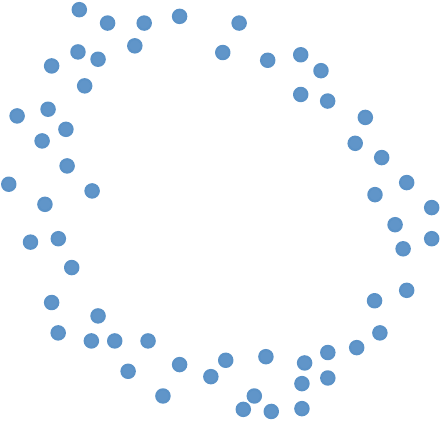}
	\caption{Point Cloud Data}
	\label{fig:point_cloud_data}
\end{figure}
Consider a finite set of points $\{x_i\}\subset \R^n$. How do we describe the perceived shape of such a set of points, such as the set depicted in Figure \ref{fig:point_cloud_data}? The human brain is a pattern-making machine that connects the dots and returns the knee-jerk response that the points in Figure \ref{fig:point_cloud_data} appears to form a circle. However, what do we mean by this and how can we automate this process so as to remove human subjectivity? 

Mathematics abounds with rigorous formulations of what does or does not look like a circle. However, only topology provides a definition that is robust with respect to perturbation and noise. One universal definition of a circle is given by the Eilenberg-Mac Lane space $K(\Z,1)$, which is homotopy equivalent to $S^1$. However, homology provides a shape descriptor based on linear algebra, which can be efficiently computed. But there is still the problem that the homology of the set of points in Figure \ref{fig:point_cloud_data} is not isomorphic to the homology of the circle. To get around this, we consider the union of Euclidean balls of some radius $r$ $\cup B(x_i,r)=:X_r$, which mirrors the ``connecting the dots'' procedure that the brain applies. Then one observes that there are natural inclusions
\[
	X_{r_0}\hookrightarrow X_{r_1} \hookrightarrow X_{r_2} \hookrightarrow X_{r_3} \cdots
\]
whenever $r_0\leq r_1\leq r_2 \leq \cdots$ and so on. Applying the $i$th homology functor $H_i(-;k)$ turns this diagram of spaces into a diagram of vector spaces, which defines a persistence module, cf. Definition \ref{defn:persistence_module}.
\begin{equation*}
	H_i(X_{r_0};k) \to H_i(X_{r_1};k) \to H_i(X_{r_2};k) \to H_i(X_{r_3};k) \to \cdots 
\end{equation*}
By applying the structure theorem~\ref{thm:crawleyboevey}, we can determine the barcodes of the collection of points. Long bars are considered to be robust topological signals in the data set. For Figure \ref{fig:point_cloud_data}, there would be one long bar in the persistence module corresponding to $H_0$, indicating that after a certain radius the the space $X_r$ is connected, and another long bar in the module corresponding to $H_1$, indicating the apparent circle in the data set. 

To summarize, we have the following prototypical pipeline of topological data analysis.
\begin{defn}[Point Cloud Persistence]\label{defn:point-cloud-persistence}\index{persistent homology!sublevel set}
The \textbf{point cloud persistence pipeline} consists of the following ingredients and operations:
\begin{enumerate}
\item Let $X$ denote a point cloud, i.e. the union of a finite set of points $\{x_i\}\subset \R^n$.
\item The union of balls $X_r:=\cup_{x_i\in X} B(x,r)$ (or the Vietoris-Rips complex~\cite{ghrist-barcodes}) defines a functor from the real line, viewed as poset, to the category of topological spaces (simplicial complexes) and maps, i.e. 
\[
G:(\R,\leq)\to\Top \qquad r \rightsquigarrow X_r \qquad r\leq r' \rightsquigarrow X_r\hookrightarrow X_{r'}
\]
\item Postcomposing this functor with homology $H_{\ast}$, defines a graded representation of the real line, which is equivalently a graded sheaf on the Alexandrov topology or a graded persistence module:
\begin{equation*}
	H_{\ast}(X_{r_0};k) \to H_{\ast}(X_{r_1};k) \to H_{\ast}(X_{r_2};k) \to H_{\ast}(X_{r_3};k) \to \cdots 
\end{equation*}
\item Applying Theorem \ref{thm:crawleyboevey} produces a Remak decomposition of this representation into a multiset of interval modules, which is visualized as a barcode by the end user. 
\end{enumerate}
\end{defn}

The first and second steps in this pipeline offer the chance for endless modification and application. Instead of considering a collection of points, one can start with a space $X$ and a function $f:X\to\R$ and consider the family of sub-level sets $X_r:=f^{-1}(-\infty,r]$. As long as the function and space are sufficiently nice, we can use Theorem \ref{thm:crawleyboevey} to produce a barcode.

\begin{exr}
Determine the barcodes associated to the function $f(x)=x^3-x$.
\end{exr}

\begin{exr}
Recast the basic theorems of Morse theory in terms of barcodes. See Figure \ref{fig:persistence_torus} for inspiration.
\end{exr}

\subsection{Level Set and Zigzag Persistence}

Despite their successes, persistence modules are not the end-all, be-all of topological data analysis. In~\cite{zigzag} Gunnar Carlsson and Vin de Silva gave three examples where diagrams of vector spaces and maps of the form
\[
	V_1 \leftrightarrow V_2 \leftrightarrow V_3 \leftrightarrow  \cdots \leftrightarrow V_{n-2} \leftrightarrow V_{n-1} \leftrightarrow V_n
\]
are of interest. One example comes from estimating the probability density function from which a point-cloud is drawn. One can try to smooth the data by defining a function $\rho_r(x)$ that counts the number of points within a radius $r$ of the point $x$. If one then tries to take the 25\% densest points measured according to a sequence of radii $r_1<r_2<\cdots <r_n$, the the only way of comparing features comes from a \textbf{zigzag} diagram of the form
\[
X_{r_1}^p \rightarrow X_{r_1}^p\cup X_{r_2}^p \leftarrow X_{r_2}^p \rightarrow \cdots \leftarrow X_{r_n}^p  
\]
where $X^p_{r_k}$ indicates the $p\%$ densest points measured according to the function $\rho_{r_k}(x)$. Similarly, if one has a function $f:X\to\R$, but not the computational power to investigate the entire sub-level set $f^{-1}(-\infty,r]$, one could choose a mesh $t_0<t_1<t_2<\cdots <t_n$ and consider the zigzag of pre-images given by
\[
f^{-1}(t_0)\rightarrow f^{-1}[t_0,t_1]\leftarrow f^{-1}(t_1) \rightarrow \cdots \leftarrow f^{-1}(t_n).
\]
Applying homology in some degree $i$ gives the traditional definition of \textbf{level set persistent} homology. This complicates the usual TDA pipeline because, a priori, the structure theorem \ref{thm:crawleyboevey} fails to apply. Fortunately, Gabriel's Theorem \ref{thm:gabriel} tells us that the direction of the arrows doesn't matter for the representations of $A_n$ type quivers, so the decomposition into interval modules still applies; we still have barcodes~\cite{zigzag}.

The fully general definition of level set persistence usually given adheres to this perspective that the assignment of homology to closed intervals is fundamental.
\begin{defn}\index{category!of intervals}
The \textbf{interval category} of $\R$, written $\Int(\R)$, is the category whose objects are closed intervals $[x,y]\subset \R$ and whose morphisms are inclusions $[x,y]\hookrightarrow [x',y']$.
\end{defn}
\begin{defn}[Level Set Persistent Homology]\index{persistent homology!level set/zigzag}
Suppose $f:X\to \R$ is a function, not necessarily continuous. The \textbf{$i$th level set persistence of $f$} $L_i$ is the representation of the interval category given by
\[
L_i:\Int(\R)\to\Vect \qquad [x,y]\rightsquigarrow H_i(f^{-1}([x,y]);k).
\]
\end{defn}

\subsubsection{Critique of the Definition of Level Set Persistence}

\begin{figure}
\centering
\includegraphics[width=.9\textwidth]{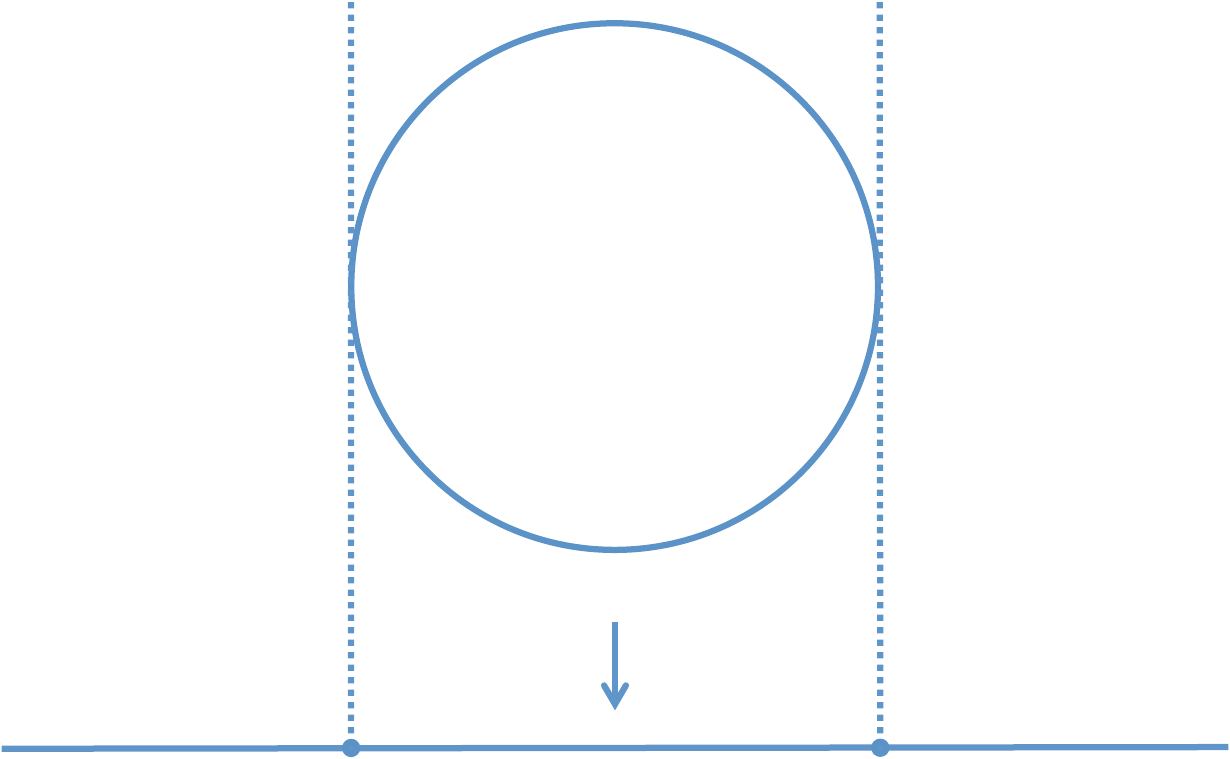}
\caption{Height Function on the Circle}
\label{fig:height_circle}
\end{figure}

The definition of level set persistence suffers in one crucial respect. The definition is \emph{non-local} and consequently requires the storage of an infinite amount of data. Consider the map $f$ depicted in Figure \ref{fig:height_circle}. Level set persistence for $H_1$ will assign the zero vector space to every interval of diameter less than $y-x$. Only after inspecting intervals large enough will the topological feature of the circle appear.

A sheaf-theoretic approach to level-set persistence offers the advantage of being local. However, the analogous sheaf version of the level set $H_1$ just examined is the zero sheaf. The trade-off appears to be too great. However, the apparent disadvantage is remedied via the use of sheaf-cohomology, which preserves all the information of the domain space while simultaneously being local.

\section{Approaching Persistence with Sheaves and Cosheaves}

\subsection{Cellular Maps and Absolute Homology Cosheaves}\label{subsec:abs_homology_cosheaves}

The contravariant nature of sheaves should remind us that that they are most naturally associated with cohomology of a space. Cosheaves, being covariant with respect to the inclusion of opens, are naturally associated with homology of a space. Since we are working over a field, we can pass back and forth between these perspectives. However, to coincide with the traditional use of homology in persistence, we first introduce our version of persistence using cosheaves.

\begin{defn}\label{defn:abs-hom-cosheaf}\index{absolute homology cosheaf}\index{cosheaf!absolute homology}
Suppose $X$ and $Y$ are cell complexes and $f:Y\to X$ is a \emph{proper} cellular map, see Definition \ref{defn:cell_map} for a reminder, then for each natural number $i\geq 0$ we have the following $i$th \textbf{absolute homology cosheaf} $\hF_i$, which assigns to a cell $\sigma$ in $X$ the $i$th homology of the pre-image $f^{-1}(\st(\sigma))$, i.e.
\[
	\hF_i(\sigma):=H_i(f^{-1}(\st(\sigma));k).
\]
This is clearly a cellular cosheaf since if $\sigma\leq \tau$, then $\st(\tau)\subseteq \st(\sigma)$ and thus we have a map 
\[
r_{\sigma,\tau}:H_i(f^{-1}(\st(\tau));k) \to H_i(f^{-1}(\st(\sigma));k).
\]
\end{defn}

\begin{figure}
\begin{center}
\includegraphics[width=\textwidth]{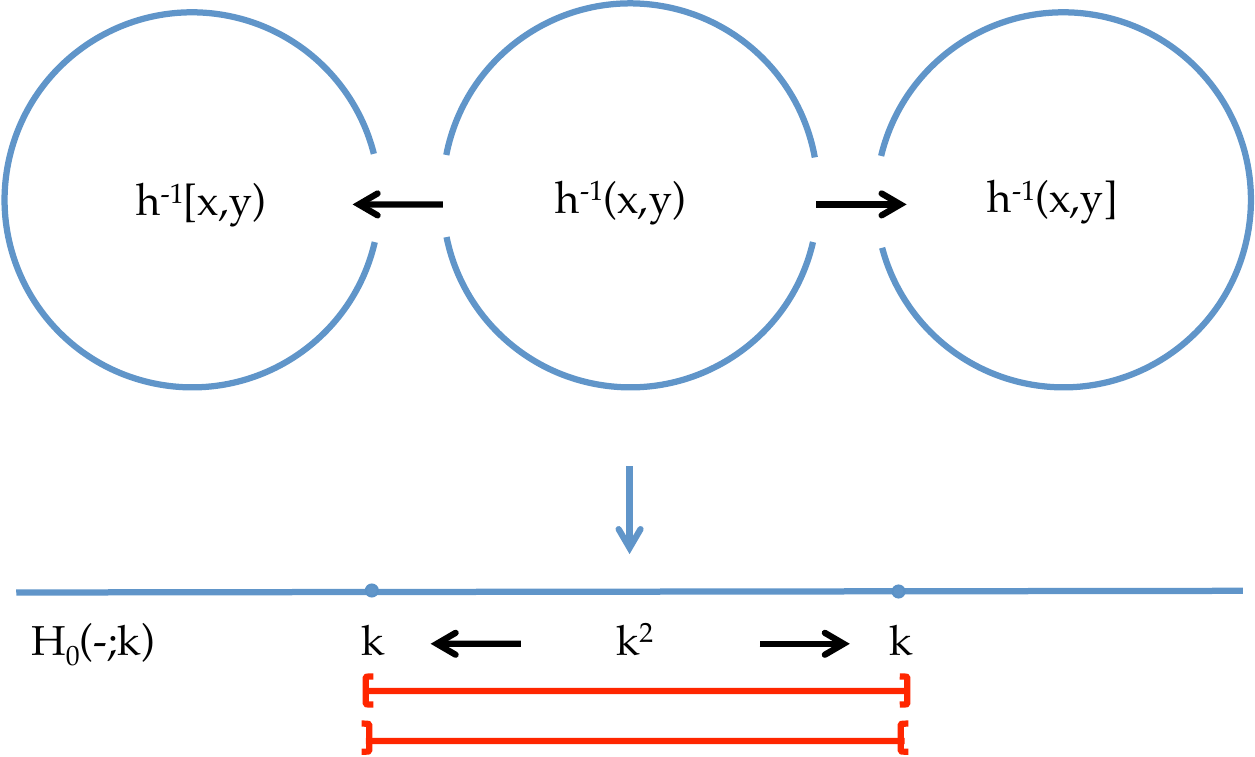}
\end{center}
\caption{Homology Cosheaf for a circle, with barcodes}
\label{fig:circle_bc}
\end{figure}

By Corollary \ref{cor:sheaves_remak} we know that every absolute homology cosheaf can be decomposed into indecomposable sheaves. For cellular maps $f:Y\to X$ where $X$ is a compact subset of $\R$, the absolute homology cosheaves have the following form:
\[
H_i(f^{-1}(\st(x_0));k) \leftarrow H_i(f^{-1}((x_0,x_1));k) \rightarrow H_i(f^{-1}(\st(x_1));k) \leftarrow \cdots
\]
Hence, by Gabriel's Theorem \ref{thm:gabriel} absolute homology cosheaves over $X\subset \R$ can be assigned barcodes. However, by inspecting the support of these indecomposable cosheaves we observe that there are four types of bars that make up any barcode:
\[
	[\textrm{---}] \qquad (\textrm{---}) \qquad [\textrm{---}) \qquad (\textrm{---}]
\]

\begin{ex}
Let $h:S^1\to\R$ be the standard height function on the circle, drawn in Figure \ref{fig:height_circle}. The only absolute homology cosheaf of interest is $\hF_0$, since the fibers have no higher homology. The associated pre-images, values of the homology cosheaf and barcode are drawn in Figure \ref{fig:circle_bc}.
\end{ex}

\begin{rmk}[Barcode Notation]
We will use intervals to represent barcodes and these will be sensitive to whether the first or last vector space in an indecomposable representation falls on a vertex or an open interval in some stratification of $[0,1]$ or $\RR$. For visual clarity, we will adopt the convention that a turned around square bracket is equivalent to a round one, i.e. $(x_i,x_{i+1})\rightsquigarrow ]x_{i},x_{i+1}[$ and $[x_i,x_{i+1})\rightsquigarrow [x_i,x_{i+1}[$ and so on.
\end{rmk}

By viewing these barcodes as cosheaves, which have a homology theory, we can compute barcode homology.

\begin{lem}\label{clm:BM_barcodes}\index{homology!of a cellular cosheaf!via barcodes}\index{cosheaf!cellular!homology via barcodes}
	Suppose $X$ is a compact subset of $\R$, equipped with some cell structure. The cosheaf homology of the four types of indecomposable cosheaves coincides with the Borel-Moore homology of the underlying barcode, i.e.
	\[
		H_0^{BM}([\textrm{---}])=k \qquad H_1^{BM}(]\textrm{---}[)=k \qquad H_i^{BM}([\textrm{---}[)=H_i^{BM}(]\textrm{---}])=0
	\]
	with all other Borel-Moore homology groups being zero. Moreover, since cosheaf homology commutes with finite direct sums, cellular cosheaf homology of $\hF$ on $X$ can be computed using the barcode $B_{\hF}$ associated to the Remak decomposition of $\hF$.
	\[
		H_i(X;\hF)\cong \oplus H_i^{BM}(B_{\hF}) \qquad i=0,1
	\]
	In short,
	\begin{quote}
	\centering
	{\em $H_0(X;\hF)$ counts \textbf{closed bars} and $H_1(X;\hF)$ counts \textbf{open bars}.}
	\end{quote}
	
\end{lem}
\begin{proof}
The proof of the claim uses simple computations, illustrated in the examples, and an invariance under subdivision argument, which is proved in Theorem \ref{thm:subdivision}. The true sticking point is why ordinary cosheaf homology becomes Borel-Moore (cosheaf) homology, which is developed in section \ref{subsubsec:BM_cosheaf_homology}.

If we imagine that we are extending the constant cosheaf supported on a barcode to a cosheaf defined on all of $[0,1]$, then the natural way of extending by zero is to use the pushforward with open supports functor $j_{\dd}$, where
\[
	j: B \hookrightarrow [0,1] \qquad \mathrm{and} \qquad \hat{k}_B \mapsto j_{\dd}\hat{k}_B.
\]
The process of taking cosheaf homology is to then push forward this cosheaf to a point. However, using Lemma \ref{lem:BM_cosheaf_homology} we see the following is true: 
\[
	p_!\hat{k}_B \cong p_* j_{\dd}\hat{k}_B \qquad H^{BM}_i(B;\hat{k}_B)\cong H_i([0,1];j_{\dd}\hat{k}_B)
\]
When it is clear that we are working on $[0,1]$ we may write $\hat{k}_B$ instead of $j_{\dd}\hat{k}_B$.
\end{proof}

One of the advantages of absolute homology cosheaves is that over the real line they can be used to compute the homology of the domain.

\begin{cor}[``The Barcode Trick'']\label{cor:barcode_homology}\index{barcode!trick}
	Assume $Y$ is compact and $f:Y\to X$ is a cellular map with $f(Y)=X\subset\RR$. For each $i$, let $B_i$ denote the barcode associated to the $i$th absolute homology cosheaf. The following is true:
	 \[
	  H_i(Y;k)\cong H_0(X;\hF_i)\oplus H_1(X;\hF_{i-1}) \cong H_0^{BM}(B_i)\oplus H_1^{BM}(B_{i-1})
	 \]
\end{cor}
\begin{proof}
This is an immediate corollary of Theorem \ref{thm:leraysheaf} and Lemma \ref{clm:BM_barcodes}.
\end{proof}

\begin{ex}[Height function on the Two Sphere]

\begin{figure}
\begin{center}
\includegraphics[width=.7\textwidth]{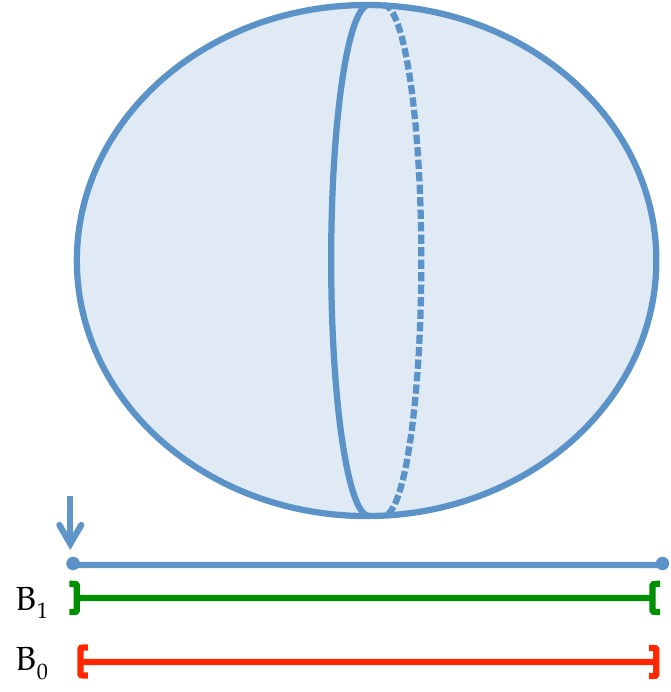}
\end{center}
\caption{Barcodes and the Two Sphere}
\label{fig:sphere_bc}
\end{figure}

Let $h:S^2\to\RR$ be the standard height function on the two sphere. In Figure \ref{fig:sphere_bc} we have drawn the map and the associated barcodes. The barcode decomposition for the cosheaf associated to taking $H_0$ of the fiber is trivial because it is already an indecomposable cosheaf.
\[
	\hF_0: \qquad \xymatrix{k & \ar[l] k \ar[r] & k}
\]
Similarly, taking $H_1$ of the fiber also yields an indecomposable cosheaf
\[
	\hF_1: \qquad \xymatrix{0 & \ar[l] k \ar[r] & 0}
\]

Now let us compute cosheaf homology. Since the space $X=[0,1]$ is compact, ordinary and compactly supported cosheaf homology agree. We label our cells as $x=0$, $a=(0,1)$ and $y=1$. To get an ordered basis and matrix representatives for our homology computation, we choose the local orientation pointing to the right and use the lexicographic ordering on the cells. For $\hF_0$ we get the following boundary matrix and homology groups:
\[
  \partial_1=\begin{bmatrix}-1 \\ 1\end{bmatrix} : k_a \to k_x\oplus k_y \qquad \Rightarrow \qquad H_1(X;\hF_0)=0 \quad H_0(X;\hF_0)=k.
\]
For $\hF_1$ the computation is even easier:
\[
  \partial_1: k_a \to 0 \qquad \Rightarrow \qquad H_1(X;\hF_1)=k \quad H_0(X;\hF_1)=0
\]
One can then check \ref{cor:barcode_homology} simply as follows:
\[
\xymatrix{H_0(X;\hF_1)=0 \ar@{.}[rd] & H_1(X;\hF_1)=k \\
H_0(X;\hF_0)=k & H_1(X;\hF_0)=0}
\]

\[
H_0(S^2)=k \qquad H_1(S^2)=0 \qquad H_2(S^2)=k
\]
\end{ex}

\begin{ex}[Height function on a Cone]

\begin{figure}[ht]
\centering
\includegraphics[width=.7\textwidth]{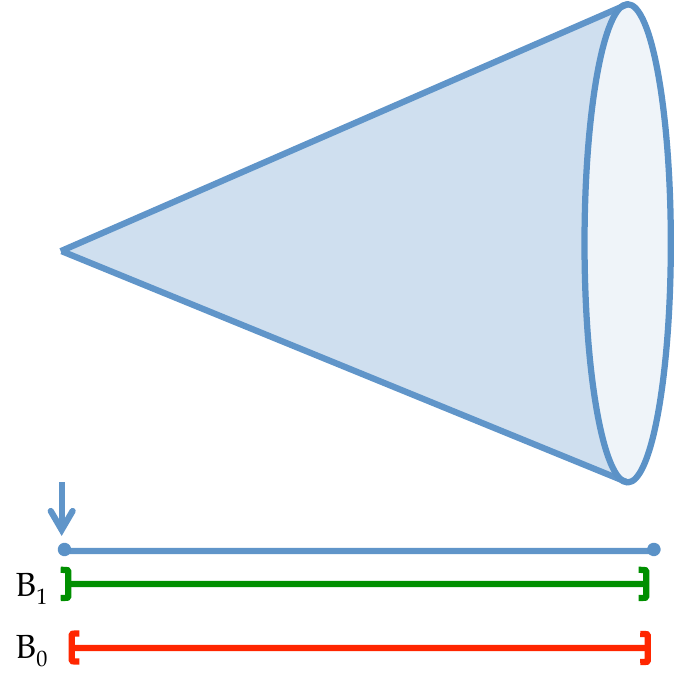}
\caption{Barcodes for the Cone}
\label{fig:cone_bc}
\end{figure}

The height function on a cone is not a Morse function because differentiability breaks down at the cone point; see Figure \ref{fig:cone_bc}. One could use stratified Morse theory as a substitute, but we'll use cosheaf homology. Here the cosheaf $\hF_0$ is the same as the previous example; we will not repeat the computation. The cosheaf $\hF_1$
\[
	\hF_1: \qquad \xymatrix{0 & \ar[l] k \ar[r] & k}
\]
exhibits different behavior. The cosheaf homology computation for this cosheaf reveals that the half-open barcode, embedded inside a compact interval, has no non-zero homology groups.
\[
  \partial_1=\id: k_a \to k_y \qquad \Rightarrow \qquad H_1(X;\hF_1)=0 \quad H_0(X;\hF_1)=0
\]
Checking \ref{cor:barcode_homology} again gives
\[
\xymatrix{H_0(X;\hF_1)=0 \ar@{.}[rd] & H_1(X;\hF_1)=0 \\
H_0(X;\hF_0)=k & H_1(X;\hF_0)=0}
\]

\[
H_0(C)=k \qquad H_1(C)=0 \qquad H_2(C)=0
\]
\end{ex}

Let's illustrate the utility of the barcode trick by computing cosheaf homology over $X$ using two different methods:
\begin{itemize}
	\item Using the computational formulae of section \ref{subsec:comp_homology}
	\item Determining the barcode decomposition and applying claim \ref{clm:BM_barcodes}.
\end{itemize}

\begin{ex}[Height function on the Torus]

\begin{figure}[ht]
\centering
\includegraphics[width=.7\textwidth]{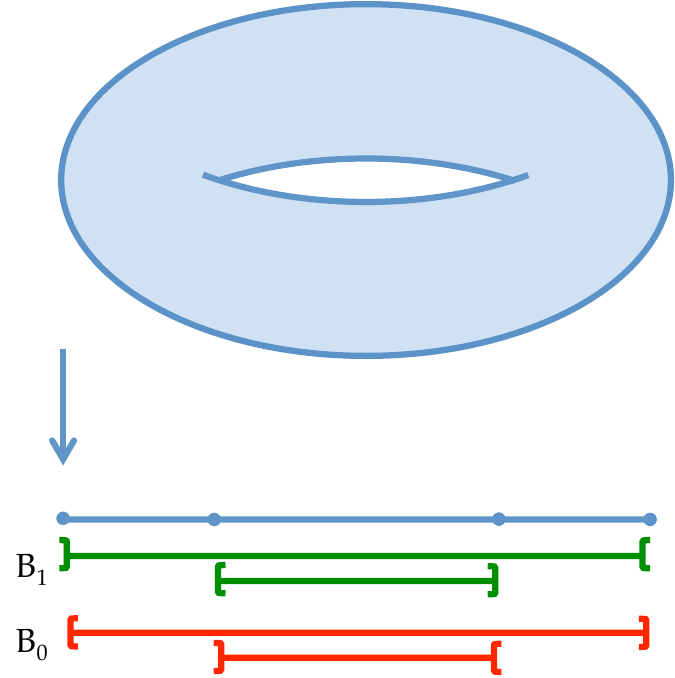}
\caption{Barcodes for Bott's torus}
\label{fig:torus_bc}
\end{figure}

The standard introductory example of Morse theory, first popularized by Raoul Bott, is the height function on the torus. In Figure \ref{fig:torus_bc} we have drawn the behavior of the fibers over the critical values and the non-critical intervals. For the sake of brevity, let us write out only the cosheaf $\hF_1$:
\[
\xymatrix{0 & \ar[l] k_a \ar[r] & k_y^2 & \ar[l] k_b^2 \ar[r] & k_z^2 & \ar[l] k_c \ar[r] & 0}
\]
Here the maps from $k_a$ to $k_y^2$ and $k_c$ to $k_z^2$ are the diagonal maps
\[
 r_{z,a}=\begin{bmatrix}1 \\ 1\end{bmatrix}=r_{z,c}
\]
and the other maps are the identity. Choosing the orientation that points to the right, we get the follow matrix representation for the boundary map:
\[
	\partial_1=\begin{bmatrix}1 & -1 & 0 & 0 \\
					1 & 0 & 1 & 0 \\
					0 & 1 & 0 & -1 \\
					0 & 0 & 1 & -1
	\end{bmatrix}
	\qquad H_1(X;\hF_1)=<\begin{bmatrix} 1 \\ 1 \\1 \\1\end{bmatrix}> \qquad H_0(X;\hF_1)\cong k
\]
However, if we change our bases as follows
\[
	\begin{bmatrix} y_1'= y_1 \\ y'_2=y_1+y_2\end{bmatrix} \qquad \begin{bmatrix} b_1'= b_1 \\ b'_2=b_1+b_2\end{bmatrix} \qquad \begin{bmatrix} z_1'= z_1 \\ z'_2=z_1+z_2\end{bmatrix}
\]
then our cosheaf $\hF_1$ can then be written as the direct sum of two indecomposable cosheaves:
\[
\xymatrix{0 & \ar[l] 0 \ar[r] & k_{y'_1} & \ar[l] k_{b'_1} \ar[r] & k_{z'_1} & \ar[l] 0 \ar[r] & 0}
\]
\[
\xymatrix{0 & \ar[l] k_a \ar[r] & k_{y'_2} & \ar[l] k_{b'_2} \ar[r] & k_{z'_2} & \ar[l] k_c \ar[r] & 0}
\]
Hence it is apparent that
\[
	H_i(X;\hF_1)\cong H_i^{BM}([\textrm{---}])\oplus H_i^{BM}(]\textrm{---}[).
\]
From which the homology of the torus can be directly observed.
\[
\xymatrix{H_0(X;\hF_1)=k \ar@{.}[rd] & H_1(X;\hF_1)=k \\
H_0(X;\hF_0)=k & H_1(X;\hF_0)=k}
\]
\[
H_0(T)=k \,\, H_1(T)=k^2 \,\, H_2(T)=k
\]
\end{ex}

\subsection{Local-to-Global Computations via Cellular Sheaves}
\label{subsec:local-to-global}

Aside from the decoration of cosheaves, section \ref{subsec:abs_homology_cosheaves} is completely classical and can be stated much more generally using non-cellular sheaves and cosheaves. However, one can still use cellular versions to compute homological information of non-cellular spaces as already noted by the author in~\cite{DMT_sheaves}.\index{discrete Morse theory for sheaves}

\subsubsection{The \v{C}ech Approach}

Recall that any topological space $X$ equipped with a cover $\cU$ has an associated simplicial approximation $N_{\cU}$ given by the nerve construction considered in Definition \ref{defn:nerve}. The nerve theorem tells us when this approximation is ``good enough'' for the purposes of cohomology.

\begin{thm}[Nerve Theorem \cite{leray45,borsuk-nerve}]\index{nerve theorem}
Given a topological space $X$ and a cover $\cU$, if the support $U_\sigma \subset X$ of each $\sigma \in N_\cU$ is acyclic (i.e., the reduced cohomology $\tilde{H}^\bullet(U_\sigma;R) = 0$ vanishes),
then $H^\bullet(N_\cU;R)\cong H^\bullet(X;R)$.
\end{thm}

Typically, the coarsest covers do not satisfy the acyclicity assumption. One achieves this by refinement of the cover, with the additional cost of more simplices in $N_{\cU}$. However, one can dodge this refinement by recording the cohomology of the intersection as a cellular sheaf on $N_{\cU}$.

\begin{defn}[\v{C}ech Sheaves]\index{sheaf!C@\v{C}ech}\index{Cech@\v{C}ech sheaves}

The {\bf \v{C}ech cellular sheaves} $\Cech^n$ associated to the cover $\cU$ of a space $X$ are defined on the nerve $N_\cU$ by the following data. Each $\sigma \in N_\cU$ is assigned the vector space $\Cech^n(\sigma) = H^n(U_{\sigma};k)$ and each face relation $\sigma \subset \tau$ is assigned the linear map $\Cech^n_{\sigma\tau}: H^n(U_{\sigma};R)\to H^n(U_{\tau};R)$ arising from the inclusion of supports $U_\tau \hookrightarrow U_\sigma$.

\end{defn}

If all simplex supports are acyclic, then $\Cech^0$ reduces to the constant sheaf on $N_{\cU}$ and all other $\Cech^n$s are trivial; in the absence of acyclicity assumptions, the following result yields a simple correction.

\begin{prop}
\label{prop:cechsheaf}
Let $X$ be a topological space and $\cU$ a cover whose nerve $N_\cU$ is at most one-dimensional. Then, for each $n \in \N$,
\begin{equation*}
		H^n(X;k) \cong H^0(N_\cU;\Cech^n)\oplus H^1(N_\cU;\Cech^{n-1}).
\end{equation*}
\end{prop}

We note that similar results have been obtained by Burghelea and Dey~\cite{dey_circle}, as well as Carlsson, de Silva, and Morozov~\cite{pyramid} in the context of zig-zag persistence. The difference between their results and ours is that their results depend on the decomposition of zig-zag persistence modules into indecomposable modules (barcodes). Our result makes the recognition that these modules are rightly conceived as sheaves over a linear nerve with a cohomology that can be quickly computed using discrete Morse theory~\cite{DMT_sheaves}.

Proposition \ref{prop:cechsheaf} generalizes the familiar Mayer-Vietoris long exact sequence, as the next example shows.
\begin{ex}[Mayer-Vietoris]\index{Mayer-Vietoris!via cellular sheaf cohomology}
The Mayer-Vietoris Theorem states that given an open cover of $X$ by two open sets $\cU=\{A,B\}$ we have the following exact sequence of $R$-modules and maps:

\[
\xymatrix{
	            0\ar@{->}[r] & H^0(X)\ar@{->}[r] & H^0(A)\oplus H^0(B) \ar@{->}[r]
	                   & H^0(A\cap B) \ar `r[d] `_l[lll] ^{\psi^0} `^d[dlll] `^r[dll] [dll] \\
	            & H^1(X)\ar@{->}[r] & \hspace{12pt}\cdots\hspace{12pt}\ar@{->}[r]
	                   & H^{n-1}(A\cap B) \ar `r[d] `_l[lll] ^{\psi^{n-1}} `^d[dlll] `^r[dll] [dll] \\
	            &H^{n}(X)\ar@{->}[r] & H^{n}(A)\oplus H^{n}(B) \ar@{->}[r]
	                   & H^{n}(A\cap B) \ar `r[d] `_l[lll] ^{\psi^{n}} `^d[dlll] `^r[dll] [dll] \\
	            & & \cdots &
	}
\]
Ideally one can determine the unknown cohomology of $X=A\cup B$ by inspecting the terms on either side. More formally, one uses universal constructions to force $H^n(X)$ into a more manageable short exact sequence:
\[
	\xymatrix{ & H^{n-1}(A)\oplus H^{n-1}(B) \ar[r]^-{\delta^0_{n-1}} & H^{n-1}(A\cap B) \ar[d] \ar[ld] & 0 \ar[d] &  \\
	0 \ar[r] & \cok(\delta^0_{n-1}) \ar[d] \ar@{.>}[r] & H^n(X) \ar@{.>}[r] \ar[d] & \ker(\delta^0_n) \ar[r] \ar[ld] & 0 \\
	& 0 & H^n(A)\oplus H^n(B) \ar[r]^-{\delta^0_{n}} & H^n(A\cap B) &  }
\]
The dotted maps are defined by the universal property of the cokernel and kernel, respectively. Over an arbitrary coefficient ring $R$ we would have to solve an extension problem in order to infer $H^n(X)$. If we take $R$ to be a field $k$, then every short exact sequence splits and we can deduce that
\[
	H^n(X)\cong \ker(\delta^0_n)\oplus \cok(\delta^0_{n-1})\cong H^0(N_{A,B};\Cech^n)\oplus H^1(N_{A,B};\Cech^{n-1})
\]
where $N_{A,B}$ is the unit interval, viewed as the nerve of a two-element cover.
\end{ex}

\subsubsection{The Leray Approach}

As Section~\ref{subsec:abs_homology_cosheaves} already indicated with examples, one can compute homological information of a space $Y$ with a suitably nice map $f:Y\to X$. We view this construction from a different perspective. By assuming that $X$ comes with a cover $\cV$, having nerve $N_\cV$, one can pull-back the associated \v{C}ech sheaf on $N_\cV$ along $f$ to yield local information about $Y$.

\begin{defn}[Leray Sheaves]\index{sheaf!Leray}\index{Leray sheaves}
The {\bf Leray cellular sheaves} $\Leray^n$ associated to a map $f:Y\to X$ and a cover $\cV$ of $f(Y) \subset X$ are defined over the nerve $N_\cV$ as follows. Each simplex $\sigma \in N_\cV$ is assigned the cohomology of the preimage of its support, i.e., $\Leray^n(\sigma) = H^n(f^{-1}(V_\sigma);k)$; furthermore, each face relation $\sigma \subset \tau$ is assigned the map induced on cohomology by the inclusion $f^{-1}(V_\tau) \hookrightarrow f^{-1}(V_\sigma)$.
\end{defn}
\begin{rmk}
\begin{itemize}
\item We will sometimes use the notation $F^n$ in place of $\Leray^n$ to emphasize the association to $f:Y\to X$. 
\item The absolute homology cosheaf $\hF^n$ of Definition \ref{defn:abs-hom-cosheaf} is clearly the cosheaf-theoretic version of $\Leray^n$ when the cover is given by the open stars of the cells.
\item A More general version of the Leray sheaf is given by the right derived pushforward of the constant sheaf along the map $f$.
\end{itemize}
\end{rmk}

In the special case where $X = Y$ and $f$ is the identity map, the Leray sheaves clearly coincide with the \v{C}ech sheaves associated to the cover $\cV$ of $X$. Thus, the following result generalizes Proposition \ref{prop:cechsheaf}.

\begin{thm}
\label{thm:leraysheaf}
	Let $f:Y \to X$ be continuous. Assume a cover $\cV$ of the image $f(Y) \subset X$ whose nerve $N_\cV$ is at most one-dimensional. Then, for each $n \in \N$,
\begin{equation*}
		H^{n}(Y;k) \cong H^0(N_\cV;\Leray^n) \oplus H^1(N_\cV;\Leray^{n-1}) .
\end{equation*}
\end{thm}
\begin{proof}
		The theorem is a simple consequence of the Leray spectral sequence which packages the cohomology of $Y$ into a coefficient system over the space $X$ from a map $f:Y\to X$~\cite{mccleary}. The restriction to a one-dimensional nerve forces the spectral sequence to collapse on the second page and hence yield the desired isomorphisms. More precisely, for each open $V\subset f(Y)$, let $C^n(V;R)$ denote the vector space freely generated by the set of all cochains defined on $V$. Clearly if $V \subset U$, then there is a surjection $C^n(U;R) \to C^n(V;R)$ defined by restriction of cochains. The sheaf $\widetilde{C}^n$ associated to this presheaf of singular cochains is consequently {\em flabby} (see \cite[p.~97]{ramanan2005global}).
		
Consider the following double complex of vector spaces:
	\[
	 \xymatrix{ \vdots & \vdots & \vdots & \vdots\\
	C^2(Y) \ar[r] \ar[u] & \bigoplus_{\dim \sigma = 0} \widetilde{C}^2(f^{-1}(V_\sigma)) \ar[r] \ar[u] & \bigoplus_{{\dim \tau = 1}} \widetilde{C}^2(f^{-1}(V_\tau) \ar[r] \ar[u] & 0 \\
	C^1(Y) \ar[r] \ar[u] &\bigoplus_{\dim \sigma = 0} \widetilde{C}^1(f^{-1}(V_\sigma)) \ar[r] \ar[u] & \bigoplus_{\dim \tau = 1} \widetilde{C}^1(f^{-1}(V_\tau) \ar[r] \ar[u]  & 0 \\
	C^0(Y) \ar[r] \ar[u] & \bigoplus_{\dim \sigma = 0} \widetilde{C}^0(f^{-1}(V_\sigma)) \ar[r] \ar[u] & \bigoplus_{\dim \tau = 1} \widetilde{C}^0(f^{-1}(V_\tau)) \ar[r] \ar[u] & 0 }
	\]
	It follows from standard results \cite[Thm~II.5.5, Thm~III.4.13]{Bredon}) that the rows are exact. By the acyclic assembly lemma \cite{weibel}, the spectral sequence converges to the cohomology of the leftmost column, i.e., $H^\bullet(Y;k)$. If one takes cohomology in the vertical direction, one obtains the defined cochain groups associated to the Leray cellular sheaves $\Leray^n$:
	\[
	 \xymatrix{  \vdots & \vdots & \vdots\\
	\bigoplus_{\dim \sigma = 0} H^2(f^{-1}(V_\sigma)) \ar[r] & \bigoplus_{\dim \tau = 1} H^2(f^{-1}(V_\tau)) \ar[r] & 0 \\
	\bigoplus_{\dim \sigma = 0} H^1(f^{-1}(V_\sigma)) \ar[r] & \bigoplus_{\dim \tau = 1} H^1(f^{-1}(V_\tau)) \ar[r]  & 0 \\
	\bigoplus_{\dim \sigma = 0} H^0(f^{-1}(V_\sigma)) \ar[r] & \bigoplus_{\dim \tau = 1} H^0(f^{-1}(V_\tau)) \ar[r] & 0}
	\]
Taking cohomology horizontally corresponds precisely to computing separately (in parallel, if one wishes) the cohomology of the Leray sheaves $\Leray^n$ over $N_\cV$, thus producing the final stable page of the spectral sequence.
	\[
	 \xymatrix{  \vdots & \vdots & \vdots\\
	 H^0(N_\cV; \Leray^2) & H^1(N_\cV;\Leray^2) & 0 \\
	 H^0(N_\cV; \Leray^1) \ar[rru] & H^1(N_\cV;\Leray^1)  & 0 \\
	 H^0(N_\cV; \Leray^0) \ar[rru] & H^1(N_\cV;\Leray^0) & 0}
	\]
Over a general ring $R$, these terms prescribe a filtration of the cohomology, giving rise to extension problems; however, over a field $k$ one can read off the cohomology directly.
\end{proof}

Note that the proof indicates precisely where we require the one-dimensional nerve restriction. Without this assumption in place, the second page of the spectral sequence may not be stable and the conclusion of the theorem need not hold.

\subsubsection{A Unifying Perspective}

There is a more sophisticated version of the nerve described originally by Segal~\cite{segal_CSSS} which is homotopically faithful to the underlying space independent of the particulars of the cover. This notion has been used in recent applications~\cite{zomorodian2008localized} and parallelizations for homology computation~\cite{lewis2012multicore}.

\begin{defn}[Mayer Vietoris Blowup]\label{defn:MV_blowup}\index{Mayer-Vietoris!blowup complex}\index{nerve!generalization}
Let $X$ be a topological space equipped with a cover $\cU$ with nerve $N_\cU$. The {\bf Mayer Vietoris blowup} $M_\cU$ associated to $\cU$ is a subset of the product $X \times N_\cU$ defined as follows. The pair $(x,s)$ lies in $M_\cU$ if and only if there is some simplex $\sigma \in N_\cU$ for which $x \in U_\sigma$ and $s \in \sigma$.
\end{defn}

\begin{rmk}
The original description given by Segal is that (up to barycentric subdivision) $M_{\cU}$ is the classifying space of a topological category, whose objects are pairs $(x,U_{\sigma})$ with $x\in U_{\sigma}$ and whose morphisms are pointed inclusions $(x,U_{\sigma})\to (y,U_{\tau})$ where $\tau\subset\sigma\subset I$ and $I$ is the indexing set of the cover $\cU=\{U_i\}_{i\in I}$.
\end{rmk}

Segal provides an updated version of the nerve theorem using this construction~\cite[Prop. 4.1]{segal_CSSS}

\begin{lem}[Generalized Nerve Theorem]\label{lem:generalized_nerve_thm}\index{nerve theorem!generalization}
	If $X$ is a paracompact Hausdorff space and $\cU=\{U_i\}_{i\in I}$ is an open cover of $X$, then $M_{\cU}$ is homotopic to $X$.
\end{lem}
\begin{proof}
An explicit proof is provided by Segal using linear homotopies. Here we take a slightly higher-brow approach.

Being a subset of the product, $M_\cU$ is equipped with natural surjective projection maps
\[
\xymatrix{ ~ & M_\cU \ar[ld]_{\rho_1} \ar[rd]^{\rho_2} & ~\\
X & ~ & N_\cU}
\]
The map $\rho_1$ has contractible fibers: for any $x \in X$, we have $\rho_1^{-1}(x) = \setof{x} \times \sigma_x$ where $\sigma_x$ is the unique simplex of maximal dimension whose support contains $x$. Thus, by Quillen's Theorem A~\cite{quillen1973higher}, the Mayer-Vietoris blowup is homotopy-equivalent to $X$ via $\rho_1$ in full generality.
\end{proof}

On the other hand, it is easy to see that the map $\rho_2$ fails to have contractible fibers precisely when the simplex supports are not contractible. In fact, given $s \in N_\cU$, the fiber $\rho_2^{-1}(s)$ has the homotopy type of the support of $\sigma_s$, which is the unique simplex of maximal dimension whose realization contains $s$. Since cohomology is a homotopy invariant, this leads to the following observation which unifies the \v{C}ech and Leray approaches.

\begin{prop}
The Leray cellular sheaves $\Leray^n$ associated to the map $\rho_2: M_\cU \to N_\cU$, where $N_\cU$ is covered by (small neighborhoods of the topological) simplices $\setof{\sigma}_{\sigma \in N_\cU}$, are isomorphic to the \v{C}ech cellular sheaves $\Cech^n$ associated to the cover $\cU$.
\end{prop}

\begin{rem}
\begin{itemize}
\item The commonality between the \v{C}ech and Leray approaches comes as no surprise to anyone sufficiently familiar with spectral sequences (and would have surprised neither \v{C}ech nor Leray).
\item Both strategies are examples of {\em distributed} cohomology computation because in order to determine the sheaf $\Cech^n$ or $\Leray^n$, one only needs to compute cohomology locally: of a non-trivial intersection of covering sets in the former case, or of a small neighborhood of the fiber $f^{-1}(x)$ in the  latter case. In principle, one can assign each local computation to a different processor, compute the appropriate sheaf cohomology over a decidedly nicer space (either $N_\cU$ or $Y$ depending on the circumstances), and aggregate this information to compute the desired cohomology of $X$.
\item By taking the appropriate linear duals and working with cosheaves, all of our results transform to computations of homology rather than cohomology.
\end{itemize}
\end{rem}

\subsection{Level Set Persistence Determines Sub-level set Persistence}
\label{subsec:level2sub}

\begin{figure}
\begin{center}
\includegraphics[width=\textwidth]{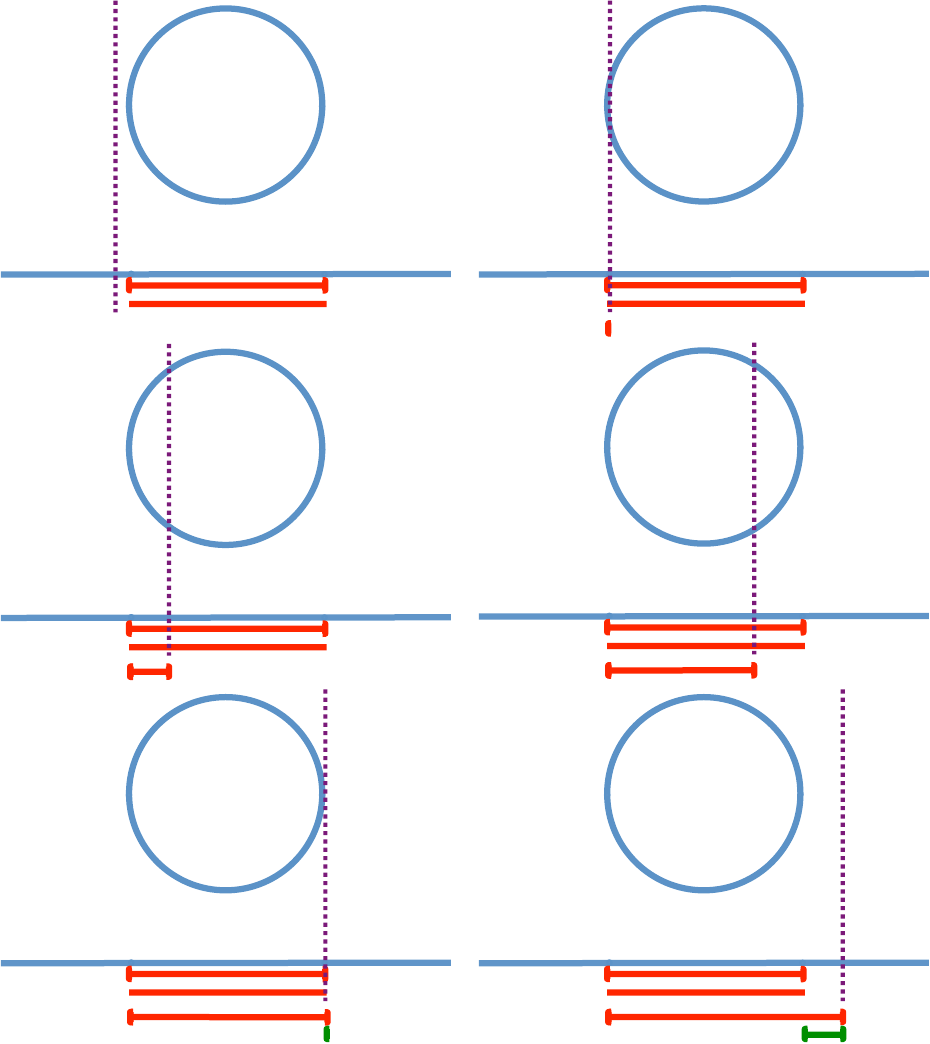}
\end{center}
\caption{Determining Sub-level Set from Level Set Persistence}
\label{fig:level_to_sub}
\end{figure}

We will use Theorem \ref{thm:leraysheaf} to obtain a non-obvious theorem in persistence, namely that level set persistence determines sub-level set persistence. By making use of \ref{clm:BM_barcodes} we illustrate how one can take the absolute homology cosheaves (or Leray sheaves) equipped with a barcode decomposition and sweep from left to right, applying \ref{cor:barcode_homology} to determine the barcodes of the associated sublevel set persistence modules. An example is drawn in Figure \ref{fig:level_to_sub}. Stated formally, we have the following theorem.

\begin{thm}[Level Set to Sublevel Set Persistence]\label{thm:level-to-sub}\index{persistent homology!level set determines sublevel set}
Let $F^k$ denote the Leray sheaf associated to a proper map $f:X\to\R$, whose stalk at $x$ is the cohomology of the fiber $H^k(f^{-1}(x))$. We can define a functor
\[
	S^k:(\RR,\leq)^{op}\to\Vect \qquad S^k(t):=H^0((-\infty,t];F^k)\oplus H^1((-\infty,t];F^{k-1})\cong H^k(f^{-1}(-\infty,t])
\]
whose value records the cohomology of the entire sublevel set. The maps
\[
	S^k(t')\to S^k(t) \qquad t\leq t'
\]
are defined sheaf theoretically by observing that we have maps
\[
	H^0((-\infty,t'];F^k)\to H^0((-\infty,t];F^k) \qquad H^1((-\infty,t'];F^{k-1})\to H^1((-\infty,t];F^{k-1})
\]
the sum of which define the desired map.
\end{thm} 

Observe that since $f:X\to\RR$ is a proper map, the restriction of $f$ to $X_{\leq t}:=f^{-1}((-\infty,t])$ is also a proper map. Consequently, we can apply Theorem \ref{thm:leraysheaf} to the space $X_{\leq t}$ instead, but we have to restrict the sheaves to the subspace $(-\infty,t]$. Fortunately, restriction of a sheaf to a subspace is a standard operation in the Grothendieck six-functor formalism, presented in~\ref{subsec:six_ops}: If $\iota:(-\infty,t]\hookrightarrow \RR$ is the inclusion, then the application of Theorem \ref{thm:leraysheaf} to the restriction reads
\[
	H^k(X_{\leq t};k)\cong H^0((-\infty,t];\iota^*F^k)\oplus H^1((-\infty,t];\iota^*F^{k-1})
\]

The upshot of this formula is that we can define a family of vector spaces, one for each $t\in\RR$ that records the cohomology of the sublevel set $X_{\leq t}$
\[
	S(t):=H^0((-\infty,t];\iota^*F^k)\oplus H_1((-\infty,t];\iota^*F^{k-1})
\]
given by computing sheaf cohomology of the restriction of the Leray sheaves to the subspace $(-\infty,t]$. What remains to be shown is that there are maps
\[
	S(t')\to S(t) \qquad t\leq t'
\]
that can be defined purely sheaf-theoretically. To do this, we will make use of some standard adjunctions in sheaf theory.

\begin{thm}
Let $f:Y\to Z$ be a continuous map. The functors $f^*:\Shv(Z)\to \Shv(Y)$ and $f_*:\Shv(Y)\to\Shv(Z)$ form an adjoint pair $(f^*,f_*)$ and thus
\[
\Hom_{\Shv(Y)}(f^*G,F)\cong \Hom_{\Shv(Z)}(G,f_*F).
\]
\end{thm}

In the above adjunction for sheaves, let $Y=(-\infty,t]$, $Z=(-\infty,t']$ and $f=j$ be the inclusion of $Y$ as a closed subspace of $Z$. Observe that if we set $F=j^*G$ in the above adjunction, then we get an isomorphism
\[
\Hom_{\Shv(Y)}(j^*G,j^*G)\cong \Hom_{\Shv(Z)}(G,j_*j^*G).
\]
that is natural in $G$. This defines the unit of the adjunction:
\[
\id_{\Shv(Y)}\to j_*j^*.
\]
We recall a basic theorem about the pushforward sheaf along a closed immersion~\cite[II.5 p. 102]{iversen}.

\begin{prop}
	For $j:Y\hookrightarrow Z$ the inclusion of a closed subspace, the functor $j_*:\Shv(Y)\to\Shv(Z)$ is exact, i.e. it sends exact sequences of sheaves to exact sequences of sheaves. Moreover, $j^*$ is always exact.
\end{prop}

\begin{lem}
	Suppose $j:Y\hookrightarrow Z$ is the inclusion of a closed subspace and $F$ is a sheaf on $Z$, then there is an induced map from the cohomology of $F$ on $Z$ to the cohomology of $j^*F$ on $Y$.
	\[
		H^i(Z;F)\to H^i(Y;j^*F)
	\]
\end{lem}
\begin{proof}
	We can read off the proof from the following diagram of spaces.
	\[
		\xymatrix{Y \ar@{^{(}->}[rr]^j \ar[dr]_{p_Y} & & Z \ar[ld]^{p_Z} \\ & \ast &}
	\]
	Sheaf cohomology is defined as the right derived functor of pushforward to a point. If we want to compute sheaf cohomology of $F$, one takes an injective resolution of $F$ 
	\[
		0\to F\to I^{\bullet}
	\]
	and applies $p_{Z*}$ to the injective resolution.
	\[
		Rp_{Z*}F:= p_{Z*}I^{\bullet}
	\]
	This results in a chain complex of vector spaces, whose cohomology is the sheaf cohomology of $F$. We usually save this step for last as it takes us out of the category of chain complexes of vector spaces and into the category of graded vector spaces. This is written as follows.
	\[
		R^ip_{Z*}F:=H^i(p_{Z*}I^{\bullet})=:H^i(Z;F)
	\]
	Since $j^*$ is exact we will consider an injective resolution of $F$ and pull that back to an injective resolution of $j^*F$. Observe the following string of identities.
	\begin{eqnarray*}
		Rp_{Y*}j^* F &:=& p_{Y*}j^*I^{\bullet} \\
		&=& (p_{Z}\circ j)_* j^*I^{\bullet} \\
		&=& p_{Z*}j_*j^* I^{\bullet}
	\end{eqnarray*}
	The unit of the adjunction defines a map of sheaves
	\[
		F \to j_*j^*F,
	\]
	which also defines a map on complexes of sheaves and hence injective resolutions.
	\[
		I^{\bullet}\to j_*j^*I^{\bullet}
	\]
	Because $j^*$ and $j_*$ are exact and preserve injectives, $j_*j^*I^{\bullet}$ is an injective resolution of $j_*j^*F$, thus
	\[
		Rp_{Z*}j_* j^* F= p_{Z*}j_*j^* I^{\bullet}= Rp_{Y*}j^* F
	\]
	and hence the unit of the adjunction defines a map
	\[
		Rp_{Z*} F\to Rp_{Y*}j^*F \quad \Rightarrow \quad H^i(Z;F)\to H^i(Y;j^*F.)
	\]
\end{proof}
\begin{rmk}[Abuse of Notation]
	Common practice in the sheaf literature is to suppress the notation $j^*F$ and to just write
	\[
		H^i(Y;F):=H^i(Y;j^*F).
	\]
	The reasoning is that $F$ is a sheaf on $Z$ and hence the only way to parse the formula on the left is to realize that the sheaf must be restricted to the subspace $Y$.
\end{rmk}

As a corollary we obtain our desired result, Theorem \ref{thm:level-to-sub}.

\section{Multidimensional Persistence}

One of the single greatest theoretical challenges to topological data analysis is a foundation for multi-dimensional persistence~\cite{lesnick-thesis}. To consider why data analysts might want such a thing, consider the following example.

\begin{figure}[ht]
\begin{center}
\includegraphics[width=.7\textwidth]{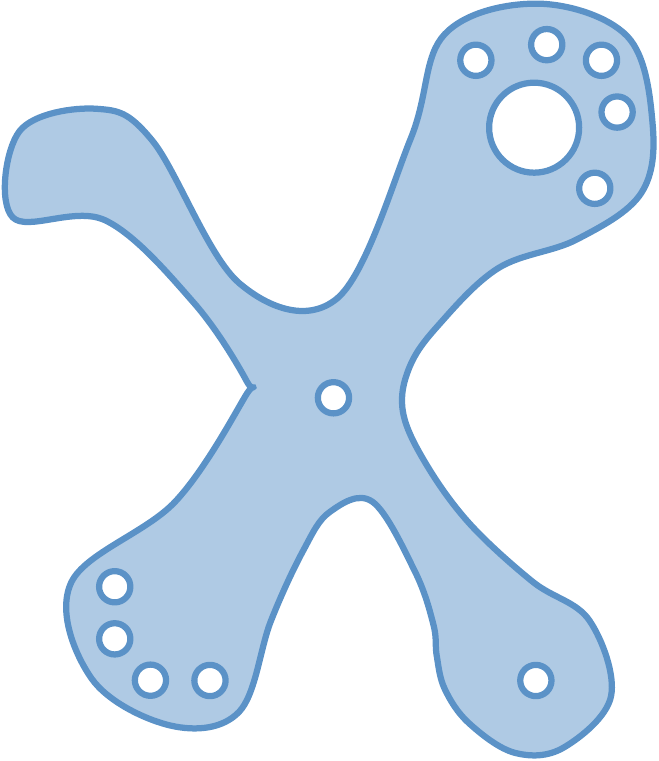}
\end{center}
\caption{A Shape Described with Multi-D Persistence}
\label{fig:flare_holes}
\end{figure}

Suppose $X$ is the shape depicted in Figure \ref{fig:flare_holes}. A common feature of interest in applications~\cite{lum2013extracting} is the presence of \emph{flares} or \emph{tendrils}. Sublevel set persistence provides a method for detecting such features. Consider the \textbf{$p$th eccentricity functional} on $X$:
\[
E^p(x):=\left(\int_{y\in X} d(x,y)^p dy\right)^p.
\]

If we filter by superlevel sets, the four endpoints of the perceived flares in Figure \ref{fig:flare_holes} will come into view. Said using homology, there are a suitable large range of values $t$ for which $E^p_{\geq t}:=\{x\in X\,|\, E^p(x)\geq t\}$ will have
\[
H_0(E^p_{\geq t};k) \cong k^4.
\]
This formally expresses the four flare-like features we see in the space $X$.

Now suppose that we are not just interested in the number of eccentric features, but rather we are interested in holes with high eccentricity value, i.e. the persistence module
\[
H_1(E^p_{\geq t};k)
\]
is of interest. However, what size of hole is of interest, and what can be regarded as noise? In other words, what is the behavior of the two-parameter family of vector spaces
\[
MP_1(t,r):=H_1((E^p_{\geq t})^r;k)
\]
where $Y^r$ denotes the set of points within distance $r$ of a subspace $Y$. This family of vector spaces defines a functor
\[
MP_1:(\R^2,\leq)\to \Vect \qquad \mathrm{where} \qquad MP_1(t,r):=H_1((E^p_{\geq t})^r;k)
\]
where $\R^2$ is viewed as a poset under the relation $(t,r)\leq (t',r')$ if and only if $t'-t$ and $r'-r$ are non-negative numbers. 

This gives the general definition of a multidimensional persistence module as introduced in~\cite{carlsson2009theory}.

\begin{defn}[Multi-dimensional Sub-Level Set Persistence]\index{persistent homology!multi-dimensional!sublevel set}
Suppose we are given a map $f:X\to\R^n$, with coordinate functions $f(x)=(f_1(x),\ldots,f_n(x))$. The \textbf{$i$th multidimensional persistence module} is defined to be the functor
\[
MP_i:(\R^n,\leq) \to \Vect \qquad (t_1,\ldots,t_n)\rightsquigarrow H_i(\{x\in X\,|\,f_j(x)\leq t_j,\, 1\leq j\leq n\};k)
\]
\end{defn}

\subsubsection{Critique of the Definition of Multi-D Persistence}

There are a few problems that such a definition suffers from:
\begin{itemize}
\item Such a definition is strongly dependent on the particular choice of basis given to $\R^n$. If one is given an abstract function $f$, valued in say a manifold $M$, then the definition fails.
\item The definition is not local.
\item There is no decomposition theorem akin to Theorem \ref{thm:crawleyboevey}. There are no barcodes. This is fundamentally due to Gabriel's Theorem \ref{thm:gabriel}, but an explicit example is provided in~\cite[Sec.5.2]{carlsson2009theory}, which suggests other invariants as well.
\end{itemize}

Sheaves and cosheaves overcome the first problem by putting level set persistence as the primary object of interest. Locality is also provided by using sheaves. Our definition of a multi-dimensional persistence module is simply the Leray sheaf.

\begin{defn}[Multi-D Persistence as the Derived Pushforward]\index{persistent homology!multi-dimensional!via derived pushforward}
Suppose $f:Y\to X$ is a continuous map. The right derived pushforward of the constant sheaf $Rf_*k_Y$ is gotten by taking a resolution of $k_Y$ by the complex of singular cochains $\widetilde{C}^{\bullet}$ and pushing forward along $f$. The Leray sheaves are the $n$th cohomology sheaves of this complex, i.e.
\[
F^n:=\cH^n(f_*\widetilde{C}^{\bullet})=:R^nf_*k_Y.
\]
If $f$ is proper, then the stalks of $F^n$ record the cohomology of the level set, i.e. $F^n_x\cong H^n(f^{-1}(x);k)$. We define \textbf{$n$th multi-dimensional level set persistence} of a \emph{proper} map $f:Y\to X$ to be the Leray sheaf $F^n$.
\end{defn}

\begin{rmk}
 If $f$ is not proper, then one can use the pushforward with compact supports to encode the compactly supported cohomology of the fiber.
\end{rmk}

One can always obtain the traditional sub-level set definition from this definition by using a multi-dimensional version of Theorem \ref{thm:level-to-sub}, the Leray spectral sequence.

Sheaves also suffer from the lack of a nice set of indecomposables, but the Grothendieck operations provide one possible approach.

\subsection{Generalized Barcodes}
\label{subsec:barcode_ss}

The need for a generalized notion of a barcode comes from the need to communicate topological summaries to non-topologists. A scientist can easily understand a histogram and the barcode is only subtly different from a histogram.

\begin{defn}[Generalized Barcodes]\label{defn:generalized_bc}\index{barcode!generalized}\index{decomposition!generalized barcode}
	Suppose $F^n$ is the Leray sheaf associated to a proper map $f:Y\to X$. A \textbf{generalized barcode} for $F^n$ is the expression of $F$ as follows.
	\[
	F^n\cong \bigoplus_{b\in B} (j_b)_! k_{Z_b}
	\]
	where $j_b:Z_b\to X$ are maps indexed by the barcode $B$ and where each $(j_b)_! k_{Z_b}$ is assumed to be indecomposable.
\end{defn}
\begin{rmk}
It is not clear when or if such a generalized barcode exists for a given $f:Y\to X$ and $n$.
\end{rmk}

\begin{figure}[ht]
\begin{center}
\includegraphics[width=.7\textwidth]{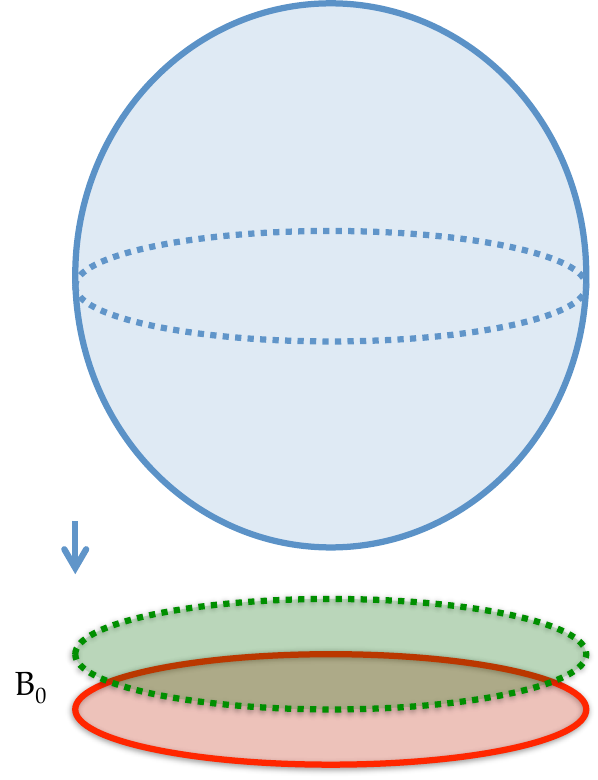}
\end{center}
\caption{Two Dimensional Barcodes for the Sphere}
\label{fig:sphere2d_bc}
\end{figure}

An example barcode is provided in Figure \ref{fig:sphere2d_bc} and discussed in the next example.

\begin{ex}[Shadow of the Sphere]
	Consider the standard Euclidean sphere $S^2$ embedded in $\RR^3$. Let $f:S^2\to\RR^2$ be the projection onto the first two factors of $\RR^3$. The image $X=f(S^2)$ is the closed unit disk.

	One can use a multi-D analog of Corollary~\ref{cor:barcode_homology} to determine the homology of the two-sphere. 
	\[
		H_0(S^2)\cong H_0(X;\hF_0)\cong k \qquad H_1(S^2)\cong H_1(X;\hF_0)\cong 0 \qquad H_2(S^2)\cong H_2(X;\hF_0)\cong k
	\]
\end{ex}

%
%

\chapter{Network Coding and Routing Sheaves}
\label{sec:nc_coding}

In the previous section, we introduced the language of barcodes and integrated them with a cosheaf-theoretic perspective on Morse theory and persistent homology. The fundamental idea was that by decomposing a cosheaf into indecomposables, we were able to understand cosheaf homology via the Borel-Moore homology of the barcode. In this section, we attempt to do the same thing for cellular sheaves on graphs. We apply the barcode perspective, wherever possible, to a class of sheaves introduced by Robert Ghrist and Yasuaki Hiraoka~\cite{GH-ncs}. These sheaves were specifically designed to model the flow of information over graphs and the generalized barcode decomposition can aid in visualizing this flow. 

First, we review some basic definitions for graphs.

\begin{defn}\index{quiver}\index{graph}
	Let $X$ be a directed graph consists of a pair of sets $E$ and $V$ of edges and vertices and a pair of functions $h,t:E\to V$ that return the \textbf{head} and \textbf{tail} of an edge respectively. A directed edge goes from its tail to its head. The set of \textbf{incoming edges} to a vertex $v$, written $in(v)$, is the set of edges whose heads are $v$, i.e. $h^{-1}(v)$. The set of \textbf{outgoing edges} at $v$ is the set of edges whose tails are at $v$, i.e. $t^{-1}(v)=v$.
\end{defn}

\begin{defn}\index{capacity}
	Let $X$ be a directed graph with vertex set $V$ and edge set $E$. A \textbf{capacity function} is a function $c$ from the edge set to either the non-negative reals $\RR_{\geq 0}$ or the non-negative integers $\ZZ_{\geq 0}$.
\end{defn}

\begin{defn}[Network Coding Cell Sheaf]\index{network coding sheaf}\index{sheaf!network coding}
	Suppose $X$ is a directed graph with a capacity function $c$. A \textbf{network coding sheaf} on $X$ is a cellular sheaf $F:X\to\Vect$ constructed as follows: 
	\begin{itemize}
		\item To an edge $e\in X$ we let $F(e)=k^{c(e)}$, a vector space of dimension equal to the capacity.
		\item To a vertex $v$ we let $F(v)=k^{c(v)}\cong \oplus_{e_i\in in(v)} k^{c(e)}$.
		\item The restriction maps are given by ordinary projections for the incoming edges, i.e. $\rho_{e_i,v}:=proj_{F(e_i)}$ for $e_i\in in(v)$, but for the outgoing edges some non-trivial coding may be performed, i.e. any linear map $\Phi_{e_k,v}:F(v)\to F(e_k)$ for $e_k\in out(v)$ will do. We write $\Phi_v=\oplus_{e_k\in out(v)} \Phi_{v,e_k}$ for the \textbf{total coding} through $v$.
	\end{itemize}
\end{defn}

\begin{rmk}
	It should be noted that in~\cite{GH-ncs}, they do not use cellular sheaves. This was primarily due to the lack of a good reference.
\end{rmk}

\begin{figure}[ht]
\centering
\includegraphics[width=.5\textwidth]{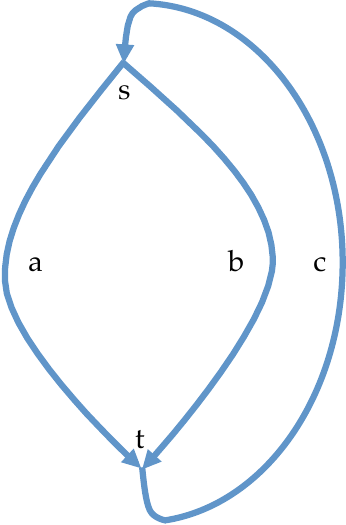}
\caption{Graph with Decoding Wire}
\label{fig:nc_graph_dir}
\end{figure}

In~\cite{GH-ncs} they do not define network coding sheaves for arbitrary directed graphs. Instead, they consider a graph with a distinguished set of sources and targets (senders and receivers) and they augment the graph by adding \textbf{decoding wires} directed to go from a target vertex back to a subset of source vertices. Heuristically for Ghrist and Hiraoka, the purpose of these edges is to make global sections of this sheaf correspond to closed loops through the graph. This topological reasoning is correct, but oversimplifies how network codings can produce counterintuitive weavings and splittings of data.

\begin{ex}
	Consider the graph in Figure~\ref{fig:nc_graph_dir} with constant capacity function $c = 1$. Consequently, all edges and vertices get a one dimensional vector space $k=\RR$ with the exception of $F(t)\cong k^2$. Define the coding maps $\rho_{a,s}=\id=\rho_{b,s}$ and
	\[
		\rho_{c,t}=\begin{bmatrix} \frac{1}{2} & \frac{1}{2}\end{bmatrix}.
	\]
	We pick a local orientation implied by the source and target vertices. The one and only coboundary matrix can be written as follows:
	\[
		\delta^0:=\begin{bmatrix}
				-1 & 1 & 0 \\
				-1 & 0 & 1 \\
				1 & -1/2 & -1/2
		\end{bmatrix}
	\]
	Consequently, $H^0(X;F)\cong H^1(X;F)\cong k$. The one global section is supported over the entire graph; it is not simply a loop through the graph.
\end{ex}

The previous example of a network coding sheaf is an example of an indecomposable sheaf that is not a generalized barcode in the sense of Definition~\ref{defn:generalized_bc}. To better understand the flow of data over graphs, as well as the utility of the barcode perspective, we consider a simpler class of network coding sheaves.

\section{Duality and Routing Sheaves}

\begin{defn}[Routing Sheaf]\index{routing sheaf}\index{sheaf!routing}
	A particular type of network coding sheaf is a \textbf{routing sheaf}. Here we assume that the capacity function is constant $c = 1$, and the coding maps $\Phi_v$ can be written as a binary matrix with at most one $1$ in each column and row. Said another way, at a vertex $v$ the total coding map maps to zero as many incoming edges as desired, so long as there is a bijection of the remaining incoming edges and a subset of the outgoing edges. The total coding map through $v$ is then a matrix representation of this bijection.
\end{defn}

The advantage of routing sheaves is that they are simple to visualize: Start at a source and use a color pen to track how an edge emanating from that source gets bounced around under the routing directions at each subsequent vertex. If at any point in your drawing you run into a vertex that sends your edge's data to zero, stop on that vertex with your pen. This argument essentially establishes the following proposition. 

\begin{prop}[Structure Theorem for Routing Sheaves]\label{prop:rout_bc}
	Suppose $X=(V,E,h,t)$ is a directed graph, then every routing sheaf $F:X\to\Vect$ can be realized as the pushforward with compact support of a disjoint union of half-open intervals $[\textrm{---}[$ or circles, whose images can intersect only at vertices.
\end{prop}

A consequence of this result combined with Poincar\'e duality is the following corollary.

\begin{cor}[Duality]\index{duality!for routing sheaves}\index{routing sheaf!duality}
	For any routing sheaf $F$ one has
	\[
		H^0(X;F)\cong H^1(X;F).
	\]
\end{cor}
\begin{proof}
	By Proposition~\ref{prop:rout_bc}, every network coding sheaf can be written as a direct sum of constant sheaves supported on half-open intervals or circles. Half-open intervals embedded into compact spaces (extending by zero using $j_!$) have trivial cohomology in both degrees. Poincar\'e duality for $S^1$ establishes the corollary.
\end{proof}

One can also prove this duality in the more general setting of network coding sheaves via a simple combinatorial argument.

\begin{prop}
	For an any network coding sheaf $F$, we have the following isomorphisms:
	\[
		\bigoplus_v F(v) \cong \bigoplus_e F(e) \qquad  H^0(X;F)\cong H^1(X;F)
	\]
\end{prop}
\begin{proof}\index{duality!network coding sheaf}\index{network coding sheaf!duality}
	By construction of a network coding sheaf there is a bijection between the sum of the vector spaces over the edges $e\in in(v)$ and the vector space over the vertex $v$.
	\[
		\bigoplus_{e\in in(v)} F(e) = F(v).
	\]
	By definition of a graph, every edge is the incoming edge for a unique vertex. Thus, by summing over all vertices, we sum over all edges without double-counting. This proves the first isomorphism. The second isomorphism follows by the rank-nullity theorem.  
\end{proof}

Ideally, one could interpret these cohomology groups as something meaningful to obtain a useful duality result, but this is still missing. Ghrist and Hiraoka interpret $H^0$ as a vector space spanned by independent information flows, but in the case of routing sheaves, $H^1$ is the group that counts closed trajectories of information flow. For routing sheaves, one could say that the Poincar\'e dual of the fundamental class of an information loop would yield a point whose removal would cease the flow of information. One might call this a ``cut equals flow'' theorem. This is only a pale shade of the greater ``Max-Cut Min-Flow'' theorem~\cite{efs-mfmc, ff-mfmc, seymour-matroids}, which compares the maximum possible flow with the minimum capacity cut required to disconnect the graph.

\section{Counting Paths Cohomologically, or Failures thereof}

Regardless of network coding sheaves connections with duality, one would like to know what information cellular sheaf cohomology can capture over a graph. Is it possible, for example, to build a sheaf that encodes source-to-target paths cohomologically? Suppose we allow the source to have its own independent capacity, without regard for the number of incoming edges. As the next example shows such a ``pseudo network coding'' (NC) sheaf cannot encode source-to-target paths cohomologically.

\begin{figure}[!ht]
\begin{minipage}[b]{0.47\linewidth}
\centering
\includegraphics[scale=1]{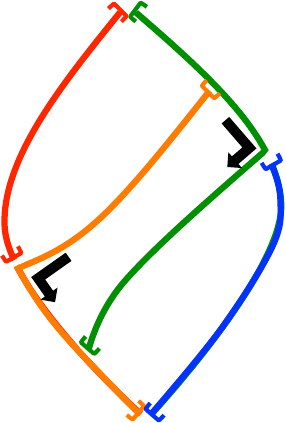}
\caption{No Decoding Edge}
\label{fig:ncs_no_decode}
\end{minipage}
\hfill
\begin{minipage}[b]{0.47\linewidth}
\centering
\includegraphics[scale=1]{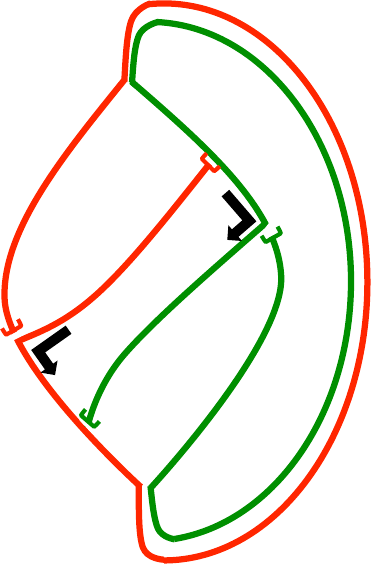}
\caption{Decoding Edge}
\label{fig:ncs_decode}
\end{minipage}
\end{figure}

\begin{ex}[Decoding Edge and Barcodes]\index{decoding edge}
In Figures~\ref{fig:ncs_no_decode} and~\ref{fig:ncs_decode}, the barcode decomposition of a network coding sheaf is drawn with and without a decoding edge.
With the particular choices made there is no flow from source to target. Yet the sheaf in Figure~\ref{fig:ncs_no_decode} decomposes as the constant sheaf on two half open intervals and two closed intervals:
\[
F_{no}\cong (j_o)_! k_{[0,1)}\oplus(j_b)_! k_{[0,1)}\oplus (i_r)_* k_{[0,1]}\oplus (i_g)_*k_{[0,1]} \Rightarrow H^0(X;F_{no})\cong k^2.
\]
This is bad if we want our sheaf to encode cohomologically the presence of source-to-target information paths.

However, with the use of a decoding edge (decoding edges) the network decomposes into only two half open intervals:
\[
F_{de}\cong (j_r)_!k_{[0,1)}\oplus(j_g)_!S k_{[0,1)} \Rightarrow H^0(X;F_{de})\cong 0
\]
which was wanted.
\end{ex}

\section{Network Coding Sheaf Homology}
\label{subsubsec:sheaf_homology_graphs}
\index{sheaf!homology, cellular!examples}
One of the virtues of the network coding sheaves is that they are easy to construct, have interesting sheaf cohomology, and provide lots of examples. As such, we will use them as a testing ground for the new theory of sheaf homology developed in Section~\ref{subsec:sheaf_homology}.

\begin{figure}[ht]
\centering
\includegraphics[width=.5\textwidth]{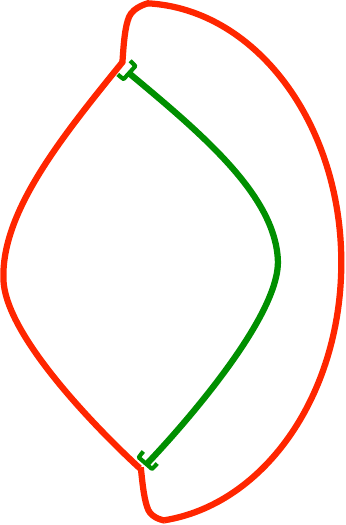}
\caption{Network Coding Sheaf}
\label{fig:nc_1}
\end{figure}

\begin{ex}
Consider the network coding sheaf implied by Figure~\ref{fig:nc_1}. Viewed as a diagram of vector spaces, it takes the following form:
\[
\xymatrix{ & \ar[ld]^1 k_s \ar[rd]_0 \ar[rrd] & & \\
		k_a & & k_b & k_c \\
		& k_t^2 \ar[ul]_{\pi_1} \ar[ur]^{\pi_2} \ar[urr]_{\pi_1} & & }
\]
Since the category of complexes of sheaves is additive, we can consider each indecomposable sheaf separately and compute its sheaf homology. If one considers just the red loop as a constant sheaf (barcode) $R$, it takes the following form:
\[
\xymatrix{ & \ar[ld]^1 k_s \ar[rd] \ar[rrd] & & \\
	k_a & & 0 & k_c \\
        & k_t \ar[ul]_{1} \ar[ur]^{} \ar[urr]_{1} & & }
\]
A projective sheaf that surjects onto $R$ is supported on the star of $s$ and $t$ respectively, i.e. $P_0:=\{s\}\oplus\{t\}$:
\[
\xymatrix{ & \ar[ld]^1 k_s \ar[rd] \ar[rrd] & & \\
	k_s\oplus k_t & & k_s\oplus k_t & k_s\oplus k_t \\
        & k_t \ar[ul]_{1} \ar[ur]^{} \ar[urr]_{1} & & }
\]
The kernel sheaf of the natural transformation $P_0\Rightarrow R$ is also projective, which we call $P_1$ and finishes the projective replacement of the sheaf $R$.
\[
\xymatrix{ & 0 \ar[ld] \ar[rd] \ar[rrd] & & \\
	[1\, -1] & & k_s\oplus k_t & [1\, -1] \\
        & 0 \ar[ul] \ar[ur]^{} \ar[urr] & & }
\]
If we take the colimit of $P_1$ and $P_0$ separately, the sheaf map $P_1\to P_0$ induces a map on colimits that defines the boundary in the chain complex computing sheaf homology:
\[
\partial_1: k^4\to k^2 \qquad 
\begin{bmatrix} 1 & 1 & 0 & 1 \\
-1 & 0 & 1 & -1\end{bmatrix}
\quad \Rightarrow \quad H_0(X;R)=0 \quad H_1(X;R)\cong k^2
\]
 Repeating the same reasoning for the green barcode $G$ yields homology groups $H_0(X;G)=0$ and $H_1(X;G)\cong k$. Since our original network coding sheaf $F$ is a direct sum $R\oplus G$ we obtain that the sheaf homology of the sheaf in Figure~\ref{fig:nc_1} is
\[
H_0(X;F)=0 \qquad H_1(X;F)\cong k^3.
\]
\end{ex}

\begin{figure}[ht]
\centering
\includegraphics[width=.5\textwidth]{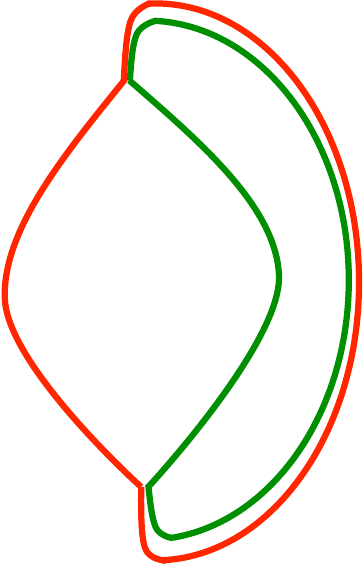}
\caption{Network Coding Sheaf with Two Decoding Wires}
\label{fig:nc_2}
\end{figure}

\begin{exr}
As an exercise, and to indicate the sensitivity of sheaf homology to its embedding, we ask the reader to verify that the sheaf homology groups of the sheaf in Figure~\ref{fig:nc_2} are
\[
H_0(X;F)=0 \qquad H_1(X;F)\cong k^8.
\]
\end{exr}

%
%

\chapter{Sheaves and Cosheaves in Sensor Networks}
\label{sec:sensors}

\begin{quote}
{\em``What strength, what art can then suffice, or what evasion bear him safe through the strict sentries and stations thick of angels watching round?''}
\begin{flushright} --- John Milton, Paradise Lost, Book II, Line 410~\cite{milton2008paradise} \end{flushright}
\end{quote}

In this section, we consider a candidate application of sheaves and cosheaves to problems in sensor networks. Section \ref{subsec:sensors} outlines some real-world  sensors as well as their mathematical abstraction. With this abstraction in hand, we consider in Section \ref{subsec:coverage} the classic problem of determining when a sensor network has completely covered a region. The introduction of time-dependent sensor networks necessitates the sheaf-theoretic approach, despite the fact that it is unwieldy in its most general form.

In Section \ref{subsec:intruder_bcs} we attempt to ``linearize'' the sheaves and cosheaves used in studying sensor networks in the hope that sheaf cohomology and cosheaf homology will give us an obstruction-theoretic approach to sensing. An approach of Henry Adams is considered in Section \ref{subsubsec:ttt_bcs}, as well as his counter-example to that approach. By using cosheaf-theoretic reasoning, we give a principled explanation for why this approach fails in Proposition \ref{prop:long_bcs}. An approach of the author and Robert Ghrist is then considered in Section~\ref{subsubsec:lss_bcs}. This approach succeeds where the previous approach fails, but it too suffers from giving false positives, as the example in Proposition \ref{prop:mobile_d4} shows. The example constructed there, which is joint with David Lipsky, uses one of the 12 indecomposable representations of the Dynkin diagram $D_4$.

Finally, a linear model for multi-modal sensing is presented in Section \ref{subsec:multi_modal}. It was there that the author realized the necessity of using indecomposables to interpret sheaf cohomology computations. A delightful examination of the act of sensing in Section \ref{subsubsec:deep_sense} shows how sheaves and cosheaves work in tandem. Theorem \ref{thm:les_sensing} uses a long exact sequence in sheaf cohomology to obtain a forcing result in multi-modal sensing. Finally, the role of higher-dimensional barcodes in multi-modal sensing is considered in Section \ref{subsubsec:evasion_bcs}. 

\section{A Brief Introduction to Sensors}
\label{subsec:sensors}\index{sensors}

Sensors are devices with delimited purview. They can measure certain properties and interact with occupants of a particular part of space-time. Examples abound in our world and they operate via differing modalities. Here are a few examples:

\begin{ex}[Sight]
	Our eyes are highly tuned sensors that can detect photons with certain frequencies (visible light and colors) and their spatial range can be on the order of kilometers. Some man-made satellites orbiting the Earth have cameras with a greater spatial resolution and frequency response --- they help us navigate by providing detailed pictures of roads, weather and climate. Eyes and satellites have a large scale and are very expensive. Cheaper sensors which can read only very coarse changes in light levels are found in our traffic lights, door ways and bathrooms.
\end{ex}

\begin{ex}[Weight and Pressure]
	Buried in roads or placed under door mats are sensors designed to respond to pressure. These open doors or gates or initiate changes in traffic signals. Some are more passive and merely collect data. A cable as thick as a thumb can be laid across a road and will record when something heavy (like a car) drives over it. Two spikes in pressure close in time indicate when a car's front and back wheels respectively drove over the cable. From this city officials can measure how fast cars are going as well as density and total volume of traffic.
\end{ex}

\begin{ex}[Radio Frequency ID]
	Some readers probably have a university card, or building card, that grants them access through locked doors merely by tapping on a sensor. Commuters drive cars equipped with sensors that allow them to pass through tolls without stopping. Some scientists tag animals to study a species' habits and movements. In all these cases, the sensor or the tag emits an electromagnetic field with limited spatial range (a few centimeters, meters, or kilometers) and only when inside this range is a tuned circuit thereby completed, connecting the sensors (card reader, toll booth, etc.) with the things being sensed (ID card, tag, etc.).
\end{ex}

Although the physical mechanisms that allow each of these sensors to sense is different, there are some broad commonalities: spatially localized sensors return data in the presence of certain occupants, which we call \textbf{intruders}. 

\section{The Coverage Problem: Static and Mobile}
\label{subsec:coverage}

The way we model sensors is to first identify the physical domain where the sensing is taking place --- a two-dimensional Euclidean plane could represent the floor of a building --- and we represent the sensors spatially via their support --- a door mat with pressure sensors would be a rectangle in the plane. Or we could think of the field of view of a camera in a ceiling pointed directly down as a disk in the plane. 

For the moment we ignore the type of data a sensor reports (we'll take that up later when we work with sheaves and cosheaves) and instead we consider the \textbf{coverage problem}: Given a collection of sensors distributed in a physical domain $D$, can we monitor the entire region without gaps?

\begin{figure}
\centering
\includegraphics[width=\textwidth]{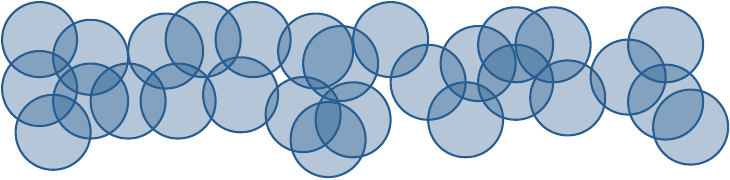}
\caption{Sensors Distributed in a Plane}
\label{fig:sensors5in}
\end{figure}

If we have good knowledge of the sensors which live on the boundary of our region, then we can, following the work of Vin de Silva and Robert Ghrist~\cite{VdSRG}, give a certificate of coverage using relative homology. However, we frame this question using sheaves of sets instead, so as to better handle the time-dependent scenario.\footnote{The author would like to thank Gunnar Carlsson and Rob Ghrist for their insights here.}

\begin{defn}\index{intruder problem}
Let $D$ be a spatial region of interest and denote by $D\times[0,1]$ a region of space-time. This carries with it a map that keeps track of time via projection onto the second factor, i.e. $\pi_2:D\times [0,1]\to[0,1]$. We assume that $D$ can be given a cell structure so that the sensors' \textbf{coverage region} $S\subset D\times [0,1]$ and the \textbf{evasion region} $E:=S^c$ can be written as the union of cells. To study the \textbf{intruder problem} is to analyze the associated sheaf of sections of the map $\pi:=\pi_2|_E :E\to [0,1]$, which we assume can be made cellular. Saying that there is an \textbf{evasion path} is to say there is a global section of this map, i.e. a $s:[0,1]\to E$ such that $\pi\circ s=\id$.
\end{defn}

\begin{figure}
\centering
\includegraphics[width=\textwidth]{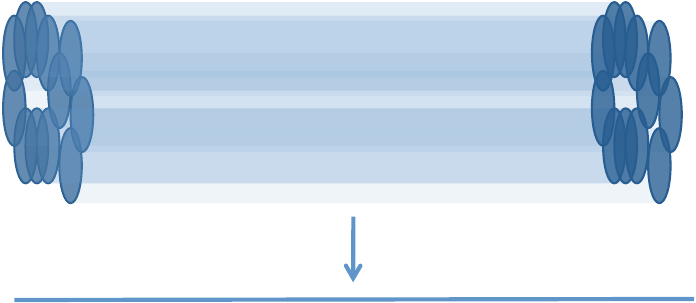}
\caption{The Space-Time Perspective}
\label{fig:coverage_section}
\end{figure}

\begin{ex}
For the situation depicted in Figure \ref{fig:coverage_section}, the intruder problem has a clear answer. An intruder can evade detection by residing in either one of the two holes present. Picking a point and then resting there for all time determines a global section of the time projection map. 
\end{ex}

It should be clear that our sheaf-theoretic question is equivalent to a much simpler one: ``Is the complement of the sensed region (the uncovered region) in $D$ non-empty?'' Thinking in terms of sheaves, at this point, buys us nothing. 

Where sheaves begin to offer a hint of leverage is in the time-dependent scenario. Here we imagine the sensors can move around in our domain $D$. Now it is possible that the sensed region $S$ does not look like a product of space and time.
 
\begin{figure}
\begin{center}
\includegraphics[width=.7\textwidth]{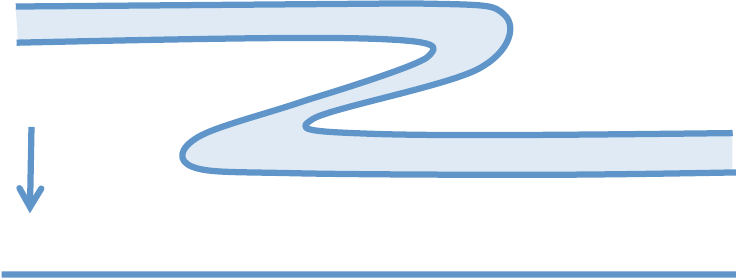}
\caption{Mobile Sensor Network}
\label{fig:mobile_s}
\end{center}
\end{figure}

\begin{ex}
In Figure \ref{fig:mobile_s} we imagine that there is a one-dimensional environment of interest that sits vertically over each point on the time axis. Between the black lines is a region that is currently being unmonitored. To begin there is only one connected component of the unmonitored region. As time marches forward to the right a second connected component of the unmonitored region opens up, followed shortly by a third. Two of these three merge and then disappear leaving only one component of unmonitored territory.

In this case the non-existence of an evasion path is clear: no intruder could have gone undetected without time-traveling. This corresponds to the ready-seen fact that this map has no section, i.e. there does not exist a map $s:[0,1]\to E$ such that $\pi\circ s=\id$.
\end{ex}

What is the purpose of considering sheaves at all? If we can stare at the drawing and detect whether a section exists or not, why bother with high-flown machinery? However, what is easily seen in toy examples, can quickly become unmanageable. The only mathematics that formalizes intuition about sections is sheaf theory and moreover, once formalized using cellular sheaves, it can be programmed on a computer. 

However, there is a disadvantage with using sheaves of sets. We'd like to be able to calculate an obstruction that would certify whether a global section exists or not. One of the stated purposes of using sheaf cohomology is to provide such a calculable obstruction. Unfortunately, cohomology requires the linear structure of vector spaces, which we do not have here. In the next section we consider what happens when we na\"ively ``linearize'' the sheaf of sections of a map.

\section{Intruders and Barcodes}
\label{subsec:intruder_bcs}

In this section, we use cellular sheaves and cosheaves to analyze the intruder problem in the time-dependent case. We assume that the time projection map $\pi$ is cellular in order to take advantage of the functors in Section \ref{sec:maps}.
By putting a sheaf or cosheaf on the evasion region and pushing forward along $\pi$, we reduce the intruder problem to one dimension where we can use the barcode perspective of Section~\ref{sec:barcodes}. There are two main approaches, both of which have their drawbacks:

\begin{itemize}\index{intruder problem!approaches to}
	\item One approach is to study the homology of the evasion region at each moment in time $\pi^{-1}(t)$. By Theorem~\ref{thm:strat_maps}, this determines a cellular cosheaf.
	\item The second approach is to linearize the space of sections of the map $\pi$. To make the space of sections finite, we pass to the Reeb graph of the evasion region. This determines a cellular sheaf and stays true to the original intruder problem.
\end{itemize}

\subsection{Tracking the Topology over Time}
\label{subsubsec:ttt_bcs}

\begin{figure}
\begin{center}
\includegraphics[width=.7\textwidth]{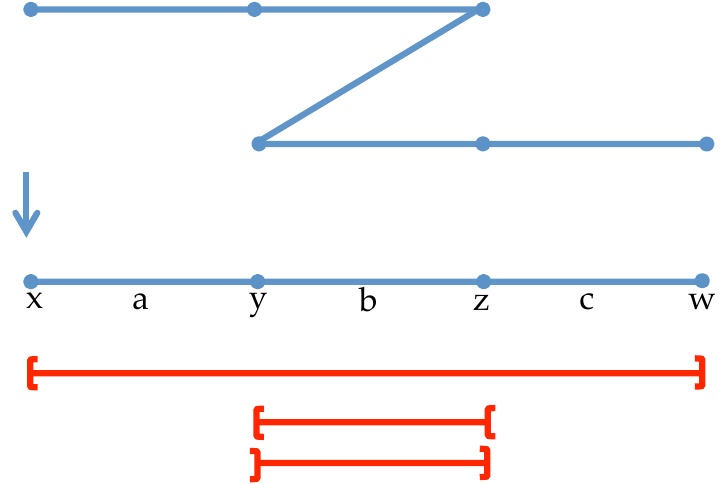}
\caption{Mobile Sensor Network}
\label{fig:mobile_reeb_bc}
\end{center}
\end{figure}

To simplify the topology, we focus on the Reeb graph version of Figure \ref{fig:mobile_s}. This is drawn and labeled in Figure \ref{fig:mobile_reeb_bc}. Since everything is occurring in two-dimensional space-time, the only interesting homological invariant of the fiber is $H_0$. Studying this is equivalent to studying the pushforward cosheaf $\hF:=\pi_*\hat{k}_E$. In the parlance of~\cite{zigzag}, this is simply a zigzag module of the following form:
\[
	\xymatrix{\hF(x) & \ar[l]_{r_{x,a}} \hF(a) \ar[r]^{r_{y,a}} & \hF(y) & \ar[l]_{r_{y,b}} \hF(b) \ar[r]^{r_{z,b}} & \hF(z) & \ar[l]_{r_{z,c}} \hF(c) \ar[r]^{r_{w,c}} & \hF(w) \\
	k_x & \ar[l] k_a \ar[r]^{} & k^2_y & \ar[l] k^3_b \ar[r] & k^2_z & \ar[l] k_c \ar[r] & k_w}
\]

If we choose for each cell in $[0,1]$ the ordered basis given by the top down ordering on the page of the cells in the fiber we get the following matrix representations of the extension maps:
\[
	r_{y,a}=\begin{bmatrix} 1\\ 0\end{bmatrix} \quad r_{y,b}=\begin{bmatrix} 1 & 0 & 0\\ 0 & 1 & 1\end{bmatrix} \quad r_{z,b}=\begin{bmatrix} 1 & 1 & 0\\ 0 & 0 & 1\end{bmatrix} \quad r_{x,c}=\begin{bmatrix} 0\\ 1\end{bmatrix}
\]
We can decompose this cosheaf into indecomposables simply by performing the correct change of basis:
\[
	\begin{bmatrix} y_1'= y_1 \\ y'_2=y_1-y_2\end{bmatrix} \qquad \begin{bmatrix} b_1'= b_1-b_2+b_3 \\ b'_2=b_1-b_2 \\ b_3'=b_2-b_3 \end{bmatrix} \qquad \begin{bmatrix} z_1'= z_2 \\ z'_2=z_1-z_2\end{bmatrix}
\]
The reader should match the resulting indecomposables with the barcodes drawn in Figure \ref{fig:mobile_reeb_bc}.
\[
\xymatrix{k_x & \ar[l] k_a \ar[r]^{} & k_{y_1'} & \ar[l] k_{b_1'} \ar[r] & k_{z_1'} & \ar[l] k_c \ar[r] & k_w \\
0 & \ar[l] 0 \ar[r] & k_{y_2'} & \ar[l] k_{b_2'} \ar[r] & 0 & \ar[l] 0 \ar[r] & 0 \\
0 & \ar[l] 0 \ar[r] & 0 & \ar[l] k_{b_3'} \ar[r] & k_{z_2'} & \ar[l] 0 \ar[r] & 0
}
\]

The presence of a long barcode may seem surprising. It indicates that there is a connected component of the evasion region that persists for all time. The following proposition explains why this long barcode must exist.

\begin{prop}\label{prop:long_bcs}
	Suppose $E\subset D\times [0,1]$ is a compact connected evasion region such that $\pi=\pi_2|E$ is surjective, i.e. there is at each point in time somewhere an intruder can evade detection, then the Remak decomposition of $\pi_*\hat{k}_E$ must have a barcode that stretches the length of $[0,1]$.
\end{prop}
\begin{proof}
	The proof starts with the easy observation that if $f:Y\to X$ is a continuous map and $\hG$ is a cosheaf on $Y$, then we have that $H_0(Y;\hG)\cong H_0(X;f_*\hG)$. This follows from the commutativity of the following diagrams and functoriality of pushforward.
	\[
		\xymatrix{Y \ar[rr]^f \ar[rd]^p & & X \ar[ld]_p \\ & \star &} \quad \xymatrix{\hG \ar[rr]^{f_*} \ar[rd]^{p_*} & & f_*\hG \ar[ld]_{p_*} \\ & p_*\hG\cong (p\circ f)_*\hG &}
	\]
	Setting $Y=E$, $X=[0,1]$, $f=\pi$ and $\hG=\hat{k}_E$, we can use the fact that $E$ is connected to get that $p_*\hat{k}_E\cong H_0(E;\hat{k}_E)\cong k$. We know that any (co)sheaf over $[0,1]$ can be written as a direct sum of constant (co)sheaves supported on barcodes.
	\[
		\pi_*\hat{k}_E\cong \hat{k}_{B_1}\oplus \cdots \oplus \hat{k}_{B_n} \qquad \qquad \mathrm{and}
	\]
	 Now we combine this with the fact that homology commutes with direct sums.
	\[
		k\cong H_0(E;\hat{k}_E)\cong H_0([0,1];\pi_*\hat{k}_E)\cong \bigoplus_i  H_0([0,1];\hat{k}_{B_i})\cong \bigoplus_i  H^{BM}_0(B_i) .
	\]
	Consequently, there can be only one closed barcode. We argue that this unique closed barcode must have support on all of $[0,1]$. Since we know that the constant section $1_E\in \Gamma(E;\hat{k}_X)$ has support on all of $E$, the pushforward section $\pi_*1_E$ that generates the closed barcode must have support on all of $[0,1]$, since $\pi$ is surjective.
\end{proof}

\begin{rmk}
	We have implicitly used sheaf-theoretic reasoning with $H^0$ taking the place of $H_0$. The argument about the support of the section is better expressed using stalks.
\end{rmk}

 As a consequence, we obtain a negative result, which is almost identical to a result of Henry Adams.

\begin{cor}
	Having a barcode associated to $\pi_*\hat{k}_E$ whose support is all of $[0,1]$ does not indicate the existence of an evasion path.
\end{cor}

\begin{rmk}
The above proof gives a cosheaf-theoretic explanation of why we shouldn't expect barcodes to detect the existence of an evasion path. Homology of the evasion region is not sensitive to its embedding, thus a long barcode will appear even if it is embedded in a way that would require an intruder to time travel. In this sense, Corollary \ref{cor:barcode_homology} can be interpreted as a stability result: although half-open barcodes can pop in and out of existence, based on the embedding, there must always be one and only one closed barcode.
\end{rmk}

\subsection{Linearizing the Sheaf of Sections}
\label{subsubsec:lss_bcs}

\begin{figure}[ht]
\begin{center}
\includegraphics[width=.5\textwidth]{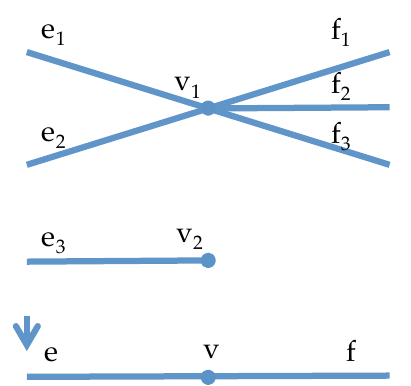}
\caption{Sheaf of All Possible Evasion Paths}
\label{fig:route_stack}
\end{center}
\end{figure}

In light of the inability of the pushforward cosheaf $\pi_*\hat{k}_E$ to distinguish when an evasion path exists or not, we return to the original sheaf-theoretic formulation of the intruder problem. To make the sheaf of sections finite enough to work with, we take the Reeb graph of the map $\pi:E\to[0,1]$. From this setup, we can extract a cellular map of 1-dimensional cell complexes, normally called $\tilde{\pi}:R(\pi)\to [0,1]$, but we will abuse notation and assume that our input $\pi: E\to [0,1]$ is already a Reeb graph.

By picking a directionality of $[0,1]$ we can endow $E$ with the structure of a directed graph. On this directed graph we can define the following cellular sheaf, which is meant to pushforward to a linear model of the sheaf of sections. It is very closely related\footnote{In a sense, the sheaf defined here gives all possible codings. It approximates a ``stack'' of network coding sheaves.} to the network coding sheaves defined in Section \ref{sec:nc_coding}.

\begin{defn}
	Let $X$ be an acyclic directed graph. We define a cellular sheaf $G$ that assigns to an edge the one-dimensional vector space $k$ and assigns to a vertex the space freely generated by all possible directed routings through that vertex. 
	
	We allow special treatment to a subset of sources $S$ and sinks $T$, where we allow $G(v)=k^{\out(v)}$ for $v\in S$ and $G(v)=k^{\inn(v)}$ for $v\in T$. All other sources and sinks get the zero vector space. The restriction mappings send a routing to the edges that participate in that routing.	
\end{defn}

\begin{ex}
For a concrete example, where we focus on a small part of a graph, consider the graph in Figure \ref{fig:route_stack}. The definition of the sheaf $G$ makes
\[
	G(v_1)=<e_i\otimes f_j \, | \, i=1,2; \,j=1,2,3>\cong k^6 \qquad \rho_{e_i,v_1}(e_j\otimes f_k)=\delta_{ij} \quad \rho_{f_j,v_1}(e_i\otimes f_k)=\delta_{jk} 
\]
and we choose to make $G(v_2)=0$. The reason we have decided to set $G(v_2)=0$ comes from the extra information of the projection map to $[0,1]$. We call such a vertex an \textbf{internal} source or sink. In the context of the intruder problem, an internal source represents an impossible entry point for an intruder. If we push the sheaf $G$ along the projection map $\pi$ we then get the following assignments of data:
\[
	\pi_*G(e)\cong<e_1,e_2,e_3> \quad \pi_*G\cong G(v_1) \quad \pi_*G(f)=<f_1,f_2,f_3>
\]
\end{ex}

\begin{figure}
\begin{center}
\includegraphics[width=.7\textwidth]{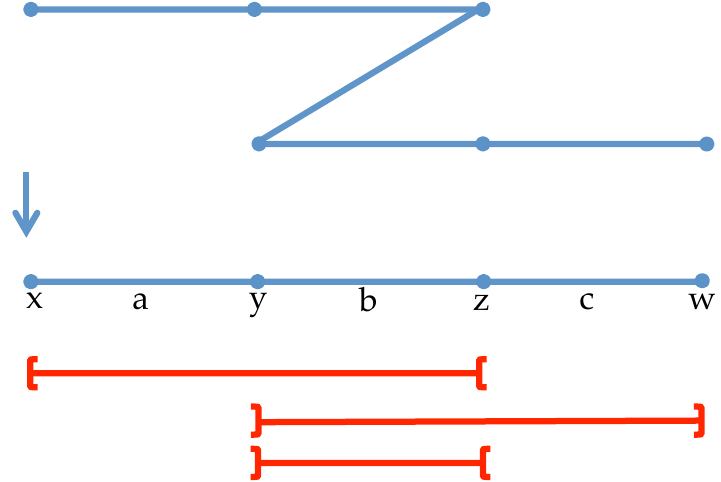}
\caption{Linearized Sheaf of Sections}
\label{fig:evade_sec_bc}
\end{center}
\end{figure}

\begin{ex}
Let us consider the example drawn in Figure \ref{fig:evade_sec_bc}, but now with the sheaf just defined. We set $F=\pi_*G$, whose values are below:

\[
	\xymatrix{F(x) \ar[r]^{\rho_{a,x}} & F(a) & \ar[l]_{\rho_{a,y}}  F(y) \ar[r]^{\rho_{b,y}} & F(b) & \ar[l]_{\rho_{b,z}} F(z) \ar[r]^{\rho_{c,z}} & F(c) & \ar[l]_{\rho_{c,w}} F(w) \\
	k_x  \ar[r] & k_a & \ar[l] k_y \ar[r] &  k^3_b & \ar[l] k_z \ar[r] &  k_c & \ar[l]  k_w}
\]
The two restriction maps of any interest include into the top section and the bottom section, respectively.
\[
	\rho_{b,y}=\begin{bmatrix}1 \\0\\0\end{bmatrix} \qquad \rho_{b,z}=\begin{bmatrix}0 \\0\\1\end{bmatrix}
\]
Without a change of basis one can see that this sheaf splits as the direct sum of indecomposables, whose barcodes are drawn in Figure \ref{fig:evade_sec_bc}.
\end{ex}

The previous example offers a glimmer of hope. No intruder can evade detection and the absence of a long barcode reflects that. Moreover, the sheaf cohomology computation shows $H^0([0,1];F)\cong 0$, which would be a promising shortcut to computing barcodes. Alas, the linearized sheaf of sections fairs no better than the cosheaf of components. Here we provide a counterexample, joint with Dave Lipsky, to either of the hopes that non-zero $H^0([0,1];F)$ or a long barcode provides an if and only if criterion for the existence of an evasion path.

\begin{figure}
\begin{center}
\includegraphics[width=.7\textwidth]{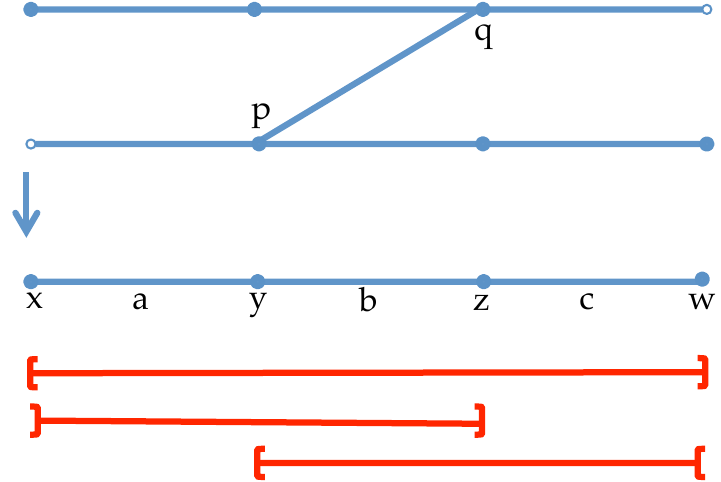}
\caption{Counterexample for the Linearized Sheaf of Sections}
\label{fig:mobile_d4}
\end{center}
\end{figure}

\begin{prop}\label{prop:mobile_d4}
	Although it is true that the existence of an evasion path implies the existence of a long barcode (and thus $H^0([0,1];F)\neq 0$) it is not true that having a long barcode (or $H^0([0,1];F)\neq 0$, which is a weaker condition) implies the existence of an evasion path.
\end{prop}
\begin{proof}
	In Figure \ref{fig:mobile_d4} we have drawn the counter-example, which we now explain. The component coming into $p$ appears immediately after time $0$, so it is impossible for an intruder to enter there. Similarly, there is a component leaving from $q$ that closes up right before time $1$. The pushforward sheaf then takes the following form
	\[
		\xymatrix{	k_x  \ar[r] & k^2_a & \ar[l] k^3_y \ar[r] &  k^3_b & \ar[l] k^3_z \ar[r] &  k^2_c & \ar[l]  k_w}
	\]
	The maps from $F(y)$ and $F(z)$ to $F(b)$ are the identity maps. The two maps that require some inspection are built out of a projection and a trace map.
	\[
		\rho_{a,y}=\begin{bmatrix}1 & 0 & 0 \\0 & 1 & 1\end{bmatrix} \qquad \rho_{c,z}=\begin{bmatrix}1 & 1 & 0 \\0 & 0 & 1\end{bmatrix}.
	\]
	The change of basis required to obtain the desired Remak decomposition indicated by the barcodes is not so easily seen. The interval decomposition algorithm outlined in~\cite{zigzag} provides a sure-fire method of obtaining it. It is left to the reader to verify the the barcodes in Figure \ref{fig:mobile_d4}.
	
	Instead, we give a sheaf-theoretic justification for the existence of a long barcode. There is a unique non-zero global section of $G$ and it is supported everywhere except on the prong incoming to $p$ and outgoing from $q$. Explicitly, it comes from choosing a compatible kernel for the restriction matrices $\rho_{a,y}$ and $\rho_{c,z}$. At $q$ the routing through the top path is annihilated by the ``negative'' of the routing through the bottom path; it is as if two intruders are traveling with opposite charges. As a consequence, its support surjects onto all of $[0,1]$. Since $H^0([0,1];F)\cong k$ we can infer the existence of one closed barcode, and because this section has global support, the barcode must be long.
\end{proof}

\begin{rmk}[Dynkin Diagrams and Stalks]
	Recall that $F:=\pi_* G$. Consider the sheaf $G$ implied by Figure \ref{fig:mobile_d4}. When restricted to the open stars at $p$ and $q$ separately $G$ is equivalent to one of the 12 indecomposable representations of the Dynkin diagram $D_4$; see~\cite{etingof_rep}, p. 83. Since the open stars intersect, one can show that the entire sheaf $G$ on $E$ is indecomposable. This cannot be used directly to show that a long barcode must exist. The pushforward of an indecomposable representation is not necessarily indecomposable. However, the argument using stalks indicates that some sections (subrepresentations), must have global support.
\end{rmk}

\section{Multi-Modal Sensing}
\label{subsec:multi_modal}\index{sensors!multi-modal}

In this section we will explore the following cartoon for multi-modal sensing:

\begin{itemize}
	\item We have a region $W$ thought of as a topological space that is tame enough to be triangulated. This space is populated by agents of interest and sensors.
	\item There is a vector space of properties $k^n$, usually $\RR^n$ or $\CC^n$, and every intruder is tagged with an unchanging \textbf{property vector} $v\in k^n$. These property vectors might record colors (which we pretend has a linear structure), sounds, thermal signatures or, in the context of wireless network data, a unique wireless SSID (we imagine scaling corresponds to the strength of the signal). In future applications, $k$ may be a ring that stores data, just as $\ZZ$ is used to record counts and $\ZZ^n$ records counts of different types of targets.
	\item There are sensors who monitor subspaces of $k^n$ and subspaces of $X$. ``Monitors'' means explicitly that a sensor $i$ with support $V_i\subset X$ has attached to it a subspace of the vector space dual to property space, i.e. $S_i\subset k^{n*}$. For simplicity, we assume that $S_i=\mathrm{span}\{\xi_i\}=:<\xi_i>$. The act of sensing corresponds to taking a property vector $v\in k^n$ and returning a number $\xi_i(v)$ that records the strength of the detection. Outside of the sensor's support $V_i$, the sensor must return zero on every vector. In the overlap of two sensors' supports, the vector space that is sensed is the internal direct sum.
\end{itemize}

This cartoon specifically suggests the use of constructible sheaves and cosheaves as a model. Because the roles of sensors and intruders are formally dual, we will have to use \emph{both} sheaves and cosheaves. Understanding the formal properties of sensing and evasion will lead us naturally to some long-exact sequences in cohomology, which will necessitate the introduction of barcodes to understand these results. 

\begin{figure}
\centering
\includegraphics[width=\textwidth]{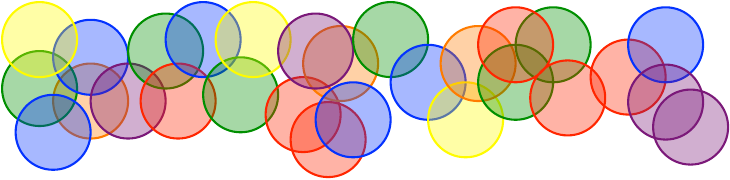}
\caption{Multi-Modal Sensors Distributed in a Plane}
\label{fig:sensors5in_color}
\end{figure}

We are going to work with a simplified version of the above cartoon. To detach ourselves from an embedding of the sensors into $W$, we will use the \v{C}ech nerve associated to the sensors supports. This will provide us with a simplicial complex $X$ and this where we will define sheaves and cosheaves. Since we can only analyze the intruder problem inside sensor's support, we call this a \textbf{relative intruder problem}. Working strictly inside the coverage region will introduce counter-intuitive results, such as Claim \ref{clm:no_sections}. Nevertheless, this setup is a prototype for future applications of sheaves and cosheaves to multi-modal sensing.

\subsection{A Deeper Look at Sensing}
\label{subsubsec:deep_sense}

\begin{figure}
\centering
\includegraphics[width=.5\textwidth]{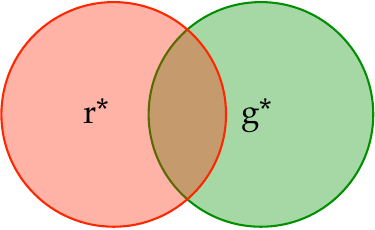}
\caption{Two Multi-Modal Sensors}
\label{fig:rg_sense}
\end{figure}

Let us investigate a little more deeply the picture of multi-modal sensing presented to us in the above cartoon. In Figure \ref{fig:rg_sense}, we consider a situation where we have a sensor capable of detecting ``red'' properties and a sensor capable of detecting ``green'' properties.\footnote{We use scare quotes to indicate that the terms can be substituted for whatever application is of interest.}

On the nerve of the sensor cover, the organizing diagram of vector spaces is clear.
\[
	<r^*>\hookrightarrow <r^*,g^*> \hookleftarrow <g^*> \qquad k \hookrightarrow k^2 \hookleftarrow k
\]
The direction of the arrows indicates that a cellular sheaf is best used to collate sensing abilities. However, the diagram of abstract vector spaces on the right has no way of telling whether an individual copy of $k$ should correspond to $<r^*>$ or $<g^*>$ or $<b^*>$. Such a distinction requires that we embed our sensing sheaf into a global system of coordinates $(k^n)^*_X$. This motivates the following definition.

\begin{defn}[Sensing Sheaf]\index{sensing sheaf}\index{sheaf!sensing}
	Suppose we have a multi-modal sensor network distributed in a space $W$. Form the nerve given by the intersections of the sensors' supports and call this simplicial complex $X$. We define a \textbf{sensing sheaf} $F$ by assigning to each vertex $v$ in $X$ the subspace $S_v\subset (k^n)^{*}$. Over higher simplices $\sigma$ we assign the following vector spaces and use the natural inclusions for the maps internal to the sheaf:
	\[
		F(\sigma)=S_{v_0}+\cdots + S_{v_n} \qquad F(\sigma)\hookrightarrow F(\tau) \quad \sigma\leq \tau.
	\]
	Here we have used the internal sum of subspaces to reflect the fact there may be dependencies. The internal sum is only defined in the presence of an ambient space, thus part of the data of a sensing sheaf is an embedding into the constant sheaf of all sensing abilities:
	\[
		\iota_F: F \hookrightarrow k^{n*}_X.
	\]
\end{defn}

Now suppose we have an intruder, which we imagine as a point in the union of the red and green sensors in Figure \ref{fig:rg_sense}. The intruder has a property vector $v\in k^n$ that lists its various attributes, its colors in this example. What number does the sensor return while the intruder is in the red sensor's domain? By design, it is $r^*(v)$, the contraction of the red co-vector and the property vector $v$. If $k^n=k^3=<v_r,v_g,v_b>$ is a three dimensional property space spanned by the attributes ``red,'' ``green,'' and ``blue,'' equipping it with the standard Euclidean inner product allows us represent this measurement by the matrix product
\[
	r^*(v)=\begin{bmatrix} 1 & 0 & 0\end{bmatrix}\begin{bmatrix} v_r \\ v_g \\ v_b\end{bmatrix}=v_r
\]

However, if the sensors can collaborate and share information, then we can store together the observations when the intruder is in the intersection of the red and green sensors' support.
\[
	\begin{bmatrix} r^*(v) \\ g^*(v)\end{bmatrix}=\begin{bmatrix} 1 & 0 & 0 \\ 0 & 1 & 0\end{bmatrix}\begin{bmatrix} v_r \\ v_g \\ v_b\end{bmatrix}=\begin{bmatrix} v_r \\ v_g\end{bmatrix}
\]
We can package these measurements into a cellular cosheaf, where two observations are the same modulo the properties unobserved by the sensors.
\[
	<r> \leftarrow <r,g> \rightarrow <g> \qquad \frac{k^3}{<g,b>} \leftarrow \frac{k^3}{<b>} \rightarrow \frac{k^3}{<r,b>}
\]

One should note that the right side gives an equivalent formulation for the \textbf{measurement cosheaf} of the Figure \ref{fig:rg_sense}. We have over each cell passed to the quotient space where the properties that are invisible to each of the sensors is treated as zero. In other words, the vectors produced by the process of measurement must naturally be considered modulo the unknown.

\begin{defn}[Evasion Co-Sheaf]\index{evasion cosheaf}\index{cosheaf!evasion}
	Given a sheaf $F$ whose restriction maps are inclusions, along with a fixed embedding into a locally constant sheaf of vector spaces $G$ (we take $G=k^{n*}_X$), we define the \textbf{annihilator cosheaf} $\widehat{Ann}(F)$ as follows:
	\begin{itemize}
		\item $\widehat{Ann}(F)(\sigma)=\{v^*\in G(\sigma)^{*}| v^*(\iota(w))=0\quad \forall w\in F(\sigma)\}$
		\item If $\sigma\subset\bar{\tau}$, then $r_{\sigma,\tau}:\widehat{Ann}(F)(\tau)\to \widehat{Ann}(F)(\sigma)$ is the inclusion.
	\end{itemize}
	When using the language of sensing sheaves, we will call $\widehat{Ann}(F)=:\hat{E}$ the \textbf{evasion cosheaf}.
\end{defn}

\begin{lem}\label{lem:evade_cok}
	Let $F$ be a sensing sheaf on $X$, then the evasion cosheaf is canonically identified as the linear dual of the cokernel of the embedding, that is to say that $\hat{E}\cong \hat{V}(\cok(\iota))$ in the diagram below.
	\[
		\xymatrix{0 \ar[r] & F \ar[r]^{\iota} \ar@{~>}[d] & G \ar[r]^{q} \ar@{~>}[d] & \cok(\iota) \ar[r] \ar@{~>}[d] & 0 \\
		0 & \ar[l] \hat{V}(F) & \ar[l] \hat{V}(G) & \ar[l] \hat{V}(\cok) & \ar[l] 0 }
	\]
\end{lem}
\begin{proof}
	Here we make use of the fact that for cellular sheaves, the cell-by-cell cokernel of the maps $\iota(\sigma):F(\sigma)\to G(\sigma)$ defines a sheaf. This is not always true for general sheaves. Reducing the argument to a cell-by-cell one, we have a short exact sequence of vector spaces
\[
	\xymatrix{0 \ar[r] & V \ar[r]^{\iota} & W \ar[r]^{q} & W/V \ar[r] & 0}
\]
where we can identify 
\[
Ann_W(V)=\{\varphi:W\to k\,|\,\varphi(v)=0\,\forall v\in V\}\cong (W/V)^{*}
\]
and of course all the restriction maps get sent to restriction maps
\[
\xymatrix{V_2 \ar[rd] & & (W/V_2)^{*} \ar[dd]\\
& W \ar[ru]\ar[rd] & \\
V_1 \ar[uu]\ar[ru] & & (W/V_1)^{*} .}
\]
\end{proof}

This identification of evasion cosheaves with the linear dual of a cokernel means that we can leverage a classical technique in studying the relative intruder problem. After all, to every short exact sequence of sheaves we get an induced long exact sequence of sheaf cohomology. In the context of multi-modal sensing this relates in a precise way the topology of the total covered region and the cohomology of the sensing and evasion sheaves. 

\begin{thm}[Sensing-Evasion Decomposition]\label{thm:les_sensing}\index{decomposition!sensing and evasion regions}
	Given a sensing sheaf of vector spaces $\iota:F\to G=\tilde{k}_X^{n^{*}}$ we obtain a long exact sequence of sheaf cohomology groups
\[
\xymatrix{
	            0\ar@{->}[r] & H^0(X;F)\ar@{->}[r] & H^0(X;k)^{\oplus n}\ar@{->}[r]
	                   & H^0(X;\cok(\iota)) \ar `r[d] `_l[lll] ^{\delta^0} `^d[dlll] `^r[dll] [dll] \\
	            & H^1(X;F)\ar@{->}[r] & \hspace{12pt}\cdots\hspace{12pt}\ar@{->}[r]
	                   & H^k(X;\cok(\iota)) \ar `r[d] `_l[lll] ^{\delta^k} `^d[dlll] `^r[dll] [dll] \\
	            &H^{k+1}(X;F)\ar@{->}[r] & H^{k+1}(X;k)^{\oplus n}\ar@{->}[r]
	                   & H^{n+1}(X;\cok(\iota)) \ar `r[d] `_l[lll] ^{\delta^{k+1}} `^d[dlll] `^r[dll] [dll] \\
	            & & \cdots &
	}
\]
Where $H^k(X;\cok(\iota))$ gets identified with the evasion co-sheaf's homology $H_k(X;\hat{E})$ via the linear duality functor, i.e. $V: \cok(\iota)) \rightsquigarrow E$.
\end{thm}
\begin{proof}
	The proof is immediate from standard homological algebra techniques.
\end{proof}

\subsection{Indecomposables, Evasion Sets, Generalized Barcodes}
\label{subsubsec:evasion_bcs}
\begin{figure}
\centering
\includegraphics[width=.3\textwidth]{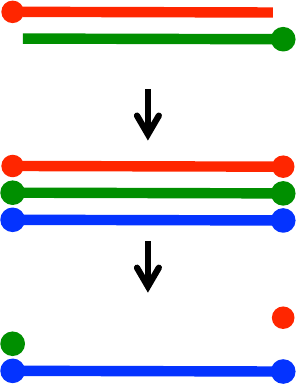}
\caption{Examining the Short Exact Sequence}
\label{fig:sense_bars}
\end{figure}

One of the drawbacks of Theorem \ref{thm:les_sensing} is that we have no good interpretation of what the sheaf cohomology groups mean. Let's consider again Figure \ref{fig:rg_sense}, but this time let us focus only on the \v{C}ech complex and each of the three sheaves that appear in the short exact sequence. This is depicted in Figure \ref{fig:sense_bars}.

As can be clearly seen each sheaf appearing in the sequence is already written as a direct sum of indecomposables, which because the nerve is a one-simplex, look like barcodes. By using the observation $H^i([0,1];F)\cong \oplus H^i_c(B_i)$, which we have already made heavy use of, we can determine all the sheaf cohomology of interest for this example.

\begin{ex}[Red-Green Sensors]
	By inspection of the indecomposable presentations of the three sheaves $F$, $k^3_X$ and $\cok(\iota)$ in Figure \ref{fig:sense_bars} we see that
	\[
	H^i(X;F)\cong 0 \, i=0,1; \qquad H^0(X;k^3_X)\cong k^3 \qquad H^0(X;\cok(\iota))\cong H_0(X;\hat{E})\cong k^3 
 	\]
The interpretation of each of the three generators in the evasion cosheaf homology is that there is a connected component where red, green and blue can separately evade.
\end{ex}

\begin{defn}[Evasion and Detection Sets]
	Let $v\in k^n$ be a property vector and $F$ a sensing sheaf on $X$. Define the \textbf{evasion set} $E_v$ to be the set of points in $X$ where an intruder with property vector $v$ can go without being detected. Dually, call the set of points where $v$ can be detected the \textbf{detection set} $D_v$.
\end{defn}

Since we are working with cellular sheaves where individual sensors have support equal to the open star of their designated vertex in the simplicial complex $X$, thus the detection set $D_v$ is equal to the union of all the stars of the sensors that can see $v$, hence $D_v$ is an open union of cells. This proves the following lemma.

\begin{lem}
	For any property vector $v$, the evasion and detection sets form an open-closed decomposition of $X$, that is
	\[
		X=E_v\cup D_v \qquad E_v\cap D_v=\emptyset, \qquad E_v \,\, \mathrm{open}.
	\]
	When $X$ is compact this means that $E_v$ is compact as well.
\end{lem}

We record another easy lemma, connected to our desire to get an indecomposable presentation for our evasion cosheaves.

\begin{lem}
	Suppose that all sensors must pull their sensing capabilities from a fixed orthonormal basis of $k^{n*}$, say $v_1^*,\ldots, v_n^*$, then the evasion cosheaf splits as a direct sum decomposition of constant cosheaves supported on the evasion sets for $v_1,\ldots,v_n$
	\[
		\hat{E}\cong \hat{k}_{E_{v_1}}\oplus \cdots \oplus \hat{k}_{E_{v_n}}
	\]
	with the further observation that each $\hat{k}_{E_{v_i}}$ has a Remak decomposition as a sum of constant cosheaves supported on the components of $E_{v_i}$.
\end{lem}
\begin{proof}
	The fact that the sensor capabilities can only be chosen from a fixed orthonormal basis, implies that we can write the constant sheaf $k^{n*}_X$ as a direct sum of $k_X\oplus\cdots\oplus k_X$ where we think of each copy of $k_X$ as being the constant sheaf generated by $<v_i^*>$. As a consequence we get the following diagram
	\[
		\xymatrix{ & k_{D_{v_1}} \ar@{^{(}->}[r] & k_X \ar[dr]^{i_{v^*_1}} & \\
		F \ar[ur]^{\pi_{v^*_1}} \ar[dr]_{\pi_{v^*_n}} & \vdots & \vdots & k^{n*}_X \\
		  & k_{D_{v_n}} \ar@{^{(}->}[r] & k_X \ar[ur]_{i_{v^*_n}} & }
	\]
	Now we can use for each factor the following standard short exact sequence of sheaves
	\[
		\xymatrix{0 \ar[r] & k_{D_{v_i}} \ar[r] & k_X \ar[r] & k_{E_{v_i}} \ar[r] & 0}
	\]
	and thus the cokernel splits as a direct sum $\oplus_{i=1}^n k_{E_{v_i}}$.
\end{proof}	

The above lemma implies that in certain cases we can interpret the homology of the evasion cosheaf in terms of the topology of the evasion sets.

\begin{figure}
\centering
\includegraphics[width=.7\textwidth]{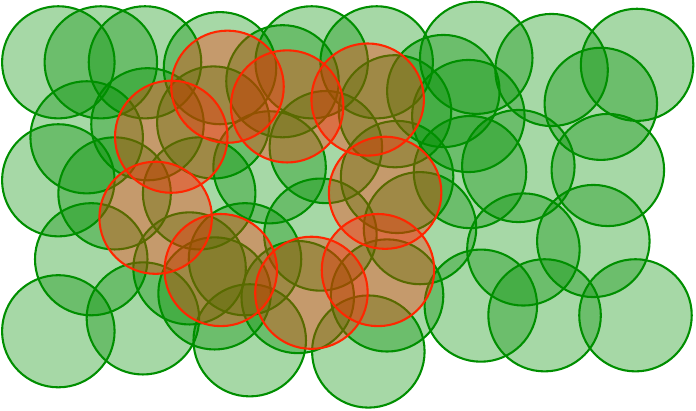}
\caption{Short Exact Sequence}
\label{fig:red_circle}
\end{figure}

We have one more observation we'd like to leverage.

\begin{clm}\label{clm:no_sections}
	If sensor's abilities are pulled from a fixed orthonormal basis $v_1^*,\ldots, v_n^*$ and moreover the detection sets are not pairwise disjoint, then the sensing sheaf has no global sections.
\end{clm}
\begin{proof}
	This follows from the fact that if an edge is common to two different detection sets, then there can be no global sections since the following sheaf has no non-zero global sections
	\[
		<v_i^*> \hookrightarrow <v_i^*,v^*_j> \hookleftarrow <v_j^*>.
	\]
\end{proof}

For the example considered in Figure \ref{fig:red_circle} assume that the space of properties is two dimensional, spanned by red and green. Then the Theorem \ref{thm:les_sensing} provides the following forcing result
\[
 0 \to H^0_c(X;F)\cong 0 \to H^0_c(X;k^{2*}_X)\cong k^2 \to H^0_c(X;\cok) \to H^1_c(X;F)\cong k \to 0 
\]
which upon careful inspection reveals that the red evasion set must be disconnected.

\ctparttext{This part represents the mathematical heart of the thesis, although many of its results were motivated by the applications considered in Part~\ref{part:applications}.

Chapter~\ref{subsec:strat} is by far the most technically demanding part of the thesis. It takes up and proves an equivalence between constructible cosheaves and representations of MacPherson's entrance path category, which hinges on a proof of the Van Kampen theorem for this category. The full machinery of stratification theory is then used to construct representations of the (definable) entrance path category from a stratified (definable) map. This part also rests on proving a codimension-criterion under which Thom's condition $a_f$ always holds.

Chapter~\ref{sec:duality} proves that Verdier duality is rightly conceived as an exchange of sheaves and cosheaves. An explicit formula for the derived equivalence of cellular sheaves and cosheaves is presented.

Chapter~\ref{sec:valuations} uses the formula of Chapter~\ref{sec:duality} to prove that compactly supported cellular sheaf cohomology can be viewed as taking a (derived) coend with the image of the constant sheaf under this formula.

Chapter~\ref{sec:graded} proves that the derived category of cellular sheaves over a one-dimensional base space is equivalent to a graded category. This formalizes the intuition of why spectral sequences over graphs always collapse on the $E_2$ page.

Chapter~\ref{sec:metric} introduces the interleaving distance for sheaves defined on a metric space. Although officially an extended pseudo-metric on the category of pre-sheaves, we prove it is an extended metric on the category of sheaves. One of the most fundamental properties of this extended metric is that global sections places sheaves into distinct connected components. To illustrate the theory more concretely, we take up an explicit description of the space of constructible sheaves over the real line.
}
\part{Novel Mathematical Contributions}
\label{part:math_contributions}

%
%

\chapter{The Definable Entrance Path Category}
\label{subsec:strat}

\begin{quote}
{\em ``Facilis descensus Averno;\\
noctes atque dies patet atri ianua Ditis;\\
sed revocare gradum superasque evadere ad auras,
hoc opus, hic labor est.''}
\begin{flushright} --- Virgil's Aeneid, Book 6, Lines 124-9 \end{flushright}
\end{quote}

Fundamentally, one-dimensional persistent homology tries to understand topological changes in a one-parameter family of spaces. Multi-dimensional persistence tries to understand topological changes in a multi-parameter family of spaces; the leap in complexity from one dimension to two can not be overstated. The model problem of interest is to describe how the homology of the fiber of a map $f:Y\to X$ changes as one queries points or subsets in $X$. For general maps, this problem is entirely too unwieldy. 

In this chapter we focus on a broad class of maps where this problem has an interesting answer: definably stratified maps. Informally, stratified maps are glued together fiber bundles. Definable maps are ones that can be defined with finitely many logical operations. Every definable map is stratified so we study simply definable maps.

The upshot of this chapter is that the homology of the fibers of a definable map give rise to a representation of a particular quiver with relations --- a category in other words --- called the definable entrance path category, whose general version was introduced by MacPherson to study general stratified maps. If one considers the opposite of the entrance path category, i.e. the exit path category, one obtains a constructible sheaf, which we now define.

\begin{defn}[Constructible Sheaves and Cosheaves]\index{sheaf!constructible}\index{sheaf!restricted sheaf}\index{cosheaf!constructible}
	Let $F$ be a sheaf valued on a topological space $X$. One says that $F$ is \textbf{constructible} if there exists a filtration by closed subsets $$\emptyset=X_{-1}\subset X_0 \subset X_1 \subset \cdots \subset X_n=X$$ 
	such that on the each connected component of the space $X^k=X_k- X_{k-1}$, the restricted sheaf $F|_{X^k}$ (the pullback of $F$ along the inclusion $X^k\hookrightarrow X$) is locally constant. Alternatively, instead of asking for a filtration one can ask for a decomposition of $X$ into disjoint pieces $X_{\sigma}$ over which the restricted sheaf is locally constant.

	Dually, we will call a cosheaf, with costalks valued in $\vect$, constructible if its linear dual is a constructible sheaf.
\end{defn}

One of the purposes of this chapter is to impose further conditions on the nature of the filtration so that we get nice properties. As stated, there is nothing to prevent us from using a one step filtration of the Cantor set. Expressing precisely these extra conditions will require the introduction of stratification theory.

\section{Stratification Theory and Tame Topology}
\label{subsec:strat_thy_tame_topology}

As wonderful as fiber bundles and local systems may be, they still fail to capture the sort of structure we are interested in because the topology of the fiber can never change. In order to bring cosheaves into contact into a larger realm of mathematics, we will need to consider stratified maps. To whet the appetite, stratified maps will allow us to describe in one language:
\begin{description}
	\item[Morse theory] --- Morse functions are just particular instances of stratified maps $f:M\to \RR$.
	\item[Picard-Lefschetz theory] --- the complex analog of Morse theory studies algebraic maps $\pi:X\to\CC$, which are necessarily stratified.
	\item[Point Cloud Data and Persistence] --- Semialgebraic families are described by semialgebraic maps, which are stratified.
\end{description}
In other words, stratification theory gives a system of geometry for exploring a wealth of examples, appearing in mathematics and nature. Stratification theory does this by breaking up a space or map into regions, over which the usual analysis of manifolds and fiber bundles apply.

\begin{defn}[Decomposition]\index{decomposition!into strata}\index{strata}
	A \textbf{decomposition} of a space $X$ is a locally finite partition of $X$ into locally closed subsets (sets of the form $U\cap Z$ for $U$ open and $Z$ closed) $\{X_{\sigma}\}_{\sigma\in P_X}$ called \textbf{pieces}, which satisfy the axiom of the frontier. Consequently, $P_X$ is a poset. When the pieces have the additional structure of being manifolds, we call them \textbf{strata}.
\end{defn}

\begin{rmk}
	A \textbf{stratum} is sometimes used to mean either a union of strata of a fixed dimension or a single connected component in a decomposition. We usually prefer the latter meaning.
\end{rmk}

We have already encountered an example of a decomposition of a space $X$, namely a cell complex. Here each piece is homeomorphic to $\RR^k$ for some $k$, which can vary from stratum to stratum. A graph is naturally decomposed into its vertices and open edges. For a decomposition that is not a cell complex, consider the complex numbers $\CC$ partitioned into the sets $\{0\}$ and $\CC-\,\{0\}$.

\begin{defn}
	Suppose $(X,P_X)$ and $(Y,P_Y)$ are decomposed spaces, then a \textbf{decomposition-preserving} map is a continuous map $f:X\to Y$ that sends pieces to pieces, i.e. we have a commutative square
	\[
	\xymatrix{X \ar[r]^f \ar[d] & Y \ar[d] \\ P_X \ar[r]^{P_f} & P_Y}
	\]
	In the case where the pieces are strata we call such a map a \textbf{stratum-preserving} map.
\end{defn}

Much like how the notion of a category emerged through the study of functors, in some sense the necessity for decompositions more general then simplicial or cell complexes came about because not all maps preserved the pieces of those decompositions. We give an example of such a map.

\begin{figure}
\centering
\includegraphics[width=.7\textwidth]{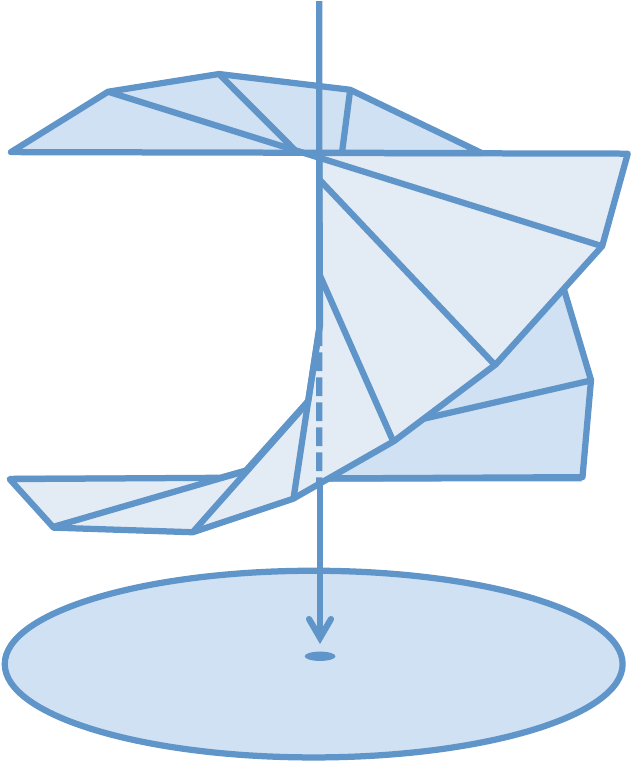}
\caption{Blowing up at a Point}
\label{fig:blowup}
\end{figure}

\begin{ex}[Blow-Ups]\label{ex:blowup}\index{blow-up}
	Consider the map
	\[
		f:\RR^2\to\RR^2 \qquad f(x,y)=(x,xy).
	\]
	This map is not triangulable, see~\cite{shiota} page 305. This map is related to the operation in algebraic geometry known as ``blowing up at a point.'' The blow-up map is an endless source of interesting geometry and counter-examples, so it is worth describing. Recall that the space of lines in $\RR^2$, written $\RR\mathbb{P}^1$ is defined to be the quotient of $\RR^2-\,\{(0,0)\}$ by the relation that $(x,y)\sim (\lambda x,\lambda y)$ for any $\lambda\neq 0$. Topologically, this quotient is the circle $S^1$. Tracing the image of the top arc of a circle from $0$ to $\pi$ through the quotient map one gets the complete circle in $\RR\mathbb{P}^1$.
	
	The blow-up $B$ of $\RR^2$ at the origin is defined to be the closure of the image of the map
	\[
		\RR^2-\,\{(0,0)\} \hookrightarrow \RR^2\times \RR\mathbb{P}^1
	\]
	where the map to the first coordinate is the inclusion and the map to the second coordinate is the quotient map. The blow-up map $\pi:B\to\RR^2$ is the projection back from the closure of the graph of this map to the closure of the domain, i.e. $\RR^2$. Thus the fiber over $(0,0)$ is a circle, but the fiber over any other point is a single point. One can visualize this by restricting the map to a closed disk centered at the origin. The image is contained in a solid torus and the closure of the image will assign the core circle to the origin. The image of $\mathbb{D}^2-\,\{(0,0)\}$ is commonly visualized as a spiral staircase as in Figure~\ref{fig:blowup} whose boundary traces out a torus knot. See \cite{arone-blowup} for a treatment of different constructions of real blow-ups and their functorial properties.
\end{ex}

\begin{figure}
\centering
\includegraphics[width=.5\textwidth]{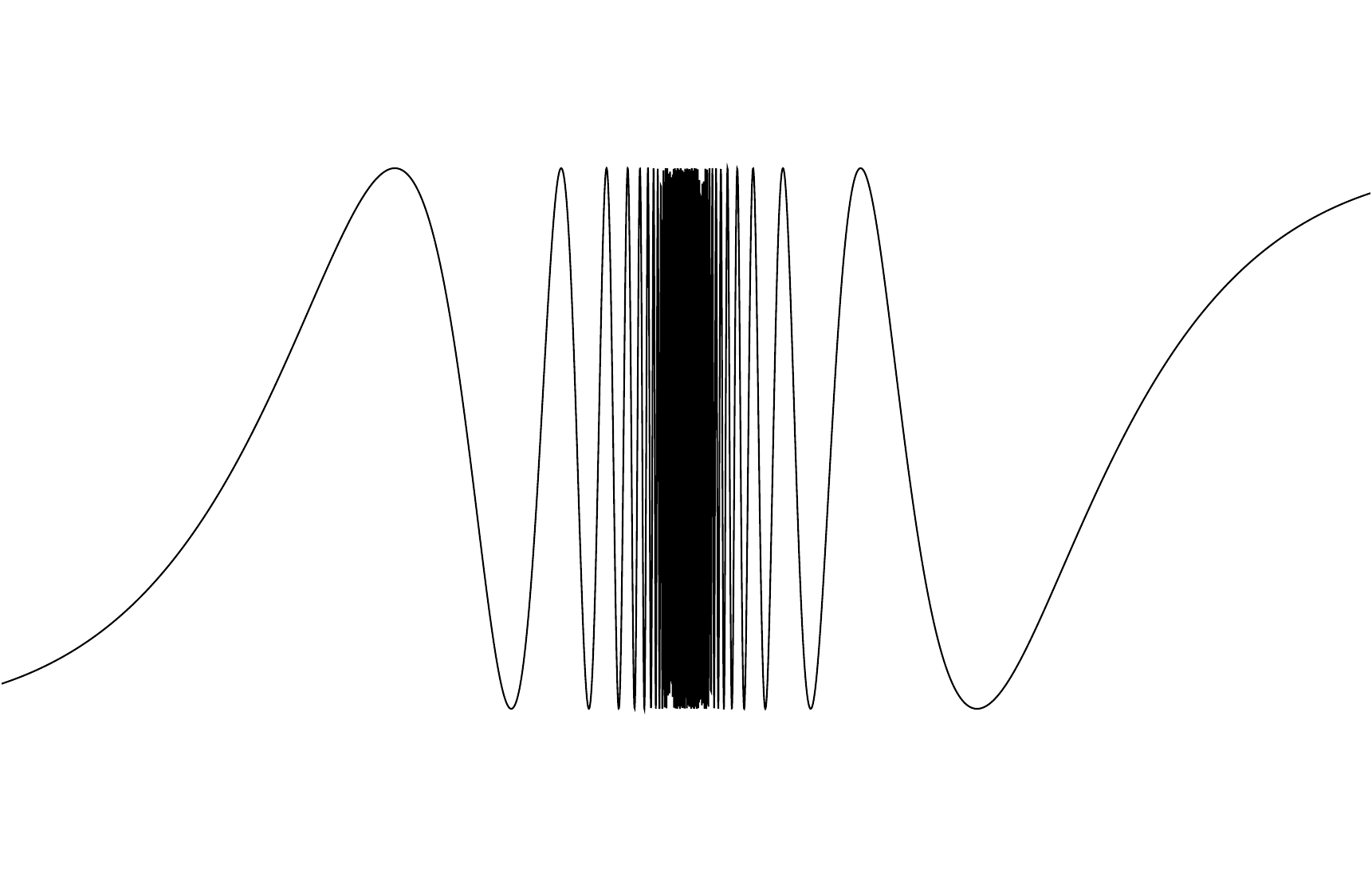}
\caption{Topologist's Sine Curve}
\label{fig:top_sine}
\end{figure}

Decomposing spaces and maps gives some control over how these things are built up out of pieces, but it is not quite strong enough to tame the geometry of interest. In particular, the topologist's sine curve drawn in Figure \ref{fig:top_sine} can be decomposed into the two pieces
\[
	X_{\tau}:=\{(x,\sin(1/x))\, |\, x\neq 0 \}\cup \{(0,y)\,|\, y\in[-1,1]\}=:X_{\sigma}
\]
that satisfy the axiom of the frontier $X_{\sigma}\subset \bar{X}_{\tau}$, but it does not have the intuitively desired property that~\cite[p.131]{lu-ctst}.
\[
	\dim X_{\sigma} < \dim X_{\tau}
\]

Further regularity conditions must be imposed to capture this property and other desired features that hold for piecewise-linear, algebraic, semi-algebraic, sub-analytic and other geometries. Systematic overviews of these different regularity conditions are overwhelming and highly technical. For a taste, one should consult J\"org Sch\"urmann's remarkable service in writing down 14 different regularity conditions and their corresponding implications in~\cite[Rmk. 4.1.9]{schurmann}. To keep the exposition light we focus on a geometric condition and its topological generalization as they have historically had a strong influence on stratification theory.

\subsection{Whitney Stratified Spaces}
\label{subsubsec:whitney}

In this section we relay two ways of fusing manifold pieces into non-manifold wholes. The champions of this section are Hassler Whitney and Ren\'e Thom.\footnote{For more historical context of these approaches, written by experts, we recommend the recent article~\cite{goresky-bull} and Part I, Section 1 of~\cite{GM}.} In 1965, Whitney, whose approach relies on the geometry of tangent planes and secant lines, defined two properties that a stratified space should possess~\cite{whitney-local, whitney1965annals}. Thom, who proposed in a 1962 paper~\cite{thom-poly} a definition of a stratified space using tubular neighborhoods, later extracted the topological consequences of Whitney's definition and outlined a more general definition of a stratified space~\cite{thom-strat}. Thom's definition was first articulated carefully by John Mather in his famous 1970 Harvard ``Notes on Topological Stability''~\cite{mather}, which went unpublished for 42 years and are to this day an excellent resource for learning the theory. 

Proving that any Whitney stratified space admits the structure of a Thom-Mather stratified space requires substantial work. Thus, we present them below as separate definitions, beginning with Whitney's. We outline the properties that make Whitney stratified spaces nice as motivation for Thom's definition. After introducing both of these definitions we will present a proof\footnote{A proof appears in Mark Goresky's thesis~\cite{goresky-thesis} that was never published and which he graciously provided to the author. We have since modified that proof to suit our purposes.} that says that closed unions of strata in a Thom-Mather space have regular neighborhoods, i.e. an open neighborhood and a weak deformation retraction. This result plays a key role in Theorem \ref{thm:strat_maps}.

\begin{figure}
	\centering
	\includegraphics[width=.6\textwidth]{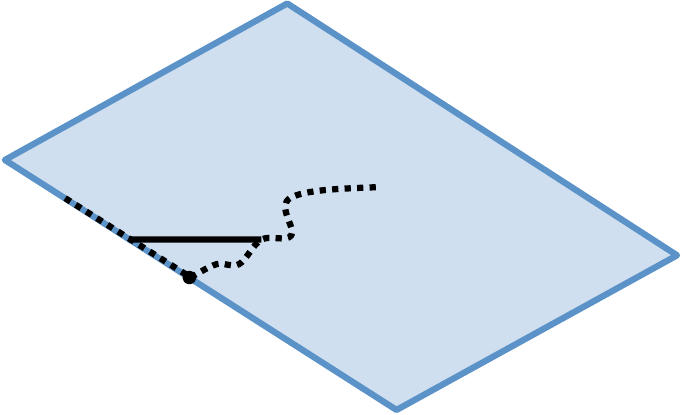}
	\caption{Diagram for Whitney Condition (b)}
	\label{fig:whitney_b}
\end{figure}

\begin{defn}[Whitney Stratified Spaces]\index{stratified space!Whitney}
	A \textbf{Whitney stratified space} is a tuple $(X,M,\{X_{\sigma}\}_{\sigma\in P_X})$ where $X$ is a closed subset of a smooth manifold $M$ along with a decomposition into pieces $\{X_{\sigma}\}_{\sigma\in P_X}$ such that
	\begin{itemize}
		\item each piece $X_{\sigma}$ is a locally closed smooth submanifold of $M$, and
		\item whenever $X_{\sigma}\leq X_{\tau}$ the pair satisfies \textbf{condition (b)}. This condition says if $\{y_i\}$ is a sequence in $X_{\tau}$ and $\{x_i\}$ is a sequence in $X_{\sigma}$ converging to $p\in X_{\sigma}$ and the tangent spaces $T_{y_i}X_{\tau}$ converges to some plane $T$ at $p$, and the secant lines $\ell_i$ connecting $x_i$ and $y_i$ converge to some line $\ell$ at $p$, then $\ell\subseteq T$.
	\end{itemize}
\end{defn}
\begin{rmk}
	We have omitted \textbf{condition (a)} because it is implied by condition (b)~\cite[Prop. 2.4]{mather}. Condition (a) states that if we only consider a sequence $y_i$ in $X_{\tau}$ converging to $p$ such that the tangent planes $T_{y_i}X_{\tau}$ converge to some plane $T$, then the tangent plane to $p$ in $X_{\sigma}$ must be contained inside $T$.
\end{rmk}

The Whitney conditions are important because so many types of spaces admit Whitney stratifications, the most important being semi-algebraic and sub-analytic spaces.  Remarkably, these conditions about limits of tangent spaces and secant lines imply strong structural properties of the space. To give the reader a taste for the properties enjoyed by Whitney stratified spaces, we provide a brief list:\footnote{Here we follow part of MacPherson's summary in the appendix of his 1991 Colloquium notes~\cite{macpherson-ih-notes}.}
\begin{itemize}\index{stratified space!general properties of}
	\item[-] \textbf{Dimension is Well-Behaved:} If $X_{\sigma}\subseteq \fr(X_{\tau}):=\bar{X}_{\tau}-X_{\tau}$, then $\dim X_{\sigma} < \dim X_{\tau}$. See Proposition 2.7 of~\cite{mather} for a proof. This rules out the topologist's sine curve in Figure \ref{fig:top_sine} from being Whitney stratified.
	\item[-] \textbf{Good Group of Self-Homeomorphisms:} If $x$ and $y$ belong to the same connected component of a stratum $X_{\sigma}$, then there is a homeomorphism $h:M\to M$ preserving $X$ and other strata such that $h(x)=y$ (\cite{mather} pp. 480-481).
	\item[-] \textbf{Local Bundle Structure:} Every stratum $X_{\sigma}$ has an open tubular neighborhood $T_{\sigma}$ and a projection map $\pi_{\sigma}:T_{\sigma}\to X_{\sigma}$ making it into a fiber bundle. This bundle is equipped with a ``distance from the stratum'' function $d_{\sigma}:T_{\sigma}\to\RR_{\geq 0}$. If we define $S_{\sigma}(\epsilon)$ to be $d^{-1}_{\sigma}(\epsilon)$, then we can identify the map $\pi_{\sigma}:T_{\sigma}\to X_{\sigma}$ with the mapping cylinder of the restricted map $\pi:S_{\sigma}(\epsilon)\to X_{\sigma}$ (\cite{goresky-fol} p. 194). Moreover, the fiber of the bundle has the stratification of a cone on a link.
	\item[-] \textbf{Triangulability:} Every Whitney stratified space can be triangulated~\cite{goresky-fol}.
\end{itemize} 

The third property is historically the most important. It guarantees that a Whitney stratification ``looks the same'' along all points in a stratum. The tubular neighborhoods exhibit this local triviality. This condition will be taken as primary when considering Thom-Mather stratifications.

\subsection{Stratified Maps and a Counterexample}
\label{subsubsec:strat_maps}

Our main purpose for considering Whitney (and hence Thom-Mather) stratified spaces is to understand stratified maps. Such maps include Morse functions as a special case and are a good model for understanding moduli problems that commonly arise in applications. Over a given stratum, a stratified map looks like a fiber bundle and all fibers are homeomorphic in a stratum-preserving way. However, as we try to compare a fiber over one stratum with a fiber over that stratum's frontier, the blow-up map of Example \ref{ex:blowup} frustrates our intuition. Thus, we introduce a more restrictive class of stratified maps called Thom maps. Finally, we illustrate that such general stratified maps are not necessarily closed under pullback. This motivates the move to tame topology in Section \ref{subsubsec:tame_topology}.

\begin{defn}[Whitney Stratified Map]\index{stratified map!Whitney}
	Suppose $f:M\to N$ is a smooth map between manifolds that contain stratified spaces $(X,\{X_{\sigma}\}_{\sigma\in P_X})$ and $(Y,\{Y_{\sigma}\}_{\sigma\in P_Y})$ such that $f(X)\subset Y$ with $f|_X$ \emph{proper}. We say $f$ is a Whitney \textbf{stratified map} if the pre-image of each stratum $Y_{\sigma}$ is a union of connected components of strata of $X$ and $f$ takes these components submersively onto $Y_{\sigma}$.
\end{defn}
\begin{rmk}
	To say that a map is (Whitney) \textbf{stratifiable} is to say there exists a stratification of $X$ and $Y$ such that the map is stratified. Often we will neglect to include the ambient manifolds and will say ``Let $f:X\to Y$ be a stratified map.'' 
\end{rmk}
\begin{rmk}\index{stratified submersion}
	When $N=Y$ is stratified as a single stratum, we say that $f$ is a \textbf{stratified submersion}, i.e. $f|_X$ is proper and for each stratum $X_{\sigma}$ $f|_{X_{\sigma}}$ is a submersion.
\end{rmk}

Recall that Ehresmann's theorem states that proper submersions are fiber bundles. Thus, over each stratum a stratified map is a fiber bundle. However, Ehresmann's theorem does not say that the local trivializations can be chosen to respect the stratification. This stratified analog of Ehresmann's theorem is expressed in Thom's first isotopy lemma~\cite[p.~41]{GM}.

\begin{lem}[Thom's First Isotopy Lemma]\label{lem:thom_1}\index{Thom's first isotopy lemma}
	Let $f:M\to \RR^n$ be a (proper) stratified submersion for $X\subseteq M$ a Whitney stratified subset. Then there is a stratum-preserving homeomorphism
	\[
	h: X \to \RR^n\times (f^{-1}(0)\cap X)
	\]
	which is smooth on each stratum and commutes with the projection to $\RR^n$. In particular, the fibers of $f|_X$ are homeomorphic by a stratum preserving homeomorphism.
\end{lem}
\begin{rmk}
	Of course, this implies that for a general stratified map, for every stratum of the codomain $Y_{\sigma}$, the fibers of $f|_{f^{-1}(Y_{\sigma})}:f^{-1}(Y_{\sigma}) \to Y_{\sigma}$ are homeomorphic in a stratum-preserving way. This lemma will be used implicitly throughout the section. It expresses the idea that stratified maps are ``glued together fiber bundles.''
\end{rmk}

As one can imagine, there is a second isotopy lemma, which applies to a more restrictive class of stratified maps. We will not state the second isotopy lemma, rather we will use some of the theory leading up to it.

\subsubsection{Counterexamples Creep In}

One would like to say that given a general (not necessarily Thom) stratified map $f:X\to Y$, one could take a path $\gamma:[0,1]\to Y$ so that the pullback $\gamma^* f:f^{-1}(\gamma)\to I$ is stratified and hence, by the above corollary, a Thom mapping. However, as the next example shows, the pullback need not be stratifiable, so the hypothesis for the corollary fails.\footnote{We are indebted to Mark Goresky for suggesting the key ideas of this example.}

\begin{figure}
\centering
\includegraphics[width=.7\textwidth]{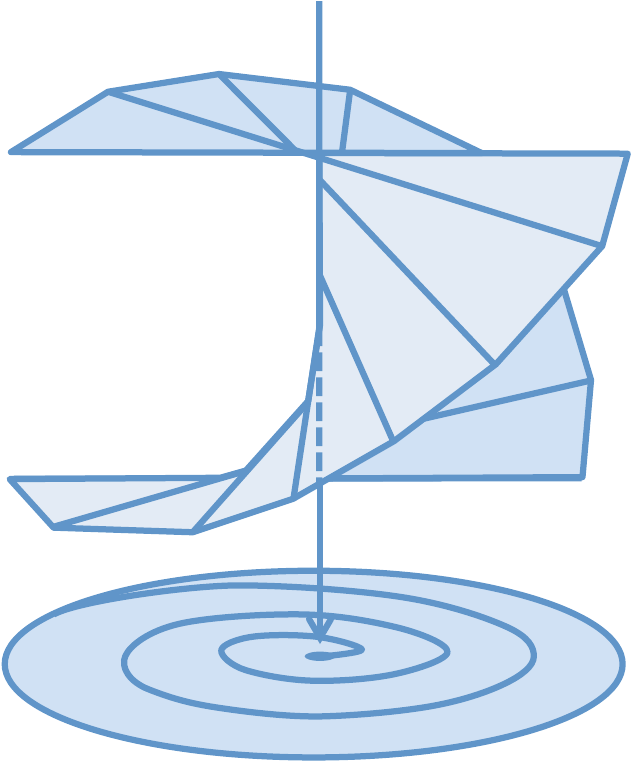}
\caption{Preimage of the Spiral is Not Stratifiable}
\label{fig:blowup_cex}
\end{figure}

\begin{ex}\label{ex:pullback_spiral}
	The blow-up map $\pi:B\to\RR^2$ is a Whitney stratified map that is not a Thom mapping. The closure $S$ of the ``quick spiral''
	\[
		S:=\mathrm{cl}\{(r,\theta)\in\RR^2 \,|\, r=e^{-\theta^2}\}
	\]
	is also Whitney stratified despite wrapping around the origin infinitely many times (\cite{pflaum} Example 1.4.8). However, the inverse image $\pi^{-1}(S)$ cannot be stratified because the inverse image of $(0,0)$ is $S^1$, which is of the same dimension as the inverse image of the spiral, despite the fact that the former is in the frontier of the latter; see Figure~\ref{fig:blowup_cex}. Since being Whitney stratified implies a drop in dimension of the frontier, contraposition shows that the inverse image cannot be Whitney stratified.
\end{ex}

In Theorem \ref{thm:strat_maps} we will give a direct geometric construction of several cosheaves associated to a stratified map. To do so we will need to consider a class of subsets and maps that have all the geometric properties of stratified spaces as well as being preserved under inverse images. This is provided in Section \ref{subsubsec:tame_topology}.

\subsection{O-minimal Structures}
\label{subsubsec:tame_topology}\index{tame topology}

Although stratification theory provides a first pass at taming geometry, it is unsuitable from our perspective because pathologies can still creep in via the inverse image, as Example \ref{ex:pullback_spiral} showed. General stratified spaces and maps are still not tame enough. However, most sets and maps encountered in nature have extra structure. For instance, computer scientists commonly work with piecewise-linear (PL) spaces, which are describable in terms of affine spaces and matrix inequalities. Some algebraic geometers work with semialgebraic spaces, which use zeros and inequalities of polynomials to define their spaces. Analysts tend to use analytic or subanalytic spaces, because the theory is well behaved. Traditionally, one has had to make a choice, once and for all, to speak only of PL geometry, or only of algebraic geometry, or only of analytic geometry. The curse of Babel has confused and separated these domains for a hundred years.

In 1984, Grothendieck declared that an axiomatic ``tame topology'' or ``topologie mod\'er\'ee'' should be developed by extracting out precisely those properties that make these classes of spaces good ones~\cite{grothendieck-sketch}. MacPherson put forth in his lecture notes for the 1991 AMS colloquium lectures a definition of what should constitute a ``good'' class of subsets of a manifold $M$~\cite{macpherson-ih-notes}. Namely, a subset $S$ is good if there is a Whitney stratification of $M$ such that $S$ is a union of strata. These subsets should be closed under the finite set-theoretic operations of unions, intersections and differences. Additionally, the closure of any good subset should be good.

In 1996, Lou van den Dries and his student Chris Miller set forth a most satisfactory definition in their paper ``Geometric Categories and O-minimal Structures''~\cite{vdd-geocat}. Taking requests from sheaf theorists~\cite{schmid-cc} and other working geometers, their paper is a valuable service to the community. It globalized a local solution to Grothendieck's program known as \textbf{o-minimal topology}. The theory of o-minimal topology is grounded in model theory and logic, but it has left almost no trace from those fields. All the logical operations of $\forall,\exists,\lor,\land$ are converted into familiar operations in geometry. Each of the above languages (PL, semialgebraic, subanalytic) are instances of an o-minimal structure. The common fundamental theorems, each expressed in their own language, can be reduced to universal logical operations, and hence geometric ones. We will start by examining o-minimal structures as they form the local models of Miller and van den Dries definition. The reader is urged to consult the textbook ``Tame Topology and O-minimal Structures''~\cite{vdd-ttos} as it is an excellent introduction that requires virtually no pre-requisites.

\begin{defn}[\cite{vdd-ttos}, p. 2]\index{o-minimal structure}
 An \textbf{o-minimal structure on $\RR$} is a sequence $\OO=\{\OO_n\}_{n\geq 0}$ satisfying
 \begin{enumerate}
  \item $\OO_n$ is a boolean algebra of subsets of $\RR^n$, i.e. it is a collection of subsets of $\RR^n$ closed under unions and complements, with $\emptyset \in\OO_n$;
  \item If $A\in\OO_n$ then $A\times\RR$ and $\RR\times A$ are both in $\OO_{n+1}$;
  \item The sets $\{(x_1,\ldots,x_n)\in\RR^n|x_i=x_j\}$ for varying $i\leq j$ are in $\OO_n$;
  \item If $A\in\OO_{n+1}$ then $\pi(A)\in\OO_{n}$ where $\pi:\RR^{n+1}\to\RR^n$ is projection onto the first $n$ factors;
  \item For each $x\in\RR$ we require $\{x\}\in\OO_1$ and $\{(x,y)\in\RR^2|x<y\}\in\OO_2$;
  \item The only sets in $\OO_1$ are the finite unions of open intervals and points.
 \end{enumerate}
When working with a fixed o-minimal structure $\OO$ on $\RR$ we say a subset of $\RR^n$ is \textbf{definable} if it belongs to $\OO_n$. A map is definable if its graph is definable.
\end{defn}

\begin{rmk}
	One should note that the third and sixth property together prohibit any spiral that wraps infinitely many times around the origin from being part of an o-minimal structure. Thus, the quick spiral in Example \ref{ex:pullback_spiral} is not definable.
\end{rmk}

Now we prove that definable sets and maps are closed under pullbacks.
\begin{lem}\label{lem:def_pb}
Suppose $f:X\to Z$ and $g:Y\to Z$ are definable maps, then the pullback $X\times_Z Y:=\{(x,y)\in X\times Y\,|\, f(x)=g(y)\}$ is a definable set and the restrictions of the projection maps are definable as well.
\end{lem}
\begin{proof}
 First note that if $X\in \OO_n$ and $Y\in \OO_m$, then $X\times Y=(X\times \RR^m)\cap (\RR^n\times Y)$ is in $\OO_{n+m}$. Since $\Gamma_f$ and $\Gamma_g$ are definable, we know that $\Gamma_f\times Y=\{(x,y,f(x))\}$ and $\Gamma_g\times X=\{(x',y',g(y')\}$ are both definable subsets of $X\times Y\times Z$. Since the intersection is definable, and a point in the intersection has $(x,y,f(x))=(x',y',g(y'))$, the image of the projection to $X\times Y$ is the pullback. One can then use B.3 of~\cite{vdd-geocat} to conclude that the restriction to the pullback of the projection maps to $X$ and $Y$ is definable as well.
\end{proof}

There are surprising facts that follow from the axioms of an o-minimal structure. For example, if $A\in\OO$, then the closure $\bar{A}$ is in $\OO$ (\cite{vdd-ttos} Ch. 1, 3.4). Another surprising fact is that definable sets can be Whitney stratified~\cite{loi-verdier}. Thus, these sets meet the requirements of MacPherson to form a good class of subsets. Perhaps even better than MacPherson's sets, definable sets can be given finite cell decompositions, where ``cell'' has its own special meaning (\cite{vdd-ttos} Ch. 3).

The prototypical o-minimal structure is the class of semialgebraic sets, which has become increasingly relevant in applied mathematics.

\begin{defn}\index{semialgebraic subset}
	A \textbf{semialgebraic} subset of $\RR^n$ is a subset of the form
	\[
		X=\bigcup_{i=1}^p\bigcap_{j=1}^q X_{ij}
	\]
	where the sets $X_{ij}$ are of the form $\{f_{ij}(x)=0\}$ or $\{f_{ij}>0\}$ with $f_{ij}$ a polynomial in $n$ variables.
\end{defn}

	The only semi-algebraic subsets of $\RR$ are finite unions of points and open intervals. From the definition, one sees that the class of semialgebraic sets is closed under finite unions and complements. The \textbf{Tarski-Seidenberg} theorem states that the projection onto the first $m$ factors $\RR^{m+n}\to\RR^m$  sends semialgebraic subsets to semialgebraic subsets~\cite{coste-sag}. We can deduce from this theorem all of the conditions of o-minimality.
	
	Semialgebraic maps are defined to be those maps $f:\RR^k\to \RR^n$ whose graphs are semialgebraic subsets of the product. It is a fact that semi-algebraic sets and maps can be Whitney stratified~\cite{shiota-geometry}. This allows us to consider the following example of a semi-algebraic family of sets:
	
	\begin{figure}
		\centering
		\includegraphics[width=.9\textwidth]{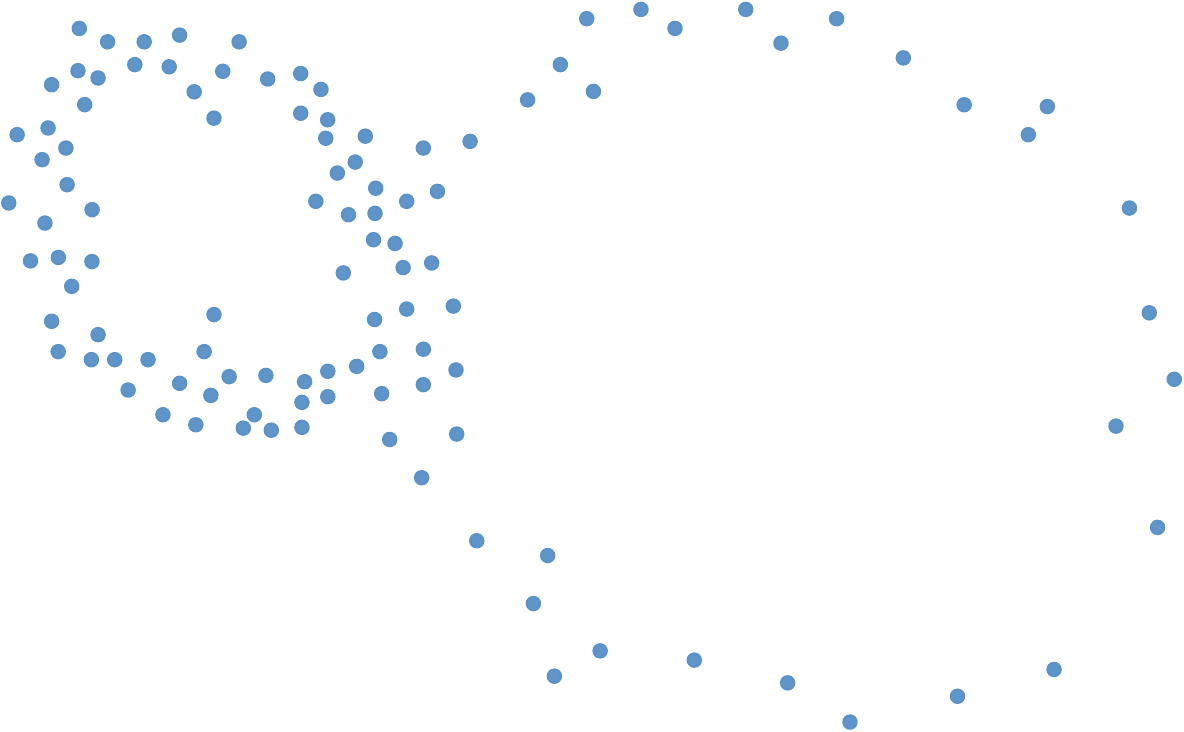}
		\caption{Point Cloud Data}
		\label{fig:point_cloud}
	\end{figure}

\begin{ex}[Point-Cloud Data]\label{ex:pcd}\index{point cloud}
	Suppose $Z$ is a finite set of points in $\RR^n$. For each $z\in Z$, consider the square of the distance function
	\[f_z(x_1,\ldots,x_n)=\sum_{i=1}^n (x_i-z_i)^2.\]
	By the previously stated facts we know that the sets
	\[
		B_z:=\{x\in\RR^{n+1} \, | \, f_z(x_1,\ldots,x_n)\leq x_{n+1}\}
	\]
	are semialgebraic along with their unions and intersections. Denote by $X$ the union of the $B_z$. The Tarski-Seidenberg theorem implies that the map
	\[
		f:X\to\RR \qquad f^{-1}(r):=\cup_{z\in Z} B(z,\sqrt{r})=\{x\in\RR^n \, | \, \exists z\in Z \, \mathrm{s.t.} \, f_z(x)\leq r\}	\]
	is semialgebraic. In particular the topology of the fiber (of the union of the closed balls) can only change finitely many times.
\end{ex}

We conclude with the definition Miller and van den Dries proposed in section 1 of~\cite{vdd-geocat}. This definition allows us to verify definability locally, and allows us to work inside manifolds other than $\RR^n$.

\begin{defn}[Analytic-Geometric Categories]\index{analytic-geometric category}\index{o-minimal structure!analytic-geometric category}
	A \textbf{analytic-geometric category} $\calG$ is given by assigning to each analytic manifold $M$ a collection of subsets $\calG(M)$ such that following conditions are satisfied:
	\begin{enumerate}
		\item $\calG(M)$ is a boolean algebra of subsets of $M$, with $M\in\calG(M)$.
		\item If $A\in\calG(M)$, then $A\times\RR\in\calG(M)$.
		\item If $f:M\to N$ is a proper analytic map and $A\in\calG(M)$, then $f(A)\in\calG(N)$.
		\item If $A\subseteq M$ and $\{U_i\}_{i\in\Lambda}$ is an open covering of $M$, then $A\in\calG(M)$ if and only if $A\cap U_i\in\calG(U_i)$ for all $i\in\Lambda$.
		\item Every bounded set in $\calG(\RR)$ has finite boundary.
	\end{enumerate}
\end{defn}
\begin{rmk}
This defines a category in the usual sense. An object of $\calG$ is a pair $(A,M)$ with $A\in\calG(M)$. A morphism $f:(A,M)\to (B,N)$ is a continuous map $f:A\to B$ whose graph
\[
	\Gamma(f):=\{(a,f(a))\in M\times N\,|\, a\in A\}
\]
is an element of $\calG(M\times N)$.
\end{rmk}

The category of $\calG$-sets and $\calG$-maps, although we will prefer to use the term ``definable,'' has all the properties one could desire, including being closed under inverse images~\cite[D.7]{vdd-geocat} (as long as the domain is closed) and Whitney stratifiability~\cite[D.16]{vdd-geocat}.\footnote{The authors of~\cite{vdd-geocat} acknowledge that there is a gap in the proof of Whitney stratifiability of $\calG$-maps, but Ta L\^{e} Loi~\cite{loi-omin} and others~\cite{definably-stratified} have since filled in this gap.}

\subsection{Thom-Mather Stratifications}

\begin{defn}[Control Data]\index{control data}
	Let $(X,M,\{X_{\sigma}\}_{\sigma\in P_X})$ be a Whitney stratified space and $\{(T_{\sigma},\pi_{\sigma},d_{\sigma})\}$ a family of tubular neighborhoods. We call this family a system of \textbf{control data} if the following commutation relations are satisfied: if $X_{\sigma}\leq X_{\tau}$, then
	\begin{eqnarray*}
		\pi_{\sigma}\circ \pi_{\tau} &=& \pi_{\sigma} \\
		d_{\sigma}\circ \pi_{\tau}&=& d_{\sigma}
	\end{eqnarray*}
whenever both sides of the equations are defined.
\end{defn}

\begin{figure}
	\centering
	\includegraphics[width=.9\textwidth]{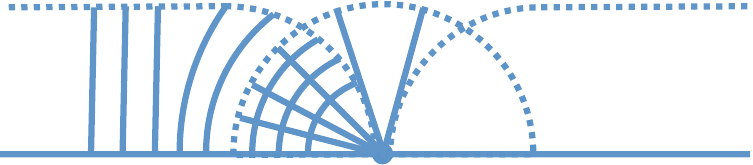}
	\caption{A System of Control Data}
	\label{fig:control}
\end{figure}
\begin{rmk}
In Figure~\ref{fig:control} we have drawn some fibers of the retraction maps for two incident strata. Notice how the fibers must bend in order for the second compatibility condition to hold.
\end{rmk}

Mather proves that every Whitney stratified space admits a system of control data. The following definition axiomatizes the properties enjoyed by a Whitney space with a system of control data.

\begin{defn}[Thom-Mather Stratified Spaces]\index{stratified space!Thom-Mather}
	A \textbf{Thom-Mather stratified space} consists of a Hausdorff, locally compact topological space $X$ with countable basis for its topology with some smooth structure; a decomposition into topological manifolds $\{X_{\sigma}\}_{\sigma\in P_X}$; and a family of control data $\{(T_{\sigma},\pi_{\sigma},d_{\sigma})\}_{\sigma\in P_X}$, where $T_{\sigma}$ is an open tubular neighborhood of $X_{\sigma}$, $\pi_{\sigma}:T_{\sigma}\to X_{\sigma}$ is a continuous retraction, and $d_{\sigma}:X_{\sigma}\to [0,\infty)$ is a continuous distance function. We require that the following conditions hold:
	\begin{itemize}
		\item $X_{\sigma}=d_{\sigma}^{-1}(0)$ for all $\sigma$.
		\item For any pair of strata $X_{\sigma}, X_{\tau}$, define $T_{\sigma,\tau}:=T_{\sigma}\cap X_{\tau}$, $\pi_{\sigma,\tau}:=\pi_{\sigma}|_{T_{\sigma,\tau}}$ and $d_{\sigma,\tau}:=d_{\sigma}|_{T_{\sigma,\tau}}$. We require that
		\[
			(\pi_{\sigma,\tau},d_{\sigma,\tau}): T_{\sigma,\tau}\to X_{\sigma}\times (0,\infty)
		\]
		is a smooth submersion. When $T_{\sigma,\tau}\neq\emptyset$, i.e. when $X_{\sigma}\leq X_{\tau}$, this implies $\dim X_{\sigma} < \dim X_{\tau}$.
		\item For any trio of strata $X_{\sigma}, X_{\tau}$ and $X_{\lambda}$ we have
		\begin{eqnarray*}
			\pi_{\sigma,\tau}\circ \pi_{\tau,\lambda}&=&\pi_{\sigma,\lambda} \\
			d_{\sigma,\tau}\circ \pi_{\tau,\lambda}&=& d_{\sigma,\lambda}
		\end{eqnarray*}
		whenever both sides of the equation are defined.
	\end{itemize}
\end{defn}
\begin{rmk}
	One should observe that the definition does not require an embedding into an ambient space. Thus Thom-Mather stratified spaces allow us to treat Whitney stratified spaces intrinsically. Any Whitney stratified space $(X,M)$ equipped with a system of control data $\{(T_{\sigma},\pi_{\sigma},d_{\sigma})\}_{\sigma\in P_X}$ defines a Thom-Mather stratified space by intersecting each $T_{\sigma}$, which is open in $M$, with $X$.
\end{rmk}

Thom-Mather stratified spaces exhibit most of the good properties of Whitney stratified spaces. The proof that Thom-Mather spaces can be triangulated was carried out by Goresky~\cite{goresky-fol}, among others. His proof views the lines of Figure \ref{fig:control} not as fibers of the retraction map $\pi_{\sigma}$, but rather\footnote{This description does not hold in higher dimensions.} as fibers of a radial projection map to the boundary of a tubular neighborhood.

\begin{defn}[Family of Lines]\label{defn:fol}\index{family of lines}
	A \textbf{family of lines} on a Thom-Mather stratified space is a system of radial projections
\[
	r_{\sigma}(\epsilon):T_{\sigma}-X_{\sigma}\to S_{\sigma}(\epsilon):=d_{\sigma}^{-1}(\epsilon)
\]
one for each stratum $X_{\sigma}$, and a positive number $\delta$, such that whenever $0<\epsilon <\delta$ and $X_{\sigma}\leq X_{\tau}$, the following commutation relations hold:
\begin{enumerate}
	\item $r_{\sigma}(\epsilon)\circ r_{\tau}(\epsilon')=r_{\tau}(\epsilon')\circ r_{\sigma}(\epsilon)\in S_{\sigma}(\epsilon)\cap S_{\tau}(\epsilon')$
	\item $d_{\sigma}\circ r_{\tau}(\epsilon)=d_{\sigma}$
	\item $d_{\tau}\circ r_{\sigma}(\epsilon)=d_{\tau}$
	\item $\pi_{\sigma}\circ r_{\tau}(\epsilon)=\pi_{\sigma}$
	\item If $0<\epsilon <\epsilon'<\delta$, then $r_{\sigma}(\epsilon')\circ r_{\sigma}(\epsilon)=r_{\sigma}(\epsilon')$
	\item $\pi_{\sigma}\circ r_{\sigma}(\epsilon)=\pi_{\sigma}$
	\item $r_{\sigma}(\epsilon)|_{T_{\sigma}(\epsilon)\cap X_{\tau}}:T_{\sigma}(\epsilon)\cap X_{\tau}\to S_{\sigma}(\epsilon)\cap X_{\tau}$ is smooth
\end{enumerate}
\end{defn}
\begin{rmk}
Every Thom-Mather stratified space admits a family of lines. This is the first proposition of~\cite{goresky-fol}.
\end{rmk}

Any family of lines can be used to identify a tubular neighborhood $T_{\sigma}$ as a mapping cylinder for the restricted projection map $\pi_{\sigma}:S_{\sigma}(\epsilon)\to X_{\sigma}$. To do so, one defines a stratum-preserving homeomorphism
\[
	h_{\sigma}:T_{\sigma}- X_{\sigma}\to S_{\sigma}(\epsilon)\times (0,\infty) \qquad h_{\sigma}(p):= (r_{\sigma}(\epsilon)(p),d_{\sigma}(p))
\]
and then extends the map in a suitable way, i.e. one takes $S_{\sigma}(\epsilon)\times [0,\infty)\sqcup X_{\sigma}$ and identifies $S_{\sigma}(\epsilon)\times\{0\}$ with its image under $\pi_{\sigma}:S_{\sigma}(\epsilon)\to X_{\sigma}$. One can check that this allows us to extend our map $h_{\sigma}$ to $T_{\sigma}$, which we do so without changing notation. One should interpret this extended homeomorphism $h_{\sigma}$ as providing a system of coordinates that is convenient for analyzing neighborhoods of strata. We use this system of coordinates in the following theorem.

\begin{figure}
	\centering
	\includegraphics[width=.9\textwidth]{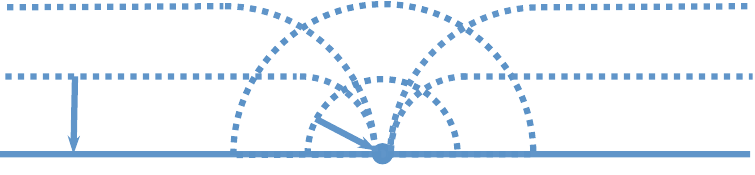}
	\caption{Two Sets of Tubular Neighborhoods}
	\label{fig:reg_nbhd}
\end{figure}

\begin{prop}[Regular Neighborhoods of Closed Unions of Strata]\label{prop:reg_nbhd}\index{strata!closed union, has regular neighborhood}
	Let $X$ be a Thom-Mather stratified space and $W=\cup_{\sigma\in P_W} X_{\sigma}$ be a closed union of strata of $X$. The inclusion
	\[
		W\hookrightarrow U_W(\epsilon/2):=\bigcup_{\sigma\in P_W} T_{\sigma}(\epsilon/2)\subseteq X
	\]
	is a homotopy equivalence.
\end{prop}
\begin{proof}
	Given a Thom-Mather stratified space, we can equip it with a family of lines~\cite{goresky-fol}. We are going to use the family of lines to construct a weak deformation retraction of $U_W(\epsilon/2)$ inside a larger open neighborhood $U_W(\epsilon):=\cup_{\sigma\in P_W} T_{\sigma}(\epsilon)$. The idea is to shrink each tubular neighborhood $T_{\sigma}(\epsilon/2)$ to $X_{\sigma}$ in such a way that a line connecting a point $p\in S_{\sigma}(\epsilon/2)$ and $r_{\sigma}(\epsilon)(p)\in S_{\sigma}(\epsilon)$ is stretched to connect $\pi_{\sigma}(p)$ and $r_{\sigma}(\epsilon)(p)$ after the homotopy. Figure \ref{fig:reg_nbhd} indicates which neighborhoods are to be collapsed.
	
	To accomplish this stretching, let $f:\RR\to [0,1]$ be any smooth function with the following properties:
	\[
		\begin{array}{r c l} 
			f(x)=0 & \mathrm{if} & x\leq \frac{1}{2} \\[\medskipamount]
						f(x)=1 & \mathrm{if} & x\geq \frac{3}{4} \\[\medskipamount]
						f'(x)>0 & \mathrm{if} & x\in(\frac{1}{2},\frac{3}{4})
		\end{array}
	\]
	The homotopy $H_{\sigma}:U\times[0,1]\to U$ defined below shrinks $T_{\sigma}(\epsilon/2)$ to $X_{\sigma}$:
	\[
		H_{\sigma}(p,t):=\left\{\begin{array}{lrl} 
		p & \mathrm{if} & p\notin T_{\sigma}(\epsilon) \\
		h_{\sigma}^{-1}(r_{\sigma}(\epsilon)(p),d_{\sigma}(p)[(1-t) f(d_{\sigma}(p)/\epsilon) + t]) & \mathrm{if} & p\in T_{\sigma}(\epsilon) \end{array}\right.
	\]	
	The homotopy is just a straight-line homotopy between the usual distance function $d_{\sigma}$ and the shrunken one $d_{\sigma}(p)f(d_{\sigma}(p)/\epsilon)$. Moreover, since the homotopy only affects the distance coordinate, properties two and three of Definition \ref{defn:fol} imply that if $X_{\sigma}\leq X_{\tau}$ then
	\[
		d_{\tau}(H_{\sigma}(p,t))=d_{\tau}(p) \qquad \mathrm{and} \qquad d_{\sigma}(H_{\tau}(p,t))=d_{\sigma}(p).
	\]
	As such, the shrinking homotopies can be applied in any order, i.e.
	\[
		H_{\tau}(H_{\sigma}(p,t),s)=H_{\sigma}(H_{\tau}(p,s),t).
	\]
	Observe that in a Thom-Mather stratified space, if $X_{\sigma}\neq X_{\sigma'}$ are two strata of the same dimension, then $T_{\sigma}\cap T_{\sigma'}=\emptyset$. Consequently, the definition for $H_{\sigma}$ extends to a homotopy $H_i$ that shrinks all the neighborhoods of strata of dimension $i$ at the same time; one just defines $H_i(p,t)=H_{\sigma}(p,t)$ if $p\in T_{\sigma}(\epsilon)$. The commutation relation now extends to the statement that for any $i$ and $j$
	\[
		H_{i}(H_{j}(p,t),s)=H_{j}(H_{i}(p,s),t).
	\]
	Thus, our desired homotopy can be defined to be
	\[
		H(p,t):=H_0(H_1(H_2(\cdots (H_m(p,t),t)\cdots,t),t),t)
	\]
	where the order of the composition doesn't matter and $m$ is the maximum dimension of a stratum appearing in $W$. If we let $r_t(p)=H|_{U(\epsilon/2)\times I}$, then $r_t$ defines a weak deformation retract of $U_W(\epsilon/2)$ to $W$, that is, $r_t(W)\subseteq W$ for all $t$, $r_0(U_W(\epsilon/2))\subseteq W$ and $r_1=\id$. It is easy to show that this implies that $W\hookrightarrow U_W(\epsilon/2)$ is a homotopy equivalence.
\end{proof}
\begin{rmk}
	One could imagine performing these homotopies at separate times by letting the homotopy parameter in dimension $i$ be a function $s_i(t)=f(t-i)$ where the shrinking homotopy in dimension $i$ is performed in the interval $(i+1/2,i+3/4)$. This is how it is done in Goresky's thesis~\cite{goresky-thesis}. This makes his homotopy
	\[
		H^G_{i}(p,t):=\left\{\begin{array}{lrl} 
		p & \mathrm{if} & p\notin T_{\sigma}(\epsilon) \\
		h_{\sigma}^{-1}(r_{\sigma}(\epsilon)(p),d_{\sigma}(p)[(1-s_i(t)) f(d_{\sigma}(p)/\epsilon) + s_i(t)]) & \mathrm{if} & p\in T_{\sigma}(\epsilon) \end{array}\right.
	\] 
	easier to visualize. However, the advantage of choosing $s_i(t)=t$ is that the homotopy is stratum preserving up until $t=0$.
\end{rmk}
\begin{rmk}
	Of course, for a given stratum $X_{\sigma}$, away from its frontier the retraction map $r_0$ coincides with the tubular projection $\pi_{\sigma}$.
\end{rmk}

As the above proposition shows, control data is essential for providing Whitney stratified spaces with good neighborhoods. Not only do they endow Whitney stratified spaces with the structure of a Thom-Mather stratified space, they allow us to construct Goresky's family of lines to carry out these retractions. These retractions are instrumental to the cosheaves that we will construct in Lemma \ref{lem:strat_maps_1d} and Theorem \ref{thm:strat_maps}. There is another technical tool that we need that can only be developed in the presence of control data.

\begin{defn}\index{stratified vector field}
A \textbf{stratified vector field} $\eta$ on $(X,\{X_{\sigma}\}_{\sigma\in P_X})$ is a collection of vector fields $\{\eta_{\sigma}\}_{\sigma\in P_X}$ with one smooth vector field on each stratum.
\end{defn}

When it is meaningful to compare these vector fields, it is remarkable to note that this collection need not be continuous. Nevertheless, in the presence of control data, the flow generated by such a discontinuous vector field is continuous.

\begin{defn}
A stratified vector field $\eta$ on $X$ is said to be \textbf{controlled} by $\{T_{\sigma},\pi_{\sigma},d_{\sigma}\}$ if the following compatibility conditions are satisfied for any pair of strata $X_{\sigma}\leq X_{\tau}$:
\begin{eqnarray*}
 \eta_{\tau}(d_{\sigma,\tau}(p)) & = & 0 \\
 d(\pi_{\sigma,\tau})(\eta_{\tau}(p)) & = & \eta_{\sigma}(\pi_{\sigma,\tau}(p)) 
\end{eqnarray*}
where ever both sides of the equation are defined.
\end{defn}

\subsection{Thom Mappings}

\begin{defn}\index{stratified map!Thom mapping}
	A \textbf{Thom mapping} is a stratified map $f:(X,M)\to (Y,N)$ that satisfies \textbf{condition $a_f$} for every pair of strata $X_{\tau}\geq X_{\sigma}$: let $x_i$ be a sequence of points in $X_{\tau}$ converging to a point $p\in X_{\sigma}$. Suppose $\ker d(f|_{X_{\tau}})_{x_i}\subseteq T_{x_i} M$ converges to a plane $K\subseteq T_{p} M$, then $\ker d(f|_{X_{\sigma}})_p\subseteq K$. 
\end{defn}

\begin{figure}
	\centering
	\includegraphics[width=.8\textwidth]{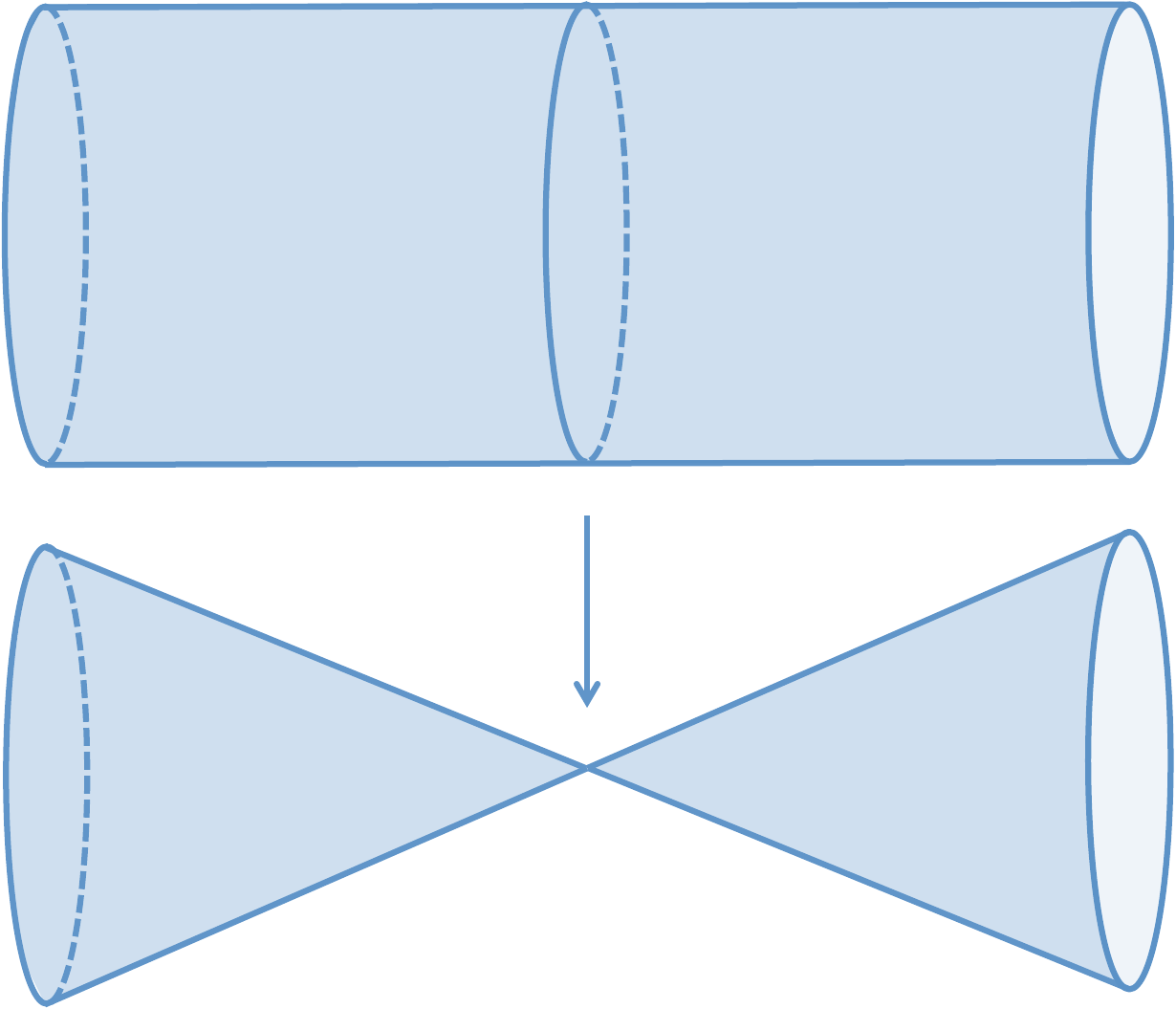}
	\caption{Not a Thom Mapping}
	\label{fig:not_thom}
\end{figure}

In Figure \ref{fig:not_thom}, we have drawn an example of a mapping that is not a Thom mapping.\footnote{This example is borrowed from~\cite{lu-ctst}.} Other non-examples include the blow-up map discussed in Example \ref{ex:blowup}. Any map that is triangulable satisfies Thom's condition $a_f$ for that triangulation viewed as a stratification. It has been a long standing conjecture that every smooth Thom mapping is triangulable. Masahiro Shiota appears to have proven this conjecture in the $C^{\infty}$ case~\cite{shiota-thom}, but we have chosen not to rely on this conjecture. Instead, we only need the following proposition of Mather's (Proposition 11.3 of~\cite{mather}).

\begin{prop}\label{prop:thom_data}
	Suppose $f:X \to Y$ is a Thom mapping and a system of control data $\{T\}$ for $Y$ is given. There exists a family of tubular neighborhoods $\{T'\}$ for $X$ over $\{T\}$, which satisfies the following compatibility conditions:
	\begin{enumerate}
		\item[(a)] If $X_{\sigma}\leq X_{\tau}$, then $\pi'_{\sigma} \circ \pi'_{\tau}=\pi'_{\sigma}$ for points in $T_{\sigma}'\cap T_{\tau}'$ in $M$. Furthermore, if $f(X_{\sigma})$ and $f(X_{\tau})$ lie in the same stratum of $Y$, then $d'_{\sigma}\circ \pi'_{\tau}=d'_{\sigma}$ where both sides are defined.
		\item[(b)] If $Y_{\sigma}$ is a stratum that contains $X_{\sigma}$, then
		\[
			f(\pi'_{\sigma}(p))=\pi_{\sigma}(f(p))
		\]
		for all $p\in T'_{\sigma}\cap f^{-1}(T_{\sigma})$.
	\end{enumerate}
\end{prop}
\begin{rmk}
	The first condition is weaker than the usual definition of control data when the strata are not mapped to the same stratum. Consequently, the above notion of a \textbf{system $\{T'\}$ of control data over $\{T\}$} is not the same as two systems of control data.
\end{rmk}

Just as the notion of control data generalizes to control data over control data, controlled vector fields generalize to controlled vector fields over controlled vector fields.

\begin{defn}
Suppose $f:X\to Y$ is a Thom mapping and $\{T'\}$ is a system of control data over $\{T\}$. If $\eta=\{\eta_{\sigma}\}$ is a controlled vector field on $\{Y_{\sigma}\}$ controlled by $\{T\}$, then there exists a stratified vector field $\eta'=\{\eta'_{\sigma}\}$ on $\{X_{\sigma}\}$ satisfying the following compatibility conditions:
\begin{enumerate}
\item[(a)] For any $X_{\sigma}$ and $p\in X_{\sigma}$, we have
\[
(df|_{X_{\sigma}})(\eta'_{\sigma}(p)) = \eta_{\sigma}(f(p))
\]
where $Y_{\sigma}$ is the stratum of $Y$ that contains $f(p)$.
\item[(b)] For any $X_{\sigma}\leq X_{\tau}$, there is a neighborhood $N'_{\sigma}$ in $T'_{\sigma}$ such that for $p\in T'_{\sigma}\cap X_{\tau}$ we have
\[
d(\pi'_{\sigma,\tau})(\eta'_{\tau}(p)) = \eta'_{\sigma}(\pi'_{\sigma,\tau}(p))
\]
and if $X_{\sigma}$ and $X_{\tau}$ are carried to the same stratum of $Y$, then we have further the condition that
\[
\eta'_{\tau}(d'_{\sigma,\tau}(p)) = 0.
\]
\end{enumerate}
Thus, the notion of a \textbf{controlled vector field $\eta'$ over $\eta$} is a weaker one than a pair of controlled vector fields on $X$ and $Y$ that commute with the Thom mapping $f$.
\end{defn}

The following result, proven with help from Mark Goresky, gives a useful criterion for determining when a stratified map is a Thom mapping, so as to make the above constructions possible there. It rests on the observation that all the classical examples of stratified maps $f:X\to Y$ that aren't Thom maps require considering a pair of strata $Y_{\sigma}<Y_{\tau}$ in $Y$ whose codimension is at least two. Combinatorially, this allows us to have a pair of strata $X_{\sigma}<X_{\tau}$ in $X$ such that $\dim X_{\sigma}\cap f^{-1}(p) > \dim X_{\tau}\cap f^{-1}(x_i)$ even though $\dim X_{\sigma}<\dim X_{\tau}$. In the following lemma we show that if the codomain only has strata of codimension 1, then the map is a Thom mapping.

\begin{lem}\label{lem:thom_codim}\index{stratified map!Thom mapping}
	Suppose $f:(X,M)\to (Y,N)$ is a Whitney stratified map that is $C^1$ on the ambient manifold $M$. Let $Y'=Y_{\sigma}\cup Y_{\tau}$ be the union of two strata whose difference in dimension is one. The restricted map $f':(X',M)\to (Y',N)$ where $X':=f^{-1}(Y')$ is a Thom map.
\end{lem} 
\begin{proof}
	The proof is local, so we consider the following setup instead: Suppose $f:(X,M) \to (Y,\RR^{k+1})$ is a Whitney stratified map where $Y$ is the upper half plane in $\RR^{k+1}$, i.e. $Y:=\{(y_1,\cdots,y_{k+1})\,|\, y_{k+1}\geq 0\}$. We assume that the stratification of the map stratifies $Y$ as $Y_{\tau}:=\{y_{k+1}>0\}\cong \RR^{k+1}$ and $Y_{\sigma}:=\{y_{k+1}=0\}\cong \RR^k$. Let $X_{\sigma}$ be a stratum of $X$ that is mapped to $Y_{\sigma}$ and $X_{\tau}$ a stratum mapped to $Y_{\tau}$. Suppose $\{x_i\}$ is a sequence in $X_{\tau}$ and $\ker df|_{X_{\tau}}(x_i)=:K_i$ converges to a subspace $K_{\infty}\subseteq T_p M$ where $p\in X_{\sigma}$. We want to show that $K_p:=\ker df|_{X_{\sigma}}(p)\subseteq K_{\infty}$. By passing to a subsequence we can further assume that the tangent planes $T_{x_i} X_{\tau}=:T_i$ converges to $T_{\infty}\subseteq T_p M$. By Whitney's condition (a), $T_p X_{\sigma}\subset T_{\infty}$.
	
	Denote by $\rho_Y(y):=\pi_{k+1}(y_1,\ldots,y_{k+1})=y_{k+1}$ the ``distance from the stratum'' function on $Y$. By pre-composing with $f$, this defines a function $\rho_X(x):=\rho_Y(f(x))$. Any vector $v\in T_i$ with $d\rho_X(x_i)(v)\neq 0$ must also have $df|_{X_{\tau}}(x_i)(v)\neq 0$ since the chain rule implies that $d\rho_X(x_i)=d\rho_Y(f(x_i))\circ df|_{X_{\tau}}(x_i)$ and thus $v\notin K_i$.
	
	Let $\pi_{\sigma}:\RR^{k+1}\to Y_{\sigma}$ be the projection onto the first $k$ coordinates. The restriction of $\pi_{\sigma}$ to $Y_{\tau}$, written $\pi_{\sigma,\tau}$, is a submersion. By virtue of $\pi_{\sigma,\tau}\circ f|_{X_{\tau}}$ being a submersion, any vector $w\in T_{\pi_{\sigma}(f(x_i))} Y_{\sigma}$ has a lift $w_{f(x_i)}\in T_{f(x_i)}Y_{\tau}$ so that $w_{f(x_i)}\in \ker d\rho_Y(f(x_i))$, which in turn has a lift $\tilde{w}_i\in T_i$. Consequently, $df|_{X_{\tau}}(x_i)(\tilde{w}_i)\neq 0$ and thus $\tilde{w}_i\notin K_i$. Moreover, $\tilde{w}_i$ is orthogonal to $\nabla \rho_{X}(x_i)$ since any lift of $w$ is chosen to factor through the kernel of $d\rho_Y(f(x_i))=\pi_{k+1}$.
	
	Thus, each $T_i$ can be written as $\tilde{T}_{\pi_{\sigma}(f(x_i))} Y_{\sigma}\oplus K_i \oplus \nabla \rho_{X}(x_i)$. Since $T_p X_{\sigma}\subset T_{\infty}$ the isomorphism $T_{\infty}\cong T_p X_{\sigma}\oplus (T_p X_{\sigma})^{\perp}$ can be further refined as $T_{\infty}\cong \tilde{T}_{f(p)} Y_{\sigma} \oplus K_p \oplus (T_p X_{\sigma})^{\perp}$. We have assumed that $f$ is $C^1$ on the ambient manifold $M$ so that the lifts $\tilde{T}_{\pi_{\sigma}(f(x_i))} Y_{\sigma}$ must converge (perhaps after passing again to a subsequence) to $\tilde{T}_{f(p)} Y_{\sigma}$. Additionally, $\nabla \rho_X(x_i)$ converges to a subspace of $(T_p X_{\sigma})^{\perp}$. Finally, since $\dim X_{\sigma}< \dim X_{\tau}$, dimension constraints force $K_p\subseteq K_{\infty}$. This proves the lemma.
\end{proof}

This lemma is instrumental for our proof of Theorem \ref{thm:strat_maps}. On it's own, it has a useful corollary.

\begin{cor}
	Any stratified map $f:(X,M)\to (Y,\RR)$ that is $C^1$ on the ambient manifold is a Thom map.
\end{cor}

\subsection{Stratified Maps to the Real Line}

\begin{lem}\label{lem:strat_maps_1d}
	Any stratified map $f:X\to\RR$ defines, for each $i$, a cellular cosheaf.
\end{lem}
\begin{proof}
The map $f:X\to \RR$ defined above has as fibers the spaces $X_r$. Because it is stratifiable with finitely many strata, we have the following decomposition of the codomain:
\[
 	(-\infty,0) \leftarrow \{0\} \rightarrow (0,t_1) \leftarrow \{t_1\} \rightarrow (t_1,t_2) \leftarrow \{t_2\} \rightarrow (t_2,t_3) \cdots
\]
The points $t_i$ indicate the radii (the ``times'') where the topology of the union of the balls changes. Since the fiber $X_{t_i}:=f^{-1}(t_i)$ is a closed union of strata, proposition \ref{prop:reg_nbhd} implies (after first choosing a system of control data and then regarding $X$ as Thom-Mather stratified) that we can fix an $\epsilon>0$ such that the neighborhood $U_{t_i}(\epsilon)=\cup_{\sigma\in X_{t_i}} T_{\sigma}(\epsilon/2)$ contains $X_{t_i}$ as a weak deformation retract. Since $f$ is proper, we claim that there exists a point $s_i^-\in (t_{i-1},t_i)$ such that $X_{s_i^-}$ is contained in $U_{t_i}(\epsilon)$. Suppose for contradiction that for all $n>>0$ there exists a point $x_n\in f^{-1}([t_i-\frac{1}{n},t_i-\frac{1}{n+1}))\cap U_{t_i}(\epsilon)^c$. If this is possible, then $\{x_n\}$ defines a sequence with no convergent subsequence, which contradicts the fact that $f^{-1}([t_{i-1},t_i])$ is compact. Consequently, there exists an $n$ such that if $s_i^-:=t_i-\frac{1}{n}$, then $f^{-1}([s_i^-,t_i])\subseteq U_{t_i}(\epsilon)$. The composition of the inclusion followed by the retraction
\[
	\xymatrix{ & U_{t_i}(\epsilon) & \\ 
	X_{s_i^-} \ar@{^{(}->}[ru] & & X_{t_i} \ar@{_{(}->}[lu]_{\cong}}
\]
allows us to define maps between the homology of the typical fiber over $(t_{i-1},t_i)$ to the homology of the fiber $X_{t_i}$.
\[
	H_i(X_{s_i^-};k)\rightarrow H_i(X_{t_i};k)
\]
An analogous argument allows us to find an $s_i^+\in (t_i,t_{i+1})$ such that $X_{s_i^+}\subset U_{t_i}(\epsilon)$. We can construct a vector field on $(t_{i-1},t_i)$ that flows from the point $s_{i-1}^+$ to $s_i^-$. Lifting this vector field to a controlled one over this one, allows us to flow the fiber over $s_{i-1}^+$ to the fiber over $s_i^-$, thus realizing the homeomorphisms $X_{s_{i-1}^+}\cong X_{s_i^-}$ explicitly. For convenience, we drop the decorations and choose any point $s_i\in (t_{i-1},t_i)$ to get our modified version of the persistence module introduced in Section~\ref{subsec:pointclouds}.
\[
	\cdots \leftarrow H_i(X_{s_i};k)\rightarrow H_i(X_{t_i};k) \leftarrow H_i(X_{s_{i+1}};k)\rightarrow H_i(X_{t_{i+1}}) \leftarrow \cdots
\]
One should note that this diagram is contravariant with respect to the poset indexing the stratification of $\RR$, thus we have constructed geometrically a cellular cosheaf.
\end{proof}

\begin{cor}\index{point cloud!determines a stratified map}
	The semialgebraic function $f:X\to\RR$ in example \ref{ex:pcd} defines, for each $i$, a cellular cosheaf.
\end{cor}

Although the above construction may appear convoluted, it is geometrically natural. Instead of using the order on $\RR$ to get a diagram of vector spaces and maps, we have a diagram indexed by the pieces of a stratification of $\RR$. This new diagram is specifically adapted to the topological changes in the family $\{X_r\}$. 

In multi-dimensional persistence we imagine the need for more than one parameter to distinguish features in a point cloud. The traditional story of persistence no longer applies since $\RR^n$ for $n\geq 2$ has no natural (partial) order. In contrast, every situation where multi-dimensional persistence can be treated as a stratified map (which is effectively always), the partial order of the pieces in a stratification presents itself as a most natural candidate. 

However, the geometry of stratified spaces in more than one dimension is subtle and a poset will not always suffice. In Section \ref{subsubsec:path_cat}, we will introduce a small category (usually equivalent to a finite one) that allows us to track persistent features in a more careful way. The proof of Lemma \ref{lem:strat_maps_1d} contains the essential ideas of this more general picture. By considering certain definable paths in the parameter space, and analyzing their inverse images, which will be definable, we can try to reduce a multi-dimensional problem to a one-dimensional one. This is the high-level outline of how Theorem \ref{thm:strat_maps} associates a constructible cosheaf to a general definable map.

\section{Representations of the Entrance Path Category}
\label{subsubsec:path_cat}

Given a Whitney (or Thom-Mather) stratified space, the entrance path category looks very much like the fundamental groupoid. It has objects that are points and morphisms that are paths. However, the paths and homotopies must respect the stratification. A path may wind around in a given stratum and it may enter deeper levels of the stratification, but upon doing so, it may never return to its higher level.

\begin{figure}
	\centering
	\includegraphics[width=.5\textwidth]{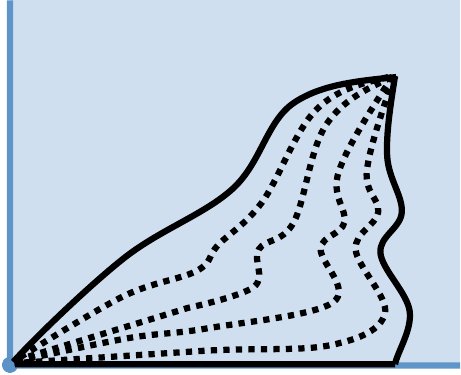}
	\caption{Two Entrance Paths and a Homotopy Between Them}
	\label{fig:entr_path}
\end{figure}

\begin{defn}\index{Entr@$\Entr(X)$ entrance path category}\index{Exit@$\Exit(X)$ exit path category}
	Let $(X,\{X_{\sigma}\}_{\sigma\in P_X})$ be a Whitney (or Thom-Mather) stratified space. We define the \textbf{entrance path category} $\Entr(X,\{X_{\sigma}\})$ to be the category whose objects are points of $X$ and whose morphisms are homotopy classes of entrance paths. An \textbf{entrance path} is a path $\gamma(t)$ whose ambient dimension (the pure dimension of the containing stratum) is non-increasing with $t$. Moreover we require the homotopies $h(s,t)$ to be entrance paths for every fixed $s$. We write $\Entr(X)$ when a given stratification is understood.
	
	Opposite to the entrance path category is the \textbf{exit path category}, written $\Exit(X)$ whose objects are the same, but whose paths ascend into higher dimensional strata.
\end{defn}
\begin{rmk}[``Tame'' Homotopies]
	David Treumann's thesis~\cite{treumann-stacks}, which was written under MacPherson's direction, contains one of the first published accounts of the exit path category. However, he added an additional hypothesis that the homotopies should be ``tame,'' which he defines by saying that $h:[0,1]^2\to X$ should admit a triangulation of $[0,1]^2$ such that the interior of each simplex in the triangulation is contained in some stratum of $X$. Jon Woolf~\cite{woolf} uses a version of the exit and entrance path category based on Quinn's theory of homotopically stratified spaces and does not require Treumann's tameness assumption. Homotopically stratified spaces are more general than Whitney or Thom-Mather stratified spaces, so we may invoke some of Woolf's results. Nevertheless, Treumann's modification foreshadows our own.
\end{rmk}

\begin{defn}[Definable Entrance Path Category]\index{Entr@$\Entr(X)$ entrance path category!definable}\index{definable entrance path category}
	For a fixed analytic-geometric category $\calG$ we can consider the \textbf{definable entrance path category} to have the same objects as before, but whose morphisms are definable entrance paths, where identify entrance paths related by definable homotopies $h:I^2\to X$. There should be a triangulation of $I^2$, so that the image of every open cell is contained in some stratum of $X$. This category will be written $\Entr_{\calG}(X,\{X_{\sigma}\}).$ Dually, we have a definable exit path category $\Exit_{\calG}(X)$.
\end{defn}
\begin{rmk}
	We will not need to use the definable entrance path category until Theorem \ref{thm:strat_maps}, so one may temporarily ignore this restrictive definition. 
\end{rmk}

From the perspective of a computer, the entrance path category definition is entirely too unwieldy to be useful. Storing the points of any space we are accustomed to thinking about (circles, tori, Klein bottles, etc.) is simply too much data to consider. Fortunately, these categories are equivalent to much simpler subcategories by choosing a single point from each connected component in the stratification and passing to a skeletal subcategory.

\begin{ex}[Entrance Path Category for $S^1$]
 Now consider the circle $S^1$ stratified as a single pure stratum. The argument above shows that we can view the entrance path category of $S^1$ as equivalent to the fundamental group $\pijuan(S^1,x_0)$. This is a category with a single object $\star$ whose Hom-set corresponds to a loop for each homotopy class of path, i.e. $\Hom(\star,\star)\cong\ZZ$.
\end{ex}

\begin{ex}[Manifolds]
More generally, if the space $X$ is a manifold, stratified as a single pure stratum, then the entrance path category is equivalent to the fundamental group.
\end{ex}

If we believe MacPherson's characterization of constructible (co)sheaves, then we can reach our much sought after explanation of why cellular sheaves and cosheaves are actually sheaves and cosheaves. Part of the explanation rests on the following characterization of the entrance path category for cell complexes.

\begin{prop}[Entrance Path Category for Cell Complexes]\label{prop:ep_for_cells}\index{Entr@$\Entr(X)$ entrance path category!specializes to $\Cell(X)$}
	If $(X,\{X_{\sigma}\}_{\sigma\in P_X})$ is stratified as a cell complex, then each stratum is contractible and there is only one homotopy class of entrance paths between any two incident cells. As such
	\[
		\Entr(X)\cong\Cell(X)^{op}=P_X^{op} \qquad \mathrm{and} \qquad \Exit(X)\cong\Cell(X)=P_X.
	\]
\end{prop}

To prove this proposition, we need a better understanding of the entrance path category. To do so, we pick out a distinguished class of entrance paths.

\subsection{Homotopy Links}

\begin{defn}[Homotopy Link]\index{homotopy link}
	Suppose $X$ is a decomposed space and $X_{\sigma}\leq X_{\tau}$ are two incident pieces. The \textbf{homotopy link} of $X_{\sigma}$ in $X_{\tau}$ is defined to be the space of paths $\gamma:I\to X_{\sigma}\cup X_{\tau}$ such that $\gamma([0,1))\subset X_{\tau}$ and $\gamma(1)\in X_{\sigma}$, i.e. it is the space of paths that enter $X_{\sigma}$ at the last possible moment.
\end{defn}

We now adapt a proof of Jon Woolf's (\cite{woolf}, Lemma 3.2) to our situation. 

\begin{lem}
	Let $(X,\{X_{\sigma}\})$ be a Thom-Mather stratified space. Any entrance path is homotopic through entrance paths to an element of the homotopy link.
\end{lem}
\begin{proof}
	Suppose $\gamma:[0,1]\to X$ is an entrance path. By compactness, it can only intersect finitely many pieces in the stratification of $X$. We write $X^j$ to denote the union of all dimension $j$ pieces. For any $i\leq j$, we have that $X^i\leq X^j$. 
	
	We claim that one can show that every entrance path $\gamma$ starting in a stratum $X^k$ and ending in a stratum $X^i$ that intersects potentially every stratum in between
	\[
		X^k\geq X^{k-1} \geq \cdots \geq X^i
	\]
	is homotopic to a path $\gamma'(t)$ which sends every $t\in [0,1)$ to $X^k$ and then enters $X^i$ at the last possible moment.
	
	To define the homotopy, one focuses on pulling the path off the last stratum $X^j$ that $\gamma$ enters before entering $X^i$, i.e. there is a partition $0<t_1<\cdots < t_n< 1$ of $[0,1]$ such that $\gamma(0)\in X^k$, $\gamma([t_n,1])\subseteq X^i$ and $\gamma([t_{n-1},t_n))\subseteq X^j$. First, we show that we can pull the path off $X^i$ into $X^j$ so that it enters only at $t=1$. The schematic uses the fundamental observation that stratified spaces can be treated locally as a system of fiber bundles.

Pick a point $x_j:=\gamma(t_n-\epsilon)\in X^j$ and consider its homeomorphic image (which we call $x_j'$) in the fiber over $\gamma(t_n)$. There is a homotopy from the path $\gamma$ relative the end points $x_j=\gamma(t_n-\epsilon)$ and $\gamma(t_n)$ to the piece-wise path that is constant in the fiber, connects $x_j$ to $x_j'$, and then heads straight to $\gamma(t_n)$ while staying in the fiber over that point. By the path lifting property for fiber bundles, we can consider a lift of $\gamma([t_n,1])$ starting with $x_j'$ which ends at $x_j''$ in the fiber over $\gamma(1)$. Repeating the same argument, we can then consider a path that heads from $x''_j$ to $\gamma(1)$ while staying in the fiber. Now the path enters the stratum $X^i$ at the last possible moment.

Repeating this argument and using the conical structure of the fiber, allows us to lift the path out of the $X^j$ stratum and into higher ones. 
\end{proof}

This result allows us to take representative entrance paths that are easy to understand. Every element of the homotopy link is an entrance path, but not every entrance path is an element of the homotopy link. Moreover, it is not clear that two paths that are homotopic as entrance paths are homotopic as entrance paths (after we have moved them into the link as in the above proof). Fortunately, David Miller has recently shown this is the case~\cite{miller-popaths}. At a high level, this provides a proof of Proposition \ref{prop:ep_for_cells}, which can also be seen using easier methods.

\begin{proof}[Proof of \ref{prop:ep_for_cells}]
	Since the pieces in a cell structure on $X$ are all contractible, each cell $X_{\sigma}$ has a single path component in its homotopy link in $X_{\tau}$. Thus the skeleton of the entrance path category for a cell complex is
	\[
		\Entr(X)\cong X^{op},
	\]
which was wanted.
\end{proof}

\subsection{Van Kampen Theorem for Entrance Paths}

\begin{figure}
	\centering
	\includegraphics[width=.8\textwidth]{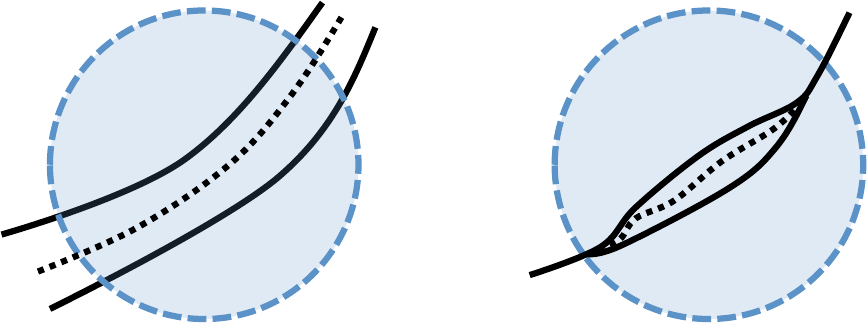}
	\caption{Modifying a Homotopy to Stay Inside an Open Set}
	\label{fig:cover_htpy}
\end{figure}

If we can show that the entrance path category can be built up locally, then we can prove that representations of this category define cosheaves. The ability to build up locally the entrance path category is the van Kampen theorem adapted to stratified spaces. Ostensibly, David Treumann's published version of his thesis~\cite{treumann-stacks} proves the van Kampen theorem for the exit path 2-category, but the elegant inductive argument in proposition 5.9 appears to have an error.\footnote{The argument inducts on the number of triangles in a triangulation of $I^2$. The statement that the closure of the complement of a single triangle is homeomorphic to $I^2$ is not true if, for example, the triangle has two vertices on one side of the square and the third on another side.} Jacob Lurie has a proof for the $\infty$-category case~\cite{lurie-dag6}. Jon Woolf has outlined another argument~\cite{woolf-email} based on his classification of $\Set$-valued representations of the entrance path category as branched covers. The following proof, joint with Dave Lipsky, is more direct and algorithmic, but less elegant in many respects.

The main difficulty in proving the van Kampen theorem is that given a cover, a homotopy of entrance paths restricts to a free homotopy between entrance paths and not a homotopy relative endpoints; see Figure \ref{fig:cover_htpy}. In contrast to the fundamental groupoid, we cannot freely add paths to make this homotopy respect endpoints, the entrance path property must be preserved and this significantly complicates the proof. We borrow Treumann's idea of using a triangulation $T$ of $I^2$ such that $h:I^2\to X$ sends open cells of $T$ to strata of $X$. We then, after sufficient refinement, define a homeomorphism of $I^2$ that allows us to treat the triangulation as a piecewise linear one. For a piecewise linear triangulation, we outline an explicit algorithm for replacing the homotopy $h$ with a composition of homotopies preserving endpoints, each of which is supported on a triangle in $I^2$. Let us state our desired version of the van Kampen theorem and give the first step of the proof.

\begin{thm}[Van Kampen Theorem for Entrance Paths]\label{thm:vkt_ep}\index{van Kampen theorem!for Entr@$\Entr(X)$}\index{Entr@$\Entr(X)$ entrance path category!van Kampen theorem}
	If $X$ is a Whitney stratified space and $\covU=\{U_i\}$ is a cover, then
	\[
		\Entr(X)\cong \varinjlim_{I\in N(\covU)}\Entr(U_I).
	\]
	Each open set is given the induced stratification from the whole space. We assume that every homotopy $h:I^2\to X$ admits a triangulation of the domain so that for each open cell in the triangulation there is a stratum of $X$ that contains its image. Moreover, the same result holds for the definable entrance path category.
\end{thm}
\begin{proof}
	The colimit is an ordinary colimit in the category of all categories. The diagram that sends each $I\in N(\covU)$ to $\Entr(U_I)$ we will call $V$. We already know that the inclusions of the open sets $U_I\hookrightarrow X$ induce functors $\phi_I:\Entr(U_I)\to \Entr(X)$ and that these define a cocone $\phi:V\Rightarrow \Entr(X)$, i.e. a natural transformation from $V$ to the constant diagram on $N(\covU)$ with value $\Entr(X)$. Now suppose $\phi':V\Rightarrow \cat$ is another cocone. We need to check that there exists a unique map $u:\Entr(X)\to\cat$ that makes all the functors commute, i.e. $u\circ \phi_I=\phi_I'$.
	
	On objects, $u(x):=\phi_i'(x)$ for whatever open set $U_i$ contains $x$. The choice doesn't matter since if $U_j$ also contains $x$, then the functor defined on the intersection causes $\phi_i\circ \phi_{ij}(x)=\phi_j\circ\phi_{ij}(x)$. Now we must define $u(\gamma)$ for $\gamma$ an entrance path in $X$. By compactness, we can pass to a finite subcover of $\{U_i\}$ to cover the path $\gamma$. We can break up $\gamma$ into shorter paths $\gamma_{i_1},\ldots,\gamma_{i_n}$, each of which lie in some element of the cover. We define $u(\gamma):=\phi'_{i_n}(\gamma_{i_n})\circ \cdots \circ \phi'_{i_1}(\gamma_{i_1})$. We must show that this definition is invariant under homotopy to complete the proof. This is accomplished by Lemma \ref{lem:vkt_alg} together with Proposition \ref{prop:PL}.
\end{proof}
	
\begin{defn}
	Call a homotopy \textbf{$U$-elementary} if there is an interval $[a,b]\subset I$ such that $h(s,t)$ is independent of $s$ so long as $t\notin [a,b]$ and the image of $I\times [a,b]$ under $h$ is contained in $U$. See Figure \ref{fig:cover_htpy} for an illustrative cartoon.
\end{defn}

\begin{rmk}
	We will use a slightly strange way of orienting the unit square $I^2=[0,1]\times [0,1]\ni (s,t)$. The ``top edge'' is the edge where $t=0$ and the ``bottom edge'' is the edge $t=1$. We will use this language because an entrance path enters ``deeper'' levels of a stratification.
\end{rmk}

\begin{lem}\label{lem:vkt_alg}
	Let $X$ be a Whitney stratified space along with a cover $\covU$ and let $\alpha(t)$ and $\beta(t)$ be entrance paths with the same start and end points. Let $h:I^2\to X$ be a homotopy (relative endpoints) through entrance paths connecting $\alpha(t)=h(0,t)$ to $\beta(t)=h(1,t)$. If $I^2$ admits a piecewise-linear triangulation $T$ such that every open cell in $T$ is mapped to a stratum of $X$, then we may define a sequence of new homotopies $h_1,\ldots, h_n:I^2\to X$, each of which are elementary for some element of the cover, so that the composite connects $\alpha\simeq \beta$. Informally speaking, each homotopy $h_i$ will be supported on a single triangle in the barycentric subdivision of the triangulation $T$.
\end{lem}
\begin{proof}
	Since the image of $I^2$ is compact, a finite subcover of $\covU$ will do. After sufficient refinement, we can assume that each triangle in $T$ is contained in some element of the subcover. By taking the barycentric subdivision $T'$, we can refer to the vertices of any triangle in $T$ via barycentric labels $v,e,f$ depending on whether the vertex is at the barycenter of a vertex, edge or face in the original triangulation. Since each open cell in $T$ is mapped to a stratum of $X$, the triangles in $T'$ satisfy the following fundamental property: $h(\sigma_f)$, where $\sigma_f:=[v,e,f]-[v,e]$, is contained in some stratum $X_{f}$; $h(\sigma_e)$, where $[v,e]-v=\sigma_e$, is contained in $X_{e}$; $h(v)$ is contained in $X_{v}$ and $X_{f}\geq X_{e} \geq X_{v}$. We will refer to the dimension of these containing strata as the ``dimension'' of $f,v$ and $e$, respectively.
	
	By the fundamental property of triangles in $T'$ we know, for example, that the path parameterized by going from $f$ to $e$ to $v$ along the boundary of a triangle is a valid entrance path and this is homotopic through entrance paths to one that goes from $f$ directly to $v$. This is the prototypical ``move'' that we will use to define a given $h_i$ in our new homotopy between $\alpha$ to $\beta$. By reparameterizing the triangle, this move defines an elementary homotopy of entrance paths.
	
	As a preparatory step we replace the entrance path $\alpha(t):=h(0,t)$ with the path that starts at $(s,t)=(1,0)$ and goes along the top edge of the square to $(0,0)$, then to $(0,1)$ and finally to $(1,1)$. Because the homotopy is constant along the top and bottom edges, this only affects the parameterization of the path, but now our modified path and $\beta(t):=h(1,t)$ share the same endpoints in $I^2$. We will now refer to our intermediate paths $\gamma$ by a sequences of vertices in $T'$, written $w_1\cdots w_n$, which taken two at a time define edges $\gamma_i=w_{a} w_{b}$ labeled by a pair of letters $fv$ or $vf$, $fe$ or $ef$, $ev$ or $ve$. Observe that if the image of $vf$ under $h$ is a valid entrance path, then this implies that $\dim v =\dim e = \dim f$ for the triangle containing that particular edge $v_if_i$. 
	
	If an entrance path $\gamma$ ever has a vertex appear twice in its list, then this indicates a loop that must be contained in the same stratum. By virtue of the fact that $I^2$ is simply connected, the portion of the path between the repeated vertices can be homotopically reduced to the constant path via the argument used to prove the van Kampen theorem for the fundamental groupoid. We will avail ourselves of this operation, which we call the \textbf{fundamental groupoid sweep} $F$. For example, if $\gamma$ contains $\cdots e_iv_ie_i\cdots $ in its list of visited vertices, then $F(\gamma)$ will replace the portion $e_iv_ie_i$ with just $e_i$. Of course, $F^2=F$.
	
	We retain the $(s,t)$ coordinates to determine valid moves in our homotopy. We do this because, by assumption, for all $s$, $h(s,t)$ is an entrance path in $t$ and thus the dimension decreases in that direction. Now we can describe our algorithm:
	
	If $F(\gamma)=\gamma$ and $s(w_i)=s(w_j)$ for all $i\neq j$, then we are done. Otherwise, apply $F$ and starting with $\gamma_1$, ask of $\gamma_i$ if there is a triangle to the left (with respect to the induced orientation of following the path) and apply one of following rules:
	\begin{enumerate}
		\item[(1a)] If $\gamma_i=vf$ or $fv$, then replace $\gamma_i$ with $\gamma_i':=vef$ or $fev$, where $e$ belongs to the triangle to the left.
		\item[(1b)] If $\gamma_i=ev$ and $s(e)>s(v)$ or if $\gamma_i=ve$ and $s(v)>s(e)$, then $\gamma_i':=efv$ or $\gamma_i':=vfe$ where $f$ belongs to the triangle to the left.
		\item[(1c)] If $\gamma_i=fe$ and $s(f)\leq s(v)\leq s(e)$ or if $\gamma_i=ef$ and $s(e)\leq s(v) \leq s(f)$ where $v$ belongs to the triangle to the left, then $\gamma_i':=fve$ or $\gamma_i':=evf$.
	\end{enumerate}
	If none of the above apply, consider adjacent paths $\gamma_i\ast\gamma_{i+1}$ two at a time and ask if the following rule is applicable:
	\begin{enumerate}
		\item[(2)] If $\gamma_i\ast \gamma_{i+1}=fev$ or $vef$ where the interior of the triangle is kept to the left, then $(\gamma_i\ast\gamma_{i+1})'=fv$ or $vf$.
	\end{enumerate}
	After each application of a rule, one must check whether $F(\gamma)=\gamma$ and $s(w_i)=s(w_j)$ and repeat as many times as necessary. The algorithm must terminate by virtue of the fact that each step reduces the number of triangles to the left.
	
	\begin{figure}
		\centering
		\includegraphics[width=.5\textwidth]{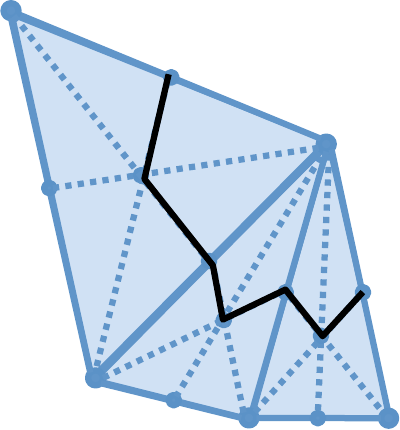}
		\caption{Forcing Move (1c) to Apply}
		\label{fig:vkt_1c}
	\end{figure}
	
	Observe that the only way for a path $\gamma_i$ not to have a triangle to its left is if it lies on the boundary of $I^2$ and it is following the boundary clockwise. If $\gamma_i$ does not belong to the $s=1$ edge, then that contradicts the assumption that $F(\gamma)=\gamma$ as the total path must return to the point $(s,t)=(1,1)$. If the edge does lie on $s=1$, then part of the desired homotopy has been achieved and it need not be moved. If there are no triangles to the left and $F(\gamma)=\gamma$, then the algorithm has finished.
	
	The rationale for rule (1b) is that any point $p$ on the edge $ev$ determines an entrance path $h(s(p),t)$, which drops into the interior $\sigma_f$ of the triangle to the left, thus bounding the dimension of $f$ by the dimension of $e$. The rule (1c) uses similar reasoning. If $s(f)\leq s(v)\leq s(e)$ where $v$ is the triangle to the left, then the entrance path determined by $v$ $h(s(v),t)$ flows into $\sigma_f$ or $\sigma_e$ thus bounding the dimension of $e$ by the dimension of $v$. Let us now prove the correctness of the algorithm.
	
	Suppose $\gamma$ has a $vf$ of $fv$ in sequence. Since a $f$ vertex cannot belong to the boundary of $I^2$, this implies that there is a triangle to the left and that rule (1a) can be applied. Thus, to show that at least one move can be applied up until the algorithm finishes, we assume that no $vf$'s or $fv$'s appear in $\gamma$. Suppose $\gamma$ consists of only $e$'s and $v$'s. Since the start of $\gamma$ has $s(v)=1$, having $s$ non-decreasing would imply that $\gamma$ is contained in $s=1$ and the algorithm would be finished. Otherwise there is a pair such that (1b) can be applied. Now assume that our path has $e$'s, $v$'s and $f$'s, with no $fv/vf$ pairs and such that for all $ev/ve$ pairs, $s$ is increasing.
	
	Because the value of $s$ must go from $1$ back to $1$, if $s$ is not constant along $\gamma$, then there must be at least one $s$ decreasing to increasing turning point. Because $\gamma$ is piecewise linear, by turning point we mean the shortest adjacent collection of edges $\gamma_i\ast\cdots\ast\gamma_{k}$ where the $s$ value goes from strictly decreasing to strictly increasing, i.e. there is an edge along which $s$ is strictly decreasing, then potentially several edges where $s$ is constant and finally an edge which increases in $s$. To determine the ``handedness'' of these turning points we must further specify the $t$ behavior. If the turning point consists of only two edges, then we can ask if the difference in $t$ of the first and last vertex is positive or negative. If the turning point has at least one constant $s$ value edge, then we can use the difference in $t$ along the edge to determine if the turning point is $t$ positive or $t$ negative.
	
	Suppose we have a $t$ negative $s$ decreasing-to-increasing turning point. If the minimal $s$ value is obtained on this turning point, then since $t$ must go from $0$ to $1$, we can conclude that the path must intersect itself at some point, contradicting the assumption that $F(\gamma)=\gamma$. To avoid self-intersection, there must be at least one $t$ positive turning point. Since there are no decreasing $ev/ve$ pairs, the decreasing edge must be either $ef$ or $fe$. If the next term is a $v$, then either rule (1a) or (2) would have to apply, respectively. Thus, we can conclude that the next term is either an $f$ or an $e$. Inducting on the length of the $s$ constant portion of the turning point and using the fact that $f,e$ and $v$ cannot be collinear, we can show so long as rules (1a), (1b) or (2) cannot be applied, that the $s$ increasing edge has to be an $fe/ef$ edge. Consequently the last two edges in such a turning point is either $fef$ or $efe$.
	
	Now we visit the last $t$ positive $s$ decreasing-to-increasing turning point. Since the reasoning is so similar, assume that the last two edges are $efe$. We aim to show that if no other rules are applicable, then the rule (1c) must be applicable for some $ef/fe$ edge. Let us refer to the vertices in the original triangulation $T$ containing these barycenters as $v_1,v_2$ and $v_3$, whose $s$ coordinates are $s_1,s_2$ and $s_3$ respectively. Since all the vertices cannot be collinear, we let $s_2$ have the largest $s$ value. The vertex $f$ is the centroid of $\{v_1,v_2,v_3\}\subset I^2$, the first $e=e_{12}$ is the centroid of $\{v_1,v_2\}$ and the next $e=e_{23}$ is the centroid of $\{v_2,v_3\}$. By assumption $1/3(s_1+s_2+s_3)=s(f)< s(e_{23})=1/2(s_2+s_3)$, thus if the next vertex visited is $v_3$, then we can apply rule (1b) to $e_{23}v_3$. If the next vertex visited is $v_2$, then we can apply rule (2). Thus, we must assume that the next vertex is $f'=1/3(v_2+v_3+v_4)$. Now we reason on $e_{23}f'$. The next vertex cannot be a $v$, otherwise (1a) could be applied. If the next edge visited is $e_{34}$, then there are two possibilities. Either $s(e_{34})< s(f')$, which would contradict the fact that we are at the last $t$ positive turning point, or $s(f')\leq s(e_{34})$, which would imply that $2s_2\leq s_3+s_4$, but this would imply that $s_2=s(v_2)\leq s(f')$ and consequently rule (1c) would apply. Thus, assuming $s(e_{34})<s(f')$, we must remain in the link of $v_2$ and proceed to $e_{24}$. Repeating inductively, and using the fact that there are only finitely many triangles, there must be a point where rule (1c) is applicable; see Figure \ref{fig:vkt_1c}. This completes the proof.
\end{proof}

Now we must select out a class of triangulations of the unit square $I^2$ that can be deformed in an entrance-path preserving way to a PL triangulation.

\begin{prop}\label{prop:PL}
	Suppose that a (definable) triangulation $\varphi:|K|\to I^2$ by a finite simplicial complex is $C^2$ when restricted to the edges of $|K|$. There is a (definable) homeomorphism $g$ of $I^2$ so that after suitable refinement, the triangulation is piecewise-linear.
\end{prop}
\begin{proof}
	The strategy of the proof is to add additional vertices $\{w_i\}$ to the image $\varphi(e)$ of each edge in $I^2$ so that the line segment connecting any two adjacent vertices $w_0,w_1$ is to one side of the curve $\varphi(e)$ between $s(w_0)$ and $s(w_1)$. We will then locally scale the $t$ value in such a way as to push that part of $\varphi(e)$ to the line segment.
	
	To add in these vertices in a principled way, we first consider the critical set of the $s$ value of $\varphi(e)$ for every edge $e$ in $|K|$. We remove the entire critical set from $e$. If the critical set contains an interval, then we know that portion of the edge is already linear and need not consider it. Now we use the implicit function theorem to write the remainder of the edge $\varphi(e)-\{ds|_{\varphi(e)}=0\}$ as a function of $s$. For each of these functions we find the critical set of its first derivative (``inflection points'') and remove these as well. What is remaining of $\varphi(e)$ is a collection of open concave and convex arcs, each of which have boundary points $w_i,w_{i+1}$ in the various critical sets we have removed. Write $\ell_i(s)$ for the equation of the $t$ coordinate of the line connecting $w_i$ to $w_{i+1}$, i.e. the graph of the line is $(s,\ell_i(s))$. We also write the portion of $\varphi(e)$ between $w_i$ and $w_{i+1}$ as $\varphi_i(s)$.
	
	Possibly after further removal of points, we assume that each arc $\varphi_i(s)$ has a tubular neighborhood $T_i$ that contains $\ell_i(s)$ and each of these neighborhoods are pairwise disjoint. We are now going to define a homeomorphism that is the identity outside of $T_i$. To do so we need one more pair of functions.	
	\[
		\begin{array}{r c l} 
			\mu_{\pm}(x)=x & \mathrm{if} & |x|\geq 1 \\[\medskipamount]
						\mu_{\pm}(x)=2x\pm 1 & \mathrm{if} & -1 \leq x\leq \frac{-1}{2} \\[\medskipamount]
						\mu_{\pm}(x)=\frac{2}{3}x\pm \frac{1}{3} & \mathrm{if} & \frac{-1}{2} \leq 1
		\end{array}
	\]
	Now we can define a homeomorphism $g_i$ on $T_i$ using $+$ if the function $\varphi_i(s)$ is concave and $-$ if the function $\varphi_i(s)$ is convex.
	\[
		g_i(s,t):=(s,2\cdot|\varphi_i(s)-\ell_i(s)|\cdot\mu_{\pm}(\frac{t-\ell_i(s)}{2|\varphi_i(s)-\ell_i(s)|})+\ell_i(s))
	\]
	Since the domains of each $T_i$ are disjoint we can define the homeomorphism $g$ to be $g_i$ when in $T_i$ and the identity otherwise. This straightens each of the $\varphi_i(s)$. The portions of $\varphi(e)$ removed are already piecewise-linear. This makes each $g\circ \varphi(\bar{\sigma})$ into a piecewise-linear polyhedron, the interior of which is mapped via $h$ to a single stratum of $X$. By adding edges and vertices, we can refine the stratification of $g\circ \varphi(\bar{\sigma})$ to be a piece-wise linear triangulation.
\end{proof}	

\subsection{The Equivalence}
\label{subsec:equivalence_cosheaves}

We will break our proof of MacPherson's characterization into two parts. The first shows that any representation defines a constructible cosheaf. The second shows that any constructible cosheaf defines a representation of the entrance path category.

\begin{thm}[Representations are Cosheaves]\index{Entr@$\Entr(X)$ entrance path category!representations are cosheaves}\index{representation!of Entr@$\Entr(X)$}
	Let $X$ be a Thom-Mather stratified space. Any representation of the entrance path category 
	\[
		\Entr(X,\{X_{\alpha}\}) \to \Vect
	\]
	defines a constructible cosheaf.
\end{thm}
\begin{proof}
	To produce a cosheaf from a representation $\hF:\Entr(X)\to\Vect$ we take colimits over the restriction of $\hF$ to the entrance path category of $U$ (with its induced stratification):
	\[
	\hF(U):= \varinjlim_{\Entr(U)} \hF|_U
	\]
	This is clearly a pre-cosheaf since if $U\hookrightarrow V$, the colimit over $\hF|_V$ defines by restriction a cocone over $\hF|_U$ and thus a unique map $\hF(U)\to\hF(V)$. We can describe more explicitly this colimit as follows:
 
	Given a point $x\in X_{\sigma}$ in a stratum of dimension $i$, there is a basis of conical neighborhoods $U_x\cong \RR^i\times C(L)$ where $L$ is the stratified fiber of the retraction map $\pi_{\sigma}$ and $C(L)$ is its open cone. For such a neighborhood, $x$ is the terminal object in $\Entr(U_x)$, thus the colimit returns the value of $\hF(x)$. Moreover, this shows that the costalks of the pre-cosheaf defined stabilize for small contractible sets containing $x$.
	
	To show this is actually a cosheaf we use the version of the van Kampen theorem adapted to the entrance path category just proved in Theorem \ref{thm:vkt_ep}:
	\[
		\varinjlim_{I\in N(\covU)} \Entr(U_I) \cong \Entr(X)
	\]
	As a consequence of colimits commuting with colimits we get that for a representation of the entrance path category $\hF$
	\[
		 \varinjlim_{I\in N(\covU)} \hF(U_I):=\varinjlim_{I\in N(\covU)} \varinjlim_{\Entr(U_I)} \hF|_{U_I} \cong  \varinjlim_{\Entr(U_I)} \varinjlim_{I\in N(\covU)} \hF|_{U_I} \cong \varinjlim_{\Entr(U)} \hF|_U.
	\]
	This establishes the cosheaf axiom.
\end{proof}

\begin{thm}[Representations of the Entrance Path Category]\index{cosheaf!constructible!is representation of entr@$\Entr(X)$}
	Every cosheaf $\hF$ with finite-dimensional costalks that is constructible with respect to a Thom-Mather stratification of $X$ determines a representation of the entrance path category.
	\[
		\Entr(X,\{X_{\alpha}\}) \to \vect
	\]
\end{thm}
\begin{proof}
	If $X$ is a Thom-Mather stratified space, then we know that every point $x\in X_{\sigma}$ in a stratum of dimension $i$ has a neighborhood $U_x\cong \RR^i\times C(L)$, where $L$ is the fiber of $\pi_{\sigma}$, and $C(L)$ is the open cone. Now suppose $\hF$ is a constructible cosheaf, which we assume has finite-dimensional costalks. We claim that for $U_x$ suitably small we can show that
	\[
		\hF_x\cong \hF(U_x).
	\]
	This is not so easy to see and a proof would require substantial more development of cosheaf theory. Heuristically, if the value of $\hF$ on a sequence of conical neighborhoods never stabilized then this would contradict the constancy of the cosheaf on sets of the form $U_x\cap X_{\tau}$. For a rigorous proof, one dualizes a constructible cosheaf to a constructible sheaf by post-composing with $\Hom_{\vect}(-,k)$, which is an equivalence, and we can apply the proof for constructible sheaves found on p. 84 of~\cite{gm-ih2}. 
	
	Consequently, if $y\in U_x\cap X_{\tau}$ is a point in a nearby stratum, then there is an analogous neighborhood $U_y$ contained in $U_x$. Repeating the same argument, we can then use the maps present in a cosheaf to define a map from the costalk at $y$ to the costalk at $x$:
	\[
		\xymatrix{\hF_x \ar[r]^-{\cong} & \hF(U_x) & \ar[l] \hF(U_y) & \ar[l]_-{\cong} \hF_y}
	\]
	Recall that the restriction of a constructible cosheaf to any stratum defines a locally constant cosheaf. For arbitrary points $y'$ in the stratum $X_{\tau}$ we can consider a path $y'\rightsquigarrow y$ and use Theorem \ref{thm:lc_sheaf_rep} to define a map from $\hF_{y'}$ to $\hF_{y}$. Postcomposing with the above map defines the map $\hF_{y'}\to \hF_{x}$. This explains why constructible cosheaves map naturally define ways of specializing a costalk over one stratum to a costalk in its frontier.

	To show homotopy invariance, we appeal to the van Kampen Theorem \ref{thm:vkt_ep} to reduce the argument to elementary homotopies of a particular form. Assume $\alpha(t)$ goes from $z\in X_{\lambda}$ directly to $x\in X_{\sigma}$ and that $\beta(t)$ goes from $z$ to $y$ and then $x$. Since restriction to any stratum defines a locally constant cosheaf, we can appeal to the homotopy invariance of Theorem \ref{thm:lc_sheaf_rep} to position these paths and points to be inside $U_x$ and so that $U_x\supset U_y\supset U_z$. By choosing a similar set of isomorphisms, we get two commutative diagrams, which followed along the top edge corresponds to the action of $\alpha(t)$ and followed along the bottom edges corresponds to $\beta(t)$.
 	\[
		\xymatrix{\hF(U_z) \ar[rr] \ar[rd] & & \hF(U_x) \\ & \hF(U_y) \ar[ru] &} \qquad \qquad \qquad 
		\xymatrix{\hF_{z} \ar[rr] \ar[rd] & & \hF_{y} \\ & \hF_{x} \ar[ru] &}
	\]
\end{proof}

\subsection{Representations from Stratified Maps}

We want to show that stratified maps induce representations of the entrance path category, which, by the first part of our equivalence, defines a constructible cosheaf.

\begin{thm}[Cosheaves from Stratified Maps]\label{thm:strat_maps}\index{cosheaf!constructible!from a definable map}\index{stratified map!determines constructible cosheaves}\index{representation!of entr@$\Entr(X)$! from a stratified map}
	Fix an analytic-geometric category $\calG$. If $Y$ is a closed set in $\calG(N)$ and $f:(Y,N)\to (X,M)$ is a $C^1$ proper definable map, then for each $i$, the assignment
	\[
		x\in X \rightsquigarrow H_i(f^{-1}(x);k)
 	\]
defines a representation of the definable entrance path category of $X$, where the stratification is gotten by the stratification induced by $f$.
\end{thm}
\begin{proof}
	Let $\gamma:I\to X$ be a definable map that satisfies the entrance path condition, i.e. as $t$ increases the dimension of the ambient stratum is non-increasing. Thus $\gamma(0)$ is in a stratum of dimension greater than or equal to $\gamma(1)$. By Lemma \ref{lem:def_pb}, we know that the pullback $Y_{\gamma}:=I\times_X Y$ is definable, as is the pullback of $f$ along $\gamma$, written $\gamma^*f$. Since definable sets can be Whitney stratified, $Y_{\gamma}$ admits a system of control data, and may be regarded as a Thom-Mather stratified space.
	
	The argument from Lemma \ref{lem:strat_maps_1d} provides us with the prototype for getting a diagram of spaces for every path. We will repeat it here for convenience and make any necessary modifications. By definable triviality (4.11 of~\cite{vdd-geocat}), there exists a finite partition of $[0,1]$ such that over each interval the inverse image is homeomorphic to the product:
	\[
		\xymatrix{f^{-1}((t_i,t_{i+1})) \ar[rr]^{\cong} \ar[rd]_f & & F\times (t_i,t_{i+1}) \ar[ld] \\ & (t_i,t_{i+1}) &}
	\]
	By properness we can, for any fixed $\epsilon>0$, find an $s_i^+$ such that $f^{-1}([t_i,s_i^+])$ is contained in $U_{i}(\epsilon):=\cup T_{\sigma}(\epsilon/2)$ for $Y_{\sigma}\subseteq f^{-1}(t_{i}):=Y_{i}$. The retraction we constructed in Proposition \ref{prop:reg_nbhd} gives a retraction map $r_i^+:=H(p,0)$ from $U_{i}(\epsilon)\to Y_{i}$. This allows us to define a map on fibers
	\[
		f^{-1}(s_i^+) \hookrightarrow U_{i}(\epsilon) \to Y_{i}.
	\]
	Applying some homology functor $H_n(-;k)$ defines the representation locally on the path. Of course, we must show that this representation is independent of the point $s_i$ taken. If $s_i^{+'}\in [t_{i},s_i^+]$ is another point, then the composition of the trivialization with the retraction witnesses the homotopy between these two choices.
	\[
	 F\times [s_i^{+'},s_i^+]\cong f^{-1}([s_i^{+'},s_i^{+}])\hookrightarrow U_{i}(\epsilon) \to Y_{i}
	\]
	Similarly, one can find a point $s_{i+1}^-$ so that it's fiber is contained in $U_{i+1}(\epsilon)$ and the retraction $r_{i+1}^-$ defines a map from that fiber to the fiber $Y_{i+1}=f^{-1}(t_{i+1})$. 
By Thom's first isotopy lemma there is a homeomorphism $\varphi_{i+1,i}$ taking the fiber over $s_{i}^+$ to the fiber over $s_{i+1}^-$. This homeomorphism is gotten by constructing a vector field that flows from $s_i^+$ to $s_{i+1}^-$ and lifting it to a controlled vector field on $f^{-1}((t_i,t_{i+1}))$ via Proposition 9.1 of~\cite{mather}. 
	Finally, one must observe that the filtration of $X$ by strata of a given dimension or less, the restriction of $\gamma$ to the half-open interval $[t_i,t_{i+1})$ is contained inside a single stratum of $X$ and thus the retraction $r_i^+$ induces a homotopy equivalence between the fiber over $s_i^+$ and the fiber over $t_i$. Applying our homology functor to the following composition defines the total action associated to this path:
	\[
		\cdots (r_{i+1})_*^-\circ (\varphi_{i+1,i})_* \circ (r_{i}^+)_*^{-1} \cdots
	\]

	\begin{figure}
		\centering
		\includegraphics[width=.8\textwidth]{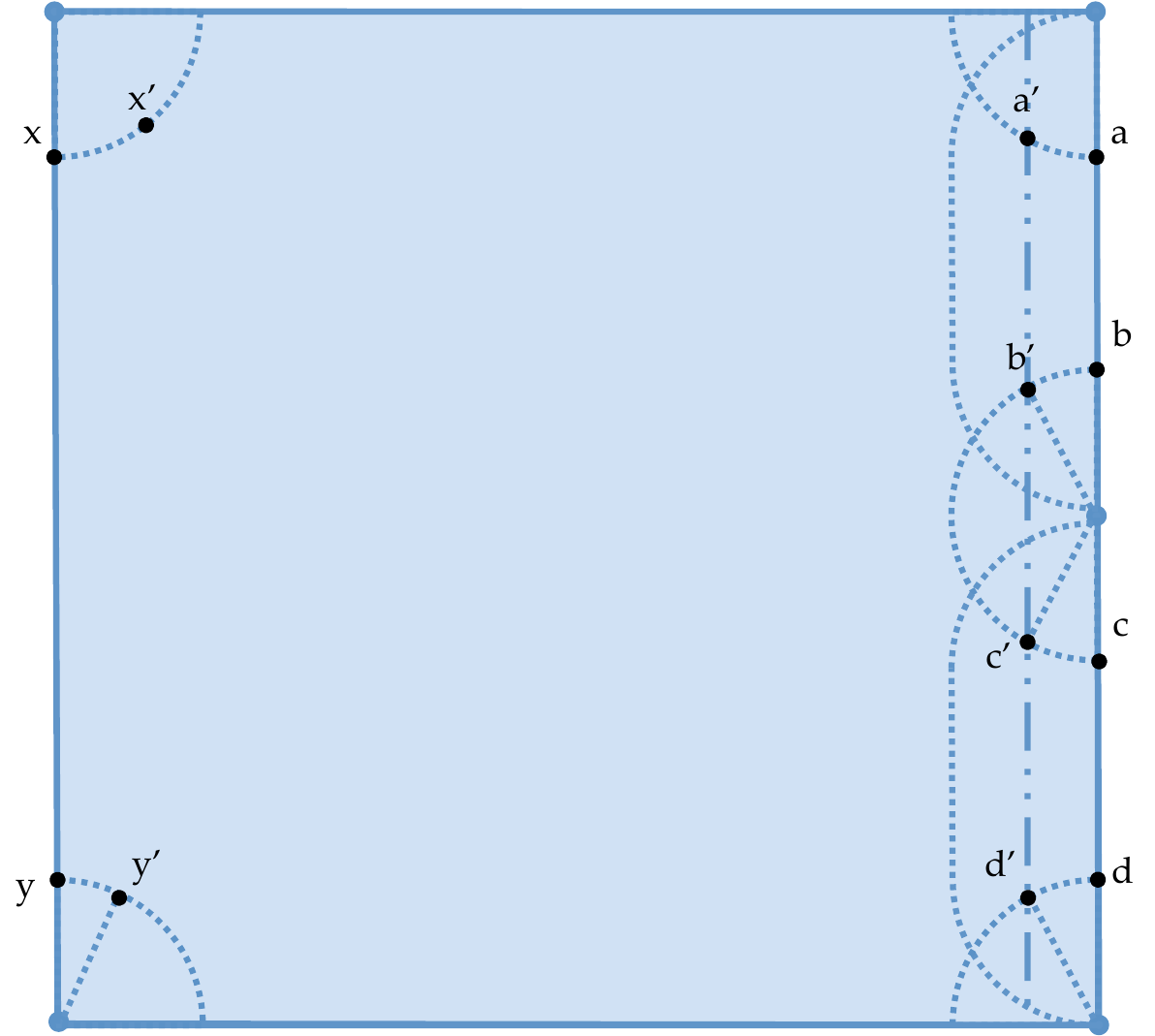}
		\caption{Argument for Homotopy Invariance}
		\label{fig:strat_maps_htpy}
	\end{figure}
	
	It remains to be seen that this map is invariant under definable homotopies of entrance paths. Suppose $h:I\times I\to X$ is a definable homotopy. Again, the pullback $Y_h:=I^2\times_X Y$ is definable, as is the map $h^*f$, and both can be stratified. Thus, we have reduced everything to considering a stratified map to the square $I^2$. By the van Kampen Theorem \ref{thm:vkt_ep}, it suffices to check homotopy invariance on an elementary homotopy, such as the one depicted in Figure \ref{fig:entr_path}. Let us assume that $h$ is a homotopy between an entrance path $\alpha(t)=h(0,t)$, which goes from a stratum $X_{\lambda}$ and enters a stratum $X_{\sigma}$ at the last possible moment $t=1$, and an entrance path $\beta(t)=h(1,t)$, which enters $X_{\tau}$ at $t=1/2$ and then goes to $X_{\sigma}$ at $t=1$. Moreover, we assume that $h$ takes the complement of $\{t=1\}\cup \{(1,t)\,|\, t\geq 1/2\}$ to the stratum $X_{\lambda}$. This guarantees that the fibers over $x,x',y,y',a,a',b,b',c'$ and $d'$ in Figure \ref{fig:strat_maps_htpy} can all be identified.

Let ${T}$ be a system of control data for $Y_h$, obtained in a specific way. By restricting to the strata over $s=0$ and $s=1$ respectively, we get control data for $Y_{s=0}$ and $Y_{s=1}$, both of which are inside $I^2\times_X Y\subset \RR^2\times N$. The spaces $Y_{s=0}$ and $Y_{s=1}$ can be identified with the inclusions of $Y_{\alpha}$ or $Y_{\beta}$, which are contained in $\RR\times N$. The manner in which Mather constructs control data in Proposition 7.1 of~\cite{mather} can be used to extend the control data for $Y_{\alpha}$ and $Y_{\beta}$ to control data for $Y_{s=0}$ and $Y_{s=1}$ inside $\RR^2\times N$ respectively. This is how we obtain those tubular neighborhoods in $\{T\}$ and the rest can be constructed to be compatible with those. This allows us to use the control data ${T}$ to meaningfully compare the construction above for $\alpha(t)$ and $\beta(t)$.
 	
	We can describe the maps associated to $\alpha(t)$ and $\beta(t)$ as follows: By properness, we assume the fiber over $x$ is contained in a regular neighborhood, which retracts via $r_x$ to the fiber over $(0,0)$. There is a homeomorphism $\varphi_{y,x}$ from the fiber over $x$ to the fiber over $y$. Finally, we can assume that the fiber over $y$ retracts via $r_y$ to the fiber over $(0,1)$. Thus the action associated to $\alpha(t)$ is the map
	\[
		(r_y)_*\circ(\varphi_{y,x})_*\circ (r_x)^{-1}_*:H_n(Y_{(0,0)})\to H_n(Y_{(0,1)})
	\]
	where we have implicitly pre-composed $r_x$ with the inclusion of the fiber.

	For $\beta(t)$, the action is similar:
	\[
		(r_d)_*\circ (\varphi_{d,c})_*\circ (r_c)_*^{-1}\circ (r_b)_*\circ (\varphi_{b,a})_*\circ (r_a)^{-1}_*: H_n(Y_{(1,0)})\to H_n(Y_{(1,1/2)}) \to H_n(Y_{(1,1)})
	\] 
	
	The strategy of the proof is to pick a path $\gamma(t)$ that interpolates $\alpha(t)$ and $\beta(t)$ and show that the associated map on homology agrees with both $\alpha(t)$ and $\beta(t)$. This path is indicated by the dotted-and-dashed line passing through $a',b',c'$ and $d'$ in Figure \ref{fig:strat_maps_htpy}. The representation associated to $\gamma(t)$ is
	\[
		(r_{1})_*\circ (\varphi_{d',c'})_*\circ (i_{c'})_*^{-1} \circ (i_{b'})_* \circ (\varphi_{b',a'})_* \circ (r_{0})^{-1}_*.
	\]
        Here the maps $i_{b'}$ and $i_{c'}$ denote the inclusion of $Y_{b'}$ and $Y_{c'}$ into the inverse image of the interval $[b',c']$. The action on homology of $(i_{c'})_*^{-1} \circ (i_{b'})_*$ agrees with an analogously constructed homeomorphism $\varphi_{c',b'}$, but we will find it easier to equate the map associated to $\beta(t)$ and $\gamma(t)$ as written above.

        Because the control data $\{T\}$ extends the control data for $Y_{\alpha}$ and $Y_{\beta}$, the retraction map $r_a$ can be taken to be the restriction of a retraction map $r_{0}:U_{Y_{(1,0)}}(\epsilon)\to Y_{(1,0)}$ constructed in Proposition \ref{prop:reg_nbhd}. This in turn can be taken to be the restriction of the tubular projections used to define a retraction map $r_{s=1}:U_{s=1}(\epsilon)\to Y_{s=1}$. The commutation relations for control data allow us to imagine first taking the fiber $Y_{a'}$ over $a'$ and retracting to the strata over the edge $e_0:=\{(1,t)\,|\, 0<t<1/2\}$, and then retracting to the fiber over $(1,0)$. This allows us to factor $r_{0}$ as
\[
r_{0}=r_{a}\circ r_{e_0},
\]
but the image of $Y_{a'}$ under $r_{e_0}$ may not be contained in $Y_{a}$ or any single fiber. This would be true if, for example, the control data defining the retraction to $Y_{e_0}$ satisfied
\[
\pi_{e_0}(f(p))=f(\pi_{\sigma}(p))
\]
for each stratum $Y_{\sigma}$ that $f$ carried to $e_0$, but in general it does not. This is what necessitates the use of the Thom properties given by Lemma \ref{lem:thom_codim} and property (b) of Proposition \ref{prop:thom_data}.

By Lemma \ref{lem:thom_codim}, we know that restricting the codomain to the complement of the vertices, the mapping $h^*f$ is a Thom mapping. Consequently, if we pick a tubular neighborhood $T_{e_{0}}$ for the edge $e_0:=\{(1,t)\,|\, 0<t<1/2\}$, there exists a system of control data $\{T'\}$ over $T_{e_0}$ and the interior of $I^2$ by Proposition \ref{prop:thom_data}. If we restrict to those tubular neighborhoods coming from strata in $Y_{e_0}$, then property (a) of Proposition \ref{prop:thom_data} implies that this restricted collection of tubular neighborhoods defines actual control data for $Y_{e_0}$, which we call $\{T'\}_{Y_{e_0}}$. A priori, the analogous restriction of $\{T\}$ to $Y_{e_0}$ defines a different system of control data. However, by Mather's uniqueness result,\footnote{Mather mentions at the bottom of page 492 of~\cite{mather} that any Whitney stratified subset $Z$ of a manifold $M$ has a unique, up to isomorphism, structure as a Thom-Mather stratified set. This is not explicitly proved, but it follows from Mather's Corollary 10.3 as explained by Goresky: Suppose $Z$ is given two different structures of control data $\{T\}$ and $\{T'\}$. If we consider $Z\times\RR$ as a Whitney stratified subset of $M\times \RR$, then $\{T\}$ and $\{T'\}$ can be extended to control data on $Z\times(-\epsilon,\epsilon)$ and $Z\times (1-\epsilon,1+\epsilon)$, respectively. Then, using the proof of prop. 7.1, one can find control data on all of $Z\times\RR$ that agrees with the $\epsilon$ extensions of $\{T\}$ and $\{T'\}$. This space, now viewed as a Thom-Mather stratified set, is then isomorphic via Corollary 10.3 to the set where just $Z\times\RR$ is given the extension of just the control data $\{T\}$.} there is a homeomorphism $\psi_{e_0}$ of $Y_{e_0}$ that takes $\{T'\}_{Y_{e_0}}$ to $\{T\}_{Y_{e_0}}$. This implies that if $Y_{\sigma}$ is a stratum that is mapped to $(1,0)$ and $Y_{\tau}$ is mapped to $e_0$, then
\[
\pi_{\sigma}=\pi_{\sigma}\circ\psi_{e_0}\circ \pi_{\tau'} \qquad \mathrm{since} \qquad \pi_{\tau}=\psi_{e_0}\circ \pi_{\sigma'}
\]
By repeating the construction of a retraction outlined in Proposition \ref{prop:reg_nbhd}, but using the control data $\{T'\}_{Y_{e_0}}$ instead to construct the family of lines, we get a map $r'_{e_0}$ that carries the fiber over $a'$ to the fiber over $a$. Post-composing $r'_{e_0}$ with $\psi_{e_0}$ gives the equality $r_{e_0}=\psi_{e_0}\circ r'_{e_0}$. This construction gives the left most triangle in the following commutative diagram:
\[
\xymatrix{Y_{(1,0)} & \ar[l]_{r_a \psi_{e_0}} Y_a \ar[r]^{\varphi_{b,a}} & Y_b \ar[r]^{r_b\psi_{e_0}} & Y_{(1,1/2)} & \ar[l]_{r_c\psi_{e_1}} Y_c \ar[r]^{\varphi_{d,c}} & Y_d \ar[r]^{r_d\psi_{e_1}} & Y_{(1,1)} \\
& \ar[lu]^{r_{0}} Y_{a'} \ar[u]_{r'_{e_0}} \ar[r]_{\varphi_{b',a'}} & Y_{b'} \ar[u]^{r'_{e_0}} \ar[ur]_{r_{1/2}} \ar[r]_{i_{b'}} & Y_{[b',c']} \ar[u]_{r_{1/2}} & \ar[l]^{i_{c'}} \ar[lu]_{r_{1/2}} Y_{c'} \ar[u]_{r'_{e_1}} \ar[r]_{\varphi_{d',c'}} & Y_{d'} \ar[u]^{r'_{e_1}} \ar[ur]_{r_{1}} & }
\]

Now we explain the other maps in this diagram. The homeomorphisms $\varphi_{b,a}$ and $\varphi_{b',a'}$ are constructed by taking a controlled vector field $\{\eta_f,\eta_{e_0},\eta_{e_1}\}$ in $I^2$ minus the vertices using the control data $\{T'\}$ over $\{T_f,T_{e_0},T_{e_1}\}$. Since $d\pi_{e_0}(\eta_f(s,t))=\eta_{e_0}(\pi_{e_0}(s,t))$ the controlled vector field over this one commutes with $f$ and gives
\[
\varphi_{b,a}\circ r_{e'_0}=r_{e_0'}\circ\varphi_{b',a'}.
\]

Again, the commutation relations in Proposition \ref{prop:reg_nbhd} allows us to, using the control data $\{T\}$ to factor $r_{1/2}=r_{b}\circ r_{e_0}$. However, the uniqueness theorem tells us that $r_{e_0}=\psi_{e_o}\circ r'_{e_0}$ where $\psi_{e_0}$. A simple diagram chase now completes the argument. Comparing the maps associated to $\alpha(t)$ and $\gamma(t)$ is much simpler and uses the same ideas. We leave it to the reader. 
\end{proof}

\begin{rmk}[Alternative Idea for a Proof]
	An alternative approach makes use of the properties of o-minimal structures. The generic triviality theorem 4.11 of~\cite{vdd-geocat} guarantees that we have a definable trivialization of the map over $(0,1)$.
	\[
		\xymatrix{f^{-1}(\gamma((0,1))) \ar[rr]^{\cong}_h \ar[rd]_{\gamma^*f} & & F\times (0,1) \ar[ld] \\ & (0,1) & }
	\]
	For each point $x$ in the fiber $F$, we get a lift $\{x\}\times (0,1)$ of the open interval. Applying the inverse homeomorphism, $h^{-1}(\{x\}\times (0,1))$ defines a definable path $\alpha_x:(0,1)\to f^{-1}(\gamma([0,1]))$. 
	
	M\'ario Edmundo and Luca Prelli, in their recent note~\cite{edmundo-six} reworking the six basic Grothendieck operations for sheaves in the o-minimal setting, have given a tantalizing reformulation of what characterizes a definable proper map. They use an idea of Ya'acov Peterzil and Charles Steinhorn~\cite{peterzil-cpt} that shows that being definably compact (equivalently, closed and bounded) is equivalent to being able to to complete curves. A map $f:Y\to X$ is definably proper if for every definable curve $\alpha:(0,1)\to Y$ and every definable map $[0,1]\to X$ there is at least one way to complete the diagram:
	\[
		\xymatrix{ (0,1) \ar[r]^{\alpha} \ar[d] & Y \ar[d]^{f} \\
		[0,1] \ar[r] \ar[ur]_{\bar{\alpha}} & X}
	\]
	If one assumes all the maps are continuous as well as definable then the completion in the diagram above is unique.\footnote{One of the unusual features of o-minimal topology is that definable maps need not always be continuous, thus the added hypothesis. Even discontinuous maps can have triangulable graphs.}
	
	In our situation, the hypotheses guarantee that for each point $x\in F$, we can complete the curve $\alpha_x:(0,1)\to I\times_X Y$ to a curve $\bar{\alpha}_x:[0,1]$. By associating endpoints over $0$ to endpoints over $1$ we define a set-theoretic map $g:f^{-1}(\gamma(0))\to f^{-1}(\gamma(1))$. The hard work is showing that this map $g$ is continuous and is invariant under homotopy.
\end{rmk}

%
%

\chapter{Duality: Exchange of Sheaves and Cosheaves}
\label{sec:duality}

\begin{quote}
{\em ``Quod superius sicut quod inferius, et quod inferius sicut quod superius, ad perpetranda miracula rei unius.''}
\begin{flushright} --- Hermes Trismegistus~\cite[p. xxix]{vaughan1888magical} \end{flushright}

\end{quote}

In this section, we are concerned with the derived equivalence of cellular sheaves and cosheaves. In Section \ref{subsec:closures}, we introduce the functor that establishes this equivalence and try to motivate it topologically via taking the ``closure'' of the data over an open cell. In the case when $X$ is a manifold, Theorem \ref{thm:mfld_duality} gives us a duality result for data that relates sheaf cohomology with our new theory of sheaf homology. Finally, the equivalence is proved in Section \ref{subsec:derived_equiv}.

\section{Taking Closures and Classical Dualities Re-Obtained}
\label{subsec:closures}

In this section we are going to explain the all-important Poincar\'e-Verdier duality as an exchange of sheaves and cosheaves. To introduce this duality, we explain an odd, but clean way of going from a cellular sheaf to a complex of cellular cosheaves. This is meant to express the idea that duality is an exchange of open and closed cells.

Suppose we start with a sheaf $F$ on the unit interval $X=[0,1]$ stratified with end points $x=0$, $y=1$, and $a=(0,1)$. Such a sheaf is just a diagram of vector spaces of the form
\[
	\xymatrix{& F(a) & \\ F(x) \ar[ur]^{\rho_{a,x}} & & \ar[lu]_{\rho_{a,y}} F(y) .}
\]
Now we are going to extend the value of the sheaf on a cell $\sigma$ to its closure $\bar{\sigma}$ by defining $\hF(\tau)=F(\sigma)$ for every cell $\tau\leq \sigma$ and using the identity maps from $\sigma$ to its faces. This in effect smears the value of the sheaf on an open cell onto all of its faces. However, what should we do to the values of the sheaf already stored on a face $\tau$? This is where we use the different slots in a complex of vector space to store independently the values:
\[
\xymatrix{F(a) & \ar[l]_{\id} F(a) \ar[r]^{\id} & F(a) \\ F(x) \ar[u]^{\rho_{a,x}} & \ar[l] 0 \ar[u] \ar[r] & \ar[u]_{\rho_{a,y}} F(y) .}
\]
For dimension reasons, it should be clear that this smearing operation defines an assignment of chain complexes to each open cell with chain maps extending to the faces:
\[
	\xymatrix{& \ar[ld]_{r^{\bullet}_{x,a}} \equivp(F)(a) \ar[rd]^{r^{\bullet}_{y,a}} & \\ \equivp(F)(x) & & \equivp(F)(y) .}
\]
This motivates the following general definition of a functor $\equivp$: to a cellular sheaf $F\in \Shv(X)$ we associate the following cosheaf of chain complexes $\equivp(F)$
\[
	\equivp(F): \qquad \sigma \qquad \rightsquigarrow \qquad F(\sigma) \to \bigoplus_{\sigma\leq_1\tau} F(\tau) \to \bigoplus_{\sigma\leq_2\gamma} F(\gamma) \to \cdots.
\]
where $F(\sigma)$ is placed in cohomological degree $\dim|\sigma|$ or homological degree $-\dim|\sigma|$
However, in order for this to be a chain complex, following two arrows in sequence should give zero. In order to guarantee this we need to use the fact that $X$ is a cell complex, and as such for any pair of cells $\sigma\leq_2\gamma$ differing in dimension by two, there are precisely two ways $\tau_1,\tau_2$ of going between $\sigma$ and $\gamma$. Using the signed incidence relations $[\sigma:\tau_i]$ and the restrictions maps internal to $F$ allows us to define the differentials in this complex by $d^{i+1}:=\oplus [\tau:\gamma]\rho^F_{\gamma,\tau}$. Now let's consider a cell $\lambda$ that is a codimension one face of $\sigma$, then the extension map $r^{\bullet}_{\lambda,\sigma}$ is defined to be the chain map
\[
	\xymatrix{0 \ar[r] \ar[d]_{r^{i-1}_{\lambda,\tau}} & F(\sigma) \ar[r]^{d^i} \ar[d]_{r^{i}_{\lambda,\tau}} & \bigoplus_{\sigma\leq_1\tau} F(\tau) \ar[r]^{d^{i+1}} \ar[d]_{r^{i+1}_{\lambda,\tau}} & \bigoplus_{\sigma\leq_2\gamma} F(\gamma) \ar[d]_{r^{i+2}_{\lambda,\tau}} \\
	F(\lambda) \ar[r]_{d^{i-1}} & \bigoplus_{\lambda\leq_1\sigma} F(\sigma) \ar[r]_{d^i} & \bigoplus_{\lambda\leq_2\tau} F(\tau) \ar[r]_{d^{i+1}} & \bigoplus_{\lambda\leq_3\gamma} F(\gamma) }.
\]
The reason it is a chain map is clear from the fact that if $\lambda\leq\sigma$ then $U_{\sigma}\subset U_{\lambda}$ and so the chain complex $\equivp(F)(\sigma)$ simply includes term by term into the chain complex $\equivp(F)(\lambda)$.

Although the idea of a cosheaf of chain complexes is perhaps easier to visualize, for actual algebraic manipulation, one uses a chain complex of cosheaves to express the same idea in a different way. 

\begin{defn}[Poincar\'e-Verdier Equivalence Functor]\index{duality!Poincar\'e-Verdier equivalence functor}
	Let $X$ be a cell complex and let $\Shv(X)$ and $\Coshv(X)$ denote the categories of cellular sheaves and cosheaves respectively. We define the \textbf{Poincar\'e-Verdier Equivalence Functor} $\equivp:D^b(\Shv(X))\to D^b(\Coshv(X))$ by the following formula: to a sheaf $F\in \Shv(X)$ we associate the following complex of projective co-sheaves, the cohomological degree corresponding to the dimension of the cell:
	\[
		\xymatrix{\cdots \ar[r] & \bigoplus_{\sigma^i\in X}[\hat{\sigma}^i]^{F(\sigma^i)} \ar[r]^-{[\sigma:\gamma]\rho^F} & \bigoplus_{\gamma^{i+1}\in X}[\hat{\gamma}^{i+1}]^{F(\gamma^{i+1})} \ar[r]^-{[\gamma:\tau]\rho^F} & \bigoplus_{\tau^{i+2}\in X}[\hat{\tau}^{i+2}]^{F(\tau^{i+2})} \ar[r] & \cdots}
	\]
	Here $\sigma^i$ denotes the $i$-cells and $[\sigma^i:\gamma^{i+1}]=\{0,\pm 1\}$ records whether the cells are incident and whether orientations agree or disagree. The maps in between are to be understood as the matrix $\oplus [\sigma^i:\gamma^{i+1}]\rho^F_{\sigma,\gamma}$.

	For a complex of sheaves
	\[
	\xymatrix{F^i \ar[d] & \rightsquigarrow & \cdots \ar[r] & \bigoplus_{\gamma^{j+1}\in X}[\hat{\gamma}^{j+1}]^{F^i(\gamma^{j+1})} \ar[r]_-{[\gamma:\tau]\rho^{F^i}} \ar[d] & \bigoplus_{\tau^{j+2}\in X}[\hat{\tau}^{j+2}]^{F^i(\tau^{j+2})} \ar[r] \ar[d] & \cdots \\
	F^{i+1} \ar[d] & \rightsquigarrow & \cdots \ar[r] & \bigoplus_{\gamma^{j+1}\in X}[\hat{\gamma}^{j+1}]^{F^{i+1}(\gamma^{j+1})} \ar[r]_-{[\gamma:\tau]\rho^{F^{i+1}}} \ar[d] & \bigoplus_{\tau^{j+2}\in X}[\hat{\tau}^{j+2}]^{F^{i+1}(\tau^{j+2})} \ar[r] \ar[d] & \cdots \\
	F^{i+2} & \rightsquigarrow & \cdots \ar[r] & \bigoplus_{\gamma^{j+1}\in X}[\hat{\gamma}^{j+1}]^{F^{i+2}(\gamma^{j+1})} \ar[r]_-{[\gamma:\tau]\rho^{F^{i+2}}}  & \bigoplus_{\tau^{j+2}\in X}[\hat{\tau}^{j+2}]^{F^{i+2}(\tau^{j+2})} \ar[r]  & \cdots}
	\]
	where we then pass to the totalization.
\end{defn}
	
Before discussing why this functor is an equivalence, let us deduce a few computational consequences of this functor.
\begin{thm}
	If $F$ is a cell sheaf on a cell complex $X$, then
	\[
		H^i_c(X;F)\cong H_{-i}(X;\equivp(F)).
	\]
\end{thm}
\begin{proof}
First we apply the equivalence functor $\equivp$ to $F$
\[
	\xymatrix{0 \ar[r] & \bigoplus_{v\in X}[\hat{v}]^{F(v)} \ar[r]_-{[v:e]\rho{e,v}} & \bigoplus_{e\in X}[\hat{e}]^{F(e)} \ar[r]_-{[e:\sigma]\rho_{\sigma,e}} & \bigoplus_{\sigma\in X}[\hat{\sigma}]^{F(\sigma)} \ar[r] & \cdots}
\]
Taking colimits (pushing forward to a point) term by term produces the complex of vector spaces
\[
	\xymatrix{0 \ar[r] & \bigoplus_{v\in X}F(v) \ar[r]_-{[v:e]\rho{e,v}} & \bigoplus_{e\in X}F(e) \ar[r]_-{[e:\sigma]\rho_{\sigma,e}} & \bigoplus_{\sigma\in X}F(\sigma) \ar[r] & \cdots}
\]
which the reader should recognize as being the computational formula for computing compactly supported sheaf cohomology.
\end{proof}

Now let us give a simple proof of the standard Poincar\'e duality statement on a manifold $X$ with coefficients in an arbitrary cell sheaf $F$, except this time the sheaf homology groups are used. 

\begin{thm}\label{thm:mfld_duality}\index{duality!over a manifold}
	Suppose $F$ is a cell sheaf on a cell complex $X$ that happens to be a compact manifold (so it has a dual cell structure $\hat{X}$), then
	\[
		H^i(X;F)\cong H_{n-i}(X;F).
	\]
	Where the group on the right is not just notational, but it indicates the left-derived functors of $p_{\dd}$ on sheaves.
\end{thm}

\begin{proof}
	We repeat the first step of the proof of the previous theorem. By feeding $F$ through the equivalence $\equivp$ we get a complex of cosheaves. Pushing forward to a point yields a complex whose (co)homology is the compactly supported cohomology of the sheaf $F$. Now we recognize that the formula yields a formula for the Borel-Moore homology for the cosheaf naturally defined on the dual cell structure.
	\[
		\xymatrix{0 \ar[r] & \ar@{~>}[d] \bigoplus_{v\in X}F(v) \ar[r]_-{[v:e]\rho{e,v}} & \ar@{~>}[d] \bigoplus_{e\in X}F(e) \ar[r]_-{[e:\sigma]\rho_{\sigma,e}} & \ar@{~>}[d] \bigoplus_{\sigma\in X}F(\sigma) \ar[r] & \cdots \\
		0 \ar[r] & \bigoplus_{\tilde{v}\in \tilde{X}}\hF(\tilde{v}) \ar[r]_-{[\tilde{v}:\tilde{e}]\rho{\tilde{e},\tilde{v}}} & \bigoplus_{\tilde{e}\in \tilde{X}}\hF(\tilde{e}) \ar[r]_-{[\tilde{e}:\tilde{\sigma}]\rho_{\tilde{\sigma},\tilde{e}}} & \bigoplus_{\tilde{\sigma}\in \tilde{X}}\hF(\tilde{\sigma}) \ar[r] & \cdots}
	\]
	Taking the homology of the bottom row is the usual formula for the Borel-Moore homology of a cellular cosheaf except the top dimensional cells are place in degree 0, the $n-1$ cells in degree -1, and so on. Everything being shifted by $n=\dim X$ we get the isomorphism
	\[
		H_{-i}(X;\equivp(F))\cong H^{BM}_{n-i}(\tilde{X};\hF).
	\]
	However, we already observed in Theorem~\ref{thm:mfld_sheaf_cosheaf} that the diagrams $\hF$ on $\tilde{X}$ and $F$ on $X$ are the same in every possible way, so in particular sheaf homology of $F$ must coincide with cosheaf homology of $\hF$. Thus using compactness to drop the Borel-Moore label and chaining together the previous theorem we get
	\[
		H^i(X;F)\cong H_{-i}(X;\equivp(F))\cong H_{n-i}(\tilde{X};\hF)\cong H_{n-i}(X;F).
	\]
\end{proof}

\section{Derived Equivalence of Sheaves and Cosheaves}
\label{subsec:derived_equiv}

Historically, the derived equivalence of cellular sheaves and cosheaves appears in a few places and is re-discovered again and again. In chronological order, the first published proof appears to be in the 1998 paper of Peter Schneider in ``Verdier Duality on the Building''~\cite{schneider-vd}, which is a follow-up of a longer paper connecting sheaves, buildings and representation theory~\cite{schneider-rep}. Unfortunately, Schneider uses the term ``local coefficient systems'' to mean what we mean by cellular cosheaves. At around the same time Maxim Vybornov made explicit mention of the relationship between sheaves and cosheaves, relating them through Koszul duality~\cite{vybornov-triang}, but it took up until 2005 for Kohji Yanagawa to explicitly state that Vybornov's work implied the derived equivalence of sheaves and cosheaves~\cite{yanagawa}.

However, the perspective presented here was arrived at independently of the above work. In early March 2012, Bob MacPherson gave a lecture (which the author attended) where he conjectured that the derived category of cellular sheaves and cosheaves should be equivalent. Within a few weeks the author produced a proof. After some truly insightful comments from David Lipsky, the equivalence was refined to its current form.

Although the ideas were foreshadowed by many sources, the use of stalk (co)sheaves appears to be a novel way of arguing.

\begin{thm}[Equivalence]\label{thm:equivalence}\index{duality!derived equivalence}\index{derived equivalence!of cellular sheaves and cosheaves}\index{sheaf!cellular!derived equivalence with cosheaves}\index{cosheaf!cellular!derived equivalence with sheaves}
	$\equivp:D^b(\Shv(X))\to D^b(\Coshv(X))$ is an equivalence.
\end{thm}
\begin{proof}
First let us point out that the functor $\equivp$ really is a functor. Indeed if $\alpha:F\to G$ is a map of sheaves then we have maps $\alpha(\sigma):F(\sigma)\to G(\sigma)$ that commute with the respective restriction maps $\rho^F$ and $\rho^G$. As a result, we get maps $[\hat{\sigma}]^{F(\sigma)}\to [\hat{\sigma}]^{G(\sigma)}$. Moreover, these maps respect the differentials in $\equivp(F)$ and $\equivp(G)$, so we get a chain map. It is clearly additive, i.e. for maps $\alpha,\beta:F\to G$ $\equivp(\alpha+\beta)=P(\alpha)+P(\beta)$. This implies that $\equivp$ preserves homotopies.

It is also clear that $\equivp$ preserves quasi-isomorphisms. Note that a sequence of cellular sheaves $A^{\bullet}$ is exact if and only if $A^{\bullet}(\sigma)$ is an exact sequence of vector spaces for every $\sigma\in X$. This implies that $\equivp(A^{\bullet})$ is a double-complex with exact rows. By standard results surrounding the theory of spectral sequences or by the acyclic assembly lemma (\cite{weibel} Lem. 2.7.3) we get that the totalization is exact.

Let us understand what this functor does to an elementary injective sheaf $[\sigma]^V$. Applying the definition we can see that
\[
\xymatrix{
 \equivp: [\sigma]^V &\rightsquigarrow& \bigoplus_{\tau^0\subset\sigma}[\hat{\tau}^0]^V \ar[r] &\cdots \ar[r] & \bigoplus_{\tau^i\subset\sigma}[\hat{\tau}^i]^V \ar[r] &\cdots \ar[r] & [\hat{\sigma}]^V
}
\]
which is nothing other than the projective cosheaf resolution of the skyscraper (or stalk) cosheaf $\hat{S}_{\sigma}^V$ supported on $\sigma$, i.e.
\[
	\hat{S}_\sigma^V(\tau)=\left\{\begin{array}{ll} V &\sigma=\tau\\ 0 & \mathrm{o.w.}\end{array}\right.
\]
Consequently, there is a quasi-isomorphism $q:\equivp([\sigma]^V)\to \hat{S}_{\sigma}^V[-\dim\sigma]$ where $\hat{S}_{\sigma}^V$ is placed in degree equal to the dimension of $\sigma$ assuming that $[\sigma]^V$ is initially in degree 0. By abusing notation slightly and letting $\equivpc$ send cosheaves to sheaves, we see that
\[
 \equivpc(q):\equivpc\equivp([\sigma]^V)\to \equivpc(S_{\sigma}^V)=[\sigma]^V
\]
and thus we can define a natural transformation from $\equivpc\equivp$ to $\id_{D^b(\Shv)}$ when restricted to elementary injectives. However, by Lemma \ref{lem:inj} we know that every injective looks like such a sum, so this works for injective sheaves concentrated in a single degree. However, it is clear that $\equivp$ sends a complex of injectives, before taking the totalization of the double complex to the projective resolutions of a complex of skyscraper cosheaves. Applying $\equivp$ to the quasi-isomorphism relating the double complex of projective cosheaves to the complex of skyscrapers, extends the natural transformation to the whole derived category. However, since $\equivp$ preserves quasi-isomorphisms, this natural transformation is in fact an equivalence. This shows $\equivpc\equivp\cong\id_{D^b(\Shv)}$. Repeating the argument starting from co-sheaves shows that
\[
 \equivp:D^b(\Shv(X))\leftrightarrow D^b(\Coshv(X)):\equivpc
\]
is an adjoint equivalence of categories.
\end{proof}

The above proof should be taken as the primary duality result from which other dualities spring. This was not always appreciated and the author's first attack on the proof was to chain together two well-known dualities, which we review in the next two sections.

\subsection{Linear Duality}
There is an endofunctor on the category of finite dimensional vector spaces $\vect$ given by sending a vector space to its dual $V\rightsquigarrow V^*$. This functor has the effect of taking a cellular sheaf $(F,\rho)$ to a cellular co-sheaf $(F^*,\rho^*)$, since the restriction maps get dualized into extension maps. It is contravariant since a sheaf morphism $F\to G$ gets sent to a co-sheaf morphism in the opposite direction $F^*\leftarrow G^*$ as one can easily check. We can promote this functor to the derived category, using a subscript $f$ to remind the reader when we restrict to the finite dimensional full subcategories.

\begin{defn}[Linear Duals]\index{duality!linear}\index{linear duality}
	Define $\hat{V}:D^b(\Shv(X))^{op}\to D^b(\Coshv(X))$ as follows
	\begin{itemize}
		\item[-] $\hat{V}(F^{\bullet})=(F^*)^{-\bullet}$, i.e. take a sheaf in slot $i$, dualize its internal restriction maps $\rho_{\sigma,\tau}^{F^i}$ to extension maps $r^{F^{i*}}_{\tau,\sigma}$ to obtain a co-sheaf and then put it in slot $-i$.
		\item[-] $\hat{V}$ sends differentials between sheaves $d^i$ to their adjoints in negative degree $\partial^{-i-1}:=(d^i)^*$
		\[\star(\xymatrix{\cdots\ar[r] & F^{i} \ar[r]^{d^i} & F^{i+1} \ar[r] & \cdots)}=\xymatrix{\cdots \ar[r] & [(F^{i+1})^*]^{-i-1} \ar[r]_{\partial^{-i-1}}& [(F^{i})^*]^{-i} \ar[r] & \cdots }\]
		We'll adopt the convention that lowering the index increases the degree $\partial^{-i-1}\to \partial_{i}$.
	\end{itemize}
	We will reserve the right to abuse notation and let $V$ map from co-sheaves to sheaves in the obvious manner, i.e. $V:D^b(\Coshv(X))^{op}\to D^b(\Shv(X))$ or formally equivalent $V:D^b(\Coshv(X))\to D^b(\Shv(X))^{op}$.
\end{defn}

\begin{lem}
	$\hat{V}_f:D^b(\Shv_f(X))\to D^b(\Coshv_f(X))^{op}$ is an equivalence of categories.
\end{lem}
\begin{proof}
	It is clear that if $\alpha:I^{\bullet}\to J^{\bullet}$ is a map in the category of complexes of sheaves homotopic to zero $\alpha\simeq 0$, i.e. there exists a map $h:I^{\bullet}\to J^{\bullet-1}$, written $h:I\to J[-1]$ such that $\alpha^n-0^n=d_J^{n-1}h^n+h^{n+1}d_I^n$. Writing out how $\hat{V}$ acts carefully we see that $\hat{V}(\alpha):\hat{V}(J)\to \hat{V}(I)$ and $\hat{V}(h):V(J[-1])=\hat{V}(J)[+1]\to \hat{V}(I)$ defines a homotopy between $\hat{V}(\alpha)$ and $\hat{V}(0)=0$ by setting $(h^*)^{\bullet}=\hat{V}(h)^{\bullet -1}$.

	$\hat{V}$ thus sends $K^b(Inj-S_f)^{op}$ to $K^b(Proj-C_f)$ and composed twice $V\hat{V}:K^b(Inj-S_f)\to K^b(Inj-S_f)$ is naturally isomorphic to the identity functor, so it is an equivalence. We can repeat the arguments for co-sheaves and use formality to put the $^{op}$ where we want.
\end{proof}

\subsection{Verdier Dual Anti-Involution}

\begin{defn}[Verdier Dual]\index{Verdier duality}\index{duality!Verdier}
The Verdier dual functor $D:D(\Shv_f(X))\to D(\Shv_f(X))^{op}$ is defined as $D:=\mathcal{H}om(-,\omega^{\bullet}_X)$. Recall that $\mathcal{H}om(F,G)$ is a sheaf whose value on a cell $\sigma$ is given by $\Hom(F|_{st(\sigma)},G|_{st(\sigma)})$, i.e. natural transformations between the restrictions to the star of $\sigma$.

The complex of injective sheaves $\omega_X^{\bullet}$ is called the dualizing complex of $X$. It has in negative degree $\omega^{-i}_X$ the sum over the one-dimensional elementary injectives concentrated on $i$-cells $[\gamma^i]$. The maps between use the orientations on cells to guarantee it is a complex.
\[
\xymatrix{\cdots \ar[r] & \oplus_{|\tau|=i+1}[\tau] \ar[r]^{\oplus[\gamma:\tau]} & \oplus_{|\gamma|=i}[\gamma] \ar[r]^{\oplus[\sigma:\gamma]} & \oplus_{|\sigma|=i-1}[\sigma] \ar[r] &\cdots }
\]
The Verdier dual of $F$ is the complex of sheaves $D^{\bullet}F:=\mathcal{H}om(F,\omega^{\bullet}_X)$. Written out explicitly it is
\[
\xymatrix{\cdots \ar[r] & \oplus_{|\tau|=i+1}[\tau]^{F(\tau)^*} \ar[r]_{\oplus[\gamma:\tau]\rho^*} & \oplus_{|\gamma|=i}[\gamma]^{F(\gamma)^*} \ar[r]_{\oplus[\sigma:\gamma]\rho^*} & \oplus_{|\sigma|=i-1}[\sigma]^{F(\sigma)^*} \ar[r] &\cdots }
\]
\end{defn}

\begin{prop}
	The functor $\equivp:D^b(\Shv_f(X))\to D^b(\Coshv_f(X))$ composed with linear duality $V:D^b(\Coshv_f(X))\to D^b(\Shv_f(X))^{op}$ gives the Verdier dual anti-equivalence, i.e. $D\cong V\equivp$.
\end{prop}
\begin{proof}
Just check by hand.
\end{proof}

\begin{rmk}
We could have used well-known facts about Verdier duality to prove a weaker version of our main theorem by restricting to finitely generated stalks.
\end{rmk}

%
%

\chapter{Cosheaves as Valuations on Sheaves}
\label{sec:valuations}

\begin{quote}
{\em``Speech is the twin of my vision....it is unequal to measure itself.''}
\begin{flushright} --- Walt Whitman's Song of Myself [25] \end{flushright}
\end{quote}

The development of cosheaves as a theory is largely fragmented. Researchers at different points in time have found a use for it here and there, at the service of different purposes and interests. The more strongly categorical and logical community have done some considerable work understanding the relationship between the topos of sheaves and cosheaves. One insight that seems very worthwhile is that cosheaves act on sheaves in a natural way. Although one can use a little bit of category theory to draw this conclusion, we use this to give some surprising reformulations of classical sheaf theory. Namely, the primary observation of this section is that the action of taking compactly supported cohomology of a sheaf can be interpreted as an action of a very particular cosheaf on the category of all sheaves.

\section{Left and Right Modules and Tensor Products}
\index{module!left and right}
Suppose $R$ is a ring with unit $1_R$. One can think of $R$ as a category with a single object $\star$ whose set of morphisms
\[
	\Hom_R(\star,\star)\cong R
\]
has the structure of an abelian group. The multiplication in the ring plays the role of a composition so $r\cdot s=r\circ s$. The abelian group structure, which corresponds to the ability to add morphisms $r+s$, reflects the fact that rings have an underlying abelian group structure. One says that $R$ is a pre-additive category, or is a category enriched in $\Ab$ --- the category of abelian groups.

An additive functor $B:R\to\Ab$ is a functor that preserves the abelian group structure, so it picks out a single abelian group, which we also call $B$, and satisfies the relation $(r+s)\cdot B=r\cdot B + s\cdot B$ and $(rs)\cdot B=r\cdot (s\cdot B)$ so such a functor is precisely the data of a \textbf{left $R$-module}. Dually, a contravariant functor $A:R^{op}\to \Ab$ prescribes the data of a \textbf{right $R$-module}. Taking the tensor product over $\ZZ$ of $A$ and $B$ allows us to define a bi-module
\[
	A\otimes B :R^{op}\times R \to \Ab \qquad (\star,\star) \mapsto A\otimes_{\ZZ} B.
\]
The latter is the group freely generated by pairs of elements from $A$ and $B$ modulo the usual relations $(a+a')\otimes b=a\otimes b + a'\otimes b$ and $a\otimes (b+b')= a\otimes b + a\otimes b'$. However, in the presence of the action of a ring $R$, there is another tensor product $A\otimes_R B$ that further quotients $A\otimes_{\ZZ} B$ by the relation $(a\cdot r)\otimes b=a\otimes (r\cdot b)$. Said using diagrams, we require that for each $r$, the following diagram commutes.
\[
	\xymatrix{A\otimes B \ar[r]^{1_A\otimes B(r)} \ar[d]_{A(r)\otimes 1_B} & A\otimes B \ar[d] \\
	A\otimes B \ar[r] & A\otimes_R B}
\]

In other words there is a coequalizer
\[
 \xymatrix{ A\otimes_{\ZZ} R \otimes_{\ZZ} B \ar@<-.5ex>[rr]_{(a,r,b)\mapsto (ar,b)} \ar@<.5ex>[rr]^{(a,r,b)\mapsto (a,rb)} & & A\otimes_{\ZZ} B \ar[r] & A\otimes_{R} B}
\]
that realizes the tensor product using purely categorical operations. This allows us to work in a greater degree of generality by making use of a special type of colimit called a \textbf{coend}, that generalizes the tensor product described above.

\begin{defn}[Tensoring Sheaves with Cosheaves]\index{coend}\index{tensoring sheaves with cosheaves}
	Let $X$ be a topological space and let $\hG$ and $F$ be a pre-cosheaf and a pre-sheaf respectively, both valued in $\Vect$. Note that for every pair of objects $U\to V$ in $\Open(X)$ we have a diagram
	\[
		\xymatrix{ & \hG(V)\otimes F(V) \\ \hG(U)\otimes F(V) \ar[ru]^{r^G_{V,U}\otimes \id} \ar[rd]_{\id\otimes\rho_{U,V}^F} & \\ & \hG(U)\otimes F(U)}
	\]
	which is the building block in defining the \textbf{coend} or \textbf{tensor product over $X$}
	\[
	\bigoplus_{U\to V}\hG(U)\otimes F(V) \rightrightarrows \bigoplus_W \hG(W)\otimes F(W) \rightarrow \int^{\Open(X)} \hG(W)\otimes F(W)=:\hG\otimes_X F.
	\]
\end{defn}

We illustrate this definition with an immediate example. 
\begin{ex}[Stalks and Skyscraper Cosheaf]
Recall that we defined the \emph{skyscraper cosheaf} at $x$ to be the cosheaf
\[
	\skycshv_x(U)=\left\{ \begin{array}{ll} k & \textrm{if $x\in U$}\\
	0 & \textrm{other wise}\end{array} \right.
\]
With some thought one can show that the tensor product of any pre-sheaf $F$ with the cosheaf $\skycshv_x$ yields
\[
	\skycshv_x\otimes_X F \cong F_x
\]
by treating $F$ as a variable which can range over all pre-sheaves, one gets, in particular, a functor
\[
	\skycshv_x\otimes_X - :\Shv(X) \to \Vect \qquad F \rightsquigarrow F_x.
\]
\end{ex}

The previous example demonstrates an important observation: 
\emph{The operation of taking stalks is equivalent to the process of tensoring with the skyscraper cosheaf.} 

To see how far this observation can be generalized, note that if we fix $\hG$ and let $F$ vary then we get a functor
\[
\hG\otimes_X -:\Shv(X)\to\Vect
\]
that is defined in terms of colimits and is thus co-continuous (it sends colimits to colimits). Now we are free to take an arbitrary cosheaf and let it act on sheaves. The ``one obvious choice'' of taking stalks at a point is run over by a veritable slew of valuations, one for each cosheaf. Moreover, it is clear that this description extends to a pairing between the symmetric monoidal categories $\Coshv(X)$ and $\Shv(X)$, i.e.
\[
	-\otimes_X-:\Coshv(X)\times \Shv(X)\to\Vect \qquad (\hG,F)\mapsto \hG\otimes_X F:=\int^{\Open(X)}\hG(U)\otimes F(U),
\]
although we haven't used the sheaf or cosheaf axiom anywhere, so the pairing is actually valid for pre-sheaves and pre-cosheaves.

\section{Compactly-Supported Cohomology}

Although the idea of using coends to tensor together sheaves and cosheaves has been independently re-discovered many times, cf. Jean-Pierre Schneider's 1987 work~\cite{schneiders-ch}, it has not been used to do any serious work. This is a shame in light of the following 1985 theorem of A.M. Pitts~\cite{pitts}.
\begin{thm}\index{Pitt's theorem}
	Let $X$ be any topological space. Every colimit-preserving functor on sheaves arises by tensoring with a cosheaf, i.e. $$\Coshv(X;\Set)\cong\Fun^{co-cts}(\Shv(X;\Set),\Set).$$
\end{thm}
This theorem is also stated in Marta Bunge and Jonathan Funk's 2006 book ``Singular Coverings of Toposes''~\cite{bungefunk} as theorem 1.4.3, which further surveys some of Lawvere's philosophy of distributions on topoi. The topos community deserves commendation for keeping the study of cosheaves alive during the past few decades, but so far work in the enriched and computable setting of vector spaces is largely missing.

We attempt to partly remedy this gap by establishing a connection between the tensor operation and the cohomology of sheaves. However, instead of establishing an enriched version of Pitt's theorem,\footnote{We delay this for another paper.} we will use it as a guide. For example, in classical sheaf theory, compactly supported cohomology is gotten by taking the constant map $p:X\to\star$ and associating to it the pushforward with compact supports functor $p_!:\Shv(X)\to\Shv(\star)\cong\Vect$. Of course, just applying $p_!$ defines only compactly supported zeroth cohomology of a sheaf $H_c^0(X;F)$. To get the higher compactly supported cohomology groups one takes an injective resolution and applies $p_!$ to the resolution. The result will be a complex of vector spaces, whose cohomology in turn produces the desired groups:
\[
	F\to I^{\bullet} \rightsquigarrow Rp_!:=p_! I^{\bullet} \qquad H^i(p_! I^{\bullet}):= H^i_c(X;F).
\] 

Historically the first fundamental duality result in sheaf theory was the statement that $Rp_!$ admits a right adjoint on the level of the derived category. This adjunction is sometimes called \textbf{global Verdier duality}:
\[
	\Hom(Rp_!F,G)\cong \Hom(F,p^!G).
\]
By applying the fact that left adjoints are co-continuous one is led to believe, in light of Pitt's theorem, that there should be a cosheaf that realizes the operation of taking derived pushforward with compact supports.

In light of the derived equivalence between cellular sheaves and cosheaves established in this paper, we provide an explicit description of the complex of cosheaves that realizes the derived pushforward.

In preparation, one should note that there are several cosheaves that realize the operation of taking stalks at a point $x$ in the cellular world. One is
\[
  \hat{\delta}_{\sigma}(\tau)=\left\{ \begin{array}{ll} k & \textrm{if $\sigma=\tau$}\\
	0 & \textrm{o.w.}\end{array} \right.
\]
The other is the correct formulation of $\skycshv_{\sigma}$ when using the Alexandrov topology
\[
  [\hat{\sigma}](\tau)=\left\{ \begin{array}{ll} k & \textrm{if $\tau\leq\sigma$}\\
	0 & \textrm{o.w.}\end{array} \right.
\]
Recall that this is also the elementary projective cosheaf concentrated on $\sigma$ with value $k$.

Observe that the first cosheaf returns the value $F(\sigma)$ because every other cell is tensored with zero. The second cosheaf works by restricting the non-zero values of $F$ to the closure of the cell $\sigma$, but this restricted diagram has a terminal object given by $F(\sigma)$, so the colimit returns $F(\sigma)$ as well.

This allows us to state the main theorem of this section.

\begin{thm}\index{cohomology!of a cellular sheaf!compactly-supported}\index{sheaf!cellular!compactly supported cohomology of}
	Let $X$ be a cell complex, then the operation $Rp_!:\Shv(X)\to\Vect$ on cellular sheaves is equivalent to tensoring with the image of the constant sheaf through the derived equivalence, i.e.
	\[
		\equivp(k_X)= \bigoplus_{v\in X} [\hat{v}] \to \bigoplus_{e\in X} [\hat{e}] \to \bigoplus_{\sigma\in X} [\hat{\sigma}] \to \cdots.
	\]
\end{thm}
\begin{proof}
	The proof is immediate given the previous description of taking stalks, i.e. one can check directly the formula
	\[
		\equivp(k_X)\otimes_X F \cong \bigoplus_{v\in X} F(v) \to \bigoplus_{e\in X} F(e) \to \bigoplus_{\sigma\in X} F(\sigma) \to \cdots
	\]
	whose cohomology is by definition the compactly supported cohomology of a cellular sheaf $F$.
\end{proof}

This perspective is especially satisfying for the following reason: it makes transparent how the underlying topology of the space $X$ is coupled with the cohomology of a sheaf $F$. Compactly supported sheaf cohomology arises by tensoring with the complex of cosheaves that computes the Borel-Moore homology of the underlying space.

\section{Sheaf Homology and Future Directions}
\index{sheaf!homology, cellular}
The perspective of tensoring sheaves and cosheaves together offers numerous directions for further research both in pure and applied sheaf theory. Just the heuristic that
\begin{quote}
\begin{center}
	\emph{each cosheaf determines a (co-)continuous valuation on the category of sheaves,}
\end{center}
\end{quote}
is suggestive of the idea that if we are going to use sheaves to model the world, then cosheaves should allow us to weight different models of the world.

After having recovered some classical operations on sheaves, we are left with many more to consider. For example the constant cosheaf $\hat{k}_X$ should act on sheaves by returning its colimit, i.e. zeroth sheaf homology
\[
	\hat{k}_X\otimes_X - :\Shv(X) \to \Vect \qquad F \rightsquigarrow H_0(X;F)=p_{\dd} F.
\]
By taking a projective resolution of the constant cosheaf once and for all, one then gets for free a way of computing \textbf{higher sheaf homology}. This yet-to-be-explored theory has only recently found its use in applications, e.g. the work of Sanjeevi Krishnan on max-flow min-cut.

Additionally, the decategorification of the pairing of the categories of constructible sheaves and cosheaves provides an alternative approach to the study of Euler integration and leads in a natural way to the study of \textbf{higher Euler calculus} through higher K-theory. More directly the operation of pairing sheaves and cosheaves is reminiscent of a convolution operation. This area is under active research in collaboration with Aaron Royer.
\index{coend!higher Euler calculus}

%
%
\chapter{Graded Descriptions of the Derived Category}
\label{sec:graded}

\begin{quote}
{\em``Mathematics is the art of giving the same names to different things.''}
\begin{flushright} --- Henri Poincar\'e~\cite{poincare1908future} \end{flushright}
\end{quote}

In Chapter~\ref{sec:derived} we introduced the derived category of cellular sheaves. The fundamental objects there are chain complexes parametrized by a cell complex. However, the derived category takes a further step by identifying objects that are ``essentially the same'' when viewed through the lens of cohomology sheaves.

\begin{figure}
\centering
\includegraphics[width=.4\textwidth]{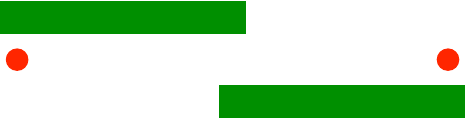}
\caption{A Chain Complex for the 2-Sphere as Barcodes}
\label{fig:bar_chain}
\end{figure}

In Section~\ref{subsec:derived_vect} we understand this principle better by demonstrating the well known fact that the derived category of chain complexes of vector spaces (sheaves over a point) is equivalent to the graded category of vector spaces. Our proof follows the standard proof in~\cite{weibel} except we use the barcode method for chain complexes introduced in Example~\ref{ex:chaincomplex_bc} to visualize explicitly what is happening. Roughly speaking, the derived category of chain complexes allows us to remove the green bars in Figure~\ref{fig:bar_chain} as they ``graph'' the chain homotopy between the identity and the map that projects onto and then includes the red dots.

This sets us up for Section~\ref{subsec:derived_graded_cellsheaves}, which culminates in a proof that the derived category of cellular sheaves over a one-dimensional base space is equivalent to a graded category of sheaves. This should seem plausible because over each cell, a chain complex is equivalent to a graded vector space. Indeed, one could repeat the proof of Lemma~\ref{lem:chaincomplex_qis} verbatim if it weren't for the pesky fact that projecting onto the cohomology cell-by-cell fails to define a sheaf map. However, the proof of the equivalence follows by constructing an explicit replacement of every object in the derived category with an object where such a na\"ive map will exist.

\section{The Derived Category for Complexes of Vector Spaces}
\label{subsec:derived_vect}

Here we warm-up with an alternative approach to the derived category of complexes of vector spaces. Recall that a chain complex simply consists of a collection of vector spaces and maps satisfying $d^2=0$.
\[
\cdots \to V^{i-1} \to V^i \to V^{i+1} \to \cdots 
\]
Cohomology defines a functor from chain complexes to graded vector spaces simply by placing the $i^{\mathrm{th}}$ cohomology in degree $i$.
\[
H^*:C^b(\Vect)\to \grVect \qquad (V^{\bullet},d_v) \rightsquigarrow \{H^*(V^{\bullet},d_V)\}
\]
A map of chain complexes $f^{\bullet}:V^{\bullet}\to W^{\bullet}$ is a quasi-isomorphism if the maps $H^i(f):H^i(V^{\bullet})\to H^i(W^{\bullet})$ is an isomorphism for every non-negative integer $i$.

The derived category of chain complexes is defined to be the category of chain complexes localized at the collection of quasi-isomorphisms. This is often simplified by saying that
\begin{quote}
\centering
{\em In the derived category quasi-isomorphisms are formally inverted.}
\end{quote}

\index{barcode!associated to chain complex}

To illustrate this slogan we will prove that every chain complex is quasi-isomorphic to a graded vector space. This is the simplest instance of a more general theorem that we prove in this chapter.

\begin{lem}\label{lem:chaincomplex_qis}\index{derived equivalence!with graded vector spaces}
A chain complex $(V,d_V)$ is quasi-isomorphic to its cohomology $(H^{\bullet}(V,d_V),0)$, viewed as complex with zero differentials.
\end{lem}
\begin{proof}
It suffices to define a chain map 
\[
\pi^{\bullet}:(V,d_V)\to (H^{\bullet}(V,d_V),0) \qquad\mathrm{or}\qquad
\iota^{\bullet}:(H^{\bullet}(V,d_V),0) \to (V,d_V)
\]
that induces the obvious isomorphism $H^i(\pi):H^i(V)\to H^i(V)$. We will use the decomposition for persistence modules from Theorem~\ref{thm:crawleyboevey} to define this map.

In Example~\ref{ex:chaincomplex_bc} we observed that any chain complex $(V,d_V)$ is isomorphic to a direct sum
\[
V\cong \bigoplus_{i\in\ZZ} S_i^{n_i}\oplus P_i^{m_i}
\]
where
\[
S_i: \qquad \cdots \to 0 \to k \to 0 \to \cdots
\]
is a length zero interval module and
\[
P_i: \qquad \cdots \to 0 \to k \to k \to 0 \to \cdots
\]
is a length one interval module, with the first non-zero term in degree $i$.
If we group terms so that 
\[
B^{i-1}:=P_i^{m_{i-1}} \qquad H^i(V):=S_i^{n_i} \qquad B^i:=P_i^{m_i}
\]
then we get an obvious map, which precomposed with the above isomorphism is our desired quasi-isomorphism.
\[
\xymatrix{\cdots \ar[r] & B^{i}\oplus H^{i}(V)\oplus B^{i-1} \ar[r] \ar[d]_{\pi} & B^{i+1}\oplus H^{i+1}(V)\oplus B^{i} \ar[r] \ar[d]_{\pi} & \cdots \\
\cdots \ar[r]_0 & H^{i}(V) \ar[r]_0 & H^{i+1}(V) \ar[r] & \cdots}
\]
One can in fact show more. Let $\pi:B^i\oplus H^{i}(V)\oplus B^{i-1} \to H^i(V)$ be the obvious projection and $\iota:H^i(V) \to B^{i-1}\oplus H^{i}(V)\oplus B^i$ be the obvious inclusion, then one can construct an explicit chain homotopy $s$ joining $\iota\circ \pi$ to the identity. 

First it is useful to observe that, in this basis, the differentials have the form
\[
d^i=\left[
\begin{array}{ccc}
0 & 0 & 0 \\
0 & 0 & 0 \\
\id & 0 & 0
\end{array}
\right]
\]
A clear candidate for the map $s^i$ is the projection onto $B^{i-1}$ that then identifies it with its isomorphic copy as the third summand in the decomposition of $V^{i-1}$, i.e. the matrix
\[
s^i=\left[
\begin{array}{ccc}
0 & 0 & \id \\
0 & 0 & 0 \\
0 & 0 & 0
\end{array}
\right]
\]
One can then check directly the equation 
\[
\id-\iota\circ\pi=d^{i-1}\circ s^i + s^{i+1}\circ d^i.
\]
\end{proof}

\begin{rmk}
Using the barcode description, the map $s^i$ simply follows length one barcodes to the left. Thus the length one barcodes can be interpreted as the ``graph'' of a chain homotopy.
\end{rmk}

\begin{cor}
For any collection of integers $m_i$ and $n_i$ we have the following isomorphisms in the derived category.
\[
\bigoplus_{i\in\ZZ} S_i^{n_i} \simeq \bigoplus_{i\in\ZZ} S_i^{n_i}\oplus P_i^{m_i}
\]
\end{cor}
The upshot of the above corollary is that 

\begin{quote}
\centering
\emph{``Indecomposables do not survive the derived category!''}~\cite{macpherson}
\end{quote}

\section{Derived Complexes of Cellular Sheaves}
\label{subsec:derived_graded_cellsheaves}

Recall from Chapter~\ref{sec:derived} that a complex of cellular sheaves $F^{\bullet}$ assigns to every cell $\sigma$ a chain complex and to every pair of incident cells $\sigma\leq\tau$ a chain map $\rho^{\bullet}_{\tau,\sigma}:F^{\bullet}(\sigma)\to F^{\bullet}(\tau)$. For each $i$ we can define the $i^{\mathrm{th}}$ cohomology sheaf as the assignment
\[
 \cohomsheaf^i(F^{\bullet}): \qquad \sigma \rightsquigarrow H^i(F^{\bullet}(\sigma)).
\]
The restriction maps being defined as the map induced on cohomology  by $\rho^{\bullet}_{\tau,\sigma}$. A quasi-isomorphism is a map of complexes $f^{\bullet}:F^{\bullet}\to G^{\bullet}$ that induces isomorphisms on each cohomology sheaf. We can also view all of the cohomology sheaves as a single graded cohomology sheaf $\cH^{*}F^{\bullet}$.

\subsection{Counterexample to the Na\"ive Approach}

\begin{figure}
\centering
\includegraphics[width=.4\textwidth]{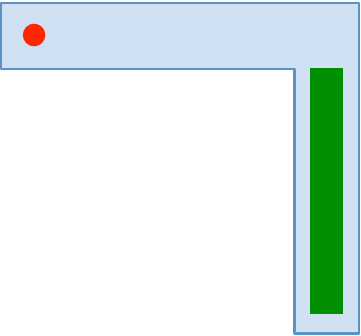}
\caption{The Counterexample}
\label{fig:bar_counterex}
\end{figure}

To begin, let us consider the simplest possible space where a complex of sheaves has interesting behavior. Let $X=[0,1)$ with cell structure $v=\{0\}$ and $e=(0,1)$. A complex of sheaves over $X$ is completely described by two chain complexes and a chain map between them:
\[
\rho^{\bullet}_{e,v}:F^{\bullet}(v)\to F^{\bullet}(e)
\]
At first glance there should be a simple adaptation of Lemma~\ref{lem:chaincomplex_qis} to a pair of chain complexes. Indeed, we can use Theorem~\ref{thm:crawleyboevey} over each cell to obtain a candidate inclusion map from the cohomology sheaf into the complex:
\[
\xymatrix{F^i(v) \ar[r]^{\rho_{e,v}} & F^i(e) \\
B^i_v \oplus H^i_v \oplus B^{i-1}_v \ar[r]^{\rho_{e,v}} \ar[u]^{\cong}  & B^i_e \oplus H^i_e \oplus B^{i-1}_e \ar[u]_{\cong} \\
 H^i_v \ar[r]_{H^i(\rho)} \ar[u] & H^i_e \ar[u]}
\]
However such a map \emph{does not commute} in general and as such fails to define a sheaf map. For example, if $F^{\bullet}$ was the following pair of complexes and chain map
\[
\xymatrix{k \ar[r]^1 & k \\
0 \ar[r] \ar[u] & k \ar[u]_1}
\]
then the map $\rho$ takes the one and only generator of $H^i_v$ to an element of the boundaries $B^{i-1}_e$, which is non-zero, but zero in cohomology, since $H^i(\rho)=0$. We have provided a barcode version of this counter example in Figure~\ref{fig:bar_counterex}.

\subsection{Using the Calculus of Fractions Formulation}

In order to address the counterexample to the na\"ive approach, we will employ a zig-zag of morphisms. In this section we briefly review why such a zigzag is natural, when viewing the derived category as the \textbf{left calculus of fractions}~\cite{gabriel1967calculus}. An alternative description of the derived category goes as follows~\cite[p.52-3]{shepard}.

\begin{defn}[Derived Category via Fractions]\index{Derived@$D^b(\aat)$ derived category!calculus of fractions}\index{category!Derived@$D^b(\aat)$ derived category!calculus of fractions}
The \textbf{bounded derived category of cellular sheaves}, $D^b(\Shv(X))$, has the same collection of objects as $\Fun(X,\Ch(\Vect))$, but with a modified class of morphisms. A morphism from $F^{\bullet}$ to $G^{\bullet}$ is a diagram of the following form
\[
\xymatrix{ & J^{\bullet} & \\
F^{\bullet} \ar[ur] & & \ar[ul]_{\simeq} G^{\bullet}}
\]
where the arrow decorated with a $\simeq$ is a quasi-isomorphism. We declare two such morphisms (diagrams) to be the same if there is a larger, commutative diagram that fits in between them:
\[
\xymatrix{ & J_1^{\bullet} \ar[d]^{\simeq} & \\
F^{\bullet} \ar[ur] \ar[r] \ar[dr] & J_{12}^{\bullet} & \ar[l]_{\simeq} \ar[ul]_{\simeq} \ar[dl]^{\simeq} G^{\bullet} \\
& J_2^{\bullet} \ar[u]_{\simeq} & }
\]
\end{defn}

\subsection{The Equivalence}\label{subsec:graded_equivalence}

\begin{figure}
\centering
\includegraphics[width=.7\textwidth]{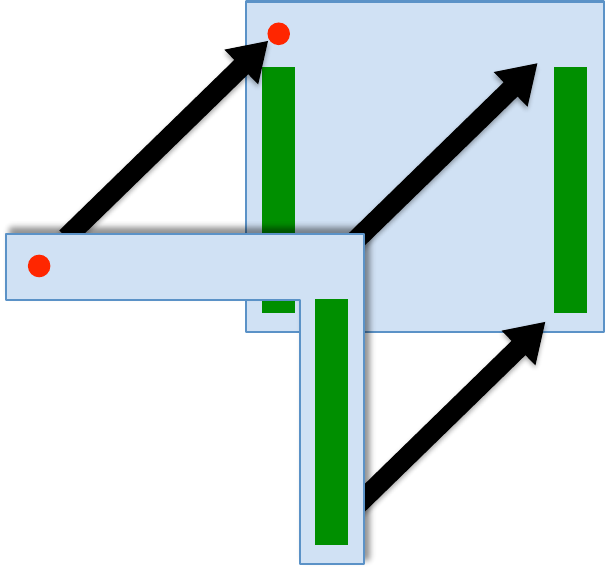}
\caption{Replacing the Counterexample with a Quasi-Isomorphic Sheaf}
\label{fig:bar_replace}
\end{figure}

We now proceed with a general method for addressing the earlier counterexample. The basic method of argument is that we will define a quasi-isomorphic complex of sheaves $J^{\bullet}$ where the na\"ive approach does work. This is visualized in Figure~\ref{fig:bar_replace}.
\[
\xymatrix{ & J^{\bullet} & \\
F^{\bullet} \ar[ur]^{\simeq} & & \ar[ul]_{\simeq} H^{*}F^{\bullet}}
\]

Let's first illustrate our method over the simple base space $X=[0,1)$.
\begin{lem}\label{lem:half-open-graded-equiv}\index{derived equivalence!with graded sheaves}\index{sheaf!cellular!derived equivalence with graded sheaves}
Over $X=[0,1)$ with cell structure $v=\{0\}$ and $e=(0,1)$, every bounded complex of cellular sheaves $F^{\bullet}$ is quasi-isomorphic to its graded cohomology sheaf $\cH^*F^{\bullet}$, i.e. $F^{\bullet}$ is \emph{isomorphic} to $\cH^*F^{\bullet}$ in the derived category.
\end{lem}
\begin{proof}
As before, we decompose the chain complex over the vertex $v$ and the edge $e$ so that $F^{\bullet}$ has the following form
\[
\xymatrix{ B^{i+1}_v\oplus H^{i+1}_v\oplus B^i_v \ar[r]^-{\rho_F^{i+1}} & B^{i+1}_e\oplus H^{i+1}_e\oplus B^i_e \\
B^{i}_v\oplus H^{i}_v\oplus B^{i-1}_v \ar[r]^-{\rho_F^{i}} \ar[u]^{d^i_v} & B^{i}_e\oplus H^{i}_e\oplus B^{i-1}_e \ar[u]_{d^i_e}
}
\]
The map $\rho^i_F$ has an easily described form in this basis.
\[
\rho^i_F=
\left[
\begin{array}{ccc}
\alpha^i & 0 & 0 \\
\gamma^i & H^i\rho & 0 \\
\delta^i & \beta^{i-1} & \alpha^{i-1}
\end{array}
\right]
\]
The na\"ive inclusion map's failure to commute is precisely described by non-zero terms in the submatrix $\beta^{i-1}$, as can be seen by inspecting the diagram below.
\[
\xymatrix{
B^i_v \oplus H^i_v \oplus B^{i-1}_v \ar[r]^-{\rho^i_F}  & B^i_e \oplus H^i_e \oplus B^{i-1}_e  \\
 H^i_v \ar[r]_-{H^i(\rho)} \ar[u] & H^i_e \ar[u]
 }
\]
However, the complex of sheaves $F^{\bullet}$ is quasi-isomorphic to the complex $J^{\bullet}$ defined as follows.
\[
\xymatrix{ B^{i+1}_v\oplus B^{i+1}_e \oplus H^{i+1}_v\oplus B^i_v \oplus B^{i}_e \ar[r]^-{\rho_J^{i+1}} & B^{i+1}_e\oplus H^{i+1}_e\oplus B^i_e \\
B^{i}_v\oplus B^{i}_e\oplus H^{i}_v\oplus B^{i-1}_v \oplus B^{i-1}_e \ar[r]^-{\rho^{i}_J} \ar[u]^{d^i_v} & B^{i}_e\oplus H^{i}_e\oplus B^{i-1}_e \ar[u]_{d^i_e}
}
\]
The matrix representation for $\rho^i_J$ has the desired form:
\[
\rho^i_J=
\left[
\begin{array}{ccccc}
0 & \id & 0 & 0 & 0\\
 \gamma^i & 0 & H^i\rho & 0 & 0 \\
\delta^i & 0 & 0 & & \id
\end{array}
\right]
\]
The quasi-isomorphism $q^{\bullet}:F^{\bullet}\to J^{\bullet}$ is defined over the vertex $v$ in any given degree $i$ as
\[
q^i_v=
\left[
\begin{array}{ccc}
\id & 0 & 0 \\
\alpha^i & 0 & 0 \\
0 & \id & 0 \\
0 & 0 & \id \\
0 & \beta^{i-1} & \alpha^{i-1} 
\end{array}
\right]
\]
and over the edge $e$ as $q^i_e=\id$. By construction, the complex $J^{\bullet}$ has a well defined sheaf map
\[
\xymatrix{
B^{i}_v\oplus B^{i}_e\oplus H^{i}_v\oplus B^{i-1}_v \oplus B^{i-1}_e \ar[r]^-{\rho^{i}_J}  & B^{i}_e\oplus H^{i}_e\oplus B^{i-1}_e \\
H^{i}_v \ar[r]_-{H^i\rho} \ar[u] & H^i_e \ar[u]
}
\]
which extends to a quasi-isomorphism $\iota:\cH^{*}F^{\bullet}\hookrightarrow J^{\bullet}$ since each of the differentials $\cH^iF^{\bullet}\to \cH^{i+1}F^{\bullet}$ are zero. This completes the proof.
\end{proof}

We can now prove the general theorem of interest.

\begin{thm}\label{thm:1D-graded-equiv}\index{derived equivalence!with graded sheaves}\index{sheaf!cellular!derived equivalence with graded sheaves}
Let $X$ be an arbitrary one dimensional cell complex and $F^{\bullet}$ a complex of sheaves over $X$. Then $F^{\bullet}$ is quasi-isomorphic to its graded cohomology sheaf. In particular, we have the equivalence of categories
\[
D^b(\Shv(X))\simeq \Shv(X;\grVect).
\]
\end{thm}
\begin{proof}
The bulk of the argument is contained in Lemma~\ref{lem:half-open-graded-equiv}, which we show extends to the desired generality. The definition of $J^i(v)$ where $v$ is a vertex with more than one incident edge is easily modified as follows:
\[
J^i(v):=B^i_v\bigoplus_{e\geq v} B^i_e \oplus H^i_v \oplus B^{i-1}_v
 \bigoplus_{e\geq v} B^{i-1}_e
\]
The restriction map to a single edge $e'$ is defined exactly as in Lemma~\ref{lem:half-open-graded-equiv} with the stipulation that factors $B^i_e$ and $B^{i-1}_e$ where $e\neq e'$ are mapped to zero.

This argument implies that we can define the desired morphism 
\[
\xymatrix{ & J^{\bullet} & \\
F^{\bullet} \ar[ur]_q^{\simeq} & & \ar[ul]_{\simeq}^{\iota} H^{*}F^{\bullet} }
\]
in an open neighborhood of a vertex. However, since over each edge $q^i_{e}:F^i(e)\to J^i(e)$ is defined to be the identity map, these locally defined maps agree on the edges and hence give a globally defined map of sheaves. This shows the equivalence with the graded cohomology sheaf.

One can check that every map of sheaves $f:F^{\bullet}\to G^{\bullet}$ extends to a map of the associated sheaves $J_f:J^{\bullet}_F\to J^{\bullet}_G$ so the construction is functorial and hence defines an equivalence of categories.
\end{proof}

%
%

\chapter{A Metric on the Category of Sheaves}
\label{sec:metric}

\begin{quote}
$A\Gamma E\Omega METPHTO\Sigma \,\, MH\Delta EI\Sigma \,\,EI\Sigma I T\Omega$\footnote{``Let no one who cannot think geometrically enter.''~\cite{plato-academy}}
\begin{flushright} --- Purported Inscription at The Academy \end{flushright}
\end{quote}

Science depends on knowing that measurements and observations should only be trusted up to an interval of uncertainty. Saying that one is traveling about 30 miles per hour depends on their being a continuum of speeds. It makes no sense to say that one is traveling at a speed that happens to be an algebraic number.~\footnote{This metaphor was used by Bob MacPherson in his opening remarks on ``Continuity and the philosophy of science.''~\cite{macpherson-seminar}} Persistent homology has addressed this issue in a rather elegant way. Given two functions $f,g:Y\to\RR$ such that
\[
||f-g||_{\infty}:=\sup_{y}|f(y)-g(y)| < \epsilon
\]
one can say that the sublevel sets obey the inclusions
\[
f^{-1}(-\infty,t] \subset g^{-1}(-\infty,t+\epsilon] \qquad \mathrm{and} \qquad g^{-1}(-\infty,t] \subset t^{-1}(-\infty,t+\epsilon]
\]
for every value of $t$. In the language of~\cite{chazal2009proximity}, one has an \textbf{interleaving} of sublevel sets:
\[
\xymatrix{ f^{-1}(-\infty,t] \ar[dr] \ar[r] & f^{-1}(-\infty,t+\epsilon] \ar[dr] \ar[r] & f^{-1}(-\infty,t+2\epsilon] \\
g^{-1}(-\infty,t] \ar[ur] \ar[r] & g^{-1}(-\infty,t+\epsilon] \ar[ur] \ar[r] & g^{-1}(-\infty,t+2\epsilon]}
\]
By functoriality of homology, one obtains a notion of interleaving of functors in $\Fun(\RR,\Vect)$, where $\RR$ with its partial order is viewed as a category. Defining the interleaving distance to be the infimum over all $\epsilon$ such that there is an $\epsilon$-interleaving gives an extended pseudo-metric on the category of such functors.

In this chapter, we will consider a generalization of this approach, first for pre-sheaves and then for sheaves. We will try to answer analogous questions such as, ``Suppose we have two maps $f,g:Y\to X$ to a metric space that are close in the supremum norm, i.e.
\[
d_{\infty}(f,g)=\sup_{y\in Y} d_X(f(y),g(y)) < \epsilon.
\]
Is there any reasonable sense where the pushforward sheaves $f_*k_Y$ and $g_*k_Y$ are close?'' It turns out that the answer is ``yes,'' but studying higher invariants of the fiber, such as $H^i$ for $i\geq 1$, is unstable to large perturbations. By studying interleavings of complexes of sheaves, one obtains a derived stability result. To conclude the chapter, we provide preliminary results towards equating the interleaving distance with a modified version of the bottleneck distance described in~\cite{cohen2007stability} for definable sheaves with finite support over the real line. The most important takeaway from this chapter is that interleavings for sheaves and pre-sheaves are obstructed by global sections.

\begin{figure}
\centering
\includegraphics[width=\textwidth]{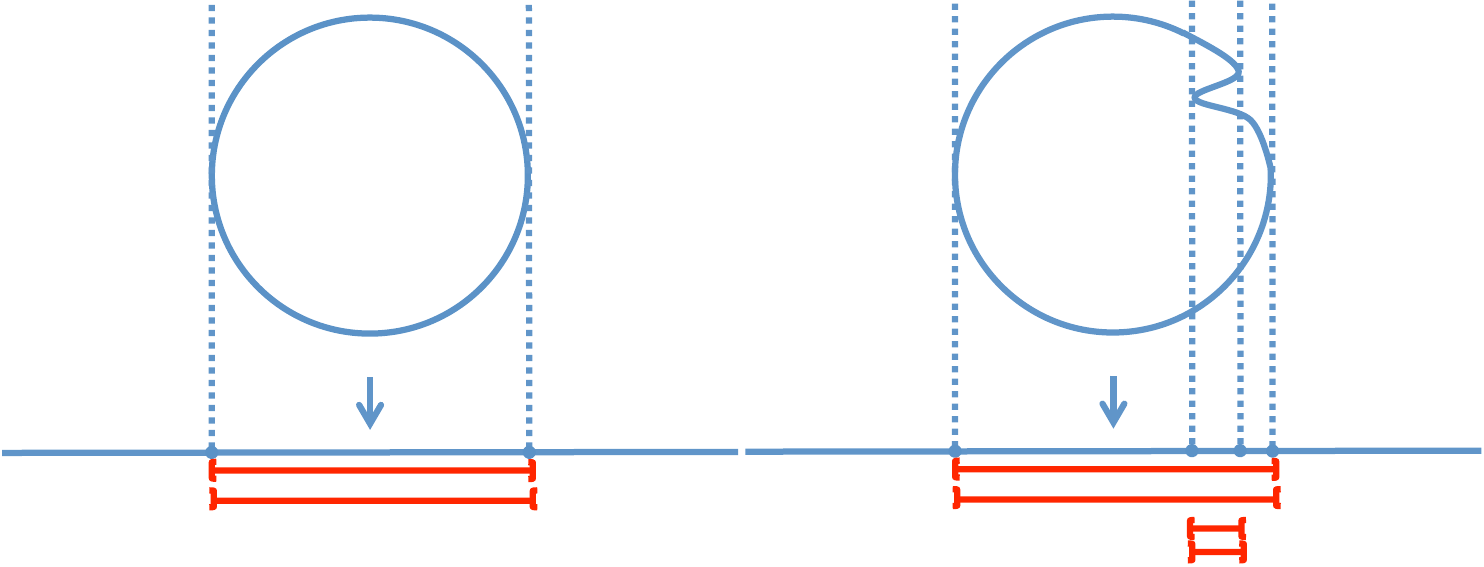}
\caption{Do Close Maps Give Rise to ``Close'' Sheaves?}
\label{fig:perturb}
\end{figure}

\section{Interleavings for Pre-Sheaves}

\begin{defn}\index{epsilon@$\epsilon$ thickening}
	Let $(X,d)$ be a metric space. The \textbf{$\epsilon$-thickening of open sets} is the map of posets
	\[
		\epsilon: \Open(X) \to \Open(X)
	\]
	given by
	\[
			U \rightsquigarrow U^{\epsilon}:=\cup_{x\in U} B(x,\epsilon).
	\]
	Of course, any map of posets dualizes to a map of posets $\epsilon:
	\Open(X)^{op}\to\Open(X)^{op}$
\end{defn}
\begin{rmk}
	For metric spaces like $\RR^n$ with the Euclidean metric, $(U^{\epsilon})^{\epsilon}=U^{2\epsilon}$. In general, the triangle inequality implies that $(U^{\epsilon})^{\epsilon}\subseteq U^{2\epsilon}$. The reverse containment is also true if $(X,d)$ is convex, for example. Convexity guarantees that intuitive results are true, but it is not strictly necessary for any of the following arguments.
\end{rmk}

\begin{defn}[Thickened Pre-Sheaf]
	Using the previous two definitions, we can define the \textbf{$\epsilon$-thickening of a pre-sheaf} $F$ via the formula
	\[
		F^{\epsilon}:=F\circ \epsilon \qquad \mathrm{i.e.} \qquad  F^{\epsilon}(U):=F(U^{\epsilon}).
	\]
	Moreover, the thickening operation is functorial. If $\varphi: F\to G$ is a natural transformation, then we get for free a natural transformation between the thickened pre-sheaves $\varphi_{\epsilon}:F^{\epsilon}\to G^{\epsilon}$. Consequently, we can define the \textbf{$\epsilon$-thickening functor} to be
	\[
		\epsilon^*:\Preshv(X)\to\Preshv(X) \qquad F \rightsquigarrow F^{\epsilon}.
	\]
\end{defn}
One of the most important observations for working with interleavings is that since $F$ is a pre-sheaf we have a canonical natural transformation
\[
	\eta^F_{\epsilon}:F^{\epsilon}\to F
\]
coming from $\rho_{U,U^{\epsilon}}:F^{\epsilon}(U)=F(U^{\epsilon})\to F(U)$. This follows by showing that for every pair $V\subseteq U$ the square
\[
	\xymatrix{ F^{\epsilon}(U) \ar[r]^{\rho_{U,U^{\epsilon}}} \ar[d]_{\rho_{V^{\epsilon},U^{\epsilon}}} & F(U) \ar[d]^{\rho_{V,U}} \\
	F^{\epsilon}(V) \ar[r]_{\rho_{V,V^{\epsilon}}} & F(V) }
\]
commutes by virtue of $F$ being a pre-sheaf:
\[
	\rho_{V,U}\circ \rho_{U,U^{\epsilon}}=\rho_{V,U^{\epsilon}}=\rho_{V,V^{\epsilon}}\circ \rho_{V^{\epsilon},U^{\epsilon}}
\]
Of course, the map
\[
	\eta^F_{2\epsilon}:F^{2\epsilon}\to F
\]
always exists as well and for metric spaces where $U^{2\epsilon}=(U^{\epsilon})^{\epsilon}$ it is equal to the composition $\eta^F_{\epsilon}\circ \epsilon^*\eta^F_{\epsilon}$.

\begin{rmk}[Notation $2\epsilon$ vs. $\epsilon\epsilon$]
To avoid cumbersome notation, we may sometimes substitute $U^{2\epsilon}$ for $(U^{\epsilon})^{\epsilon}$, $\eta^F_{2\epsilon}$ for $\epsilon^*\eta^F_{\epsilon}$, $F^{2\epsilon}$ for $\epsilon^*F^{\epsilon}$, and so on. For metric spaces such as $\RR^n$ with the Euclidean metric these differences do not exist and can be safely ignored.
\end{rmk}

A version of the following definition was communicated to the author by Amit Patel~\cite{patel}.

\begin{defn}[Interleaving of Pre-Sheaves]\index{interleaving!of presheaves}
	Let $F,G:\Open(X)^{op}\to\dat$ be two pre-sheaves on a metric space $X$. We define an \textbf{$\epsilon$-interleaving} of $F$ and $G$ to be a pair of natural transformations 
	\[
	\varphi_{\epsilon}:F^{\epsilon}\to G \qquad \psi_{\epsilon}: G^{\epsilon}\to F
	\]
	that satisfy the compatibility relations
	\[
		\eta^F_{2\epsilon}=\psi_{\epsilon}\circ \epsilon^*\varphi_{\epsilon} \qquad \eta^G_{2\epsilon}=\varphi_{\epsilon}\circ \epsilon^*\psi_{\epsilon}.
	\]
	An interleaving is better summarized via the following commutative diagram:
	\[
	\xymatrix{F^{2\epsilon} \ar[d] \ar[rrd]^(.25){\epsilon^*\varphi_{\epsilon}} & & G^{2\epsilon} \ar[d] \ar[lld]_(.25){\epsilon^*\psi_{\epsilon}} \\
	F^{\epsilon} \ar[d] \ar[rrd]^(.25){\varphi_{\epsilon}} & & G^{\epsilon} \ar[d] \ar[lld]_(.25){\psi_{\epsilon}} \\
	F & & G
	}
	\]
	Observe that we have abused notation by writing $F^{2\epsilon}$ for $\epsilon^*F^{\epsilon}$. Also observe that there is a logical dualization for two pre-cosheaves.
\end{defn}

\begin{figure}
\centering
\[
	\xymatrix{F^{2\epsilon} \ar[dd] \ar[rddd] &   & F^{2\epsilon'} \ar[dd] \ar[rddd]  \ar[ll] &   \\
	&  G^{2\epsilon} \ar[dd] \ar[ld]  &   &   G^{2\epsilon'} \ar[dd] \ar[ll] \ar[ld] \\
	F^{\epsilon} \ar[dd] \ar[rddd] &    &   F^{\epsilon'} \ar[dd] \ar[ll] \ar[rddd] &  \\
	&  G^{\epsilon} \ar[dd] \ar[ld] &    &  G^{\epsilon'} \ar[dd] \ar[ll] \ar[ld] \\
	F &   &  F \ar[ll]  &   \\
	  &  G  &   &  G \ar[ll]
	}
\]
\caption{Diagram for the Proof of Lemma~\ref{lem:closed_up_int}}
\label{fig:closed_up_int}
\end{figure}
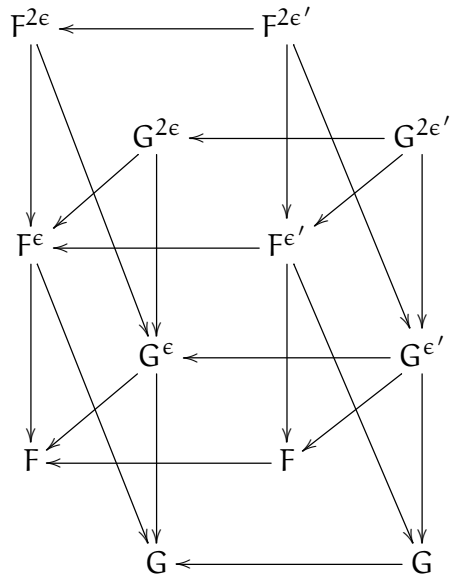

\begin{lem}\label{lem:closed_up_int}
If two presheaves $F$ and $G$ are $\epsilon$-interleaved for some $\epsilon\geq 0$, then they are $\epsilon'$-interleaved for every $\epsilon'\geq \epsilon$.
\end{lem}
\begin{proof}
Suppose that we have an $\epsilon$ interleaving, i.e. maps $\varphi_{\epsilon}:F^{\epsilon}\to G$ and $\psi_{\epsilon}:G^{\epsilon}$ that induce the correct commutative diagram. If $\epsilon'\geq \epsilon$ then the natural map $F^{\epsilon'}\to F$ factors through $F^{\epsilon}\to F$, allowing us to define $\varphi_{\epsilon'}=\varphi_{\epsilon}\eta^F_{\epsilon,\epsilon'}$ and symmetrically for $\psi_{\epsilon'}$. The map $(\epsilon')^*\varphi_{\epsilon'}$ is defined by
\[
F((U^{\epsilon'})^{\epsilon'})=F^{\epsilon'}(U^{\epsilon'})\to F^{\epsilon}(U^{\epsilon'})\to G(U^{\epsilon'})=G^{\epsilon'}(U).
\]
Observing that the map $\varphi_{\epsilon}$ is natural for a pair of open sets $U^{\epsilon}\subset U^{\epsilon'}$ proves that we have the commutative diagram in Figure~\ref{fig:closed_up_int}.
\end{proof}

\begin{rmk}
	One should note that if $F$ and $G$ are $\epsilon$-interleaved and $G$ and $H$ are $\epsilon'$-interleaved, then $F$ and $H$ are $\epsilon+\epsilon'$-interleaved.
	
	One should also note that if $F$ and $G$ are $0$-interleaved, then they are isomorphic. 
\end{rmk}

\begin{defn}
	We can define the \textbf{interleaving distance} on pre-sheaves by declaring
	\[
		d(F,G):=\inf \{ \epsilon\geq 0 \,|\, \exists \epsilon-\mathrm{interleaving}\}.
	\]
	If no interleaving exists, we define $d(F,G)=\infty$. This is what we mean by an extended metric.
\end{defn}

For pre-sheaves the interleaving distance is an extended pseudo-metric. ``Extended'' means that the distance $\infty$ is allowed and ``pseudo'' means that if $d(F,G)=0$, then it does not follow that $F=G$. For sheaves, it is true that if a map induces isomorphisms on stalks, then the map is an isomorphism of sheaves. This suggests that if for every $\epsilon>0$ there is an $\epsilon$-interleaving, then perhaps the sheaves are isomorphic. We will present an argument in this direction later on.

\subsection{Easy Stability}
\index{stability}
The following results were obtained independently of~\cite{bubenik-metrics}.

\begin{lem}
Suppose $f,g:Y\to X$ are continuous maps to a metric space that are less than $\epsilon$ distance apart in the supremum norm, then the functors
\[
\hF:U\rightsquigarrow f^{-1}(U) \qquad \mathrm{and} \qquad \hG:U\rightsquigarrow g^{-1}(U)
\]
are $\epsilon$-interleaved, when viewed as cosheaves.
\end{lem}
\begin{proof}
By hypothesis, for every open set we have the following inclusions:
\[
f^{-1}(U)\subseteq g^{-1}(U^{\epsilon}) \qquad \mathrm{and} \qquad f^{-1}(U^{\epsilon})\supseteq g^{-1}(U)
\]
This implies that we have an interleaving of pre-images:
	\[
	\xymatrix{ f^{-1}(U^{2\epsilon}) & g^{-1}(U^{2\epsilon}) \\
	f^{-1}(U^{\epsilon}) \ar[ru] \ar[u] & g^{-1}(U^{\epsilon}) \ar[lu] \ar[u]\\
	f^{-1}(U) \ar[ru] \ar[u] & g^{-1}(U) \ar[lu] \ar[u]}
	\]
Which is equivalent to saying that the functors $\hF$ and $\hG$ are $\epsilon$-interleaved.
\end{proof}

\begin{cor}\label{cor:stability}
If $f,g:Y\to X$ are continuous maps such that $d_{\infty}(f,g)<\epsilon$, then the cohomology presheaves
\[
H^iF:U\rightsquigarrow H^i(f^{-1}(U);k) \qquad \mathrm{and} \qquad H^iG:U\rightsquigarrow H^i(g^{-1}(U);k)
\]
are $\epsilon$-interleaved for every $i\geq 0$.
\end{cor}
\begin{proof}
This follows from the lemma since cohomology is a functor.
\end{proof}

\begin{cor}
If $f,g:Y\to X$ are continuous maps such that $d_{\infty}(f,g)<\epsilon$, then the pushforwards of the constant sheaf on $Y$,
\[
f_*k_Y \qquad \mathrm{and} \qquad g_*k_Y,
\]
are $\epsilon$-interleaved.
\end{cor}
\begin{proof}
The result follows by recalling the definition of the pushforward and that the constant sheaf records $H^0$ of an open set.
\end{proof}

\subsection{Global Sections Obstruct Interleavings}

In this section, we introduce the first major difference between the interleaving distance for functors $S:(\RR,\leq)\to\Vect$ and interleavings for pre-sheaves.

\begin{lem}[Obstruction to Interleavings]\label{lem:obstructs}\index{interleaving!of presheaves!obstructed by global sections}
	If $F$ and $G$ are pre-sheaves on a metric space $X$ and $F(X)\ncong G(X)$, then there is no $\epsilon$-interleaving, for any value of $\epsilon$.
\end{lem}
\begin{proof}
Clearly $X^{\epsilon}=X$. Recalling the compatibility condition for interleavings
\[	
\xymatrix{ F(X)=F^{2\epsilon}(X) \ar[dd]_{\id} \ar[rd]^{\varphi_{2\epsilon}(X)} & \\
 & G^{\epsilon}(X)=G(X) \ar[ld]^{\psi_{\epsilon}} \\
F(X) & }
\]
implies that $G(X)\to F(X)$ is a surjection. Considering the analogous triangle
\[
\xymatrix{
	& G^{2\epsilon}(X)=G(X) \ar[dd]^{\id} \ar[dl]_{\psi_{2\epsilon}} \\
	F(X)=F^{\epsilon}(X) \ar[dr]_{\varphi_{\epsilon}} & \\
	& G(X)
}
\]
implies that $F(X)\to G(X)$ is a surjection, which together proves that $F(X)\cong G(X)$. Contraposition proves the result.
\end{proof}

\begin{ex}[Skyscraper Sheaf vs. Ephemeral Module]
Suppose $\delta_{x}:\RR\to\Vect$ is the functor that assigns $k$ to the points $x\in\RR$ and $0$ everywhere else. Such a persistence module is sometimes referred to as an \textbf{ephemeral module}. Without too much trouble, one can see that this functor is interleaving distance $0$ from the zero functor.

In contrast, the skyscraper sheaf $S_x$, which assigns to every open set $U\ni x$ the vector space $k$ and $0$ to any open set not containing $x$, is infinite distance away from the zero sheaf.
\end{ex}

\section{Interleavings for Sheaves}

Now we want to emulate the above constructions for functors $F:\Open(X)^{op}\to\cD$ that satisfy the sheaf axiom. Of course, we can regard any sheaf $F$ as a pre-sheaf and apply the thickening construction to produce a pre-sheaf $F^{\epsilon}$. Is the resulting pre-sheaf also a sheaf? No, because thickening can create intersections where there shouldn't be any, as in the following example.

\begin{ex}
Let $F=k_X$ be the locally constant sheaf on $X=\R$. If $U_1=(0,1)$ and $U_2=(1,2)$, then for any $\epsilon>0$ we have that
\[
\xymatrix{ & F^{\epsilon}(U_1\cup U_2) \ar[ld] \ar[rd] & \\
F^{\epsilon}(U_1) \ar[rd] & & F^{\epsilon}(U_2) \ar[ld] \\
& F^{\epsilon}(U_1\cap U_2) & }
\qquad 
\xymatrix{ & k \ar[ld] \ar[rd] & \\
k \ar[rd] & & k \ar[ld] \\
& 0 & }
\]
is not a limit diagram. This is due to the obvious defect 
\[
(U_1\cap U_2)^{\epsilon}\neq U_1^{\epsilon}\cap U_2^{\epsilon}.
\]
\end{ex}

In light of the above example, thickening the underlying pre-sheaf must be followed by sheafification. Thus one approach to defining the thickening of a sheaf is to use the the sheafification of the thickened pre-sheaf. However, there is a nicer definition of the thickening of a sheaf that has the advantage that it can be defined in the derived setting, as well as being easier to prove theorems with. We state this definition now:

\begin{defn}[Thickening via Convolution]\index{sheaf!thickening}
Let $X$ be a metric space and let $\Bar{\Delta}^{\epsilon}_X$ be the subset of the product $\{(x,y)|d(x,y)\leq\epsilon\}$. Denote by $\epsilon_i$ the projection onto the $i^{th}$ coordinate of $\Bar{\Delta}^{\epsilon}_X$. Define the \textbf{$\epsilon$-thickened sheaf} $\widetilde{F}^{\epsilon}$ to be $\epsilon_{2*}\epsilon^*_1F$. When the context is clear we will omit the tilde.
\end{defn}

We will not show that the sheafification of the thickened pre-sheaf is the same as the above definition, as they very well could be different in certain cases.

\begin{defn}[Interleavings for Sheaves]\index{interleaving!of sheaves}
We say two sheaves $F$ and $G$ are \textbf{$\epsilon$-interleaved} if there are maps $\varphi_{\epsilon}:\widetilde{F}^{\epsilon}\to G$ and $\psi_{\epsilon}:\widetilde{G}^{\epsilon}\to F$ such that the following diagram commutes
	\[
	\xymatrix{ \widetilde{F}^{2\epsilon} \ar[d] \ar[rrd]^(.25){\epsilon^*\varphi_{\epsilon}} & & \widetilde{G}^{2\epsilon} \ar[d] \ar[lld]_(.25){\epsilon^*\psi_{\epsilon}} \\
	\widetilde{F}^{\epsilon} \ar[d] \ar[rrd]^(.25){\varphi_{\epsilon}} & & \widetilde{G}^{\epsilon} \ar[d] \ar[lld]_(.25){\psi_{\epsilon}} \\
	F & & G
	}
	\]
\end{defn}

\subsubsection{Derived Stability}

\begin{thm}\index{stability!for derived interleavings}
If $f,g:Y\to X$ are $\epsilon$-close in the sup norm, then the derived pushforwards of the constant sheaf $k_Y$ are $\epsilon$-interleaved, i.e.
\[
	d(Rf_* k_Y,Rg_* k_Y)\leq d_{\infty}(f,g)
\]
\end{thm}
\begin{proof}
\[
f^{-1}(U)\subseteq g^{-1}(U^{\epsilon}) \qquad \mathrm{and} \qquad f^{-1}(U^{\epsilon})\supseteq g^{-1}(U)
\]
The functors $\hF, \hG:\Open(X)\to\mathbfsf{Top}$ are $\epsilon$-interleaved, as already noted. This implies that the complexes of singular cochains are interleaved, which, after subdivision, defines a flabby resolution of the constant sheaf.
\end{proof}

\subsection{The Effect of Sheafification}

In this section we investigate the impact of sheafification on interleavings of pre-sheaves. After all, every sheaf is a presheaf, so we can apply the notion of interleaving to both structures. As we will show by example, two pre-sheaves can be finite distance apart, but then be infinite distance apart (not $\epsilon$-interleaved for any $\epsilon$) after sheafification. Conversely, we can produce two pre-sheaves that are infinite distance apart, whose sheafifications are interleaved.

\begin{figure}
\centering
\includegraphics[width=\textwidth]{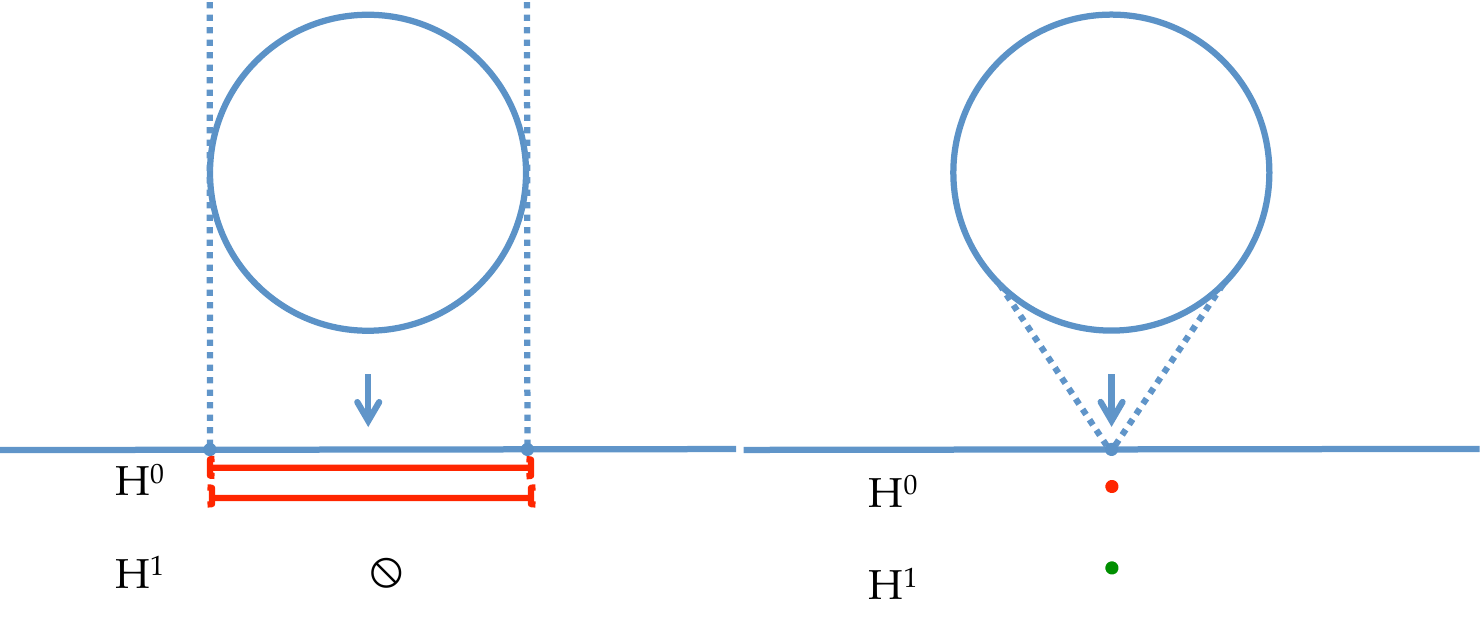}
\caption{Two Maps and their associated Leray Sheaves}
\label{fig:interleaving12}
\end{figure}

Let the map drawn on the left of Figure~\ref{fig:interleaving12} be called $f$ and let the map on the right be called $g$. As already observed in Corollary~\ref{cor:stability}, the presheaves 
\[
H^1F:U\rightsquigarrow H^1(f^{-1}(U);k) \qquad \mathrm{and} \qquad H^1G:U\rightsquigarrow H^1(g^{-1}(U);k)
\]
are $\epsilon$-interleaved for $\epsilon$ larger than the radius of the circle. However, the sheafification of both of these presheaves produce radically different sheaves. 

Recall from Definition~\ref{defn:sheafification} that the sheafification can be viewed as the sheaf of sections of the product over all the stalks mapping down to the base space.
\[
		\xymatrix{ \prod\limits_{x\in \RR} F_x \ar[d]^{\pi} \\ \RR}
\]
For $H^1F$, every stalk is the zero vector space. For $H^1G$, the stalk at $p$ is non-zero and all other ones are zero. Consequently the sheafifications are
\[
		\widetilde{F}^1 \cong 0 \qquad \mathrm{and} \qquad \widetilde{G}^1\cong S_p
\]
where $S_p$ is the skyscraper sheaf at $p$. By lemma~\ref{lem:obstructs}, these two sheaves are not interleaved.

One might conjecture in light of the above example that sheafification is a distance increasing operation. This is not the case. Consider the presheaf $H^1F$ from the example above. One can easily see that it is \emph{not} interleaved with the zero sheaf. However, the sheafification of $H^1F$ is the zero sheaf. So sheafification took two presheaves that were infinite distance apart and returned isomorphic (distance zero) sheaves.

\subsection{Thickening Global Sections}

We will now prove that global sections remain unaffected by this thickening procedure. This provides us with a sheafified version of Lemma~\ref{lem:obstructs}.

\begin{lem}
If $X$ is a metric space such that the closed $\epsilon$-ball is connected for every point of $x$ and the projection map $\epsilon_1:\Bar{\Delta}^{\epsilon}_X\to X$ is a closed map, then thickening preserves global sections, i.e.
\[
	\widetilde{F}^{\epsilon}(X)\cong F(X).
\]
\end{lem}
\begin{proof}
This follows from a simpler version of the Vietoris mapping theorem, stated as Theorem 11.7 of~\cite{Bredon}. Essentially, the Vietoris mapping theorem guarantees that
\[
H^0(X;F)\cong H^0(\Bar{\Delta}^{\epsilon}_X;\epsilon_1^*F)
\]
and since the pushforward $\epsilon_{2*}\epsilon_1^*F$ has the same global sections of $\epsilon_1^*F$ by functoriality (global sections are gotten by pushing forward to a point), then the result follows.
\end{proof}

\subsection{Metric on Sheaves}

\begin{clm}
If $F,G\in\Shv(X)$ are interleaving-distance zero apart, i.e. there is a sequence of $\{\epsilon_n\}$ converging to zero where $F$ and $G$ are $\epsilon_n$ interleaved for each $n$, then $F\cong G$.
\end{clm}
\begin{proof}[Sketch]
We are going to define maps between the \'etal\'e spaces
\[
\varphi:\prod_{x\in X}F_x \to \prod_{x\in X} G_x \qquad \mathrm{and} \qquad \psi:\prod_{x\in X} G_x \to \prod_{x\in X}F_x
\]
with the property that they are inverses of one another.

Given any element $s_x\in F_x$ there exists a $U\ni x$ and a section $s_U$ such that $(s_U)_x=s_x$. However, since $U$ is open, there exists an $r>0$ and $\epsilon_n>0$ such that $B(x,r+2\epsilon_n)\subset U$. This implies in turn that $B(x,r)^{2\epsilon_n}\subset U$. Consequently, there is a lift $s^{2\epsilon_n}_x\in \widetilde{F}^{2\epsilon_n}_x$ of $s_x$, gotten by taking the image of $s_u$ under the following composition:
\[
F(U)\to F(B(x,r+2\epsilon_n))\to F(B(x,r)^{2\epsilon_n})\to \widetilde{F}^{2\epsilon_n}_x
\]

We then define 
\[
\varphi(s_x):=\varphi^{\epsilon_n}_x\circ\widetilde{\eta_{\epsilon_n,2\epsilon_n}}^F_x(s^{2\epsilon_n}_x)
\]
Of course, by the definition of interleaving, we have a natural choice of a lift of $\varphi(s_x)$ to $\widetilde{G}^{\epsilon_n}_x$ given by 
\[
\widetilde{\eta_{\epsilon_n}}^G_x\varphi^{2\epsilon_n}_x(s^{2\epsilon_n}_x),
\]
which by construction has the property that applying $\psi^{\epsilon}_x$ yields $s_x$.

Of course it needs to be checked that such a lift is picked out by the symmetric construction of the map $\psi$. It could have happened that a lift to $\widetilde{G}^{\epsilon_n'}_x$ for $\epsilon_n'\neq \epsilon_n$. However, one only needs to observe that for any sheaf $F$ we have the following commutative diagram for any pair $0<\epsilon_n'<\epsilon_n$:
\[
\xymatrix{\widetilde{F}^{2\epsilon_n} \ar[d] \ar[r] & \widetilde{F}^{\epsilon_n} \ar[d] \ar[r] & F \\
\widetilde{F}^{2\epsilon_n'} \ar[r] & \widetilde{F}^{\epsilon_n'} \ar[ur] & }
\]
Arguing symmetrically and checking continuity of the resulting maps $\varphi$ and $\psi$ completes the proof.
\end{proof}

\begin{thm}[Skeletal Metric Space]
Let $X$ be an arbitrary metric space. The interleaving distance induces an extended metric on the skeleton of the category of sheaves $\Shv(X)$. Here extended means that the value $d(F,G)=\infty$ is allowed, i.e. there is no interleaving between $F$ and $G$ whatsoever.
\end{thm}
\begin{proof}
Recall that the skeleton of a category $\cC$ is a full, isomorphism-dense subcategory $\cS$ in which no two distinct objects are isomorphic. 

If we take $\cC=\Shv(X)$, on which the interleaving distance already defines an extended pseudo-metric, then the above result implies that if $d(F,G)=0$, then $F\cong G$ and hence in a skeletal subcategory $F=G$. This implies that the interleaving distance is an extended metric when restricted to any skeletal subcategory of $\Shv(X)$.
\end{proof}

As one can imagine, the space of sheaves viewed as a metric space can be enormously complicated. Every map $f:X\to\R^n$ has an associated sheaf on $\R^n$, simply by considering the pushforward of the constant sheaf. This includes every possible subspace $Y\subset \R^n$ with the pushforward of the constant sheaf along this inclusion serving as a sort of ``indicator function'' on it. The interleaving distance would give us one notion of distance between all these possible subspaces. In this case, there is a ready comparison to be made with the Gromov-Hausdorff distance between metric spaces, which is more refined than the interleaving distance. However, the interleaving distance also gives a distance between \emph{information} on top of a metric space, where information is encoded via a sheaf. 

\section{The Space of Constructible Sheaves over $\R$}

In this section we will give an explicit description of the space of constructible/cellular sheaves on $\R$ with the interleaving distance as a metric. It turns out that one can use the indecomposable sheaves to give a set of ``coordinates'' on this space. It will turn out that the space resembles a disjoint union of configuration spaces, where there the components are divvied up by global sections, i.e. $H^0$. A comparison to McDuff's construction of the tangent space use a configuration space of points, where points can disappear, is made.

First we recall the definition of constructible sheaf pertinent to this section. Because there are competing, more general notions of a constructible sheaf, we will use slightly different terminology.

\begin{defn}\index{category!of definable sheaves}\index{sheaf!definable}
Let $\Shv_d(\R)$ denote the \textbf{category of definable sheaves} over the real line $\R$, equipped with the usual Euclidean topology. Specifically, a sheaf $F$ will be regarded as a contravariant functor from the open set category with the necessary gluing properties. Such a sheaf $F$ is \textbf{definable} if $\R$ can be written as the finite union of open intervals and points 
\[
\R=(-\infty,a_0)\cup \{a_0\}\cup \cdots (a_i,a_{i+1}) \cdots \cup \{a_n\} \cup (a_n,\infty)
\]
such that when restricted to each interval the sheaf is locally constant. We do \textbf{not} assume that every sheaf is constructible with respect to the same set of intervals.

We will find it convenient to work with a subcategory of this category given by the definable sheaves with \textbf{finite support}, $\Shv_{d,f}(X)$, where the sheaf must restrict to zero on the two half-open intervals including $\pm\infty$.
\end{defn}

As already established in this thesis, such a sheaf is completely described via a zig-zag of vector spaces and linear maps
\[
F(a_0)\leftarrow F(x_0) \rightarrow F(a_1) \leftarrow F(x_1)\rightarrow \cdots \leftarrow F(x_n) \rightarrow F(a_n)
\]
where we have abbreviated the intervals and points by using $a$'s and $x$'s along with subscripts appropriately.

By Gabriel's theorem, we know that such a diagram amounts to a representation of an $A_n$-type quiver and can be decomposed into finitely many indecomposable representations. These indecomposables, in view of the cell structure, can be regarded as one of the following sheaves:
\begin{itemize}
\item $k_{[x_i,x_j]}$ --- the constant sheaf on the closed interval $[x_i,x_j]$
\item $k_{(x_i,x_j)}$ --- the constant sheaf on the open interval $(x_i,x_j)$
\item $k_{[x_i,x_j)}$ --- the constant sheaf on the half-open interval $[x_i,x_j)$
\item $k_{(x_i,x_j]}$ --- the constant sheaf on the half-open interval $(x_i,x_j]$
\item $k_{\R}$ --- the constant sheaf supported on the whole real line $\R$
\item $S_x$ --- the skyscraper sheaf concentrated on the vertex $x$
\end{itemize}
The skyscraper sheaf, $S_x$, is just a special instance of the constant sheaf on a closed interval $[x_i,x_j]$ where $x_i=x_j$. As such, we will not always distinguish the skyscraper sheaf from the constant sheaf on the closed interval. Similarly, one can view the constant sheaf on $\R$, $k_{\R}$, as a degenerate version of the open interval. 

To keep the notation clean and free ourselves from a particular declaration of cell structure on the real line, we will speak of the \textbf{four indecomposable sheaves} on the real line:
\[
k_{[b,d]} \qquad k_{(b,d)} \qquad k_{[b,d)} \qquad k_{(b,d]}
\]

\begin{rmk}
The constant sheaf, $k_{\R}$, and any of the others where $b$ or $d$ is $\pm\infty$ are excluded from the category of definable sheaves with finite support.
\end{rmk}

In the following sections it will be paramount to understand when there is and isn't a non-zero map of sheaves between these four types.

\begin{prop}[D\'evissage for 1D Indecomposable Sheaves]\label{prop:devissage}
We have the following explicit characterizations for the space of sheaf morphisms:
\begin{itemize}
\item For closed intervals $I_1=[b_1,d_1]$ and $I_2=[b_2,d_2]$ we have
\begin{equation*}
 \Hom_{\Shv}(k_{I_1},k_{I_2})=
\begin{cases}
 k & \text{if } I_2\subseteq I_1, \\
0 & \text{o.w.}
\end{cases}
\end{equation*}
\item For open intervals $I_1=(b_1,d_1)$ and $I_2=(b_2,d_2)$ we have
\begin{equation*}
 \Hom_{\Shv}(k_{I_1},k_{I_2})=
\begin{cases}
 k & \text{if } I_1\subseteq I_2, \\
0 & \text{o.w.}
\end{cases}
\end{equation*}
\item For half open intervals of the form $I_n=[b_n,d_n)$ we have
\begin{equation*}
 \Hom_{\Shv}(k_{I_1},k_{I_2})=
\begin{cases}
 k & \text{if } b_1\leq b_2 \text{ and } b_2<d_1\leq d_2 , \\
0 & \text{o.w.}
\end{cases}
\end{equation*}
\item For half open intervals of the form $I_n=(b_n,d_n]$ we have
\begin{equation*}
 \Hom_{\Shv}(k_{I_1},k_{I_2})=
\begin{cases}
 k & \text{if } b_2\leq b_1 \text{ and } b_1<d_2\leq d_1, \\
0 & \text{o.w.}
\end{cases}
\end{equation*}
\item For a closed interval $I_1=[b_1,d_1]$ and any non-compact interval $I_2$
\begin{equation*}
 \Hom_{\Shv}(k_{I_1},k_{I_2})= 0
\end{equation*}
\end{itemize}
\end{prop}
\begin{proof}
The calculations all follow from considering the behavior of a cellular sheaf map near the endpoints of the above indecomposables. Specifically, in a small enough neighborhood of any point in the real line, an indecomposable cellular sheaf has one of the following three forms, possibly after reflection.
\[
M=\quad 0\leftarrow k \rightarrow k
\]
\[
U=\quad k\leftarrow k \rightarrow k
\]
\[
H=\quad k\leftarrow 0 \rightarrow 0
\]
Clearly there are non-zero natural transformations $H\Rightarrow U\Rightarrow M$, but every natural transformation $M\Rightarrow U\Rightarrow H$ must be zero.
\end{proof}

\begin{rmk}[D\'evissage for Constructible Sheaves]
In David Nadler's beautiful application of constructible sheaves to the study of the Fukaya category~\cite{nadler-catmorse}, he refers to the diagram
$$
\xymatrix{
\Shv(V) \ar@/^2pc/[rr]^-{j_!}\ar@/_2pc/[rr]^-{j_*} && \ar[ll]_-{j^! \simeq j^*} \Shv(X)  \ar@/^2pc/[rr]^-{i^*} 
\ar@/_2pc/[rr]^-{i^!} 
&& \ar[ll]_-{i_! \simeq i_*} \Shv(Y)
}
$$
as the ``d\'evissage pattern for constructible sheaves'' --- an ode to Grothendieck's method for studying \emph{coherent} sheaves. Here $j:V\to X$ is the inclusion of an open set and $i:Y:=X-V\to X$ is the inclusion of the closed complement. Here the categories $\Shv(-)$ refer to the full differential graded category of constructible complexes of sheaves, $H^0$ of which is the usual derived category.
\end{rmk}

\subsection{Interleavings and Dynamics on Indecomposable Sheaves}
\index{interleaving!of sheaves!over the real line}
Since indecomposable sheaves are the underlying elements that build up a sheaf, we investigate the behavior of each of these sheaves under epsilon-thickening. We summarize the result of each of these calculations below:
\begin{itemize}
\item If $F=k_{[b,d]}$ and $\epsilon \geq 0$, then $\widetilde{F}^{\epsilon}=k_{[b-\epsilon,d+\epsilon]}$.
\item If $F=k_{(b,d)}$ and $0\leq \epsilon < d-b$, then $\widetilde{F}^{\epsilon}=k_{(b+\epsilon,d-\epsilon)}$. If $d-b\leq \epsilon$, then $\widetilde{F}^{\epsilon}=0$.
\item If $F=k_{[b,d)}$ and $\epsilon \geq 0$, then $\widetilde{F}^{\epsilon}=k_{[b-\epsilon,d-\epsilon]}$.
\item If $F=k_{(b,d]}$ and $\epsilon \geq 0$, then $\widetilde{F}^{\epsilon}=k_{[b+\epsilon,d+\epsilon]}$.
\end{itemize}
With these calculations in hand, we can then discuss distance between these sheaves.
\begin{itemize}
\item The sheaf $k_{[b,d]}$ is interleaving distance $r=(d-b)/2$ from the skyscraper sheaf $S_{m}$ where $m=(d+b)/2$ --- the midpoint --- and infinite distance from any of the other three types of indecomposables, as well as the zero sheaf.
\item The sheaf $k_{(b,d)}$ is interleaving distance $r=(d-b)/2$ from the zero sheaf $0$.
\item The sheaves $k_{[b,d)}$ and $k_{(b.d]}$ are interleaving distance $r=(d-b)/2$ from the zero sheaf $0$.
\end{itemize}

We now give the supporting arguments for these calculations.

\begin{prop}
If $F=k_{[b,d]}$ and $\epsilon\geq 0$, then $\widetilde{F}^{\epsilon}=k_{[b-\epsilon,d+\epsilon]}$.
\end{prop}
\begin{proof}
To construct $\widetilde{F}^{\epsilon}$ it suffices to consider the stalks of the thickened pre-sheaf. It suffices to consider the extreme points. Consider the point $x=b-\epsilon$, then any ball $B(x,r)$ has the property that $[b,d]\cap B(x,r)^{\epsilon}=[b,b+r)\neq\emptyset$. Since $F$ is defined as the pushforward of the constant sheaf along $j:X=[b,d]\hookrightarrow \R$, then $F^{\epsilon}(B(x,r))=k_X([b,b+r))=k$. This proves that $F_x\cong k$.
\end{proof}

\begin{prop}
If $F=k_{(b,d)}$ and $0\leq \epsilon < d-b$, then $\widetilde{F}^{\epsilon}=k_{(b+\epsilon,d-\epsilon)}$. If $d-b\leq \epsilon$, then $\widetilde{F}^{\epsilon}=0$.
\end{prop}
\begin{proof}
If $j:W=(b,d)\hookrightarrow \R$ denotes the inclusion of the open interval, then we can identify $F=j_! k_W$. Here $j_!$ denotes the pushforward with compact supports functor. For the inclusion of a locally closed subspace $W$ into a general topological space $X$ \cite{iversen} provides a precise description. For a sheaf $F$ on $W$ the sheaf $j_!F$ has sections on an open set $U$ given by
\[
\Gamma(U,j_!F):=\{s\in \Gamma(W\cap U,F) \,|\, \supp(s) \, \mathrm{closed}\,\mathrm{rel.}\, U\}
\]
The support of a section is the set of points where a section has non-vanishing stalks.

For an open set $U\subset\R$, we can describe the sheaf pertinent to us even more explicitly: $j_! k_W(U)$ is non-zero if and only if there is a closed set $Y$ such that $U\subset Y\subset W$ --- we say that $U$ is \textbf{completely contained} in $W$. Note that $W=(b,d)$ is \emph{not} completely contained in itself, so in particular $(b,d)$ is assigned the zero vector space by $j_!k_W$.

Consequently, any point $x$ within distance $\epsilon$ of the boundary of $(b,d)$ will fail to possess a ball $B(x,r)$ such that $B(x,r)^{\epsilon}$ is completely contained in $(b,d)$. Thus we say that $(b,d)$ is ``eroded'' by distance $\epsilon$.
\end{proof}

\begin{prop}
If $F=k_{[b,d)}$ and $\epsilon \geq 0$, then $\widetilde{F}^{\epsilon}=k_{[b-\epsilon,d-\epsilon]}$.
\end{prop}
\begin{proof}
This follows by considering each of the above arguments separately about each endpoint.
\end{proof}

Because the other case is obviously symmetric, we omit a separate argument.

An extremely important observation that distinguishes interleavings of sheaves from the usual context of an interleavings of persistence modules~\cite{bauer2013induced} is that the sheaf $k_{[b,d]}$ is not interleaved with the zero sheaf. This follows from the fact that the space of global sections is preserved by thickening.

\subsection{Coordinates for the Category of Sheaves}
\index{category!of sheaves!coordinates}

We can now use the decomposition theorem for constructible sheaves over the real line to give explicit coordinates for each isomorphism class of a definable sheaf with finite support.

Let $\HH=\{(x,y)\in \R^2 \,|\, y\geq 0 \}$ denote the closed upper half plane. To each of the four indecomposable sheaves with finite support --- $k_{[b,d]}$, $k_{(b,d)}$, $k_{[b,d)}$, and $k_{(b,d]}$ --- we can associate a point in $\HH$ as follows:
\[
I \rightsquigarrow (x,y)=(m(I),r(I)):=(\frac{b+d}{2},\frac{d-b}{2})
\]
Here $I$ is a stand-in for any of the four types of indecomposable sheaves. The variable names $m$ and $r$ are meant to connote the midpoint and the radius, respectively, of the underlying bar in the ``barcode''. One can then associate to any definable sheaf $F$ with finite support the following coordinates
\[
F\cong \bigoplus_{k=1}^n k_{I_k}\rightsquigarrow \{(m(I_1),r(I_1)),\cdots,(m(I_n),r(I_n))\}\in \HH^{n_1} \sqcup \HH^{n_2} \sqcup \HH^{n_3} \sqcup \HH^{n_4}
\]
where $n_1,n_2,n_3,n_4$ refers to the number of closed, open, half-open on the right and half-open on the left indecomposables occurring in the decomposition for $F$, respectively. Of course, $n=n_1+n+2+n_3+n_4$. 

As presented, the space has too many points for the simple reason that the line $r=0$ in $\HH$ must be identified with the zero sheaf for the non-closed indecomposable types. To capture the full category of sheaves we must append a distinguished basepoint $\star$ to represent the zero sheaf, quotient the upper half plane so as to identify $\{r=0\}\sim \{\star\}$ and then form a few infinite symmetric products. The first part is simple. We define
\[
Z:=(\HH\sqcup \{\star\})/\sim \qquad \mathrm{where} \qquad (m,r)\sim \star \,\,\mathrm{iff}\,\, r=0.
\]
The next construction begins by observing that $Z$ is naturally a pointed space, where $\star$ serves as the distinguished basepoint. To every pointed space we can associate a new space called the infinite symmetric product.

\begin{defn}
 Recall that the \textbf{$n$-fold symmetric product} of $X$, denoted $SP_n(X)$, is given by forming the $n$-fold Cartesian product and quotienting out by the action of the symmetric group, i.e. $SP_n(X):=X^n/\Sigma_n$. 

 Let $Z$ be a pointed topological space, whose distinguished point is $\star$. There is a system of embeddings $Z^n\hookrightarrow Z^{n+1}$ giuven by sending any point 
\[
(z_1,\ldots,z_n) \mapsto (z_1,\ldots,z_{n+1},\star).
\] 
This embedding descends to an embedding $SP_n(Z)\hookrightarrow SP_{n+1}(Z)$, which forms a directed system of spaces. The \textbf{infinite symmetric product} $SP(Z)$ is the direct limit of this system, i.e.
\[
SP(Z)=\varinjlim_{n>0}SP_n(Z).
\]
\end{defn}

We can now state a theorem.

\begin{thm}
Isomorphism classes of $\Shv_{d,f}(\R)$ are in bijective correspondence with points in the following space:
\[
\BB:=\bigsqcup_{n>0}SP_n(\HH)\times \left(SP(Z)\bigvee_{\star} SP(Z) \bigvee_{\star} SP(Z)\right)
\]
\end{thm}

\begin{rmk}
If one could show that the topology induced by the interleaving distance made $\HH$ and $Z$ into cell complexes, then the Dold-Thom theorem would tell us the singular homology of the space of definable sheaves on the real line is isomorphic to a countably infinite number of copies of $\ZZ$ in degree zero and zero in all higher degrees.
\end{rmk}

It remains to be seen what geometry is induced on $\BB$ by pulling back the interleaving distance. Conjecturally, this should be accomplished by the bottleneck distance~\cite{chazal2009proximity,bauer2013induced}, with the stipulation that only points in $SP(Z)$ can be matched with zero.

\subsection{Towards a Bottleneck Distance for Sheaves}

We now investigate simpler descriptions of the interleaving distance for definable sheaves with finite support on the real line. Our goal is to establish a connection between the interleaving distance for sheaves and the bottleneck distance, which we now define.

\begin{defn}[The Bottleneck Distance]\index{bottleneck distance}
Let $\Delta_+:=\{(x,y)\in\R^2| x\leq y\}$ be equipped with the sup norm 
\[
d_{\infty}(p,p')=\sup\{|x_1-x_2|,|y_1-y_2|\}
\]
A multiset $D$ in $\Delta_+$ is a subset $|D|$ of $\Delta_+$ equipped with a multiplicity function $\mu:|D|\to \NN$. Any multiset can be considered as a set via disjoint unions
\[
D=\bigcup_{p\in |D|}\coprod_{i=1}^{\mu(p)}\{p\}.
\]
A multi-bijection $m:D\to D'$ is a bijection between the underlying sets, where if $\mu(p)=k$, then the underlying set has $k$ elements corresponding to $p$. The \textbf{bottleneck distance} is a distance between multisets defined by the formula
\[
d_B(D,D'):=\inf_{m}\sup_{p\in D} d_{\infty}(p,p').
\]
\end{defn}

In usual sub-level set persistence, the multisets $D$ consist of the diagonal $\Delta=\{(x,x)\}$, equipped with infinite multiplicity, and the points $(b,d)$ corresponding to the interval modules $k_{[b,d)}$ making up the interval decomposition of theorem \ref{thm:crawleyboevey}. The stipulation that the diagonal has infinite multiplicity reflects the fact that for half-open intervals, the module can be interleaved with zero. For sheaves, the obstruction by global sections result implies that the number of points corresponding to indecomposables $k_{[b,d]}$ is an invariant of the multi-set --- infinite multiplicity of the diagonal cannot be used there. Nevertheless, we can describe the geometry there in our choice of coordinates.

\begin{lem}
Suppose $F=k_{[m_1-r_1,m_1+r_1]}$ and $G=k_{[m_2-r_2,m_2+r_2]}$ are two indecomposable sheaves supported over closed intervals, then the interleaving distance for $F$ and $G$ is their distance in a taxicab metric on $\HH$, i.e.
\[
d(k_{[m_1-r_1,m_1+r_1]},k_{[m_2-r_2,m_2+r_2]})=|m_1-m_2|+|r_1-r_2|.
\]
\end{lem}
\begin{rmk}
If we write $[m_1-r_1,m_1+r_1]=:[b_1,d_1]$ and $[m_2-r_2,m_2+r_2]=:[b_2,d_2]$ then the taxicab metric specializes to a sup-norm on the space $\{(b,d)\in\R^2\,|\,b\leq d\}$, that is to say
\[
|m_1-m_2|+|r_1-r_2|=\sup\{|b_1-b_2|,|d_1-d_2|\},
\]
which is a special instance of the bottleneck distance on persistence diagrams.
\end{rmk}
\begin{proof}
Without loss of generality we can assume that $m_1\leq m_2$. We know that there can only be a non-zero map from $F^{\epsilon}\to G$ if $[m_2-r_2,m_2+r_2]\subset [m_1-r_1-\epsilon,m_1+r_1+\epsilon]$, i.e. if $\epsilon\geq (m_2-m_1)+(r_2-r_1)$. Similarly, there is a non-zero map $G^{\epsilon}\to F$ only if $\epsilon\geq (m_2-m_1)+(r_1-r_2)$. Since in order for a non-zero interleaving to exist both maps must be non-zero, we conclude that
\[
\epsilon=\sup\{(m_2-m_1)+(r_2-r_1),(m_2-m_1)+(r_1-r_2)\}=m_2-m_1 + |r_1-r_2|.
\]
This proves that there is a non-zero interleaving for this value of $\epsilon$. However, the interleaving distance is the infimum over all such $\epsilon$. However, for any smaller $\epsilon$ one of the maps $\varphi_{\epsilon}$ or $\psi_{\epsilon}$ must be zero. However, global sections obstructs such a pair of maps from defining an interleaving, since for every $\epsilon$ the maps
\[
\eta^F_{2\epsilon}:F^{2\epsilon}(X) \to F(X) \qquad \mathrm{and} \qquad \eta^G_{2\epsilon}:G^{2\epsilon}(X) \to G(X)
\]
are non-zero. They are, in fact, the identity map $\id:k\to k$.
\end{proof}

The situation for constant sheaves supported on half-open intervals or the open interval is more complicated since they have no global sections and can be interleaved with the zero sheaf.

\begin{lem}\label{lem:distance_from_0}
Suppose $F$ and $G$ are each interleaved with the zero sheaf, i.e. $d(F,0),d(G,0)<\infty$. Then the interleaving distance between $F$ and $G$ is bounded above by the greater of the two distances from zero,
\[
d(F,G)\leq \sup\{d(F,0),d(G,0)\}.
\]
\end{lem}
\begin{proof}
Choose any $\epsilon\geq \sup\{d(F,0),d(G,0)\}$, then the zero maps define an $\epsilon$ interleaving, factoring through zero, between $F$ ang $G$. Since the interleaving distance is the infimum, the inequality follows.
\end{proof}

We can use the above lemma to establish what the interleaving distance between pairs of sheaves of the other three types looks like.

\begin{lem}
Let $F=k_{[m_1-r_1,m_1+r_1)}$ and $G=k_{[m_2-r_2,m_2+r_2)}$ be indecomposable sheaves supported on half-open intervals, then
\[
 d(k_{[m_1-r_1,m_1+r_1)},k_{[m_2-r_2,m_2+r_2)})=\inf\{\max(r_1,r_2),|m_1-m_2|+|r_1-r_2|\}
\]
\end{lem}
\begin{proof}
Without loss of generality, we assume that $F=k_{[-R,R)}$ and $G=k_{[m-r,m+r)}$ where $m\geq 0$. This can be done since the $\epsilon$-thickening operation simply translates sheaves of this form to the left by $\epsilon$. In view of lemma \ref{lem:distance_from_0} it suffices to consider the possible values for $m$ and $r$ such that $G^{\epsilon}$ admits a non-zero map to $F$. By proposition \ref{prop:devissage} we can determine precisely what inequalities $\epsilon$ must satisfy.

First we assume that $r\leq R$. In this case, the inequalities are
\[
-r-R\leq m-\epsilon \leq r-R \qquad \mathrm{and} \qquad R-r \leq m+\epsilon < R+r.
\]
The first set of inequalities bounds when $G^{\epsilon}\to F$ can be non-zero and the second set of inequalities bounds when $F^{\epsilon}\to G$ can be non-zero. The smallest such $\epsilon$ that satisfies the first set is $\epsilon=m+R-r$, which only satisfies the second set if $m< R$. For these values of $m$ and $r$, there is the smallest non-zero interleaving, so the distance is
\[
d(F,G)=m+R-r.
\]
Similarly for $r\geq R$, $\epsilon=m+r-R$ is the smallest possible value for $G^{\epsilon}\to F$ to be non-zero. Requiring $F^{\epsilon}\to G$ to be non-zero as well implies that $m\leq R$ and for these values of $m$ and $r$ the distance is
\[
d(F,G)=m+r-R.
\]
Consequently, wherever a non-zero interleaving is possible, the smallest such value is given in the desired taxicab form $m+|r-R|$. Now we can apply lemma \ref{lem:distance_from_0} to determine whether the non-zero interleavings are the smallest possible. We split into two cases. Suppose $r\leq R$, then $m+R-r\leq R$ ---  the largest distance from zero --- when $m\leq r$, precisely where we determined a non-zero interleaving exists. Similarly if $R\leq r$, then we have that $m+r-R\leq r$ precisely if $m\leq R$.
\end{proof}

With this evidence in hand, we believe a version of the isometry theorem~\cite{lesnick-thesis,bubenik2012categorification,bauer2013induced} should hold for definable sheaves with finite support over the real line. However, we must delay this for another time.


\clearpage

\backmatter

\bibliographystyle{alpha}
\bibliography{THESIS}

\printindex

\end{document}